%% file: DiffGeo5_0.tex
\pdfoutput=1    

\documentclass[11pt,a4paper]{article} 

\usepackage[OT2,T1]{fontenc}
\usepackage[utf8]{inputenc}
\usepackage[greek,russian,english]{babel}
\usepackage{amsmath}
\usepackage{amssymb,verbatim}
\usepackage{amsfonts}
\usepackage{amsthm}
\usepackage{pifont}
\usepackage{graphicx}
\usepackage{datetime}
\usepackage{scrdate}
\usepackage{enumitem}
\usepackage{geometry}
\usepackage{float}
\usepackage{bigints}
\usepackage{multirow}
\usepackage{bbm} 
\usepackage{dsfont} 
\usepackage{cases}
\usepackage[labelfont=bf]{caption}
\usepackage{tabularx}
\usepackage{hyperref}
\usepackage[
  backend=biber,
  natbib=true,
  ]{biblatex}
\usepackage{parskip} 
\usepackage{sidecap} 
\usepackage{esvect}
\usepackage{microtype} 
\usepackage{chngcntr} 
\counterwithin{figure}{section}
\usepackage{makeidx}
\usepackage{tikz}
\usepackage{esint} 
\usepackage[outdir=./]{epstopdf}
\usepackage{csquotes}

\hypersetup{colorlinks=true, linkcolor=blue, citecolor=blue, pdfstartview=FitH, linktocpage=false}

\makeindex

\settimeformat{hhmmsstime} 
\addto\captionsngerman{%
}

\newtheorem{thm}{Theorem}[section]
\newtheorem{lem}{Lemma}[section]
\newtheorem{cor}{Corollary}[section]

\theoremstyle{definition}
\newtheorem{defin}{Definition}[section]

\newtheorem{remark}{Remark}[section]
\newtheorem{example}{\bs\;Example}[section]

\numberwithin{equation}{section}
\addto\captionsenglish{

}
\numberwithin{table}{section}

\makeatletter 
\newcommand\mynobreakpar{\par\nobreak\@afterheading} 
\makeatother

\newcommand{\N}{\mathbb{N}}
\newcommand{\Z}{\mathbb{Z}}
\newcommand{\R}{\mathbb{R}}
\newcommand{\CC}{\mathbb{C}}
\newcommand{\beq}{\begin{equation*}}
\newcommand{\eeq}{\end{equation*}}
\newcommand{\beqn}{\begin{equation}}
\newcommand{\eeqn}{\end{equation}}
\newcommand{\dd}{\mathrm{d}}
\newcommand{\tn}{\textnormal}

\newcommand{\ph}{\varphi}
\newcommand{\rh}{\varrho}
\newcommand{\ep}{\varepsilon}
\newcommand{\om}{\omega}
\newcommand{\vt}{\vartheta}
\newcommand{\la}{\lambda}

\newcommand{\Om}{\varOmega}
\newcommand{\bs}{$\bigstar$}
\newcommand{\db}{\displaybreak[0]}
\newcommand{\ii}{\mathrm{i}}
\newcommand{\ee}{\mathrm{e}}
\newcommand{\p}{\partial}

\newcommand{\K}{\mathcal{C}}
\newcommand{\g}{^\tn{o}}

\newcommand{\fab}{f_{a,\:\!b}}

\DeclareMathOperator{\sgn}{sgn}
\DeclareMathOperator{\Rez}{Re}
\DeclareMathOperator{\Imz}{Im}

\newcommand{\entspr}{\;\widehat{=}\;} 

\begin{document}
\input{DiffGeo5_1}
\input{DiffGeo5_2}
\input{DiffGeo5_3a}
\input{DiffGeo5_3b}
\input{DiffGeo5_4a}
\input{DiffGeo5_4b}
\input{DiffGeo5_4c}
\input{DiffGeo5_4d}
\input{DiffGeo5_4e}
\input{DiffGeo5_4f}
\input{DiffGeo5_4g}
\input{DiffGeo5_5}
\input{DiffGeo5_L}
\end{document}

%% file: DiffGeo5_1.tex

\title{Application of complex-valued functions in plane\\ differential geometry and kinematics}
\author{Uwe Bäsel}
\date{} 
\maketitle
\thispagestyle{empty}
\begin{abstract}
\noindent
In this paper, we discuss some problems of elementary plane differential geometry and kinematics.
Although the results are not new, the consistent use of complex-valued functions (plane curves) of a real variable (parameter) allows to derive them ab ovo in a particularly simple, uniform and transparent way. 
A number of examples with figures complete the explanations.\\[0.2cm]
\textbf{2020 Mathematics Subject Classification:}
53-01, 
53A04, 
51M25, 
53A17, 
70B15  
\\[0.2cm]
\textbf{Keywords:} 
Scalar product of complex numbers,
quasi vector product,
parametric curve,
signed area,
arc length,
oriented curvature,
evolute,
involute,
envelope,
support function,
Gerono lemniscate,
Tait-Kneser theorem,
four-vertex theorem,
dyad,
five-bar linkage,
coupler curve,
plane motion,
Bresse circles,
zero-normal jerk circle,
cubic of stationary curvature,
four-bar linkage,
cam mechanism,
machine elements,
shaft-hub-connection,
polygon profiles,
machine tools
\end{abstract}

\tableofcontents


%% file: DiffGeo5_2.tex

\section{Some basics about complex numbers}

\begin{figure}[ht]
\begin{minipage}{0.48\textwidth}
  \centering
  \includegraphics[width=0.7\textwidth]{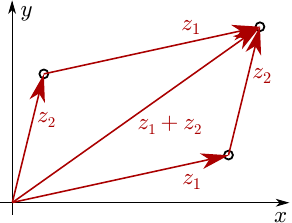}
  \caption{Addition of complex numbers $z_1$ and~$z_2$}
  \label{Abb:Addition_kompl_Zahlen}
\end{minipage}
\hfill
\begin{minipage}{0.48\textwidth}
  \centering
  \includegraphics[width=0.7\textwidth]{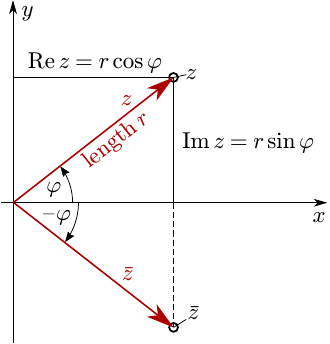}
  \caption{Complex number $z$ and complex conjugate $\bar{z}$ of $z$}
  \label{Abb:Komplexe_Zahl01} 
\end{minipage}
\end{figure}

A complex number $z$ (see e.g.\ \cite{Gamelin}, \cite{Pieper}, \cite{Remmert}) is an expression of the form
\beq
  z = x + \ii y\,,\quad
  \ii = \sqrt{-1}\,,
\eeq
where $x$ and $y$ are real numbers.
$x$ and $y$ are called {\em real part} and {\em imaginary part} of $z$, respectively; one writes for them
\beq
  \Rez z = x\,,\qquad
  \Imz z = y\,.
\eeq
A complex number $z$ is called {\em real} if $z = x$, and {\em purely imaginary} if $z = \ii y$.
The complex numbers can be visualized as points or vectors in the {\em Cartesian complex plane} with real and imaginary axis.
The unit vector on the real axis is the $1$, the unit vector on the imaginary axis is $\ii$.
The addition of two complex numbers $z_1$ and $z_2$ (see Fig.\ \ref{Abb:Addition_kompl_Zahlen}) corresponds to the ordinary vector addition:
\beq
  z_1 + z_2
= (x_1 + \ii y_1) + (x_2 + \ii y_2)
= (x_1 + x_2) + \ii (y_1 + y_2)\,.  
\eeq

The {\em complex multiplication} (the {\em complex product}\index{product!complex product@complex $ \sim $}) is defined by
\beq
  z_1 \cdot z_2
= (x_1 + \ii y_1) \cdot (x_2 + \ii y_2)
= (x_1 x_2 - y_1 y_2) + \ii (x_1 y_2 + y_1 x_2)\,.  
\eeq
For the complex multiplication the commutative law $z_1 \cdot z_2 = z_2 \cdot z_1$, the distributive laws $z_1 \cdot (z_2 + z_3) = z_1 \cdot z_2 + z_1 \cdot z_3$ and $(z_1 + z_2) \cdot z_3 = z_1 \cdot z_3 + z_2 \cdot z_3$, and the assoziative law $(z_1 \cdot z_2) \cdot z_3 = z_1 \cdot z_2 \cdot z_3 = z_1 \cdot (z_2 \cdot z_3)$ are valid.

For any complex number $z = x + \ii y$, the complex number $\bar{z} = x - \ii y$ is called the {\em complex conjugate} of $z$\index{complex conjugate of $z$}.
One obtains $\bar{z}$ by reflection of $z$ on the real axis (see Fig.\ \ref{Abb:Komplexe_Zahl01}).
Obviously, $\bar{\bar{z}} = z$.
Moreover we have
\begin{align*}
 \overline{\mathstrut z_1+z_2}
= {} & \overline{\mathstrut x_1 + \ii y_1 + x_2 + \ii y_2}
= \overline{x_1 + x_2 + \ii (y_1 + y_2)}
= x_1 + x_2 - \ii (y_1 + y_2)\\
= {} & x_1 - \ii y_1 + x_2 - \ii y_2
= \bar{z}_1 + \bar{z}_2  
\end{align*}
and
\begin{align*}
  \overline{\mathstrut z_1\cdot z_2}
= {} & \overline{(x_1 + \ii y_1)(x_2 + \ii y_2)}
= \overline{x_1 x_2 - y_1 y_2 + \ii (x_1 y_2 +y_1 x_2)}
= x_1 x_2 - y_1 y_2 - \ii (x_1 y_2 +y_1 x_2)\\
= {} & (x_1 - \ii y_1)(x_2 - \ii y_2)
= \bar{z}_1 \cdot \bar{z}_2  
\end{align*}
and
\beqn \label{Eq:Re_z_und_Im_z}
\left.
\begin{aligned}
  z + \bar{z}
= {} & x + \ii y + x - \ii y
= 2 x
= 2 \Rez z\,,\\  
  z - \bar{z}
= {} & (x + \ii y) - (x - \ii y)
= 2 \ii y
= 2 \ii \Imz z\,.  
\end{aligned}
\right\}
\eeqn
The change from the Cartesian coordinates $x$, $y$ to the complex representation can be regarded as a linear coordinate transformation
\beqn \label{Eq:Abbildung_reell->komplex}
  z = x + \ii y\,,\quad
  \bar{z} = x - \ii y
\eeqn
\autocite[p.\ 17]{Wunderlich:Komplexe_Zahlen}.
This mapping is obviously the inverse mapping of
\beqn \label{Eq:Abbildung_komplex->reell}
  x = \frac{z + \bar{z}}{2}\,,\quad
  y = \frac{z - \bar{z}}{2\ii} 
\eeqn
according to \eqref{Eq:Re_z_und_Im_z}.
The Jacobian determinant\index{Jacobian determinant} of \eqref{Eq:Abbildung_reell->komplex} is imaginary, but the mapping \eqref{Eq:Abbildung_reell->komplex} is regular, because
\beq
  \left|\begin{array}{@{\,}cc@{\,}}
  	\frac{\p z}{\p x} & \frac{\p z}{\p y}\\[0.2cm]
  	\frac{\p\bar{z}}{\p x} & \frac{\p\bar{z}}{\p y}
  \end{array}\right|
= \left|\begin{array}{@{\,}cr@{\,}}
	1 & -\ii\\ 1 & -\ii
  \end{array}\right|
= -2\ii
\ne 0\,.    
\eeq
$z$ and $\bar{z}$ are also called {\em isotropic coordinates}\index{coordinates!isotropic coordinates@isotropic $\sim $} (not to be confused with the isotropic coordinates of spacetime) or {\em minimal coordinates}\index{coordinates!minimal coordinate@minimal $\sim $}\label{p:Minimalkoordinaten}
\autocite[p.\ 17]{Wunderlich:Komplexe_Zahlen}.
The implicit function representation $F(x,y) = 0$ can be expressed by the  mapping \eqref{Eq:Abbildung_komplex->reell} in the complex form
\beq
  \widetilde{F}(z,\bar{z})
:= F\left(\frac{z+\bar{z}}{2},\frac{z-\bar{z}}{2\ii}\right) 
= 0\,.
\eeq

The absolute value $|z|$ of a complex number $z = x + \ii y$ is the length of its vector (see Fig.\ \ref{Abb:Komplexe_Zahl01}), thus
\beq
  |z|
= \sqrt{x^2 + y^2}
= \sqrt{\mathstrut z \cdot \bar{z}}\,.  
\eeq
The main properties of the absolute value are
\begin{align*}
& |z| > 0 \quad\mbox{if}\quad z \ne 0\,,\\[0.1cm]
& |z_1 \cdot z_2| = |z_1|\,|z_2|\,,\\[0.1cm]
& |z_1 + z_2| \le |z_1| + |z_2| \qquad \mbox{(triangle inequality)}\,. 
\end{align*}  
The distance between the points $z_1$ and $z_2$ is given by
\beq
  |z_1 - z_2|
= \sqrt{(x_1-x_2)^2 + (y_1-y_2)^2}\,.  
\eeq
 
\begin{SCfigure}[0.9][ht]
  \includegraphics[width=0.25\textwidth]{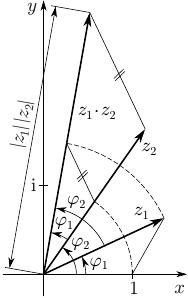}
  \caption{Rotational stretching}
  \label{Abb:Drehstreckung01}
\end{SCfigure}

A complex number $z$ can be written in polar coordinates $r$ and $\ph$ as
\beqn \label{Eq:Polarform}
  z = r(\cos\ph + \ii \sin\ph)\,,
\eeqn
where $x = r\cos\ph$ and $y = r\sin\ph$ (see Fig.\ \ref{Abb:Komplexe_Zahl01}).
$r$ is the distance to the coordinate origin (the length of the vector), as
\beq
  |z|
= \sqrt{r^2 \cos^2\ph + r^2 \sin^2\ph}
= r\,\sqrt{\cos^2\ph + \sin^2\ph}
= r   
\eeq
shows, and $\ph$ is angle with the real axis.  
$\ph$ is called {\em argument} of $z$\index{argument of $z$}.
One writes $\ph = \arg z$.
Since $\ph$ and $\ph + 2\pi k$, $k \in \mathbb{Z}$, can be assigned to the same number $z$, the polar coordinate representation is not unique.
Uniqueness in the polar coordinate representation is obtained by requiring, e.g., $-\pi < \ph \le \pi$.
Then the corresponding $\ph$ is called the {\em principal value}\index{principal value} of the argument.
Using {\em Euler's formula}\index{Euler's formula}
\beq
  \cos\ph + \ii \sin\ph
= \ee^{\ii\ph}  
\eeq
we can write \eqref{Eq:Polarform} also as
\beqn \label{Eq:Exponentialform}
  z 
= r\ee^{\ii\ph}
= |z|\,\ee^{\ii\ph}\,.
\eeqn
Especially the form \eqref{Eq:Exponentialform} offers the decisive advantages in the use of the complex numbers in plane geometry and kinematics.
For the complex conjugate also holds
\beq
  z
= r(\cos\ph - \ii\sin\ph)
= r\left[\cos(-\ph) + \ii\sin(-\ph)\right] 
= r\ee^{-\ii\ph}\,.
\eeq
For the complex multiplication one obtains by means of the form \eqref{Eq:Exponentialform}
\beqn \label{Eq:Drehstreckung}
  z_1 \cdot z_2
= r_1\ee^{\ii\ph_1} \cdot r_2\ee^{\ii\ph_2}
= r_1 r_2\,\ee^{\ii(\ph_1+\ph_2)}
= r_1 r_2\,[\cos(\ph_1+\ph_2) + \ii\sin(\ph_1+\ph_2)]\,.  
\eeqn
From this it can be seen that the multiplication of $z_1 = r_1 \ee^{\ii\ph_1}$ by $z_2 = r_2 \ee^{\ii\ph_2}$ can be interpreted as the transformation of the vector $z_1$ into the vector $z_1 \cdot z_2$, rotating it by the angle $\ph_2$ and changing its length by the factor $r_2$.
It is said that the multiplication by $z_2$ causes a {\em rotational stretching}\index{rotational stretching} of $z_1$ (see Fig.\ \ref{Abb:Drehstreckung01}), which of course can also be a ``rotational shortening''.
The vector $z_1 \cdot z_2$ can be easily constructed.
To do this, complete $z_1$ by adding the point $(1,0)$ to a triangle, rotate this triangle around the coordinate origin with angle $\ph_2$ and then scale it with center in the coordinate origin by the factor $r_2$, and from $z_1$ one has received $z_1 \cdot z_2$.
That this construction provides the correct absolut value (the correct length) of $z_1 \cdot z_2$ can be seen by
\beq
  \frac{|z_1 \cdot z_2|}{|z_2|}
= \frac{|z_1| \cdot |z_2|}{|z_2|}
= \frac{|z_1|}{1}\,.  
\eeq   

If $r_2 = |z_2| = 1$ in \eqref{Eq:Drehstreckung}, the complex multiplication causes only a rotation of $z_1$, but no change in $r_1 = |z_1|$.
Indeed, with $z_1 = r_1\ee^{\ii\ph_1}$ and $z_2 = \ee^{\ii\ph_2}$ we obtain
\beq
  \tilde{z}_1
= z_1 \cdot z_2
= r_1\ee^{\ii\ph_1}\ee^{\ii\ph_2}
= r_1\ee^{\ii(\ph_1+\ph_2)}\,,  
\eeq
where $\tilde{z}_1$ is the vector $z_1$ rotated by $\ph_2$.
For a rotation by $90\g$, it is specifically
\beq
  z_2
= \ee^{\ii\ph_2}
= \ee^{\ii\frac{\pi}{2}}
= \cos\frac{\pi}{2} + \ii\sin\frac{\pi}{2}
= \ii\,,  
\eeq
hence
\beq
  \tilde{z}_1
= \ii z_1
= \ii r \ee^{\ii\ph_1}\,.   
\eeq
For a rotation by $-90\g$ one analogously gets
\beq
  \tilde{z}_1
= -\ii z_1
= -\ii r \ee^{\ii\ph_1}\,.   
\eeq

Now, we consider the complex products $\bar{z}_1 \cdot z_2$ and $z_1 \cdot \bar{z}_2$.
For the first product we get
\beq
  \bar{z}_1 \cdot z_2
= (x_1 - \ii y_1) \cdot (x_2 + \ii y_2)
= x_1 x_2 + y_1 y_2 + \ii\left(x_1 y_2 - y_1 x_2\right). 
\eeq
We see that the real part of this product is the {\em scalar product}\index{scalar product}\index{product!scalar product@scalar $\sim $} of the vectors $z_1$ and $z_2$ and use the common notation $\langle z_1,z_2 \rangle$ for this in the following; the imaginary part is the {\em standard symplectic form on $\R^2$} \cite{Gotay&Isenberg}.\footnote{I owe the hint for this connection with symplectic geometry to Stefan Gössner, Dortmund (see also \cite{Goessner:SG-Fundamentals}).}
As this symplectic form, apart from the missing unit vector, is equal to the cross product (vector product) of two vectors $z_1,z_2\in\R^2$, we will refer to this form as the {\em quasi vector product (QVP)} of $z_1$ and $z_2$, and write $[z_1,z_2]$ for it.\footnote{For a good tutorial see the video \cite{Youtube:Dot_and_Cross_Products_of_Complex_Numbers}.}
For the second complex product we get
\beq
  z_1 \cdot \bar{z}_2
= (x_1 + \ii y_1) \cdot (x_2 - \ii y_2)
= x_1 x_2 + y_1 y_2 - \ii\left(x_1 y_2 - y_1 x_2\right), 
\eeq
and therefore
\begin{align}
  \bar{z}_1 \cdot z_2
= {} & \langle z_1,z_2\rangle + \ii\left[z_1,z_2\right],\label{Eq:conj(z_1)*z_2}\\
  z_1 \cdot \bar{z}_2
= {} & \langle z_1,z_2\rangle - \ii\left[z_1,z_2\right].\label{Eq:z_1*conj(z_2)}  
\end{align}
From \eqref{Eq:conj(z_1)*z_2} and \eqref{Eq:z_1*conj(z_2)} we immediately get the following complex expressions for the scalar product or quasi vector product
\begin{align}
  \langle z_1,z_2\rangle
= {} & \frac{\bar{z}_1 z_2 + z_1 \bar{z}_2}{2}\,,\label{Eq:Def_SP_complex}\\  
  [z_1,z_2]
= {} & \frac{\bar{z}_1 z_2 - z_1 \bar{z}_2}{2\ii}\,.\label{Eq:Def_QVP_complex}
\end{align}

{\bf Scalar product.}\index{product!scalar product@scalar $\sim $} Now, we will now look at the scalar product of two complex numbers in more detail.
This scalar product is defined by \eqref{Eq:Def_SP_complex} or, in real form, by
\beqn \label{Eq:Def_SP_real}
  \langle z_1,z_2\rangle
= x_1 x_2 + y_1 y_2\,.
\eeqn
From \eqref{Eq:Def_SP_complex} with $z_1 = r_1\ee^{\ii\ph_1}$, $z_2 = r_2\ee^{\ii\ph_2}$ we get
\begin{align} \label{Eq:Skalarprodukt_mit_Kosinus}
  \langle z_1,z_2 \rangle
= {} & \frac{r_1\ee^{-\ii\ph_1}\,r_2\ee^{\ii\ph_2} + r_1\ee^{\ii\ph_1}\,r_2\ee^{-\ii\ph_2}}{2} 
= r_1 r_2\,\frac{\ee^{\ii(\ph_2-\ph_1)} + \ee^{-\ii(\ph_2-\ph_1)}}{2}\nonumber\\
= {} & |z_1|\, |z_2| \cos(\ph_2-\ph_1)
= |z_1|\, |z_2| \cos(\ph_1-\ph_2)\,.
\end{align}
The scalar product is
\beq
\begin{array}{l@{\quad}l}
  \mbox{a) commutative:} & 
  \langle z_1,z_2 \rangle
= \langle z_2,z_1 \rangle\,,\\[0.1cm]
  \mbox{b) distributive:} &
  \langle z_1+z_2,z_3 \rangle
= \langle z_1,z_3 \rangle + \langle z_2,z_3 \rangle\,,\\[0.1cm]  
  &
  \langle z_1,z_2+z_3 \rangle
= \langle z_1,z_2 \rangle + \langle z_1,z_3 \rangle\,,\\[0.1cm]
  \mbox{c) non-associative:} &
  \langle\langle z_1,z_2 \rangle,z_3 \rangle
\ne \langle z_1,\langle z_2,z_3 \rangle\rangle\,.  
\end{array}
\eeq
The relation \eqref{Eq:Skalarprodukt_mit_Kosinus} allows to calculate the angle $\measuredangle(z_1,z_2)$ between the two complex numbers (vectors) $z_1$ and $z_2$ by means of
\beq
  \cos\measuredangle(z_1,z_2)
= \frac{\langle z_1,z_2 \rangle}{r_1r_2}
= \frac{\langle z_1,z_2 \rangle}{|z_1|\,|z_2|}
= \frac
	{\langle z_1,z_2 \rangle}
	{\sqrt{\langle z_1,z_1 \rangle}\,\sqrt{\langle z_2,z_2 \rangle}}\,.  
\eeq
(Here we used $\langle z,z \rangle = |z|^2$ from \eqref{Eq:Rechenregeln1}.)
Since the function $\cos \colon [0,\pi] \rightarrow [-1,1]$ is bijective, this uniquely defines the angle $\measuredangle(z_1,z_2)$.
We have
\beqn \label{Eq:Angle_with_inner_product}
  \measuredangle(z_1,z_2)
= \mathrm{Arccos}
	\left(\frac{\langle z_1,z_2 \rangle}{|z_1|\,|z_2|}\right),   
\eeqn
where $\mathrm{Arccos}$ denotes the principal value of $\arccos$, i.\,e.\ its restriction to the interval $[0,\pi]$.

\begin{SCfigure}[][ht]
  \includegraphics[width=0.4\textwidth]{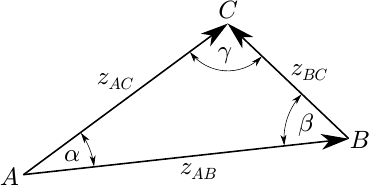}
  \caption{Triangle and law of cosines}
  \label{Abb:Law_of_cosines01}
\end{SCfigure}

Now, we consider the triangle $ABC$ in Fig.\ \ref{Abb:Law_of_cosines01}.
It holds
\beq
  z_{AC} = z_{AB} + z_{BC}\,,
  \quad\mbox{hence}\quad
  z_{BC} = z_{AC} - z_{AB}\,.
\eeq

From this we get
\begin{align*}
  z_{BC}\,\overline{z_{BC}}
= {} & {\left(z_{AC}-z_{AB}\right) \overline{\left(z_{AC}-z_{AB}\right)}}\db\\ 
= {} & {\left(z_{AC}-z_{AB}\right) \left(\overline{z_{AC}}-\overline{z_{AB}}\right)}\db\\
= {} & z_{AC}\,\overline{z_{AC}} - z_{AC}\,\overline{z_{AB}} - z_{AB}\,\overline{z_{AC}} + z_{AB}\,\overline{z_{AB}}\,,
\end{align*} 
hence
\beq
  |z_{BC}|^2
= |z_{AC}|^2 + |z_{AB}|^2 - 2 \left\langle z_{AB},z_{AC}\right\rangle, 
\eeq
and, using \eqref{Eq:Skalarprodukt_mit_Kosinus},
\beq
  |z_{BC}|^2
= |z_{AC}|^2 + |z_{AB}|^2 - 2\left|z_{AB}\right| \left|z_{AC}\right|\cos\alpha\,.
\eeq
With
\beq
  a := |z_{BC}|\,,\quad b := |z_{AC}|\,,\quad c := |z_{AB}|
\eeq
we have
\beqn \label{Eq:Law_of_cosines_a}
  a^2
= b^2 + c^2 - 2bc\cos\alpha\,.
\eeqn
Cyclic permutation gives
\beqn \label{Eq:Law_of_cosines_b}
  b^2
= a^2 + c^2 - 2ac\cos\beta\,,\quad
  c^2
= a^2 + b^2 - 2ab\cos\gamma\,.
\eeqn
\eqref{Eq:Law_of_cosines_a} together with \eqref{Eq:Law_of_cosines_b} is the {\em law of cosines}\index{cosine!law of cosines@law of $\sim $s}.

\begin{figure}[ht]
\begin{minipage}{0.48\textwidth}
  \centering
  \includegraphics[width=0.63\textwidth]{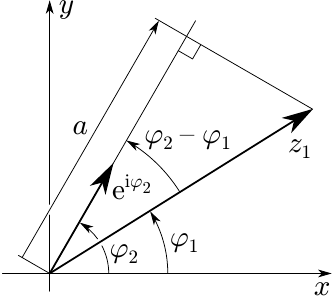}
  \caption{Scalar product\index{product!scalar product@scalar $\sim $} and projection}
  \label{Abb:Projection}
\end{minipage}
\hfill
\begin{minipage}{0.48\textwidth}
  \centering
  \includegraphics[width=0.55\textwidth]{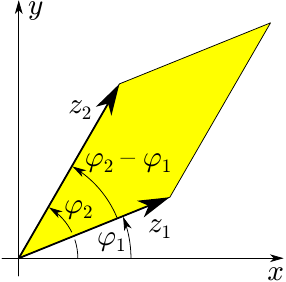}
  \caption{Geometric interpretation of the quasi vector product $[z_1,z_2]$}
  \label{Abb:Parallelogramm}
\end{minipage}
\end{figure}

We determine the length $a$ of the projection of a vector $z_1$ in the direction of a unit vector $\ee^{\ii\ph_2}$ (see Fig.\ \ref{Abb:Projection}).
If we put $z_2:=\ee^{\ii\ph_2}$ in the scalar product, then with \eqref{Eq:Skalarprodukt_mit_Kosinus} we get
\beq
  \left\langle z_1, \ee^{\ii\ph_2} \right\rangle
= \left|z_1\right| \left|\ee^{\ii\ph_2}\right| \cos(\ph_2-\ph_1)
= \left|z_1\right| \cos(\ph_2-\ph_1)\,.
\eeq
Because of
\beq
  \cos(\ph_2-\ph_1)
= \frac{a}{|z_1|}
\eeq   
we have
\beqn \label{Eq:Projection_with_inner_product}
  a
= |z_1| \cos(\ph_2-\ph_1)
= \left\langle z_1, \ee^{\ii\ph_2} \right\rangle\,.
\eeqn
For the projection as vector we get
\beq
  a\,\ee^{\ii\ph_2}
= \left\langle z_1, \ee^{\ii\ph_2} \right\rangle \ee^{\ii\ph_2}\,.  
\eeq

If $z_1(t)$ and $z_2(t)$ are differentiable functions of the real parameter $t$, then
\begin{align}
  \frac{\dd}{\dd t} \left\langle z_1(t), z_2(t) \right\rangle
= {} & \frac{\dd}{\dd t} \left(x_1(t)\,x_2(t) + y_1(t)\,y_2(t)\right)\nonumber\db\\[0.05cm] 
= {} & x_1'x_2 + x_1x_2' + y_1'y_2 + y_1y_2'\nonumber\db\\[0.05cm] 
= {} & x_1'x_2 + y_1'y_2 + x_1x_2' + y_1y_2'\nonumber\\[0.05cm]
= {} & \langle z_1'(t), z_2(t) \rangle + \langle z_1(t), z_2'(t) \rangle\,. \label{Eq:Product_rule_SP}
\end{align}
This is the {\em product rule}\index{product rule} for the scalar product\index{product!scalar product@scalar $\sim $}. 
  
{\bf Quasi vector product\index{product!quasi vector product@quasi vector $\sim $}.} The {\em quasi vector product} of two complex numbers always yields a real number and is defined by \eqref{Eq:Def_QVP_complex} or, in real form, by 
\beqn \label{Eq:Def_QVP_real}
  [z_1,z_2]
= x_1 y_2 - y_1 x_2
= \left|\begin{array}{@{\:}cc@{\:}}
	x_1 & y_1 \\ x_2 & y_2
  \end{array}\right|.
\eeqn
From this we get
\begin{align} \label{Eq:QVP_with_sine}
  [z_1,z_2]
= {} & \frac{r_1\ee^{-\ii\ph_1}\,r_2\ee^{\ii\ph_2} - r_1\ee^{\ii\ph_1}\,r_2\ee^{-\ii\ph_2}}{2\ii}
= r_1r_2 \frac{\ee^{\ii(\ph_2-\ph_1)} - \ee^{-\ii(\ph_2-\ph_1)}}{2\ii}\nonumber\\[0.05cm]
= {} & |z_1|\, |z_2| \sin(\ph_2-\ph_1)\,.       
\end{align}
Thus, the real number $[z_1,z_2]$ is equal to the (signed) area of the parallelogram spanned by the vectors $z_1$ and $z_2$ (see Fig.\ \ref{Abb:Parallelogramm}).

It holds
\begin{align*}
  \Imz(\bar{z}_1\,z_2)
= {} & \Imz((x_1-\ii y_1)(x_2+\ii y_2))
= \Imz(x_1 x_2 + y_1 y_2 + \ii (x_1 y_2 - y_1 x_2))
= x_1 y_2 - y_1 x_2\,,\\[0.1cm]  
  \Imz(z_1\,\bar{z}_2)
= {} & \Imz((x_1+\ii y_1)(x_2-\ii y_2))
= \Imz(x_1 x_2 + y_1 y_2 - \ii (x_1 y_2 - y_1 x_2))
= -(x_1 y_2 - y_1 x_2)\,,
\end{align*}
hence
\beqn \label{Eq:[z_1,z_2]_mit_Im}
  [z_1,z_2]
= \Imz(\bar{z}_1\,z_2)\,,\qquad
  [z_1,z_2]
= -\Imz(z_1\,\bar{z}_2)\,.
\eeqn

The quasi vector product is
\beqn \label{Eq:QVP_rules}
\begin{array}{l@{\quad}l}
  \mbox{a) alternating:} &
  [z_1,z_2]
= -[z_2,z_1]\,,\\[0.1cm]
  \mbox{b) distributive:} & 
  [z_1+z_2,z_3]
= [z_1,z_3] + [z_2,z_3]\,,\\[0.1cm]
  &
  [z_1,z_2+z_3]
= [z_1,z_2] + [z_1,z_3]\,,\\[0.1cm]
  \mbox{c) non-associative:} &
  [[z_1,z_2],z_3]
\ne [z_1,[z_2,z_3]]\,.
\end{array}
\eeqn

If $z_1(t)$ and $z_2(t)$ are differentiable functions of the real parameter $t$, then
\begin{align}
  \frac{\dd}{\dd t} \left[z_1(t), z_2(t)\right]
= {} & \frac{\dd}{\dd t} \left(x_1(t)\,y_2(t) - y_1(t)\,x_2(t)\right)\nonumber\db\\[0.05cm] 
= {} & x_1'y_2 + x_1y_2' - (y_1'x_2 + y_1x_2')\nonumber\db\\[0.05cm] 
= {} & x_1'y_2 - y_1'x_2 + (x_1y_2' - y_1x_2')\nonumber\\[0.05cm]
= {} & [z_1'(t), z_2(t)] + [z_1(t), z_2'(t)]\,. \label{Eq:Product_rule_QVP}
\end{align}
This is the {\em product rule}\index{product rule} for the quasi vector product.

{\bf More calculation rules for the scalar and quasi vector product.}
For the scalar and the quasi vector product also the following, easily provable calculation rules \autocite[p.~414]{Luck&Modler} apply, which can be used beneficially in concrete calculations:
\beqn \label{Eq:Rechenregeln1}
\left.
\begin{gathered}
  \langle 1,z \rangle = \Rez z\,,\qquad
  [1,z] = \Imz z\,,\qquad
  \langle z,z \rangle = z \bar{z} = |z|^2\,,\\[0.1cm]
  \langle \ii z_1,z_2 \rangle = [z_1,z_2]\,,\quad
  \langle z_1,\ii z_2 \rangle = -[z_1,z_2]\,,\quad
  [\ii z_1,z_2] = -\langle z_1,z_2\rangle\,,\quad
  [z_1,\ii z_2] = \langle z_1,z_2 \rangle\,,     
\end{gathered}
\;\right\}
\eeqn
\beqn \label{Eq:Rechenregeln2}
\left.
\renewcommand{\arraystretch}{1.1}
\begin{array}{|l|l|} \hline
  \langle z_1,xz_2 \rangle = \langle xz_1,z_2 \rangle
= x \langle z_1,z_2 \rangle\,,\; x \in \R\,, & 
  [z_1,xz_2] = [xz_1,z_2] = x [z_1,z_2]\,,\; x \in \R\,,\\[0.1cm]
  \langle z_1,z_2z_3 \rangle = \langle z_1\bar{z}_2,z_3 \rangle\,, &
  [z_1,z_2z_3] = [z_1\bar{z}_2,z_3]\,,\\[0.1cm] 
  \langle z_1z_3,z_2z_3 \rangle = \langle z_1z_3\bar{z}_3,z_2 \rangle
= z_3\bar{z}_3 \langle z_1,z_2 \rangle\,, & 
  [z_1z_3,z_2z_3] = [z_1z_3\bar{z}_3,z_2] 
= z_3\bar{z}_3 [z_1,z_2]\,,\\[0.1cm]
  \langle z_1\ee^{\ii\ph},z_2\ee^{\ii\ph} \rangle
= \langle z_1,z_2 \rangle \quad \mbox{(rotation rule)}\,, &
  [z_1\ee^{\ii\ph},z_2\ee^{\ii\ph}]
= [z_1,z_2] \quad \mbox{(rotation rule)}\,,\\[0.1cm]
  \langle z_1,\ii z_2 \rangle = -\langle \ii z_1,z_2 \rangle\,, &
  [z_1,\ii z_2] = -[\ii z_1,z_2]\,.\\ \hline          
\end{array}
\renewcommand{\arraystretch}{1}
\;\right\}
\eeqn
We show as an example that $\langle z_1z_3,z_2z_3 \rangle = z_3 \bar{z}_3 \langle z_1,z_2 \rangle$:
\begin{align*}
  \langle z_1z_3,z_2z_3 \rangle
= {} & \frac{1}{2}\left(z_1z_3\,\overline{z_2z_3}
+ \overline{z_1z_3}\,z_2z_3\right)
= \frac{1}{2}\left(z_1z_3\bar{z}_2\bar{z}_3
+ \bar{z}_1\bar{z}_3z_2z_3\right)\\
= {} & z_3\bar{z}_3\,\frac{1}{2}\left(z_1\bar{z}_2+\bar{z}_1z_2\right)
= z_3\bar{z}_3 \langle z_1,z_2 \rangle\,.  
\end{align*}
With $z_3 = \ee^{\ii\ph}$ follows also directly the rotation rule for the scalar product
\beq
  \langle z_1\ee^{\ii\ph},z_2\ee^{\ii\ph} \rangle
= \ee^{\ii\ph}\ee^{-\ii\ph} \langle{z_1,z_2} \rangle
= \langle z_1,z_2 \rangle\,.      
\eeq
    
{\bf Straight lines.}  A straight line $g$ can be represented by a fixed point $z_\tn{P}$ and a direction vector~$z_\tn{t}$ using the parameter representation
\beq
  z
= z_\tn{P} + \la z_\tn{t}\,,\qquad
  -\infty < \la < \infty\,,  
\eeq
or by
\beq
  [z-z_\tn{P},z_\tn{t}] = 0\,.
\eeq
The last equation can be written as
\beq
  [z,z_\tn{t}] = [z_\tn{P},z_\tn{t}] = c_1\,,
\eeq
where $c_1$ is a real constant.
We determine the intersection of the two straight lines defined by the equations
\beq
  [z,z_\tn{t1}] = c_1\,,\qquad
  [z,z_\tn{t2}] = c_2\,,
\eeq
respectively.
According to \eqref{Eq:Def_QVP_complex} we have
\begin{align*}
  [z,z_\tn{t1}] 
= {} & c_1 = \frac{\ii}{2}\,(z \bar{z}_\tn{t1} - \bar{z} z_\tn{t1})\,,\\[0.1cm]
  [z,z_\tn{t2}] 
= {} & c_2 = \frac{\ii}{2}\,(z \bar{z}_\tn{t2} - \bar{z} z_\tn{t2})\,,
\end{align*}
hence
\beqn \label{Eq:Lineares_Gleichungssystem_komplex}
\left.
\begin{aligned}
  z \bar{z}_\tn{t1} - \bar{z} z_\tn{t1}
= {} & {-}2 \ii c_1\,,\\[0.05cm]
  z \bar{z}_\tn{t2} - \bar{z} z_\tn{t2} 
= {} & {-}2 \ii c_2\,.
\end{aligned}
\:\right\}
\eeqn
The two equations \eqref{Eq:Lineares_Gleichungssystem_komplex} can be understood as a linear equation system for the determination of the two unknowns $z$ and $\bar{z}$.
Using Cramer's rule\index{Cramer's rule} we obtain
\begin{align*}
  z
= {} & \frac{
  \left|\begin{array}{@{\,}cc@{\,}}
	-2 \ii c_1 & -z_\tn{t1}\\
	-2 \ii c_2 & -z_\tn{t2} 	
  \end{array}\right|}
  {
  \left|\begin{array}{@{\,}cc@{\,}}
  	\bar{z}_\tn{t1} & -z_\tn{t1}\\
  	\bar{z}_\tn{t2} & -z_\tn{t2}
  \end{array}\right|}  
= -\frac{
  \left|\begin{array}{@{\,}cc@{\,}}
	2 \ii c_1 & z_\tn{t1}\\
	2 \ii c_2 & z_\tn{t2} 	
  \end{array}\right|}
  {
  \left|\begin{array}{@{\,}cc@{\,}}
  	\bar{z}_\tn{t1} & z_\tn{t1}\\
  	\bar{z}_\tn{t2} & z_\tn{t2}
  \end{array}\right|}
= -\frac
	{2\ii(c_1 z_\tn{t2}-c_2 z_\tn{t1})}
	{\bar{z}_{t1} z_\tn{t2} - \bar{z}_{t2} z_\tn{t1}}
= -\frac
	{\ii(c_2 z_\tn{t1}-c_1 z_\tn{t2})}
	{\frac{1}{2}\,(z_{t1} \bar{z}_\tn{t2} - \bar{z}_{t1} z_\tn{t2})}\\[0.1cm]
= {} & \frac
	{c_2 z_\tn{t1}-c_1 z_\tn{t2}}
	{\frac{\ii}{2}\,(z_{t1} \bar{z}_\tn{t2} - \bar{z}_{t1} z_\tn{t2})}
= \frac
	{c_2 z_\tn{t1}-c_1 z_\tn{t2}}
	{[z_\tn{t1},z_\tn{t2}]}\,.	  
\end{align*}
Three pairwise distinct points $z_1$, $z_2$, $z_3$ lie on a straight line if
\begin{gather*}
  [z_1-z_2,z_1-z_3] = 0
  \quad\Longleftrightarrow\quad
  [z_1,z_1-z_3] - [z_2,z_1-z_3] = 0
  \quad\Longleftrightarrow\\[0.1cm]
  \underbrace{[z_1,z_1]}_0 - [z_1,z_3] - ([z_2,z_1] - [z_2,z_3]) = 0
  \quad\Longleftrightarrow\quad
  [z_1,z_2] + [z_2,z_3] + [z_3,z_1] = 0\,.
\end{gather*}
If we put for $z_3$ any point $z$ of the straight line $g$, we get
\beq
  [z_1,z_2] + [z_2,z] + [z,z_1] = 0
  \qquad\mbox{or}\qquad
  [z_1-z_2,z_1-z] = 0
\eeq
or
\beq
  [z_1-z_2,z] = [z_1,z_2]\,.
\eeq
Using \eqref{Eq:Projection_with_inner_product}, we can also determine an equation for a straight line $g$.
Writing $z$ instead of $z_1$, and $\ph$ instead of $\ph_2$, we have
\beqn \label{Eq:Line_equation_with_projection}
  \left\langle z, \ee^{\ii\ph} \right\rangle = a\,.
\eeqn

\begin{figure}[ht]
\begin{minipage}[b]{0.5\textwidth}
  \centering
  \includegraphics[width=0.55\textwidth]{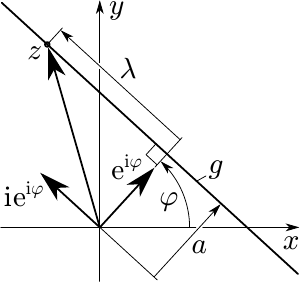}
  \caption{Straight line $g$}
  \label{Abb:Hessesche_Normalform}
\end{minipage}
\hfill
\begin{minipage}[b]{0.5\textwidth}
  \centering
  \includegraphics[width=0.55\textwidth]{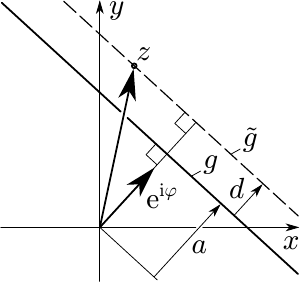}
  \caption{Distance $d$ of point $z$ from line $g$}
  \label{Abb:Hesse_Nf_Distance}
\end{minipage}
\end{figure}

The projection of the position vector $z$ onto the direction of the unit vector $\ee^{\ii\ph}$ is the same for every point $z \in g$ (see Fig.\ \ref{Abb:Hessesche_Normalform}).

It follows that \eqref{Eq:Line_equation_with_projection} is the equation of the line $g$ with direction perpendicular to $\ee^{\ii\ph}$ and distance~$a$ from the coordinate origin.
We can write \eqref{Eq:Line_equation_with_projection} as
\beqn \label{Eq:Hesse_normal_form}
  x\cos\ph + y \sin\ph = a\,,
\eeqn
which is the {\em \textsc{Hesse}\footnote{\textsc{Ludwig Otto Hesse}, 1811-1874} normal form}\index{Hesse normal form}.

If $z \notin g$, then $z$ lies on a line $\tilde{g}$ parallel to $g$ with distance $d$ (see Fig.\ \ref{Abb:Hesse_Nf_Distance}) and Hesse normal form
\beq
  \left\langle z, \ee^{\ii\ph} \right\rangle = a + d\,.
\eeq
Therefore
\beqn \label{Eq:Hesse_Nf_Distance} 
  d = \left\langle z, \ee^{\ii\ph} \right\rangle - a
\eeqn
is the distance of point $z$ from line $g$\index{distance!distance of a point from a line@$\sim $ of a point from a line}. 
The distance $d$ according to \eqref{Eq:Hesse_Nf_Distance} is positive if the coordinate origin and $z$ are on different sides of $g$.
It is negative if the coordinate origin and $z$ are on the same side of $g$.

We can also describe the line $g$ in Fig.\ \ref{Abb:Hessesche_Normalform} with the equation
\beqn \label{Eq:z(phi,a)}
  z(\ph,a)
= a\ee^{\ii\ph} + \la\ii\ee^{\ii\ph} 
= \left(a+\ii\la\right)\ee^{\ii\ph}\,,\quad
  -\infty < \la < +\infty\,.  
\eeqn
Splitting \eqref{Eq:z(phi,a)} into real and imaginary part gives
\beqn \label{Eq:x(phi,a)_and_y(phi,a)}
\begin{array}{c@{\;=\;}*{2}{c@{\:}}c}
  x(\ph,a)
& a\cos\ph & - & \la\sin\ph\,,\\[0.05cm]   
  y(\ph,a)
& a\sin\ph & + & \la\cos\ph\,.
\end{array}
\eeqn
Eliminating $\la$ from \eqref{Eq:x(phi,a)_and_y(phi,a)} yields \eqref{Eq:Hesse_normal_form}, as can be easily checked. 
On the other hand, putting \eqref{Eq:z(phi,a)} into $\langle z,\ee^{\ii\ph}\rangle$ and using the rotation rule from \eqref{Eq:Rechenregeln2} gives
\beq
  \left\langle \left(a+\ii\la\right)\ee^{\ii\ph},\,\ee^{\ii\ph} \right\rangle
= \left\langle a+\ii\la,\,1 \right\rangle
= a\,.  
\eeq
Thus, for any $\la$, $z(\ph,a)$ from \eqref{Eq:z(phi,a)} satisfies \eqref{Eq:Line_equation_with_projection} identically.  

Note that \eqref{Eq:Line_equation_with_projection}, \eqref{Eq:Hesse_normal_form}, \eqref{Eq:z(phi,a)}, and \eqref{Eq:x(phi,a)_and_y(phi,a)} even hold for $\ph = 0$ and $\ph = \pi$.

{\bf Circles.} A circle is the set of all points $z$ with equal distance $r$ from a fixed point (center point) $z_0$.
So we have
\beq
  |z-z_0|
= r\,.  
\eeq
Squaring gives
\beqn \label{Eq:squared_circular_equation}
  r^2
= |z-z_0|^2
= \left(z-z_0\right) \overline{\left(z-z_0\right)}  
= \left(z-z_0\right) \left(\bar{z}-\bar{z}_0\right)
= z\bar{z} - \bar{z}_0 z - z_0 \bar{z} + |z_0|^2\,.
\eeqn 
A frequently used parametric equation (with parameter $\ph$) is given by
\beq
  z
= z(\ph)
= z_0 + r \ee^{\ii\ph}\,,\qquad 0 \le \ph < 2\pi\,.
\eeq
We consider the scalar product
\beqn \label{Eq:Thaleskreis}
  \langle z-z_1,z-z_2 \rangle = 0\,,
\eeqn
where $z$ is a variable and $z_1$ and $z_2$ are constants.
Because of \eqref{Eq:Skalarprodukt_mit_Kosinus} the cosine of the angle between the vectors $z-z_1$ and $z-z_2$ is zero.
Consequently, the angle between these vectors is equal to $\pi/2$ or $-\pi/2$; the vectors are orthogonal.
All points $z$ for which this is true lie on the Thales circle\index{Thales circle} over the diameter line whose endpoints are $z_1$ and $z_2$.
\eqref{Eq:Thaleskreis} is thus a circle equation to two diametrical points $z_1$ and $z_2$.

\begin{SCfigure}[0.8][ht]
  \includegraphics[width=0.3\textwidth]{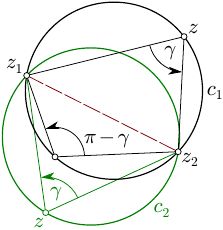}
  \caption{Circles with common chord $z_1 z_2$}
  \label{Abb:Cyclic_quadrilateral}
\end{SCfigure}

Now, we consider Fig.\ \ref{Abb:Cyclic_quadrilateral}.
In the arc $z_1 z z_2$ of a circle $c_1$, the angles $\gamma = \measuredangle\left(\vv{z z_1},\vv{z z_2}\right)$ over a chord $z_1 z_2$ are constant.

Therefore
\beq
  \left[\left(z_1-z\right)\ee^{\ii\gamma},z_2-z\right] = 0
  \quad\mbox{or, equivalently,}\quad
  \left[\left(z-z_1\right)\ee^{\ii\gamma},z-z_2\right] = 0 
\eeq
applies.
Since opposite angles in a cyclic quadrilateral\index{quadrilateral!cyclic quadrilateral@cyclic $\sim $} are supplementary, for the second arc $z_1 z z_2$ of $c_1$ (in Fig.\ \ref{Abb:Cyclic_quadrilateral} ``below'' $\overline{z_1 z_2}$) it holds
\beq
  \big[(z_2-z)\,\ee^{\ii(\pi-\gamma)},z_1-z\big] = 0
  \quad\mbox{or, equivalently,}\quad
  \big[{-(z-z_2)}\,\ee^{-\ii\gamma},z-z_1\big] = 0\,,
\eeq
and therefore
\beq
  \left[z-z_1,\left(z-z_2\right)\ee^{-\ii\gamma}\right] = 0\,,
\eeq
thus
\beqn \label{Eq:Equation_of_c_1}
  \left[\left(z-z_1\right)\ee^{\ii\gamma},z-z_2\right] = 0\,.
\eeqn
\eqref{Eq:Equation_of_c_1} is therefore the equation of the complete circle $c_1$.
For the circle $c_2$ one analogously finds
\beqn \label{Eq:Equation_of_c_2}
  \left[z-z_1,\left(z-z_2\right)\ee^{\ii\gamma}\right] = 0\,.
\eeqn
For $\gamma = \pi/2$, \eqref{Eq:Equation_of_c_1} becomes
\begin{align*}
  0
= {} & \bigl[\left(z-z_1\right)\ee^{\ii\pi/2},z-z_2\bigr]  
= \bigl[\ii\left(z-z_1\right),z-z_2\bigr]
= -\bigl[z-z_2,\ii\left(z-z_1\right)\bigr]\\
= {} & {-\langle z-z_2,z-z_1\rangle}
= -\langle z-z_1,z-z_2\rangle\,,
\end{align*}
hence we have \eqref{Eq:Thaleskreis}.
For $\gamma = \pi/2$, from \eqref{Eq:Equation_of_c_2} we again get \eqref{Eq:Thaleskreis}.
The circles $c_1$ and $c_2$ are identical if $\gamma = \pi/2$. 


The following vector identities can be useful (cf.\ \cite[p.\ 5, Tab.~2]{Goessner:SG-Fundamentals}):
\begin{gather}
  z_1[z_2,z_3]
= \ii z_2\langle z_3,z_1\rangle - \ii z_3\langle z_1,z_2\rangle,\label{Eq:Grassmann-Goessner}\db\\[0.05cm]
 [z_1,[z_2,z_3]]
= \langle z_3,\langle z_1,z_2\rangle\rangle - \langle z_2,\langle z_3,z_1\rangle\rangle,\label{Eq:Grassmann}\db\\[0.05cm]
  z_1[z_2,z_3] + z_2[z_3,z_1] + z_3[z_1,z_2]
= 0\,,\label{Eq:Jacobi-Goessner}\db\\[0.05cm]
  [z_1,[z_2,z_3] + [z_2,[z_3,z_1]] + [z_3,[z_1,z_2]]
= 0\,,\label{Eq:Jacobi}\db\\[0.05cm]
  [z_1,z_2]^2 + \langle z_1,z_2\rangle^2
= |z_1|^2\,|z_2|^2\,,\label{Eq:Lagrange}\db\\[0.05cm] 
  \left[z_1,z_2\right]\left[z_3,z_4\right]
= \langle z_1,z_3\rangle \langle z_2,z_4\rangle - \langle z_1,z_4\rangle \langle z_2,z_3\rangle.\label{Eq:Binet-Cauchy} 
\end{gather}
Proof of \eqref{Eq:Grassmann-Goessner}:
It holds
\begin{align*}
  z_1\left[z_2,z_3\right]
= {} & \frac{z_1\left(\bar{z}_2 z_3-z_2 \bar{z}_3\right)}{2\ii}
= \frac{z_1 \bar{z}_2 z_3 - z_1 z_2 \bar{z}_3 + \bar{z}_1 z_2 z_3 - \bar{z}_1 z_2 z_3}{2\ii}\db\\[0.05cm]
= {} & \frac{-z_2\left(\bar{z}_1 z_3 + z_1 \bar{z}_3\right) + z_3\left(\bar{z}_1 z_2 + z_1 \bar{z}_2\right)}{2\ii}\db\\[0.05cm] 
= {} & \frac{\ii z_2\left(\bar{z}_1 z_3 + z_1 \bar{z}_3\right) - \ii z_3\left(\bar{z}_1 z_2 + z_1 \bar{z}_2\right)}{2}\db\\[0.05cm]
= {} &  \ii z_2\left\langle z_1,z_3\right\rangle - \ii z_3\left\langle z_1,z_2\right\rangle\\[0.05cm]
= {} &  \ii z_2\left\langle z_3,z_1\right\rangle - \ii z_3\left\langle z_1,z_2\right\rangle.
\end{align*}
Proof of \eqref{Eq:Grassmann}:
\begin{align*}
  [z_1,[z_2,z_3]]
= {} & [x_1 +\ii y_1,x_2 y_3 - y_2 x_3]\db\\
= {} & {-}y_1 (x_2 y_3 - y_2 x_3)\db\\
= {} & y_1 y_2 x_3 - y_1 y_3 x_2 + x_1 x_2 x_3 - x_1 x_2 x_3\db\\
= {} & (x_1 x_2 + y_1 y_2)x_3 - (x_1 x_3 + y_1 y_3)x_2\db\\
= {} & \langle z_1,z_2\rangle x_3 - \langle z_1,z_3\rangle x_2\db\\
= {} & \langle\langle z_1,z_2\rangle, x_3 + \ii y_3\rangle - \langle\langle z_1,z_3\rangle, x_2 + \ii y_2\rangle\db\\
= {} & \langle z_3,\langle z_1,z_2\rangle\rangle - \langle z_2,\langle z_3,z_1\rangle\rangle.
\end{align*}
Proof of \eqref{Eq:Jacobi-Goessner}:
Applying \eqref{Eq:Grassmann-Goessner}, we have
\begin{align*}
& \hspace{-0.5cm} z_1\left[z_2,z_3\right] + z_2\left[z_3,z_1\right] + z_3\left[z_1,z_2\right]\db\\[0.05cm]
= {} & \ii z_2\left\langle z_3,z_1\right\rangle - \ii z_3\left\langle z_1,z_2\right\rangle
	 + \ii z_3\left\langle z_1,z_2\right\rangle - \ii z_1\left\langle z_2,z_3\right\rangle
	 + \ii z_1\left\langle z_2,z_3\right\rangle - \ii z_2\left\langle z_3,z_1\right\rangle\\[0.05cm]
= {} & 0\,. 
\end{align*}
Proof of \eqref{Eq:Jacobi}: It follow directly from \eqref{Eq:Grassmann}.

Proof of \eqref{Eq:Lagrange}: It holds 
\begin{align*}
  [z_1,z_2]^2 + \langle z_1,z_2\rangle^2
= {} & \langle \ii z_1,z_2\rangle^2 + \langle z_1,z_2\rangle^2
= \left(\frac{-\ii\bar{z}_1 z_2 + \ii z_1\bar{z}_2}{2}\right)^2 + \left(\frac{\bar{z}_1 z_2 + z_1\bar{z}_2}{2}\right)^2\db\\[0.05cm]
= {} & \frac{1}{4} \left(-\bar{z}_1^2 z_2^2 + 2\bar{z}_1 z_2 z_1 \bar{z}_2 - z_1^2 \bar{z}_2^2
						+ \bar{z}_1^2 z_2^2 + 2\bar{z}_1 z_2 z_1 \bar{z}_2 + z_1^2 \bar{z}_2^2\right)\\[0.05cm]
= {} & z_1 \bar{z}_1 z_2 \bar{z_2}					 
= |z_1|^2 |z_2|^2\,.
\end{align*}
Identity \eqref{Eq:Binet-Cauchy} can also be easily proved by straightforward calculation.

\begin{example} 
We derive a formula for the circumcenter $z_C$ of a triangle\index{triangle!circumcenter of a triangle@circumcenter of a $\sim $} (cf.\ \cite[p.\ 2]{Goessner:SG-Triangle}).
The circumcenter is the intersection point between the three perpendicular bisectors of the triangle sides.
Let the triangle be determined by the two vectors $z_1$ and $z_2$ (see Fig.\ \ref{Abb:Circumcircle}).  
\begin{figure}[ht]
\begin{minipage}[b]{0.48\textwidth}
  \centering
  \includegraphics[width=0.75\textwidth]{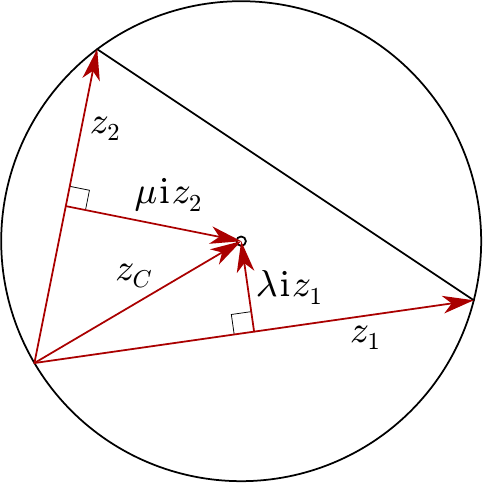}
  \caption{Circumcircle of a triangle\index{triangle!circumcircle of a triangle@circumcircle of a $\sim $}}
  \label{Abb:Circumcircle}
\end{minipage}
\hfill
\begin{minipage}[b]{0.48\textwidth}
  \centering
  \includegraphics[width=0.9\textwidth]{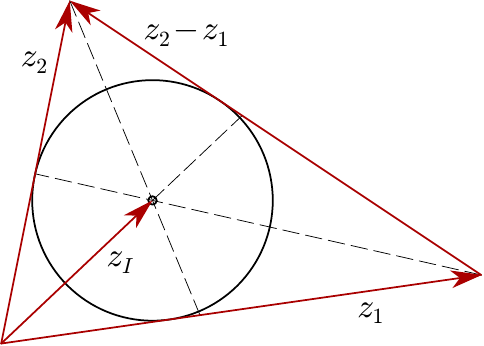}
  \caption{Incircle\index{Incircle} of a triangle\index{triangle!incircle of a triangle@incircle of a $\sim $}}
  \label{Abb:Incircle}
\end{minipage}
\end{figure}
Then we have
\beq
  z_C = \frac{1}{2}\,z_1 + \la \ii z_1 = \frac{1}{2}\,z_2 + \mu \ii z_2\,.
\eeq
Scalar multiplication by $z_2$ gives
\beq
  \left\langle\tfrac{1}{2}\,z_1 + \la \ii z_1,z_2\right\rangle
= \left\langle\tfrac{1}{2}\,z_2 + \mu \ii z_2,z_2\right\rangle,
\eeq
hence
\beq
  \tfrac{1}{2} \left\langle z_1,z_2 \right\rangle + \la \left\langle \ii z_1,z_2 \right\rangle
= \tfrac{1}{2} \left\langle z_2,z_2 \right\rangle + \mu \left\langle \ii z_2,z_2 \right\rangle  
\eeq
and therefore
\beq
  \langle z_1,z_2 \rangle + 2\la \left[z_1,z_2\right]
= |z_2|^2\,.   
\eeq
So we get
\beq
  \la
= \frac{|z_2|^2 - \langle z_1,z_2 \rangle}{2\left[z_1,z_2\right]}\,,
\eeq
and this gives
\beq
  z_C
= \frac{1}{2}\,z_1 + \ii\,\frac{|z_2|^2 - \langle z_1,z_2 \rangle}{2\left[z_1,z_2\right]}\,z_1
= \frac{z_1\left[z_1,z_2\right] + \ii z_1\,|z_2|^2 - \ii z_1\left\langle z_1,z_2 \right\rangle}{2\left[z_1,z_2\right]}\,.  
\eeq
From \eqref{Eq:Grassmann-Goessner} we get
\beq
  z_1\left[z_1,z_2\right]
= \ii z_1\left\langle z_1,z_2\right\rangle - \ii z_2\left|z_1\right|^2   
\eeq
and therefore
\begin{align*}
  z_C
= {} & \frac{\ii z_1\left\langle z_1,z_2\right\rangle - \ii z_2\left|z_1\right|^2 + \ii z_1\,|z_2|^2 - \ii z_1\left\langle z_1,z_2 \right\rangle}
			{2\left[z_1,z_2\right]}\\[0.05cm]
= {} & \ii\,\frac{z_1\,|z_2|^2 - z_2\,|z_1|^2}{2\left[z_1,z_2\right]}\,.
\end{align*}
Clearly, the radius $r_C$ of the circumcircle is given by
\beq
  r_C
= |z_C|
= \frac{1}{2} \left|\frac{z_1\,|z_2|^2 - z_2\,|z_1|^2}{\left[z_1,z_2\right]}\right|. \tag*{\bs}
\eeq
\end{example}

\begin{example} 
We determine the incircle of a triangle\index{triangle!incircle of a triangle@incircle of a $\sim $} formed by two vectors $z_1$ and $z_2$ (see Fig.\ \ref{Abb:Incircle}).
The incenter ($=$ center of the incircle)\index{triangle!incenter of a triangle@incenter of a $\sim $} $z_I$ is the intersection point of the angle bisectors,
\beq
  z_I
= \la \left(\frac{z_1}{|z_1|} + \frac{z_2}{|z_2|}\right)
= z_1 + \mu \left(\frac{z_2-z_1}{|z_2-z_1|} - \frac{z_1}{|z_1|}\right),   
\eeq
thus
\beq
  \la \left(\frac{z_1}{|z_1|} + \frac{z_2}{|z_2|}\right) - z_1
= \mu \left(\frac{z_2-z_1}{|z_2-z_1|} - \frac{z_1}{|z_1|}\right).  
\eeq
Quasi vector multiplication with the bracket expression on the right-hand side results in
\beq
  \left[\la \left(\frac{z_1}{|z_1|} + \frac{z_2}{|z_2|}\right) - z_1,\, \frac{z_2-z_1}{|z_2-z_1|} - \frac{z_1}{|z_1|}\right]
= 0\,,
\eeq
hence
\beq
  \la \left[\frac{z_1}{|z_1|} + \frac{z_2}{|z_2|},\, \frac{z_2-z_1}{|z_2-z_1|} - \frac{z_1}{|z_1|}\right]
= \left[z_1,\, \frac{z_2-z_1}{|z_2-z_1|} - \frac{z_1}{|z_1|}\right].  
\eeq
It follows that
\beq
  \la \left(\frac{[z_1,z_2]}{|z_1|\,|z_2-z_1|} + \frac{[z_2,-z_1]}{|z_2|\,|z_2-z_1|} + \frac{[z_2,-z_1]}{|z_2|\,|z_1|}\right)
= \frac{[z_1,z_2]}{|z_2-z_1|}\,,  
\eeq
thus
\beq
  \la \left(\frac{1}{|z_1|\,|z_2-z_1|} + \frac{1}{|z_2|\,|z_2-z_1|} + \frac{1}{|z_2|\,|z_1|}\right)
= \frac{1}{|z_2-z_1|}\,.  
\eeq
With $|z_3| = |z_2 - z_1|$ we get
\beq
  \la \left(\frac{1}{|z_1|\,|z_3|} + \frac{1}{|z_2|\,|z_3|} + \frac{1}{|z_1|\,|z_2|}\right)
= \la\, \frac{|z_1| + |z_2| + |z_3|}{|z_1|\,|z_2|\,|z_3|}
= \frac{1}{|z_3|}\,,  
\eeq
and consequently
\beq
  \la
= \frac{|z_1|\,|z_2|}{|z_1| + |z_2| + |z_3|}\,.  
\eeq
This gives
\beq
   z_I
= \frac{|z_1|\,|z_2|}{|z_1| + |z_2| + |z_3|} \left(\frac{z_1}{|z_1|} + \frac{z_2}{|z_2|}\right)
= \frac{z_1\,|z_2| + z_2\,|z_1|}{|z_1| + |z_2| + |z_3|}\,.
\eeq
The radius $r_I$ of the incircle is the distance of point $z_I$ from the triangle side given by $z_1$.
From \eqref{Eq:Hesse_Nf_Distance} (see also Fig.\ \ref{Abb:Hesse_Nf_Distance}) follows
\beq
  r_I
= \left|\left\langle z_I,\,\ii\,\frac{z_1}{|z_1|}\right\rangle\right|  
= \left|\left[z_I,\,\frac{z_1}{|z_1|}\right]\right|
= \left|\left[\frac{z_1\,|z_2| + z_2\,|z_1|}{|z_1| + |z_2| + |z_3|},\,\frac{z_1}{|z_1|}\right]\right|
= \frac{|[z_1,z_2]|}{|z_1| + |z_2| + |z_3|}\,. \tag*{\bs}
\eeq
\end{example}

%% file: DiffGeo5_3a.tex

\section{Plane curves} \label{Sec:Plane_curves}
\subsection{Parametric curves}

\begin{defin}
A {\em (plane) parametric $C^r$-curve}\index{curve!parametric curve@parametric $\sim $} is an $r$-times continuously differentiable complex-valued function 
\beq
  z \colon I \to \CC\,,\quad
  t \mapsto z(t) = x(t) + \ii y(t)
\eeq
defined on an open interval $I\subset\R$, where $x(t)$ and $y(t)$ are respectively the real and imaginary part of $z(t)$.
The variable $t$ is called {\em parameter} of the curve.
\end{defin}

$z$ is a $C^r$-curve if its real and imaginary part are $C^r$-functions.  
However, it is often not necessary to differentiate real and imaginary part of $z(t)$ separately in order to find the derivatives of $z(t)$, because one can differentiate $z$ like a real-valued function by treating the complex number $\ii$ as a real constant (see \cite{Amateur&Singh&saz}).  
The complex number (the vector) $z'(t)$ is called the {\em tangent vector}\index{vector!tangent vector@tangent $\sim $} of the curve $z$ at (point) $t$. 
We call the image $z(I) \subset \CC$ as {\em graph}\index{graph} of $z$.
Note that $z$ is a function while the graph of $z$ is a subset of $\CC$.\footnote{We will not use the term curve uniformly and sometimes refer to a graph as a curve. $z(t)$ is then referred to as the {\em parametric equation}\index{parametric equation} of the curve.}

\begin{defin}
A {\em regular curve}\index{curve!regular curve@regular $\sim $} is a parametric $C^1$-curve with $z'(t) \ne 0$ for every $t \in I$.
A point $t_0 \in I$ with $z'(t_0) = 0$ is called {\em singular point}\index{singular point}. 
\end{defin}

Henceforth, we will alternatively refer to $t_0$ and $z(t_0)$ as point on the curve $z$.

\begin{example}
The parametric curve
\beq
  z(t) = \ee^{\ii t} = \cos t + \ii\sin t\,,\quad 0-\ep < t < 2\pi+\ep\,,\quad \ep > 0\,, 
\eeq
is a $C^\infty$-curve and a regular curve. 
Its graph is the unit circle $x^2 + y^2 = 1$.
The tangent vector at point $t_0$ (resp.\ $\ee^{\ii t_0}$) is $z'(t_0) = \ii\ee^{\ii t_0}$.
Every tangent vector has length $\left|\ii\ee^{\ii t}\right| = \left|\ii\right| \left|\ee^{\ii t}\right| = 1$ (see Fig.\ \ref{Abb:Circle01}).

\begin{figure}[ht]
\begin{minipage}[b]{0.48\textwidth}
\centering
\includegraphics[width=0.7\textwidth]{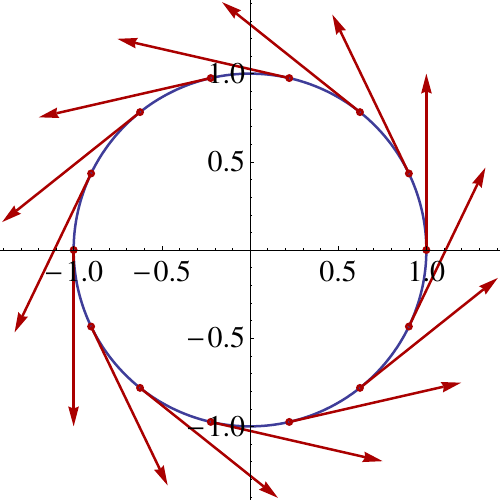}
\caption{Circle with tangent vectors}
\label{Abb:Circle01}
\end{minipage}
\hfill
\begin{minipage}[b]{0.48\textwidth}
\centering
\includegraphics[width=0.65\textwidth]{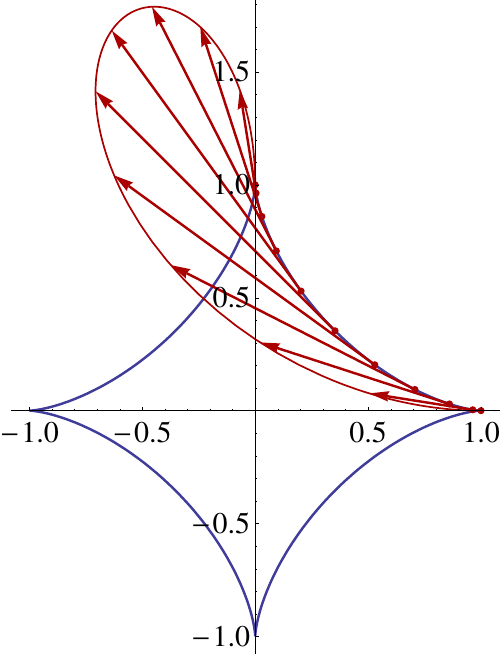}
\caption{Astroid\index{astroid} with tangent vectors}
\label{Abb:Astroid01}
\end{minipage}
\end{figure}

The parametric curve
\beq
  \tilde{z}(t) = \ee^{2\ii t} = \cos(2t) + \ii\sin(2t)\,,\quad 0-\ep < t < \pi+\ep\,,\quad \ep > 0\,, 
\eeq
is also a regular $C^\infty$-curve.
Its graph is also the unit circle, but every tangent vector has length $\left|2\ii\ee^{2\ii t}\right| = 2\left|\ii\right| \left|\ee^{2\ii t}\right| = 2$.\hfill\bs 
\end{example}

\begin{example}
The parametric curve
\beq
  z(t) = \frac{3}{4}\,\ee^{\ii t} + \frac{1}{4}\,\ee^{-3\ii t}\,,\quad 0-\ep < t < 2\pi+\ep\,,\quad \ep > 0\,, 
\eeq
is a $C^\infty$-curve.
But it is not a regular curve, because from
\beq
  z'(t) = \frac{3}{4}\,\ii\ee^{\ii t} - \frac{3}{4}\,\ii\ee^{-3\ii t}
\eeq
follows
\beq
  z'(0) = z'(\pi/2) = z'(\pi) = z'(3\pi/2) = z'(2\pi) = 0\,.
\eeq
The points $t = 0$, $\pi/2$, $\pi$, $3\pi/2$, $2\pi$ are singular points\index{singular point}.
The graph with some tangent vectors is shown in Fig.\ \ref{Abb:Astroid01}.
This curve is called {\em astroid}, and the singular points are {\em cusps}\index{cusp}.\hfill\bs   
\end{example}

\subsection{Area of the region bounded by a closed curve}

Now we derive a formula for the area\index{area} $A$ of the region bounded by a closed curve\index{curve!closed curve@closed $\sim $}.

\begin{thm} \label{Thm:Area_with_QVP}
Let $z$ be a regular closed parametric curve
\beq
  z = z(t)\,,\quad 0 \le t \le T\,,
\eeq
then the oriented area $A$ of the region bounded by $z$ is given by the line integral
\beqn \label{Eq:A}
  A
= \frac{1}{2} \oint \left[z(t),z'(t)\right] \dd t
= \frac{1}{2} \int_0^T \left[z(t),z'(t)\right] \dd t\,.\footnote{This formula with $[z(t),z'(t)] = x(t)\,y'(t)-y(t)\,x'(t)$ is  \textsc{Leibniz}'s sector formula \autocite[pp.\ 161-162]{Rothe1}, \autocite[pp.\ 52-54]{Rothe2}, \autocite{wiki:Sektorformel}, \cite[pp.\ 149-151]{Baule1}. It is a special case of \textsc{Green}'s theorem which in turn is a special case of \textsc{Stokes}' theorem \autocite[pp.\ 359-360]{Santalo}. \textsc{Gottfried Wilhelm Leibniz}, 1646-1716; \textsc{George Green}, 1793-1841; \textsc{George Gabriel Stokes}, 1819-1903}
\eeqn
\end{thm}

\begin{figure}[ht]
\begin{minipage}{0.48\textwidth}
\centering
\includegraphics[width=0.7\textwidth]{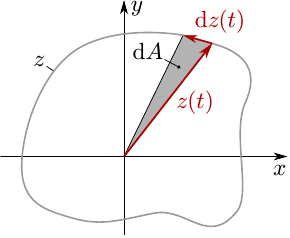}
\caption{Area element} 
\label{Abb:Area_element}
\end{minipage}
\hfill
\begin{minipage}{0.48\textwidth}
\centering
\includegraphics[width=0.6\textwidth]{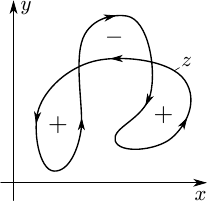}
\caption{Signed areas of oriented loops}
\label{Abb:Area_as_sum}
\end{minipage}
\end{figure}

\begin{proof}
The signed (oriented) area of the parallelogram spanned by the vectors $z(t)$ and $\dd z(t)$ is equal to $[z(t),\dd z(t)]$  (see Fig.\ \ref{Abb:Area_element}).
Thus, the area $\dd A$ (area element) of the gray triangle is given by
\beq
  \dd A
= \dd A(t)
= \frac{1}{2} \left[z(t),\dd z(t)\right]
= \frac{1}{2} \left[z(t),z'(t)\,\dd t\right]
= \frac{1}{2} \left[z(t),z'(t)\right] \dd t\,.
\eeq
The assertion follows immediately.
\end{proof}

\begin{remark} \label{Rem:A}
If a closed parametric curve $z(t)$ intersects itself, i.e., it is not a {\em simple closed curve}\index{curve!simple closed curve@simple closed $\sim $} {\em (\textsc{Jordan}\footnote{\textsc{Camille Jordan}, 1838-1922} curve)}\index{curve!Jordan curve@Jordan $\sim $}, then \eqref{Eq:A} yields the algebraic sum of the areas of all loops, where the sign of the area of each loop is given by its mathematical orientation (see Fig.\ \ref{Abb:Area_as_sum} and, for applications, Example \ref{Bsp:Gerono}, Subsubsection \ref{Sec:Example} and Eq.\ \eqref{Eq:A_profile_curve}).
\end{remark}

\begin{example} \label{Bsp:Gerono}
We calculate the area of the region bounded by the {\em \textsc{Gerono}\footnote{\textsc{Camille-Christophe Gerono}, 1799-1891} lemniscate}\index{Gerono lemniscate} {\em (eight curve)} with parametric representation
\beq
  x = x(t) = a\cos t\,,\quad
  y = y(t) = a\sin t\cos t = \tfrac{1}{2}\:\!a\sin(2t)\,,\quad
  0 \le t \le 2\pi\,,
\eeq
hence
\begin{align}
  z(t)
= {} & a\left(\cos t + \ii\,\tfrac{1}{2}\sin(2t)\right), \label{Eq:z(t)_Gerono}\\[0.05cm] 
  z'(t)
= {} & a\left(-\sin t + \ii\cos(2t)\right) \label{Eq:z'(t)_Gerono} 
\end{align}
(see Figs.\ \ref{Abb:Gerono_lemniscate1a} and \ref{Abb:Gerono_lemniscate1b}).

\begin{figure}[ht]
\begin{minipage}{0.48\textwidth}
\includegraphics[width=\textwidth]{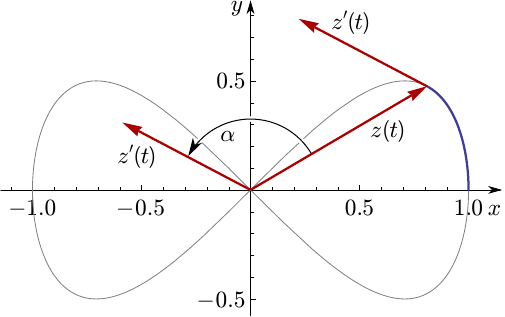}
\caption{\textsc{Gerono} lemniscate with $a = 1$; $t = \pi/5$}
\label{Abb:Gerono_lemniscate1a}
\end{minipage}
\hfill
\begin{minipage}{0.48\textwidth}
\includegraphics[width=\textwidth]{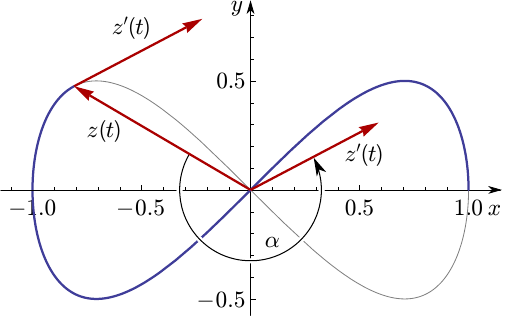}
\caption{\textsc{Gerono} lemniscate with $a = 1$; $t = 6\pi/5$}
\label{Abb:Gerono_lemniscate1b}
\end{minipage}
\end{figure}

From \eqref{Eq:[z_1,z_2]_mit_Im} we know that
\beq
  [z(t),z'(t)]
= \Imz\left(\overline{z(t)}\:z'(t)\right),  
\eeq
and with \eqref{Eq:z(t)_Gerono} and \eqref{Eq:z'(t)_Gerono} follows
\begin{align*}
  \Imz\left(\overline{z(t)}\:z'(t)\right)
= {} & a^2 \Imz\left(\left(\cos t - \ii\,\tfrac{1}{2}\sin(2t)\right) \bigl(-\sin t + \ii\cos(2t)\right)\\
= {} & a^2 \left(\cos t\cos(2t) + \tfrac{1}{2}\sin t\sin(2t)\right)\\[0.15cm]
= {} & a^2\cos^3 t\,, 
\end{align*}
thus
\beq
  A
= \frac{1}{2} \int_0^{2\pi} \left[z(t),z'(t)\right] \dd\ph
= \frac{a^2}{2} \int_0^{2\pi} \cos^3 t\, \dd t
= 0\,. 
\eeq
Considering Figs.\ \ref{Abb:Gerono_lemniscate1a} and \ref{Abb:Gerono_lemniscate1b}, this result is not surprising.
If $0 \le t < \pi/2$ or $3\pi/2 < t \le 2\pi$, then for the angle $\alpha$ between $z(t)$ and $z'(t)$ we have $0 \le \alpha \le \pi$ (see Fig.\ \ref{Abb:Gerono_lemniscate1a}), and from \eqref{Eq:QVP_with_sine} it follows that $[z(t),z'(t)] \ge 0$.
If $\pi/2 < t < 3\pi/2$, then $\pi \le \alpha \le 2\pi$ (see Fig.\ \ref{Abb:Gerono_lemniscate1b}), thus $[z(t),z'(t)] \le 0$.
The right loop of the Gerono lemniscate is positively oriented while the left loop is negatively oriented.

Denoting by $A_1$ the area of the region bounded by the right loop, we get
\beq
  A_1
= \frac{a^2}{2} \int_{-\pi/2}^{\pi/2} \cos^3 t\, \dd t
= \frac{a^2}{2} \left(\frac{3}{4}\sin t + \frac{1}{12}\sin(3t)\right)\bigg|_{-\pi/2}^{\pi/2}
= \frac{a^2}{2} \left(\frac{3}{4} - \frac{1}{12} + \frac{3}{4} - \frac{1}{12}\right)
= \frac{2a^2}{3}\,,
\eeq
hence the whole unoriented area of the region bounded by the \textsc{Gerono} lemniscate is equal to $2A_1 = 4a^2/3$. \hfill\bs
\end{example}

\subsection{Arc length}

We denote by $L(t_0,t)$ the arc length\index{arc length} of a parametric curve $z(t)$ between two points with the parameter values $t_0$ and $t$.
With the arc length element
\beq
  \dd s
= \dd s(t)
= \sqrt{(\dd x(t))^2 + (\dd y(t))^2}
= \sqrt{x'(t)^2 + y'(t)^2}\: \dd t
= \sqrt{z'(t)\,\overline{z'(t)}}\: \dd t
= |z'(t)|\, \dd t
\eeq
we have
\beqn \label{Eq:Arc_length}
  L(t_0,t)
= \int_{\tau\,=\,t_0}^t \dd s(\tau)
= \int_{ \tau\,=\,t_0}^t |z'(\tau)|\: \dd\tau\,.
\eeqn

\begin{example}
The length $L$ of the Gerono lemniscate (eight curve) in Example \ref{Bsp:Gerono} is with \eqref{Eq:z'(t)_Gerono} given by 
\begin{align*}
  L
= {} & \int_0^{2\pi} |z'(t)|\: \dd\ph
= a \int_0^{2\pi} \left|-\sin t + \ii\cos(2t)\right| \dd t
= 4a \int_0^{\pi/2} \bigl|-\sin t + \ii\cos(2t)\bigr|\: \dd t\\
= {} & 4a \int_0^{\pi/2} \sqrt{\sin^2 t + \cos^2(2t)}\; \dd t\,.
\end{align*}
Numerical evaluation of this integral with {\em Mathematica}\index{Mathematica} yields
\beq
  L = 6.09722347010491604643\ldots a \approx 6.09722 a
\eeq
(see \cite{OEIS:A118178}). \hfill\bs
\end{example}

%% file: DiffGeo5_3b.tex

\subsection{Curvature}

We will now derive a general curvature formula for a parametric curve $z(t)$.

\begin{thm} \label{Thm:Curvature}
Let the parametric curve $z$ be twice continuously differentiable at point $t$.
Then the oriented curvature\index{curvature!oriented curvature@oriented $\sim $} $\kappa(t)$ of $z$ at $t$ is given by
\beqn \label{Eq:Curvature}
  \kappa(t)
= \frac{\left[z'(t),z''(t)\right]}{\left|z'(t)\right|^3}\,,
\eeqn
and the center point $\tilde{z}(t)$ of the osculating circle\index{osculating circle}\footnote{The {\em osculating circle} is the circle that best approximates the curve at a point.} (center of curvature) by
\beqn \label{Eq:Center_of_curvature}
  \tilde{z}(t) = z(t) + \ii\, \frac{1}{\kappa(t)}\, \frac{z'(t)}{|z'(t)|}\,.
\eeqn
\end{thm}

\begin{proof}
According to the definition of the curvature of a curve, we have $\kappa = \dd\tau/\dd s$, where $s$ is the arc length, and $\tau$ is the argument of the tangent vector $z'(t)$.

We first determine a formula for the argument $\phi$ of a complex number $z = r\,\ee^{\ii\phi}$.
From $z = r\,\ee^{\ii\phi}$ follows
\[ \ln z = \ln r + \ii\phi + \ii 2k\pi\,; \]
while from $\bar{z} = r\,\ee^{-\ii\phi}$ follows
\[ \ln\bar{z} = \ln r - \ii\phi - \ii 2k\pi\,. \]
Thus we have
\[ 2\ii\phi + 4k\ii\pi
= \ln z - \ln\bar{z}
  \quad\Longrightarrow\quad
  \phi = \arg z = \frac{1}{2\ii}\, (\ln z - \ln\bar{z}) - 2k\pi\,,\quad
  k\in\Z\,. \]
For $z$ we now substitute $z'(t)$ into the last formula; then $\phi$ is to be replaced by $\tau(t)$, thus
\[ \tau(t) = \arg z'(t) = \frac{1}{2\ii} \left(\ln z'(t) - \ln\overline{z'(t)}\right) - 2k\pi\,,\quad
  k\in\Z\,. \]
According to the chain rule, we have
\[ \kappa
= \frac{\dd\tau}{\dd t} \frac{\dd t}{\dd s}\,,
  \quad\mbox{hence}\quad
  \kappa(t)
= \frac{\tau'(t)}{s'(t)}\,. \]
With
\[ \tau'(t)
= \frac{1}{2\ii} \left(\frac{z''(t)}{z'(t)} - \frac{\bar{z}''(t)}{\bar{z}'(t)}\right)
= \frac{1}{2\ii} \frac{\bar{z}'(t)\,z''(t) - z'(t)\,\bar{z}''(t)}{z'(t)\,\bar{z}'(t)}
= \frac{1}{2\ii} \frac{\bar{z}'(t)\,z''(t) - z'(t)\,\bar{z}''(t)}{\left|z'(t)\right|^2} \]
and
\[ s'(t) = \frac{\dd s(t)}{\dd t} = |z'(t)|\,, \]
taking into account \eqref{Eq:Def_QVP_complex}, it follows that
\beq
  \kappa(t)
= \frac{\bar{z}'(t)\,z''(t) - z'(t)\,\bar{z}''(t)}{2\ii\,\left|z'(t)\right|^3}
= \frac{\overline{z'(t)}\,z''(t) - z'(t)\,\overline{z''(t)}}{2\ii\,\left|z'(t)\right|^3}
= \frac{\left[z'(t),z''(t)\right]}{\left|z'(t)\right|^3}\,.
\eeq
The center of curvature (center of the osculating circle) $\tilde{z}(t)$ lies on the line normal to the oriented curve $z(t)$ at the point $t$.
It lies on the left side of the curve traveled in the direction of increasing values of $t$ when the curvature is positive and on the right side when the curvature is negative.
Since $1/|\kappa(t)|$ is the radius of curvature\index{curvature radius} and the normal vector $\ii z'(t)/|z'(t)|$ always points to the left, $\tilde{z}(t)$ is given by \eqref{Eq:Center_of_curvature}.
\end{proof}

The curvature $\kappa(t)$ is positive if the tangent vector $z'(t)$ rotates with positive direction, and negative if $z'(t)$ rotates with negative direction.\footnote{Note that the sign of the curvature depends on the orientation of the curve. Reversing the orientation causes the change of the sign.}
Due to the definition of the quasi vector product, $\left[z'(t),z''(t)\right]$ is the signed area of the parallelogram spanned by the vectors $z'(t)$ and $z''(t)$.
Therefore, $\kappa(t) > 0$ if $0 < \measuredangle(z'(t),z''(t)) < \pi$, and $\kappa(t) < 0$ if $\pi < \measuredangle(z'(t),z''(t)) < 2\pi$.

From \eqref{Eq:Curvature}, using \eqref{Eq:Def_QVP_real}, follows
\beq
  \kappa(t)
= \frac{x'(t)\,y''(t) - y'(t)\,x''(t)}{\left(x'(t)^2 + y'(t)^2\right)^{3/2}}\,, 
\eeq
the well-known formula for the curvature of a plane curve in real parametric representation $x(t), y(t)$ (see e.\,g.\ \textcite[p.\ 446]{Gellert}, \textcite[p.\ 415]{Bosch}, \textcite[p.\ 20]{Eschenburg&Jost}.).

\begin{example}
Using \eqref{Eq:Curvature}, \eqref{Eq:z'(t)_Gerono}, $z''(t) = a\left(-\cos t-2\ii\sin(2 t)\right)$ and \eqref{Eq:[z_1,z_2]_mit_Im}, the curvature of the \textsc{Gerono} lemniscate (eight curve) from Example \ref{Bsp:Gerono} is given by
\begin{align} \label{Eq:Gerono_curvature}
  \kappa(t)
= {} & \frac{\left[z'(t),z''(t)\right]}{\left|z'(t)\right|^3}
= \frac{\left[a\left(-\sin t+\ii\cos(2t)\right),\,a\left(-\cos t-2\ii\sin(2t)\right)\right]}{a^3\left|-\sin t+\ii\cos(2t)\right|^3}
  \nonumber\db\\[0.05cm]
= {} & \frac{a^2\Imz\left\{\left(-\sin t-\ii\cos(2t)\right)\left(-\cos t-2\ii\sin(2t)\right)\right\}}{a^3\left(\cos^2(2t)+\sin^2 t\right)^{3/2}}
  \nonumber\db\\[0.05cm]
= {} & \frac{2\sin t\sin(2t)+\cos t\cos(2t)}{a\left(\cos^2(2t)+\sin^2 t\right)^{3/2}}   
= \frac{3\cos t-\cos(3t)}{2a\left(\cos^2(2t)+\sin^2 t\right)^{3/2}}
\end{align}
(see Fig.\ \ref{Abb:Gerono_curvature}).
If $0 \le t < \pi/2$ or $3\pi/2 < t \le 2\pi$, then for the angle $\beta$ between $z'(t)$ and $z''(t)$ we have $0 \le \beta \le \pi$ (see Fig.\ \ref{Abb:Gerono_lemniscate2a}), and from \eqref{Eq:QVP_with_sine} it follows that $[z'(t),z''(t)] \ge 0$.
If $\pi/2 < t < 3\pi/2$, then $\pi \le \beta \le 2\pi$ (see Fig.\ \ref{Abb:Gerono_lemniscate2b}), thus $[z'(t),z''(t)] \le 0$.

\begin{SCfigure}[0.85][ht]
\includegraphics[width=0.45\textwidth]{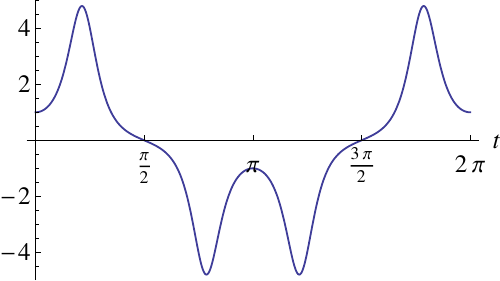}
\caption{Curvature function of the \textsc{Gerono} lemniscate with $a = 1$}
\label{Abb:Gerono_curvature}
\end{SCfigure}

\begin{figure}[ht]
\begin{minipage}{0.48\textwidth}
\includegraphics[width=\textwidth]{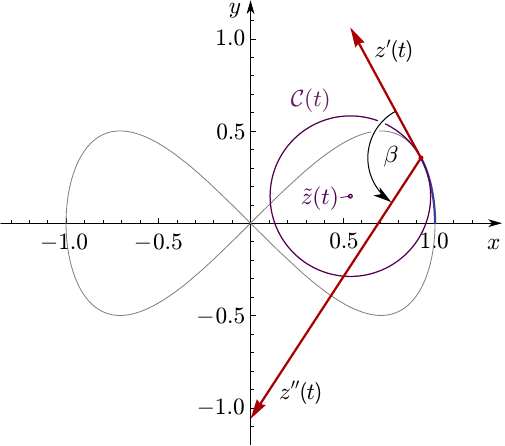}
\caption{\textsc{Gerono} lemniscate ($a = 1$) with curvature circle $\K(t)$ for $t = \pi/8$}
\label{Abb:Gerono_lemniscate2a}
\end{minipage}
\hfill
\begin{minipage}{0.48\textwidth}
\includegraphics[width=\textwidth]{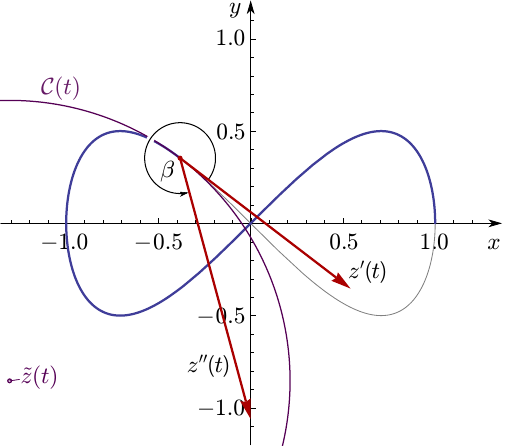}
\caption{\textsc{Gerono} lemniscate ($a = 1$) with curvature circle $\K(t)$ for $t = 11\pi/8$}
\label{Abb:Gerono_lemniscate2b}
\end{minipage}
\end{figure}

From \eqref{Eq:Center_of_curvature} with \eqref{Eq:z(t)_Gerono}, \eqref{Eq:z'(t)_Gerono} and \eqref{Eq:Curvature} we get
\beq
  \tilde{z}(t)
= a\left(\cos t+\ii\,\tfrac{1}{2}\sin(2t)\right) 
  - 2a\, \frac{\cos^2(2t)+\sin^2 t}{3\cos t-\cos(3t)}\, \bigl(\cos(2t)+\ii\sin t\bigr)
\eeq
as center point of the curvature circle.
Splitting into real and imaginary part provides
\begin{align*}
  \tilde{x}(t)
= {} & \Rez\left(\tilde{z}(t)\right)
= \frac{4-3\cos(2t)-\cos(6t)}{2\left(3\cos t-\cos(3t)\right)}\,a\,,\db\\[0.05cm]   
  \tilde{y}(t)
= {} & \Imz\left(\tilde{z}(t)\right)
= \frac{-6\sin t+7\sin(3t)-3\sin(5t)}{4\left(3\cos t-\cos(3t)\right)}\,a\,. \tag*{\bs}
\end{align*}
\end{example}

\begin{lem} \label{Lem:T'(t)}
Let $T(t)$ be the tangent unit vector at point $t$ of a parametric curve $z$.
Then the derivative $T'(t) = \frac{\dd}{\dd t}\,T(t)$ is given by
\beq
  T'(t)
= \kappa(t) \left|z'(t)\right| \ii\, T(t)\,,  
\eeq
where $\kappa(t)$ is the oriented curvature.\footnote{Since $\ii T$ is the normal unit vector $N$, the formula for $T'(t)$ is the first \textsc{Frenet} formula for plane curves (see e.g.\ \cite[p.\ 16]{doCarmo}, \cite[p.\ 20]{Eschenburg&Jost}), now in complex-valued form. Note $\frac{\dd T}{\dd s} = \frac{\dd T}{\dd t}\, \frac{\dd t}{\dd s} = \frac{\dd T}{\dd t}\, \frac{1}{|z'(t)|} = \kappa N$ for this. \textsc{Jean Fr\'ed\'eric Frenet}, 1816-1900.}
\end{lem}

\begin{proof}
We have
\beq
  \frac{\dd T(t)}{\dd t}
= \frac{\dd}{\dd t}\, \frac{z'(t)}{|z'(t)|}
= \frac{\dd}{\dd t} \left(z'(t)\left|z'(t)\right|^{-1}\right)
= z''(t) \left|z'(t)\right|^{-1} - z'(t) \left|z'(t)\right|^{-2} \left|z'(t)\right|'. 
\eeq
With
\begin{align*}
  \left|z'(t)\right|'
= {} & \frac{\dd}{\dd t} \left(z'(t)\,\overline{z'(t)}\right)^{1/2}
= \frac{1}{2} \left(z'(t)\,\overline{z'(t)}\right)^{-1/2} \left(z''(t)\,\overline{z'(t)} + z'(t)\,\overline{z''(t)}\right)\\[0.05cm]
= {} & \frac{1}{2}\, \frac{z''(t)\,\overline{z'(t)} + z'(t)\,\overline{z''(t)}}{\left|z'(t)\right|}\,,  
\end{align*}
it follows
\begin{align*}
  \frac{\dd T(t)}{\dd t}
= {} & \frac{z''(t)}{|z'(t)|}
  - \frac{z'(t)}{|z'(t)|^2}\, \frac{1}{2}\, \frac{z''(t)\,\overline{z'(t)} + z'(t)\,\overline{z''(t)}}{\left|z'(t)\right|}\db\\[0.05cm]
= {} & \frac{z''(t)\left|z'(t)\right|^2 - \frac{1}{2}\,z''(t)\,z'(t)\,\overline{z'(t)} - \frac{1}{2}\,z'(t)^2\,\overline{z''(t)}}{|z'(t)|^3}\db\\[0.05cm]   
= {} & \frac{1}{2}\, \frac{z''(t)\left|z'(t)\right|^2 - z'(t)^2\,\overline{z''(t)}}{|z'(t)|^3}
= \frac{1}{2}\, \frac{z''(t)\,z'(t)\,\overline{z'(t)} - z'(t)^2\,\overline{z''(t)}}{|z'(t)|^3}\db\\[0.05cm]
= {} & \frac{z'(t)}{|z'(t)|}\, \frac{1}{2}\, \frac{z''(t)\,\overline{z'(t)} - z'(t)\,\overline{z''(t)}}{|z'(t)|^2}
= T(t) \left|z'(t)\right| \frac{1}{2}\, \frac{z''(t)\,\overline{z'(t)} - z'(t)\,\overline{z''(t)}}{|z'(t)|^3}\db\\[0.05cm]
= {} & T(t) \left|z'(t)\right| \left(-\frac{1}{2}\, \frac{z'(t)\,\overline{z''(t)} - z''(t)\,\overline{z'(t)}}{|z'(t)|^3}\right)
= T(t) \left|z'(t)\right| \ii\, \frac{\ii}{2}\, \frac{z'(t)\,\overline{z''(t)} - z''(t)\,\overline{z'(t)}}{|z'(t)|^3}\db\\[0.05cm]
= {} & T(t) \left|z'(t)\right| \ii\, \frac{[z'(t),z''(t)]}{|z'(t)|^3}
= \kappa(t) \left|z'(t)\right| \ii\, T(t)\,. \qedhere
\end{align*} 
\end{proof}

\subsection{Evolutes and involutes (evolvents)}

\begin{defin}
The {\em evolute}\index{evolute} of a curve is the locus of all its centers of curvature.
\end{defin}

The evolute is also the envelope\index{envelope} of the family of normal lines to the curve \cite[p.\ 95]{Baule1}, \cite[p.\ 64]{Schoene}, \cite{Ghys&Tabachnikov&Timorin}.
We will prove this fact in Theorem \ref{Thm:Envelope=evolute} (see also Figs.\ \ref{Abb:Ellipse_mit_Evolute}, \ref{Abb:Evolute_of_coupler_curve} and \ref{Abb:Circle_with_involutes01}).

A parametric equation of the evolute of a curve is given by \eqref{Eq:Center_of_curvature}. 

\begin{example}
We consider the ellipse given by the parametric equation
\beq
  z(t)
= a\cos t + \ii b\sin t\,,\quad
  0 \le t \le 2\pi\,.   
\eeq
From \eqref{Eq:Curvature} with
\beq
  z'(t)
= -a\sin t + \ii b\cos t\,,\qquad
  z''(t)
= -a\cos t - \ii b\sin t
\eeq
we get
\beqn \label{Eq:Curvature_of_ellipse}
  \kappa(t)
= \frac{ab}{\left(a^2\sin^2 t + b^2\cos^2 t\right)^{3/2}}\,.  
\eeqn
Now from \eqref{Eq:Center_of_curvature} one easily gets
\beqn \label{Eq:Evolute_of_the_ellipse_1}
  \tilde{z}(t)
= \frac{a^2-b^2}{a}\cos^3 t + \ii\,\frac{b^2-a^2}{b}\sin^3 t  
\eeqn
as equation of the evolute (see Fig.\ \ref{Abb:Ellipse_mit_Evolute}). \hfill\bs

\begin{SCfigure}[][ht]
\includegraphics[width=0.5\textwidth]{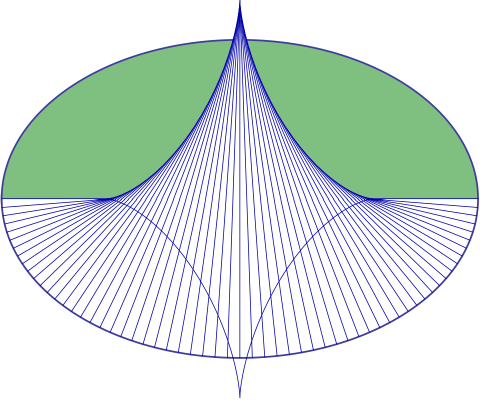}
\caption{Ellipse with $a/b = 3/2$ and its evolute (curve with four cusps)\\[0.2cm] Every line segment connects corresponding points $z(t)$ and $\tilde{z}(t)$ of the ellipse and its evolute, respectively. The length of such a line segment is equal to the curvature radius (radius of the osculating circle) $= 1/\kappa(t)$.\\[0.2cm] Minimum curvature radius $= 1/\kappa(0) = \frac{4}{9}\,a = 0.\overline{4}\,a$, maximum curvature radius $1/\kappa(\pi/2) = \frac{3}{2}\,a = 1.5\,a$\\[0.2cm] The evolute is also the envelope of the family of normals to the ellipse.}
\label{Abb:Ellipse_mit_Evolute}
\end{SCfigure}

\begin{SCfigure}[][ht]
\includegraphics[width=0.5\textwidth]{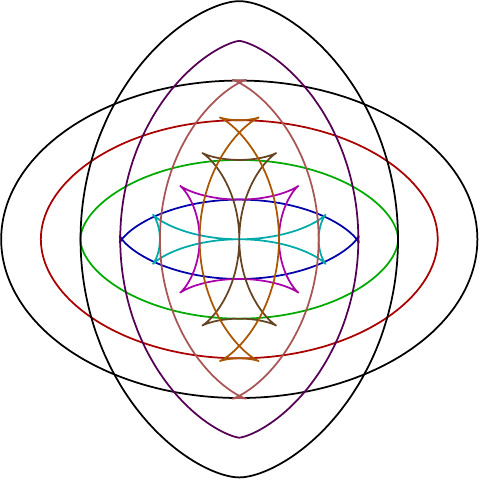}
\caption{Ellipse ($a/b = 3/2$) with ``inner'' parallel curves\index{parallel curve}\index{curve!parallel curve@parallel $\sim $} with $\la/b = k/4$, $k = 1,2,\ldots,10$}
\label{Abb:Ellipse_with_parallel_curves}
\end{SCfigure}
\end{example}

From the parametric curve $z(t)$ one easily gets the parametric equation $z_\la(t)$ of the parallel curve\index{parallel curve} at signed distance $\la \in \R$ with
\beq
  z_\la(t)
= z(t) + \ii\,\la\,\frac{z'(t)}{|z'(t)|}\,.  
\eeq
Fig.\ \ref{Abb:Ellipse_with_parallel_curves} shows the ellipse from Fig.\ \ref{Abb:Ellipse_mit_Evolute} with some of its parallel curves with $\la > 0$.

\begin{defin}
Each point of a straight line rolling on a given plane curve describes a curve which is called {\em an involute}\index{involute} (also {\em an evolvent}\index{evolvent}) of the given curve.
Instead of the straight line, one can also imagine a taut string that is unwrapped from or wrapped onto the given curve.
\end{defin}

Let $z$ be a regular plane curve with parametric equation $z(t)$ and nowhere vanishing curvature, then an involute $z_a$ of $z$ is given by the parametric equation
\beqn \label{Eq:Involute}
  z_a(t)
= z(t) - \frac{z'(t)}{\left|z'(t)\right|}\, \int_a^t \left|z'(\tau)\right| \dd\tau\,. 
\eeqn 

\begin{SCfigure}[][!htb]
\includegraphics[width=0.6\textwidth]{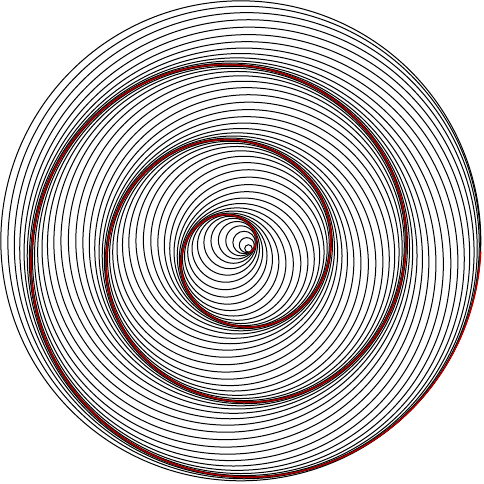}
\caption{Half involute to a circle (circle not shown) with osculating circles}
\label{Abb:Involute_with_osculating_circles}
\end{SCfigure}

\begin{example}
The equation of the/one involute\index{involute} (evolvent)\index{evolvent} of a circle with radius $r$ and center in the coordinate origin (see Figs.\ \ref{Abb:Evolvente01}, \ref{Abb:Kreisevolventen03}, and, with osculating circles according to Theorem \ref{Thm:Curvature}, Fig.\ \ref{Abb:Involute_with_osculating_circles}) is
\begin{align} \label{Eq:Circle_involute01}
  z(\ph)
= {} & r\ee^{\ii\ph} + r\ph(-\ii\ee^{\ii\ph})\nonumber\\ 
= {} & r \left(1-\ii\ph\right) \ee^{\ii\ph}\,. 
\end{align}
The involute with this equation has its cusp at point $r\ee^{\ii 0}$.

\begin{figure}[!ht]
\begin{minipage}{0.47\textwidth}
\centering
\includegraphics[width=0.7\textwidth]{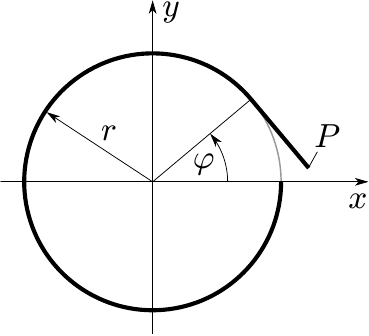}
\caption{Involute of the circle: unwinding/winding a string from/onto the circle}
\label{Abb:Evolvente01}
\vspace{0.3cm}
\end{minipage}
\hfill
\begin{minipage}{0.47\textwidth}
\centering
\includegraphics[width=0.7\textwidth]{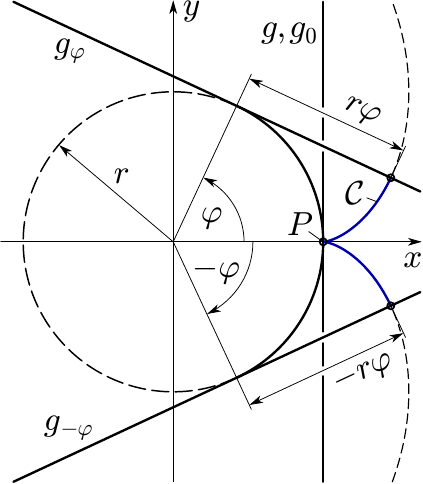}
\caption{Involute $\mathcal{C}$ of the circle: rolling of a line $g$ on the circle}
\label{Abb:Kreisevolventen03}
\vspace{0.3cm}
\end{minipage}
\end{figure}

\begin{figure}[!ht]
\begin{minipage}{\textwidth}
\begin{minipage}{0.48\textwidth}
\centering
\includegraphics[width=\textwidth]{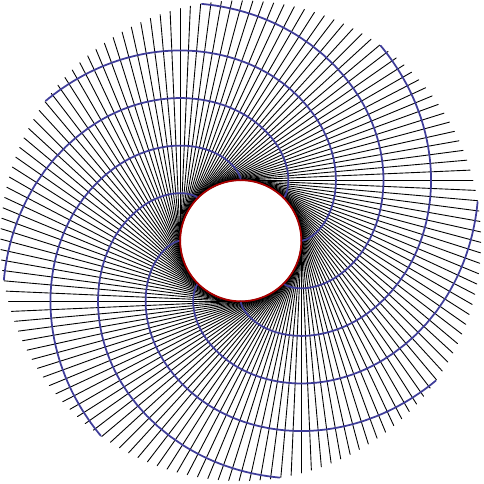}
\end{minipage}
\hfill
\begin{minipage}{0.48\textwidth}
\centering
\includegraphics[width=\textwidth]{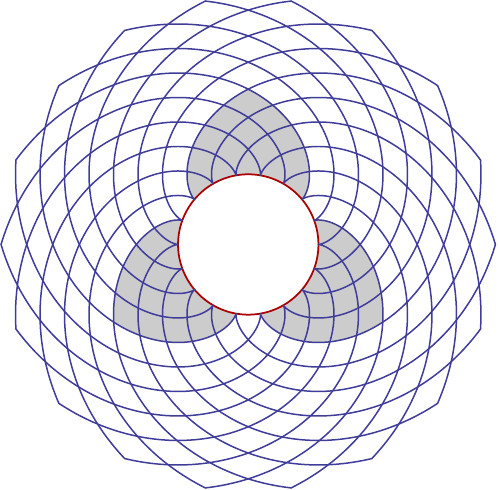}
\end{minipage}
\end{minipage}
\vspace*{0.3cm}

\begin{minipage}{\textwidth}
\begin{minipage}[t]{0.48\textwidth}
\caption{Circle with (half) tangents and (half) involutes\\[0.2cm] The circle is the evolute of the involutes. (A curve is the evolute of each of its involutes.) The circle tangents and the involutes are perpendicular to each other.\\[0.2cm] The circle (evolute) is the envelope of the family of normal lines to the involutes.}
\label{Abb:Circle_with_involutes01}
\end{minipage}
\hfill
\begin{minipage}[t]{0.48\textwidth}
\caption{Circle with involutes}
\label{Abb:Circle_with_involutes02}
\end{minipage}
\end{minipage}
\end{figure}

In order to obtain the equation of the involute with cusp at point $r\ee^{\ii\alpha}$ (see Figs.\ \ref{Abb:Circle_with_involutes01} and \ref{Abb:Circle_with_involutes02}), we multiply the right side of \eqref{Eq:Circle_involute01} by $\ee^{\ii\alpha}$,
\beqn \label{Eq:Circle_involute02}
  z(\ph)
= r \left(1-\ii\ph\right) \ee^{\ii\ph} \ee^{\ii\alpha}    
= r \left(1-\ii\ph\right) \ee^{\ii(\ph+\alpha)}\,.\footnote{The cusps always occur at $\ph = 0$.}
\eeqn
\hfill\bs
\end{example} 

\begin{thm} \label{Thm:Evolute_of_its_involutes}
Let $z(t)$ be a twice continuously differentiable parametric curve with curvature $\kappa(t) > 0$, and arc length $s$ with $s(t) > s(a)$ for every $t > a$.
Then $z$ is the evolute of each of its involutes $z_a$.
\end{thm}

\begin{proof}
We start from \eqref{Eq:Involute}.
$\frac{z'(t)}{|z'(t)|}$ is the tangent unit vector $T(t)$ at point $t$ of $z$.
Denoting by $s$ the arc length of the curve $z$, we can write \eqref{Eq:Involute} as
\beqn \label{Eq:z_a(t)}
  z_a(t)
= z(t) - T(t)\,(s(t)-s(a))\,.  
\eeqn
It follows
\begin{align*}
  z_a'(t)
= {} & z'(t) - T'(t)\,(s(t)-s(a)) - T(t)\,s'(t)\db\\
= {} & z'(t) - T'(t)\,(s(t)-s(a)) - \frac{z'(t)}{|z'(t)|}\,|z'(t)|\db\\
= {} & {-}T'(t)\,(s(t)-s(a))\,.  
\end{align*}
Using Lemma \ref{Lem:T'(t)} we get
\beqn \label{Eq:Tangent_vector_of_the_involute}
  z_a'(t)
= -\kappa(t) \left|z'(t)\right| \ii\, T(t)\, (s(t)-s(a))\,,  
\eeqn
hence the tangent unit vector of the involute $z_a$ is given by
\beq
  T_a(t)
= \frac{z_a'(t)}{|z_a'(t)|}
= \frac{-\kappa(t)\,|z'(t)|\,\ii\,T(t)\,(s(t)-s(a))}{|\kappa(t)|\,|z'(t)|\,|s(t)-s(a)|}\,,   
\eeq
thus
\beqn \label{Eq:T_a(t)}
  T_a(t)
= -\ii\,T(t) 
  \quad\mbox{for}\quad
  t > a\,. 
\eeqn
Again with Lemma \ref{Lem:T'(t)} follows
\beq
  T_a'(t)
= -\ii\,T'(t)  
= -\ii\, \kappa(t) \left|z'(t)\right| \ii\, T(t)   
= \kappa(t) \left|z'(t)\right| T(t)
\eeq
and
\beq
  T_a'(t)
= \kappa_a(t) \left|z_a'(t)\right| \ii\, T_a(t)\,,   
\eeq
and so we have
\beq
  \kappa_a(t) \left|z_a'(t)\right| \ii\, T_a(t)
= \kappa(t) \left|z'(t)\right| T(t)\,,  
\eeq
hence
\beq
  \ii\,\kappa_a(t)\,T_a(t)
= \frac{\kappa(t)\,|z'(t)|\,T(t)}{|z_a'(t)|}\,.  
\eeq
Using \eqref{Eq:Tangent_vector_of_the_involute} we get    
\beq
  \ii\,\kappa_a(t)\,T_a(t)
= \frac{\kappa(t)\,|z'(t)|\,T(t)}{|\kappa(t)|\,|z'(t)|\,|s(t)-s(a)|}
= \frac{T(t)}{s(t)-s(a)}\,,  
\eeq
and with \eqref{Eq:T_a(t)}
\beqn \label{Eq:kappa_a(t)}
  \kappa_a(t)
= \frac{1}{s(t)-s(a)}\,.  
\eeqn
The evolute of $z_a$ is given by
\begin{align}
  \tilde{z}(t)
= {} & z_a(t) + \ii\,\frac{1}{\kappa_a(t)}\,\frac{z_a'(t)}{|z_a'(t)|}\nonumber\\ 
= {} & z_a(t) + \ii\,\frac{1}{\kappa_a(t)}\,T_a(t)\,.\label{Eq:Evolute_of_z_a} 
\end{align}
Substituting \eqref{Eq:z_a(t)}, \eqref{Eq:T_a(t)} and \eqref{Eq:kappa_a(t)} into \eqref{Eq:Evolute_of_z_a} gives
\begin{align*}
  \tilde{z}(t)
= {} & z(t) - T(t)\,(s(t)-s(a)) + \ii\,(s(t)-s(a))\,(-\ii\,T(t))\db\\
= {} & z(t) - T(t)\,(s(t)-s(a)) + (s(t)-s(a))\,T(t)\db\\
= {} & z(t)\,, 
\end{align*}
thus the evolute of $z_a$ is $z$ itself. 
\end{proof}

\begin{example}
In the following, we will determine the involutes to the parametric curve \eqref{Eq:Evolute_of_the_ellipse_1}.
From Theorem \ref{Thm:Evolute_of_its_involutes} we know that a curve is the evolute of each of its involutes.
But we can not use this in our case because the curve \eqref{Eq:Evolute_of_the_ellipse_1}, having four singularities, is not regular, and its curvature is not always $> 0$.  
From \eqref{Eq:Involute} we know that for the determination of the parametric equation of the involute of a curve we need the tangent unit vector and the integral $\int|z'(\tau)|\,\dd\tau$.

By omitting the tilde we write \eqref{Eq:Evolute_of_the_ellipse_1} as
\beqn \label{Eq:Evolute_of_the_ellipse_2}
  z(t)
= \left(b^2-a^2\right) \left(-\frac{\cos^3 t}{a} + \ii\,\frac{\sin^3 t}{b}\right).  
\eeqn
It follows
\begin{gather*}
  z'(t)
= \frac{3}{2}\, \frac{b^2-a^2}{ab} \sin(2 t)\, (b\cos t + \ii\,a\sin t)\,,\db\\[0.05cm]
  |z'(t)|
= \frac{3}{2}\, \frac{\left|a^2-b^2\right|}{ab} \left|\sin(2 t)\right| \sqrt{(a\sin t)^2 + (b\sin t)^2}\,,    
\end{gather*}
and the tangent unit vector is
\beq
  T(t)
= \frac{z'(t)}{|z'(t)|}
= \sgn(b^2-a^2)\, \sgn(\sin(2t))\,
  \frac{b\cos t + \ii\,a\sin t}{\sqrt{(a\sin t)^2 + (b\cos t)^2}}\,,   
\eeq
where sgn is the sign function defined by
\beq
  \sgn x
= \left\{\begin{array}{r@{\quad\mbox{if}\quad}c}
	-1 & x < 0\,,\\
	 0 & x = 0\,,\\
	 1 & x > 0\,. 
  \end{array}\right.	  
\eeq
Under the assumption $\sin(2t) > 0$, {\em Mathematica}\index{Mathematica} with assistance finds
\beq
  \int\left|z'(t)\right|\dd t
= \sgn\left(a^2-b^2\right) \frac{((a\sin t)^2+(b\cos t)^2)^{3/2}}{ab}\,. 
\eeq
Apart from the term with the sign function, this is the radius of curvature of the ellipse (see \eqref{Eq:Curvature_of_ellipse}).

One finds that an involute $z_\la$ of $z$ is given by the parametric equation
\beq
  z_\la(t)
= z(t) + \left(\frac{((a\sin t)^2+(b\cos t)^2)^{3/2}}{ab}-\la\right) \frac{b\cos t+\ii\,a\sin t}{\sqrt{(a\sin t)^2+(b\cos t)^2}}\,,\quad
  \la\in\R\,,
\eeq
with $z(t)$ from \eqref{Eq:Evolute_of_the_ellipse_2}.

For $a = 3$ and $b=2$, Figs.\ \ref{Abb:Ellipse_as_involute1} and \ref{Abb:Ellipse_as_involute2} show the involute $z_0$, which is the ellipse with $a=3$ and $b=2$ (see Fig.\ \ref{Abb:Ellipse_mit_Evolute}).
For $\la \ne 0$ one gets parallel curves of this ellipse (see Fig.\ \ref{Abb:Ellipse_with_parallel_curves}, and Figs.\ \ref{Abb:Curve_C_with_involute1}, \ref{Abb:Curve_C_with_involute2}). \hfill\bs

\begin{figure}[H]
\begin{minipage}{0.49\textwidth}
\centering
\includegraphics[width=\textwidth]{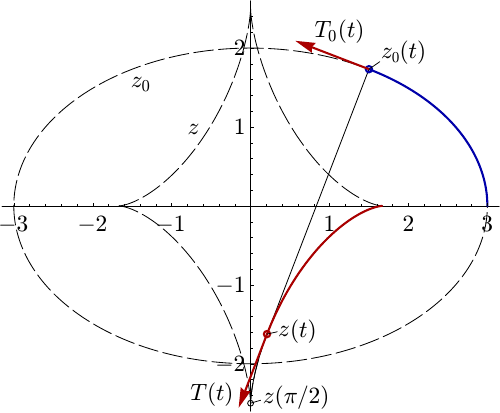}
\caption{Unwrapping a taut string, fixed at point $z(\pi/2)$, from the curve $z$; $T_0(t) = -\ii T(t)$}
\label{Abb:Ellipse_as_involute1}
\end{minipage}
\hfill
\begin{minipage}{0.49\textwidth}
\centering
\includegraphics[width=\textwidth]{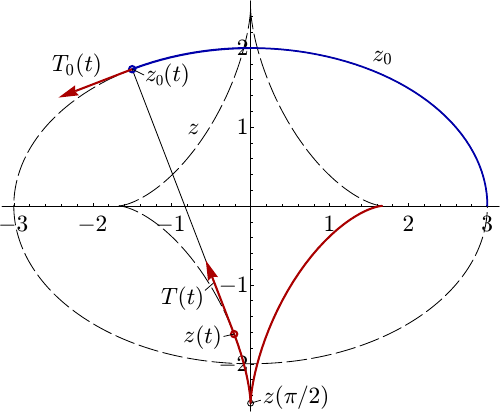}
\caption{Wrapping the taut string, fixed at point $z(\pi/2)$, onto the curve $z$; $T_0(t) = \ii T(t)$}
\label{Abb:Ellipse_as_involute2}
\end{minipage}
\end{figure}

\begin{figure}[!ht]
\begin{minipage}{0.49\textwidth}
\centering
\captionsetup{justification=centering}
\includegraphics[width=\textwidth]{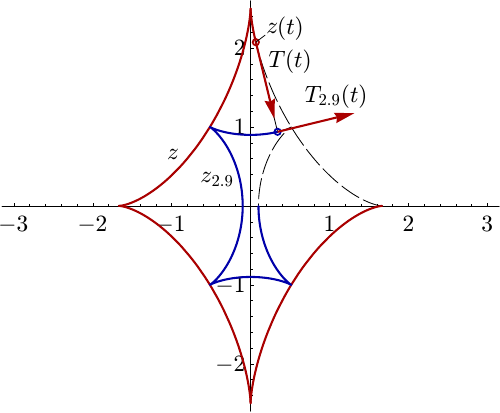}
\caption{Curve $z$ with involute $z_{2.9}$;\\ $t = 290\g\cdot\pi/180\g$; $T_{2.9}(t) = \ii T(t)$}
\label{Abb:Curve_C_with_involute1}
\end{minipage}
\hfill
\begin{minipage}{0.49\textwidth}
\centering
\captionsetup{justification=centering}
\includegraphics[width=\textwidth]{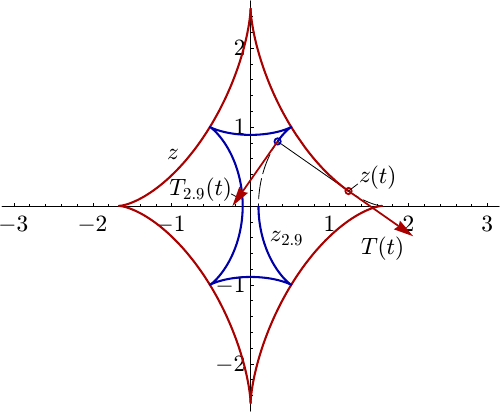}
\caption{Curve $z$ with involute $z_{2.9}$;\\ $t = 335\g\cdot\pi/180\g$; $T_{2.9}(t) = -\ii T(t)$}
\label{Abb:Curve_C_with_involute2}
\end{minipage}
\end{figure}
\end{example}

\begin{thm}[\textsc{Tait}-\textsc{Kneser} Theorem\footnote{\textsc{Peter Guthrie Tait}, 1831-1901 (see also \cite[pp.\ 189-191]{Nahin1}); \textsc{Adolf Kneser}, 1862-1930}] \label{Thm:Tait-Kneser}
Let $z(t)$ be a twice continuously differentiable parametric curve with continuously increasing or decreasing curvature $ \kappa(t) > 0$, then the osculating circles of this curve do not intersect and are nested (like russian matryoshka (\foreignlanguage{russian}{matr\"eshka}) dolls).
\end{thm}

\begin{figure}[!htb]
\begin{minipage}{0.41\textwidth}
\centering
\includegraphics[width=0.88\textwidth]{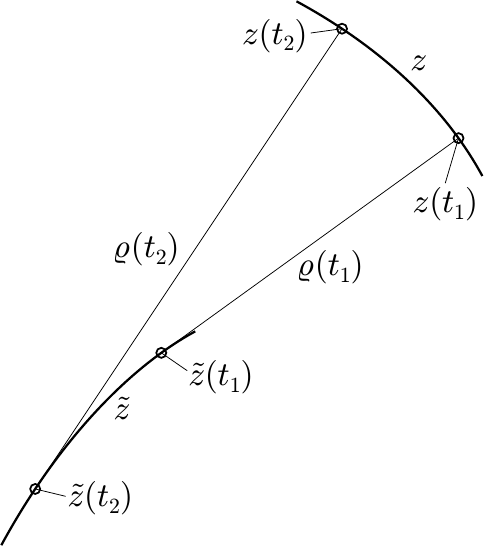}
\caption{To the proof of Theorem \ref{Thm:Tait-Kneser}}
\label{Abb:Tait_Kneser01}
\end{minipage}
\hfill
\begin{minipage}{0.55\textwidth}
\centering
\includegraphics[width=0.94\textwidth]{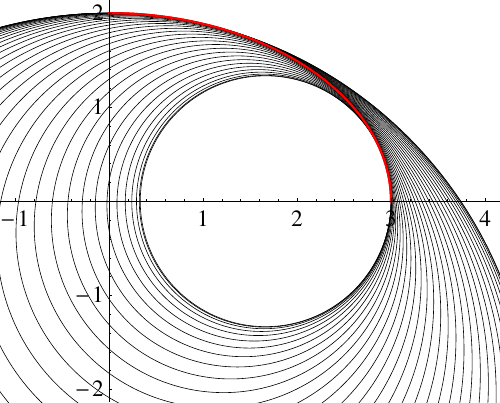}
\caption{Arc of an ellipse (red) with osculating circles}
\label{Abb:Ellipse_with_osculating_circles01}
\end{minipage}
\end{figure}

\begin{proof}
Let $z(t_1)$ and $z(t_2)$ with $t_2 > t_1$ be two points of $z$, and let $\K_1$ and $\K_2$, respectively, be the osculating circles of $z$ at this points.
Let $\tilde{z}(t_1)$ and $\tilde{z}(t_2)$ be the centers, and $\rh(t_1)$ and $\rh(t_2)$ the radii (see Fig.\ \ref{Abb:Tait_Kneser01}) of $\K_1$ and $\K_2$, respectively.
$|\rh(t_2)-\rh(t_1)|$ is equal to the length of the arc $\tilde{z}(t_1)\,\tilde{z}(t_2)$ on the evolute $\tilde{z}$ because, according to Theorem \ref{Thm:Evolute_of_its_involutes}, the curve $z$ is an involute of its evolute $\tilde{z}$.
Since the evolute $\tilde{z}$ cannot contain a straight piece, it follows that
\beq
  \left|\tilde{z}(t_2)-\tilde{z}(t_1)\right| < |\rh(t_2)-\rh(t_1)|\,.
\eeq  
This means that $\K_1$ and $\K_2$ do not intersect.
If the curvature $\kappa(t)$ decreases with increasing values of $t$ (as in Fig.\ \ref{Abb:Tait_Kneser01}), then the curvature radius $\rh(t)$ increases and $\K_1$ lies within $\K_2$. 
If the curvature $\kappa(t)$ increases with increasing values of $t$, then the curvature radius $\rh(t)$ decreases and $\K_2$ lies within $\K_1$.
\end{proof}

Concerning the \textsc{Tait}-\textsc{Kneser} theorem\index{Tait-Kneser theorem} see \cite{Tait}, \cite{Kneser}, and \cite{Ghys&Tabachnikov&Timorin}, \cite{wiki:Tait-Kneser_theorem}.
For examples see Figs.\ \ref{Abb:Ellipse_with_osculating_circles01} and \ref{Abb:Involute_with_osculating_circles}.
Fig.\ \ref{Abb:Ellipse_with_osculating_circles01} shows osculating circles of the ellipse arc $z(t) = 3\cos t + 2\ii\sin t$, $0 \le t \le \pi/2$, whereas Fig.\ \ref{Abb:Involute_with_osculating_circles} shows osculating circles of the half involute to a circle defined by \eqref{Eq:Circle_involute01} or \eqref{Eq:Circle_involute02} with $\ph \ge 0$ in each case.
Here, the radius of the osculating circle for parameter value $\ph$ is equal to $r\ph$, hence the curvature is strictly increasing.

\subsection{The four-vertex theorem}

\begin{figure}[ht]
\begin{minipage}{0.48\textwidth}
\centering
\includegraphics[width=0.88\textwidth]{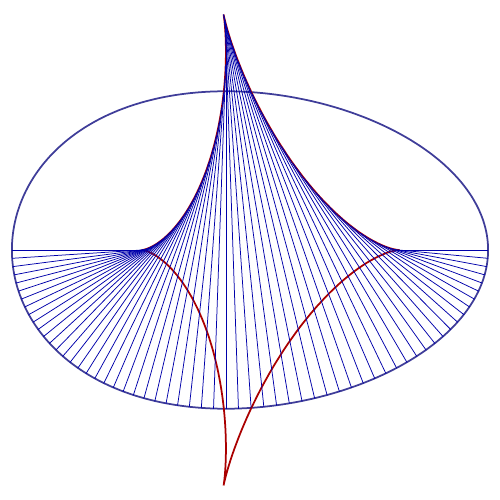}
\caption{Strictly convex simple closed curve (``Eilinie'') with evolute}
\label{Abb:Eilinie}
\end{minipage}
\hfill
\begin{minipage}{0.48\textwidth}
\centering
\includegraphics[width=0.88\textwidth]{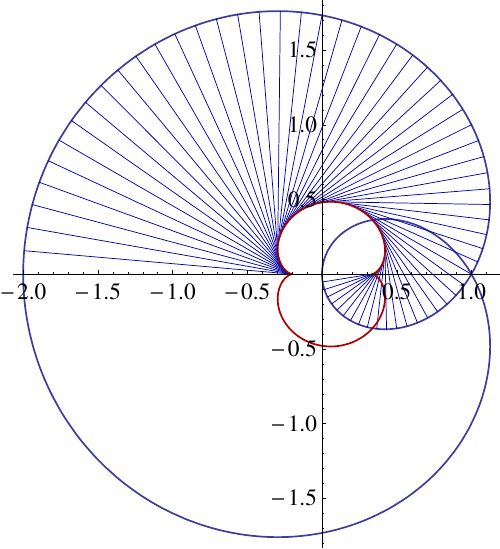}
\caption{One lima\c{c}on of Pascal with evolute}
\label{Abb:Kreiskonchoide}
\end{minipage}
\end{figure}

Clearly, an ellipse -- apart from the special case of a circle -- has exactly four points at which the curvature attains an extremum; the corresponding centers of curvature are cusps of the evolute (see Fig.\ \ref{Abb:Ellipse_mit_Evolute}).
Also, the egg-shaped curve in Fig.\ \ref{Abb:Eilinie} has exactly four points at which the curvature attains its extreme values.

Now, we consider the lima\c{c}on of \textsc{Pascal}\footnote{\textsc{\'Etienne Pascal}, 1588-1651} (see \cite[pp.\ 3-4, 9]{Krause&Carl}, \cite[pp.\ 51-62]{Wunderlich}, \cite[pp.\ 103-104; 106]{Baesel:Kolloquium}, \cite{wiki:Limacon}) with parametric equation
\beq
  z(\ph)
= \ee^{\ii\ph} - \ee^{2\ii\ph}\,,\quad
  0 \le \ph \le 2\pi\,.
\eeq
This is a closed self-intersecing plane curve (see Fig.\ \ref{Abb:Kreiskonchoide}) whose continuous curvature function
\beq
  \kappa(\ph)
= \frac{9-6\cos\ph}{(5-4\cos\ph)^{3/2}}  
\eeq
has exactly one minimum and one maximum.
(For a similar curve with evolute see \cite[p.\ 150]{Fuchs&Tabachnikov}.) 

A point of a curve with an extreme value of the curvature is called {\em vertex}.

The following theorem, called the {\em four-vertex theorem}\index{theorem!four-vertex theorem@four-vertex $\sim $} (in German: {\em Vierscheitelsatz}\index{Vierscheitelsatz}) was (probably) first proved by \textcite{Mukhopadhyaya} (see \cite[p.\ 18]{Blaschke:Differentialgeometrie}). 

\begin{thm} \label{Thm:Four-vertex_theorem}
Let $z(t)$, $a \le t \le b$, be a strictly convex simple closed parametric curve, other than a circle, whose curvature function $\kappa(t)$ and its derivative $\kappa'(t)$ are continuous.
Then this curve has at least four vertices.
\end{thm}

The idea of the following proof (by contradiction) is due to \textsc{G.\ Herglotz}\footnote{\textsc{Gustav Herglotz}, 1881-1953} and was reported by \textsc{W.\ Blaschke}\footnote{\textsc{Wilhelm Blaschke}, 1885-1962} \cite[p.\ 18]{Blaschke:Differentialgeometrie}.

\begin{proof}
From the (Weierstrass\index{Weierstrass, Karl}) extreme value theorem\index{theorem!extreme value theorem@extreme value $\sim $} \cite{wiki:Extreme_value_theorem} we know that the real-valued continuous function $\kappa$ must attain a maximum and a minimum on the closed interval $[a,b]$, each at least once.
Since there is always a minimum between two maxima, the number of extrema is even.

We assume that the curve $z$ has exactly two vertices $z_1$ and $z_2$, and consider the integral
\beq
  I
:= \oint \left[z(t)-z_1, z_2-z_1\right] \kappa'(t)\, \dd t    
\eeq
along $z$.
The straight line with equation
\beq
  \left[z-z_1, z_2-z_1\right]
= 0
\eeq
through $z_1$ and $z_2$ divides $z$ into two disjoint parts.
In one of the parts is $\kappa'(t) > 0$, while in the other is $\kappa'(t) < 0$.
If $z(t)$ is on the right side of the vector $z_2-z_1$ fixed at point $z_1$, then $\left[z(t)-z_1, z_2-z_1\right] > 0$; if $z(t)$ is on the left side of the said vector, then $\left[z(t)-z_1, z_2-z_1\right] < 0$ (see \eqref{Eq:QVP_with_sine} and Fig.\ \ref{Abb:Parallelogramm}).
The integrand therefore has no change of sign, hence $I \ne 0$.
On the other hand, with \eqref{Eq:Rechenregeln1} we have
\begin{align*}
  I
= {} & \oint \big\langle \ii\,(z(t)-z_1), z_2-z_1 \big\rangle\, \kappa'(t)\, \dd t\db\\[0.05cm]
= {} & \oint \big\langle \ii z(t), z_2-z_1 \big\rangle\, \kappa'(t)\, \dd t
		- \oint \big\langle \ii z_1, z_2-z_1 \big\rangle\, \kappa'(t)\, \dd t\db\\[0.05cm]
= {} & \oint \big\langle \ii z(t), z_2-z_1 \big\rangle\, \kappa'(t)\, \dd t
		- \big\langle \ii z_1, z_2-z_1 \big\rangle \oint \kappa'(t)\, \dd t\db\\[0.05cm]
= {} & \oint \big\langle \ii z(t), z_2-z_1 \big\rangle\, \kappa'(t)\, \dd t
		- \underbrace{\big\langle \ii z_1, z_2-z_1 \big\rangle\, \kappa(t)\,\big|_{t\,=\,a}^b}_{0}\,.		
\end{align*}
Partial integration yields
\beq
  I
= \underbrace{\big\langle \ii z(t), z_2-z_1 \big\rangle\, \kappa(t)\,\big|_{t\,=\,a}^b}_0\,
		- \oint \big\langle \ii z'(t), z_2-z_1 \big\rangle\, \kappa(t)\, \dd t\,,
\eeq
where we have applied \eqref{Eq:Product_rule_SP}.
Using the first formula for the scalar product in \eqref{Eq:Rechenregeln2} twice, we get
\beq
  I
= \oint \big\langle \ii \kappa(t) z'(t), z_1-z_2 \big\rangle\,	\dd t
= \oint \left\langle \ii \kappa(t) \left|z'(t)\right| T(t), z_1-z_2 \right\rangle \dd t\,,	
\eeq
where $T(t)$ is the tangent unit vector.
Applying Lemma \ref{Lem:T'(t)} gives
\beq
  I
= \oint \left\langle T'(t), z_1-z_2 \right\rangle \dd t
= \left\langle T(t), z_1-z_2 \right\rangle\big|_{t\,=\,a}^b
= 0\,.	
\eeq
This is a contradiction to $I \ne 0$, therefore the curve $z$ has at least four vertices.
Examples for curves in the sense of Theorem \ref{Thm:Four-vertex_theorem} with the smallest number of four vertices are ellipses, see e.g.\ Fig.\ \ref{Abb:Ellipse_mit_Evolute}, and the curve in Fig.\ \ref{Abb:Eilinie}.
\end{proof}

\textsc{Blaschke}'s proof \cite{Blaschke:Vierscheitelsatz}, \cite[pp.\ 160-161]{Blaschke:Kreis_und_Kugel} of the four-vertex theorem, which he came up with through \textsc{Carath\'eodory},\footnote{\textsc{Constantin Carath\'eodory}, 1873-1950} is also very interesting.
It is also a proof by contradiction, whereby the curvature radius function of the curve is applied as mass density to the unit circle.
It is shown that the center of gravity of such a circle lies at its center, but this is not the case for the circle associated with a curve with only two vertices.
Three further proofs can be found in \cite[pp.\ 171-175]{Fuchs&Tabachnikov}.
\textsc{Kneser} \cite{Kneser} already in 1912 proved that the four-vertex theorem even holds true for all closed non-self-intersecting plane curves having continuous curvature function.
In this context, he first proved the result of Theorem \ref{Thm:Tait-Kneser} and then assumed that the curve has exactly one minimum and exactly one maximum of the curvature function. Next, he projected the curve and its circles of curvature stereographically onto an arbitrary sphere and deduced a contradiction using a theorem by \textsc{Möbius}.\footnote{\textsc{August Ferdinand Möbius}, 1790-1868}

\textsc{Herman Gluck} proved the converse of the four-vertex theorem (in a more general context): {\em Every continuous strictly positive function on the unit circle which has at least two maxima and two minima is the curvature function of a plane simple closed curve} \cite{Gluck}, \cite{DeTurck&Gluck&Pomerleano&Vick}.
\textsc{Björn Dahlberg} proved this converse without the restriction of convexity \cite{DeTurck&Gluck&Pomerleano&Vick}.

\subsection{Envelopes and support functions}

The following theorem gives a necessary condition for a two-parameter family $z(t,\la)$ of plane curves to have an envelope\index{envelope}.

\begin{thm} \label{Thm:Envelope}
For a family of plane curves given by the complex-valued parameter representation $z = z(t,\la) \in C^1$ with real parameters $t$ and $\la$, the necessary condition for the existence of an envelope is
\beqn \label{Eq:Huellkurvengleichung1}
  \left[\frac{\p z}{\p t}, \frac{\p z}{\p \la}\right]
= 0\,.
\eeqn
\end{thm}

\begin{proof}
Let us suppose that $t$ is the parameter along the curve, and $\la$ is the parameter that determines one curve of the family.
Now, we use the fact that a point of the envelope is the intersection of two closely adjacent curves $z(t,\la)$ and $z(t, \la + \Delta\la)$ with $\Delta\la \rightarrow 0$.
If such an intersection point exists, then it can be found with
\beq
  \left.\begin{aligned}
  z(t + \Delta t, \la + \Delta\la) = z(t,\la)\,,\\
  \bar{z}(t + \Delta t, \la + \Delta\la) = \bar{z}(t,\la)\,.
\end{aligned}\,\right\}
\eeq
Taylor series expansion up to first order terms yields
\beq
\left.\begin{aligned}
  z(t,\la) + \frac{\p z(t,\la)}{\p t}\,\Delta t + \frac{\p z(t,\la)}{\p\la}\,\Delta\la = z(t,\la)\,,\\[0.05cm] 
  \bar{z}(t,\la) + \frac{\p\bar{z}(t,\la)}{\p t}\,\Delta t + \frac{\p\bar{z}(t,\la)}{\p\la}\,\Delta\la = \bar{z}(t,\la)\,,
\end{aligned}\;\right\}
\eeq
hence the homogeneous linear system of equations
\beqn \label{Eq:Homog_lin_GS}
\left.\begin{aligned}
  z_t(t,\la)\,\Delta t + z_\la(t,\la)\,\Delta\la = 0\,,\\ 
  \bar{z}_t(t,\la)\,\Delta t + \bar{z}_\la(t,\la)\,\Delta\la = 0\,\phantom{,}
\end{aligned}\;\right\}
\eeqn
for $\Delta t$ and $\Delta\la$, where
\beq
  z_t \equiv \frac{\p z}{\p t}
  \quad\mbox{and}\quad
  z_\la \equiv \frac{\p z}{\p \la}\,.
\eeq
This system has a non-trivial solution only if the coefficient determinant (Jacobian determinant)\index{Jacobian determinant} vanishes,
\beq
  \left|\begin{array}{cc}
	z_t(t,\la) & z_\la(t,\la)\\[0.1cm]
	\bar{z}_t(t,\la) & \bar{z}_\la(t,\la)
  \end{array}\right| = 0\,,	  
\eeq
hence
\beq
  \bar{z}_t(t,\la)\,z_\la(t,\la) - z_t(t,\la)\,\bar{z}_\la(t,\la)
= 0\,.  
\eeq
Apart from the irrelevant constant factor $1/(2\ii)$ (see \eqref{Eq:Def_QVP_complex}), this is Equation \eqref{Eq:Huellkurvengleichung1}.
\end{proof}

Eq.\ \eqref{Eq:Huellkurvengleichung1} provides a functional relation between $t$ and $\la$ and thus the equation of the envelope.
Furthermore, Eq.\ \eqref{Eq:Huellkurvengleichung1} states that the two vectors $z_t$ and $z_\la$ are collinear at an intersection point ($=$ point of the envelope) (see \eqref{Eq:QVP_with_sine} and Fig.\ \ref{Abb:Parallelogramm}).\footnote{One obtains the result of Theorem \ref{Thm:Envelope} directly from the observation that the collinearity of the vectors $z_t$ and $z_\la$ is the necessary condition for a point to be on the envelope (see \cite{Encyclopedia:Envelope} which was adapted from an original arcticle by \foreignlanguage{russian}{Viktor Abramovich Zalgaller}, 1920-2020). However, although our proof of Theorem \ref{Thm:Envelope} is longer, the starting point with the intersection of adjacent curves seems more accessible to us than that with the collinear vectors. Concerning the determination of the envelope for a family of parametric functions (curves) see also \cite{Icaro_Lorran&Blatter}.}

\begin{thm} \label{Thm:Envelope=evolute}
The envelope\index{envelope} of the family of normal lines to a curve is also the evolute of this curve. 
\end{thm}

\begin{proof}
Let the curve be given in parametric representation by $z(t) = x(t) + \ii y(t)$.
Then, $z'(t) = x'(t) + \ii y'(t)$ is the non-normalized tangent vector, and $\ii z'(t) = \ii x'(t) - y'(t)$ the non-normalized normal vector.
Thus, the normal at point $t$ is given by
\beqn \label{Eq:Normal_to_curve}
  z_n(t,\la)
= z(t) + \ii\la z'(t)\,,\quad
  -\infty < \la < \infty\,.  
\eeqn
Applying Theorem \ref{Thm:Envelope} to \eqref{Eq:Normal_to_curve}, using (\ref{Eq:QVP_rules}\,b), \eqref{Eq:QVP_with_sine}, the first rule for the quasi vector product in \eqref{Eq:Rechenregeln2}, the rotation rule for the quasi vector product in \eqref{Eq:Rechenregeln2}, and (\ref{Eq:QVP_rules}\,a), gives
\begin{align*}
  0
= {} & \left[\frac{\p z_n(t,\la)}{\p t},\, \frac{\p z_n(t,\la)}{\p \la}\right]\db\\[0.05cm]
= {} & \left[z'(t) + \ii\la z''(t),\, \ii z'(t)\right]\db\\[0.05cm]  
= {} & \left[z'(t), \ii z'(t)\right] + \left[\ii\la z''(t), \ii z'(t)\right]\db\\[0.05cm]
= {} & \left|z'(t)\right|^2 + \la \left[z''(t), z'(t)\right]\db\\[0.05cm]
= {} & \left|z'(t)\right|^2 - \la \left[z'(t), z''(t)\right],
\end{align*}
hence
\beqn \label{Eq:lambda_for_normals}
  \la
= \la(t)
= \frac{\left|z'(t)\right|^2}{\left[z'(t), z''(t)\right]}\,.   
\eeqn
Plugging \eqref{Eq:lambda_for_normals} into \eqref{Eq:Normal_to_curve} yields
\begin{align}
  z_n(t,\la(t))
= {} & z(t) + \ii\, \frac{\left|z'(t)\right|^2}{\left[z'(t), z''(t)\right]}\, z'(t)\nonumber\db\\[0.05cm]   
= {} & z(t) + \ii\, \frac{\left|z'(t)\right|^3}{\left[z'(t), z''(t)\right]}\, \frac{z'(t)}{\left|z'(t)\right|}\,.\label{Eq:Envelope=evolute}
\end{align}
Considering \eqref{Eq:Curvature}, we see that \eqref{Eq:Envelope=evolute} is identical to \eqref{Eq:Center_of_curvature}.
\end{proof}

For examples to Theorem \ref{Thm:Envelope=evolute} see Figs.\ \ref{Abb:Ellipse_mit_Evolute}, \ref{Abb:Evolute_of_coupler_curve} and \ref{Abb:Circle_with_involutes01}.

\begin{cor} \label{Cor:Envelope_of_line_family}
Let $g_\ph$ be a family of lines defined by
\beqn \label{Eq:Family_of_lines}
  z(\ph,\la)
= a(\ph)\,\ee^{\ii\ph} + \la\,\ii\ee^{\ii\ph}\,,\quad
  -\infty < \la < \infty\,,  
\eeqn
where $a(\ph)$ is the distance of the line $g_\ph$ with $\ph = \tn{const}$ from the origin, and $\ph$ the angle between the positive $x$-axis and the normal of this line (cf.\ \eqref{Eq:z(phi,a)} and {\em Fig.\ \ref{Abb:Hessesche_Normalform}}).
Then a parametric representation of the envelope of the family $g_\ph$ is given by
\beqn \label{Eq:Envelope_of_line_family}
  \tilde{z}(\ph)
= \left(a(\ph) + \ii a'(\ph)\right)\ee^{\ii\ph}\,.  
\eeqn 
\end{cor}

\begin{proof}
We have
\begin{align*}
  z_\ph(\ph,\la)
= \frac{\p z(\ph,\la)}{\p \ph}
= {} & a'(\ph)\,\ee^{\ii\ph} + \left(a(\ph) + \ii\la\right)\ii\ee^{\ii\ph}
= \left(a'(\ph) - \la + \ii a(\ph)\right)\ee^{\ii\ph}\,,\db\\[0.05cm]
  z_\la(\ph,\la)
= \frac{\p z(\ph,\la)}{\p \la}
= {} & \ii\ee^{\ii\ph}\,.  
\end{align*}
Using Theorem \ref{Thm:Envelope}, the rotation rule in \eqref{Eq:Rechenregeln2}, and the first formula in \eqref{Eq:[z_1,z_2]_mit_Im}, we find
\begin{align*}
  0 
= {} & \left[z_\ph(\ph,\la),z_\la(\ph,\la)\right]\db\\[0.05cm]
= {} & \left[\left(a'(\ph) - \la + \ii a(\ph)\right)\ee^{\ii\ph},\,\ii\ee^{\ii\ph}\right]\db\\[0.05cm]  
= {} & \left[a'(\ph) - \la + \ii a(\ph),\,\ii\right]\db\\[0.05cm]
= {} & \Imz\bigl\{\left(a'(\ph) - \la - \ii a(\ph)\right)\ii\bigr\}\db\\[0.05cm]
= {} & a'(\ph) - \la\,,
\end{align*}
hence
\beq
  \la = \la(\ph) = a'(\ph)\,,
\eeq
and
\beq
  \tilde{z}(\ph)
= z(\ph,\la(\ph))
= \left(a(\ph) + \ii a'(\ph)\right)\ee^{\ii\ph}\,. \qedhere  
\eeq  
\end{proof}

From \eqref{Eq:Envelope_of_line_family} one easily gets the real parametric representation
\beq
\begin{array}{c@{\;=\;}*{2}{c@{\;}}l}
  \tilde{x}(\ph)
& a(\ph)\cos\ph & - & a'(\ph)\sin\ph\,,\\[0.15cm]
  \tilde{y}(\ph)
& a(\ph)\sin\ph & + & a'(\ph)\cos\ph  
\end{array}
\eeq
of the envelope of the family $g_\ph$ of lines defined by \eqref{Eq:Family_of_lines} (cf.\ \cite[pp.\ 2-3]{Santalo}).

\begin{example}
Let $\sigma$ be a line segment (in position $\sigma_\ph$) (see Fig.\ \ref{Abb:Enveloppe_Astroide01}) of constant length $\ell$ whose endpoints $z_A(\ph)$ and $z_B(\ph)$ are sliding along the $x$-axis and $y$-axis, respectively.
We will now determine the envelope of the resulting family of line segments $\sigma_\ph$,  $0 \le \ph < 2\pi$ (see Fig.\ \ref{Abb:Enveloppe_Astroide02}) and thus of the family of straight lines $g_\ph$ with $\sigma_\ph \subset g_\ph$.
Using \eqref{Eq:Line_equation_with_projection} (see also Fig.\ \ref{Abb:Hessesche_Normalform}) and $\sin\ph = x_A(\ph)/\ell$ (see Fig.\ \ref{Abb:Enveloppe_Astroide01}), for the distance $a = a(\ph)$ of $g_\ph$ from the coordinate origin we have
\beq
  a(\ph)
= \left\langle z_A(\ph),\,\ee^{\ii\ph} \right\rangle
= \left\langle x_A(\ph),\,\ee^{\ii\ph} \right\rangle
= \bigl\langle \ell\sin\ph,\, \cos\ph+\ii\sin\ph \bigr\rangle
= \ell\sin\ph\cos\ph
= \tfrac{1}{2}\,\ell\sin(2\ph)\,.\footnote{Note that $a(\ph) < 0$ if $\pi/2 < \ph < \pi$ or $3\pi/2 < \ph < 2\pi$.}
\eeq

\begin{figure}[ht]
\begin{minipage}{0.46\textwidth}
  \centering
  \includegraphics[width=0.7\textwidth]{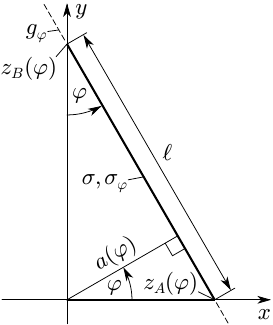}
  \caption{Line segment $\sigma$}
  \label{Abb:Enveloppe_Astroide01}
\end{minipage}
\hfill
\begin{minipage}{0.54\textwidth}
  \includegraphics[width=\textwidth]{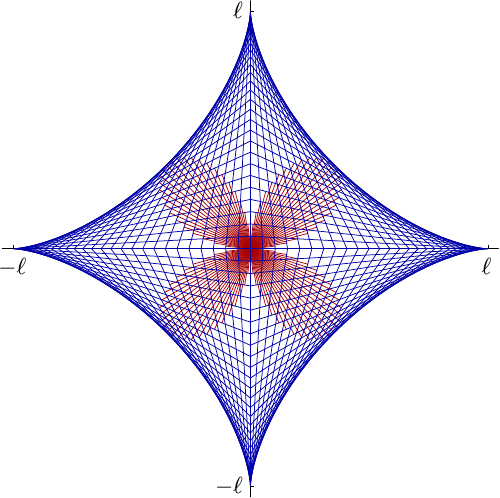}
  \caption{Families of line segments $\sigma_\ph$ (blue) and $a(\ph)\,\ee^{\ii\ph}$ (red)}
  \label{Abb:Enveloppe_Astroide02}
\end{minipage}
\end{figure}

With $a'(\ph) = \ell\cos(2\ph)$, according to \eqref{Eq:Envelope_of_line_family} in Corollary \ref{Cor:Envelope_of_line_family}, the envelope of $g_\ph$ is given by
\beq
  \tilde{z}(\ph)
= \ell \left(\tfrac{1}{2}\sin(2\ph) + \ii\cos(2\ph)\right) \ee^{\ii\ph}  
= \ell \left(\sin^3\ph + \ii\cos^3\ph\right).
\eeq
This is a parametric equation of the {\em astroid}\index{astroid}.\footnote{For more examples of envelopes see \cite[pp.\ 49-51, 394-399]{Boettcher:Algebraische_Kurven_ohne_Index}. On pages 394-399, in particular, the representation of envelopes as {\em algebraic curves}\index{curve!algebraic curve@algebraic $\sim $} is discussed.} \hfill\bs
\end{example}

\begin{defin}
If the envelope given by \eqref{Eq:Envelope_of_line_family} is the boundary $\partial K$ of a convex set $K$ and the origin is an interior point of $K$, then the function $a(\ph)$ is called the {\em support function of $K$} or the {\em support function of the convex curve $\partial K$} with reference to the origin (cf.\ \textcite[p.\ 3]{Santalo}).
According to \textcite[p.\ 600]{Schneider&Weil2008}, the {\em support function} of a convex body (set) $K \subset \R^d$ is defined by
\beqn \label{Eq:Support_function_S&W}
  h(K,u) := \max\left\{\langle x, u\rangle\, \colon x \in K\right\}
  \quad\mbox{for}\; u \in \R^d\,.\footnote{In \cite[p.\ 183]{Schneider&Weil1992}, $u \in S^{d-1}$ is assumed, where $S^{d-1}$ is the unit sphere $S^{d-1} := \{x \in \R^d \colon \|x\| = 1\}$. This means that every vector $u$ is a unit vector.}
\eeqn  
(The value of the support function for the direction $u$ is the maximum of the scalar product\index{product!scalar product@scalar $\sim $} $\langle x, u\rangle$ for all points $x \in K$, that is the maximum projection of a vector $x \in K$ onto the direction $u$ (cf.\ Fig.\ \ref{Abb:Projection} and Eq.\ \eqref{Eq:Projection_with_inner_product}).)  
\end{defin}

\begin{example}
We determine the support function\index{support function} $\tilde{a}(\ph)$ of an ellipse\index{ellipse} $E$ with boundary\index{boundary} $\partial E$ given by the parametric representation
\beq
  z(\tau) = a\cos\tau + \ii b\sin\tau\,,\quad
  0 \le \tau \le 2\pi\,.
\eeq
From \eqref{Eq:Support_function_S&W} it follows that
\begin{align*}
  \tilde{a}(\ph)
= {} & h(E,\ee^{\ii\ph}) = h(\partial E,\ee^{\ii\ph})\\[0.05cm]
= {} & \max\left\{\left\langle z(\tau),\, \ee^{\ii\ph}\right\rangle \colon 0\le\tau\le 2\pi\right\}\db\\[0.05cm]
= {} & \max\left\{\left\langle a\cos\tau + \ii b\sin\tau,\, \ee^{\ii\ph}\right\rangle \colon 0\le\tau\le 2\pi\right\}\db\\[0.05cm]
= {} & \max\left\{a\cos\tau\cos\ph + b\sin\tau\sin\ph\: \colon 0\le\tau\le 2\pi\right\}\db\\[0.05cm]
= {} & \max\left\{a\cos\ph\cos\tau + b\sin\ph\sin\tau\: \colon 0\le\tau\le 2\pi\right\}\db\\[0.05cm]
= {} & \max\left\{\left\langle a\cos\ph + \ii b\sin\ph,\, \ee^{\ii\tau}\right\rangle \colon 0\le\tau\le 2\pi\right\}.
\end{align*}
The scalar product\index{product!scalar product@scalar $\sim $} $\left\langle a\cos\ph + \ii b\sin\ph, \ee^{\ii\tau}\right\rangle$ is the projection of the vector $a\cos\ph + \ii b\sin\ph$ onto the direction $\ee^{\ii\tau}$.
The maximum projection for $0 \le \tau \le 2\pi$ is the length of the vector $a\cos\ph + \ii b\sin\ph$, thus
\beq
  \tilde{a}(\ph)
= \left|a\cos\ph + \ii b\sin\ph\right|  
= \sqrt{a^2\cos^2\ph+b^2\sin^2\ph}\;.   
\eeq
Here we have essentially followed the argumentation of user147263 from June 5, 2016 in \cite{J.Doe&user147263&gIS}. \hfill\bs
\end{example}

%% file: DiffGeo5_4a.tex

\section{Applications in plane kinematics} \label{Sec:Plane_kinematics}

\subsection{Kinematics of a dyad\index{dyad}} \label{Subsec:Kinematics_of_a_dyad}

\begin{SCfigure}[0.7][ht]
  \includegraphics[width=0.3\textwidth]{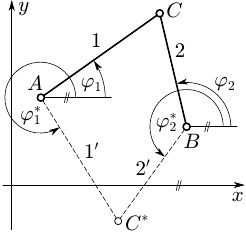}
  \caption{Dyad}
  \label{Abb:Dyad}
\end{SCfigure}

We consider two in the plane movable bars (line segments) 1 and 2 of length $\ell_1$ or $\ell_2$ which are connected by a pivot joint at point $C$ (see Fig.\ \ref{Abb:Dyad}).
Such an assembly is referred to as a {\em dyad}\index{dyad}.
We will determine the angle $\ph_1$ for given points $A$ and $B$.

It holds
\beq
  z_C = z_A + \ell_1\,\ee^{\ii\ph_1}
  \quad\mbox{and}\quad
  z_C = z_B + \ell_2\,\ee^{\ii\ph_2}\,.
\eeq
Equating gives
\beq
  z_A + \ell_1\,\ee^{\ii\ph_1}
= z_B + \ell_2\,\ee^{\ii\ph_2}\,,
\eeq
hence
\beq
  \ell_2\,\ee^{\ii\ph_2}
= \ell_1\,\ee^{\ii\ph_1} - z_{AB}
  \quad\mbox{with}\quad
  z_{AB}
:= z_B - z_A\,.  
\eeq
In order to eliminate $\ph_2$, we use the complex conjugate
\beq
  \ell_2\,\ee^{-\ii\ph_2}
= \ell_1\,\ee^{-\ii\ph_1} - \overline{z_{AB}}\,.
\eeq
With this we get
\beq
  \ell_2^2\,\ee^{\ii\ph_2}\,\ee^{-\ii\ph_2}
= \left(\ell_1\,\ee^{\ii\ph_1} - z_{AB}\right) \left(\ell_1\,\ee^{-\ii\ph_1} - \overline{z_{AB}}\right),
\eeq
hence
\beqn \label{Eq:Dyade_Uebertragungsgleichung_1}
  \ell_1\,\overline{z_{AB}}\,\ee^{\ii\ph_1} + \ell_1\,z_{AB}\,\ee^{-\ii\ph_1} + \ell_2^2 - \ell_1^2 - |z_{AB}|^2
= 0\,.
\eeqn
Putting
\beq
  h
:= 2 \ell_1 z_{AB}
  \quad\mbox{and}\quad
  c
:= \ell_2^2 - \ell_1^2 - |z_{AB}|^2
= \ell_2^2 - \ell_1^2 - \frac{|h|^2}{4\ell_1^2}\,,  
\eeq
we can write \eqref{Eq:Dyade_Uebertragungsgleichung_1} as
\beqn \label{Eq:Dyade_Uebertragungsgleichung_2}
  \overline{h}\,\ee^{\ii\ph_1} + h\,\ee^{-\ii\ph_1} + 2c
= 0\,.
\eeqn
Complex multiplication of \eqref{Eq:Dyade_Uebertragungsgleichung_2} with $\ee^{\ii\ph_1}$ yields the quadratic equation
\beq
  \overline{h}\,\ee^{2\ii\ph_1} + 2c\,\ee^{\ii\ph_1} + h
= 0  
\eeq
for $\ee^{\ii\ph_1}$ with solutions
\beqn \label{Eq:e^{i phi_1}}
  \ee^{\ii\ph_1}
= \frac{-c \pm \ii\, \sqrt{|h|^2 - c^2}}{\overline{h}}\,,
\eeqn
and therefore
\beqn \label{Eq:dyad_phi_1}
  \ph_1
= \arg\left(\ee^{\ii\ph_1}\right)
= \arg\left(\frac{-c \pm \ii\, \sqrt{|h|^2 - c^2}}{\overline{h}}\right).
\eeqn
The $\pm$-sign in  \eqref{Eq:e^{i phi_1}} and \eqref{Eq:dyad_phi_1} corresponds to the two assembly variants $ACB$ and $AC^*B$ of the dyad (see Fig.\ \ref{Abb:Dyad}). 
From \eqref{Eq:e^{i phi_1}} one easily finds
\beqn \label{Eq:Re_and_Im_of_exp(i*alpha)}
\left.
\begin{aligned}
  \Rez\big(\ee^{\ii\ph_1}\big)
= {} & \cos\ph_1
= \frac{-a\,c \mp b\,\sqrt{a^2 + b^2 - c^2}}{a^2 + b^2}\,,\\[0.15cm]
  \Imz\big(\ee^{\ii\ph_1}\big)
= {} & \sin\ph_1
= \frac{-b\,c \pm a\,\sqrt{a^2 + b^2 - c^2}}{a^2 + b^2} 
\end{aligned}
\;\right\}
\eeqn
with
\beq
  a := \Rez(h) = 2\ell_1\left(x_B-x_A\right)\,,\quad   
  b := \Imz(h) = 2\ell_1\left(y_B-y_A\right)\,.
\eeq
Using \eqref{Eq:Def_SP_complex}, from \eqref{Eq:Dyade_Uebertragungsgleichung_2} follows
\beqn \label{Eq:Dyade_Uebertragungsgleichung_3}
  \left\langle h, \ee^{\ii\ph_1} \right\rangle + c
= 0  
\eeqn
or
\beq
  \left\langle a + \ii b, \ee^{\ii\ph_1} \right\rangle + c
= 0\,,  
\eeq
thus, using \eqref{Eq:Def_SP_real},
\beq
 a\cos\ph_1 + b\sin\ph_1 + c
= 0\,.
\eeq
From this it is possible to obtain a quadratic equation for $\tan(\alpha/2)$, and hence for $\alpha$ (see \cite[p.\ 99]{Luck&Modler}).

Now, we assume $z_A$ and $z_B$ to be (complex-valued) functions of a real parameter $\ph$, hence
\beq
  h(\ph) = 2 \ell_1 \left(z_B(\ph) - z_A(\ph)\right)
  \quad\mbox{and}\quad
  c(\ph) = \ell_2^2 - \ell_1^2 - \frac{|h(\ph)|^2}{4\ell_1^2}\,.
\eeq
The angle $\ph_1$ is now a real-valued function of $\ph$.
In order to obtain the derivative $\ph_1'(\ph)$, we differentiate \eqref{Eq:Dyade_Uebertragungsgleichung_3} implicitly,
\begin{align*}
  0
= {} & \left\langle h'(\ph),\ee^{\ii\ph_1(\ph)}\right\rangle + \left\langle h(\ph),\ii\ph_1'(\ph)\ee^{\ii\ph_1(\ph)}\right\rangle + c'(\ph)\db\\[0.05cm]
= {} & \left\langle h'(\ph),\ee^{\ii\ph_1(\ph)}\right\rangle + \ph_1'(\ph)\left\langle h(\ph),\ii\ee^{\ii\ph_1(\ph)}\right\rangle + c'(\ph)\db\\[0.05cm]
= {} & \left\langle h'(\ph),\ee^{\ii\ph_1(\ph)}\right\rangle - \ph_1'(\ph)\left[h(\ph),\ee^{\ii\ph_1(\ph)}\right] + c'(\ph)\,,
\end{align*}
and get
\beqn \label{Eq:phi_1'}
  \ph_1'(\ph)
= \frac{\left\langle h'(\ph),\ee^{\ii\ph_1(\ph)}\right\rangle + c'(\ph)}{\left[h(\ph),\ee^{\ii\ph_1(\ph)}\right]}  
\eeqn
with $\ee^{\ii\ph_1(\ph)}$ from \eqref{Eq:e^{i phi_1}}, and
\beq
  h'(\ph) = 2\ell_1\left(z_B'(\ph)-z_A'(\ph)\right),\quad
  c'(\ph) = -\frac{\left\langle h(\ph),h'(\ph)\right\rangle}{2 l_1^2}\,.
\eeq
Analogously, we obtain the second derivative
\beqn \label{Eq:phi_1''}
  \ph_1''
= \frac{\left\langle h''-\ph_1'^2 h,\ee^{\ii\ph_1}\right\rangle - 2\ph_1'\left[h',\ee^{\ii\ph_1}\right] + c''}{\left[h,\ee^{\ii\ph_1}\right]}\,,  
\eeqn
where
\beq
  h'' = 2\ell_1(z_B''-z_A'')\,,\quad
  c'' = -\frac{\langle h', h'\rangle + \langle h, h''\rangle}{2 l_1^2}\,,      
\eeq
and the third derivative  
\beq
  \ph_1'''
= \frac{\left\langle h'''-3\ph_1'^2h'-3\ph_1'\ph_1''h,\ee^{\ii\ph_1}\right\rangle - \left[3\ph_1'h''+3\ph_1''h'-\ph_1'^3h,\ee^{\ii\ph_1}\right]
  + c'''}{\left[h,\ee^{\ii\ph_1}\right]}\,,  
\eeq
where
\beq
  h''' = 2\ell_1(z_B'''-z_A''')\,,\quad
  c''' = -\frac{3\langle h', h''\rangle + \langle h, h'''\rangle}{2 l_1^2}\,.      
\eeq

The following subroutine\index{subroutine} {\em DyadC} programmed in {\em Mathematica}\index{Mathematica} returns upon input of
\begin{gather*}
  \texttt{l1} = \ell_1\,,\quad
  \texttt{l2} = \ell_2\,,\quad
  \texttt{zA} = z_A(\ph)\,,\quad
  \texttt{zA1} = z_A'(\ph)\,,\quad
  \texttt{zA2} = z_A''(\ph)\,,\quad
  \texttt{zA3} = z_A'''(\ph)\,,\\
  \texttt{zB} = z_B(\ph)\,,\quad
  \texttt{zB1} = z_B'(\ph)\,,\quad
  \texttt{zB2} = z_B''(\ph)\,,\quad
  \texttt{zB3} = z_B'''(\ph)\,,\\
  \texttt{vz} = \pm 1 = \mbox{sign for assembly variant}
\end{gather*}
a list \texttt{\{exp, phi1, phi11, phi12, phi13\}} containing the vector $\ee^{\ii\ph_1(\ph)}$, and the angle $\ph_1(\ph)$ (in radians) with its first three derivatives $\ph_1'(\ph)$, $\ph_1''(\ph)$, $\ph_1'''(\ph)$: 

\label{Subroutine-code}
\begin{verbatim}
DyadC[l1_, l2_, zA_, zA1_, zA2_, zA3_, zB_, zB1_, zB2_, zB3_, vz_] :=
  
  Module[{h, h1, h2, h3, c, c1, c2, c3, exp, phi1, phi11, phi12, 
    phi13},
   h = 2*l1*(zB - zA); h1 = 2*l1*(zB1 - zA1); h2 = 2*l1*(zB2 - zA2);
   h3 = 2*l1*(zB3 - zA3); c = l2^2 - l1^2 - Abs[h]^2/(4*l1^2);
   c1 = -SP[h, h1]/(2*l1^2); c2 = -(SP[h1, h1] + SP[h, h2])/(2*l1^2);
   c3 = -(3*SP[h1, h2] + SP[h, h3])/(2*l1^2);
   exp = N[(-c + vz*I*(Abs[h]^2 - c^2)^(1/2))/h\[Conjugate]];
   phi1 = Arg[exp];
   phi11 = (SP[h1, exp] + c1)/QVP[h, exp];
   phi12 = (SP[h2 - phi11^2*h, exp] - 2*phi11*QVP[h1, exp] + c2)/
     QVP[h, exp];
   phi13 = (SP[h3 - 3*phi11^2*h1 - 3*phi11*phi12*h, exp] - 
       QVP[3*phi11*h2 + 3*phi12*h1 - phi11^3*h, exp] + c3)/QVP[h, exp];
   {exp, phi1, phi11, phi12, phi13}];
\end{verbatim}

It is $\texttt{I} \equiv \ii$.
\texttt{SP} and \texttt{QVP} denote the scalar product or quasi vector product defined by
\begin{gather*}
  \texttt{SP[z1\_, z2\_] := Re[z1]*Re[z2] + Im[z1]*Im[z2];}\\
  \texttt{QVP[z1\_, z2\_] := Re[z1]*Im[z2] - Im[z1]*Re[z2];}
\end{gather*}

The subroutine is immediately executable after copying it into {\em Mathematica}\index{Mathematica}.
If the angle $\ph_2$ or $\ph_2^*$ with unit vector and derivatives is needed, one has to call the subroutine with
\beq
  \texttt{DyadC[l2, l1, zB, zB1, zB2, zB3, zA, zA1, zA2, zA3, -vz]}
\eeq
or
\beq
  \texttt{DyadC[l2, l1, zB, zB1, zB2, zB3, zA, zA1, zA2, zA3, vz]}\,.
\eeq
Additionally, we note that the expressions \eqref{Eq:phi_1'} and \eqref{Eq:phi_1''} can respectively be written as 
\beq
  \ph_1'
= \frac{a'\cos\ph_1 + b'\sin\ph_1 + c'}{a\sin\ph_1 - b\cos\ph_1}
\eeq
and
\beq
  \ph_1''
= \frac{\left(a''+2\ph_1'b'-\ph_1'^2 a\right)\cos\ph_1 + \left(b''-2\ph_1'a'-\ph_1'^2 b\right)\sin\ph_1 + c''}{a\sin\ph_1 - b\cos\ph_1}\,.  
\eeq

An alternative way to compute $\ee^{\ii\ph_1}$ is as follows:
We have (see Fig.\ \ref{Abb:Dyad})
\beqn \label{Eq:e^{i phi_1}_with_e^{i alpha}}
  \ee^{\ii\ph_1}
= \frac{z_{AC}}{|z_{AC}|}
= \frac{z_{AC}}{\ell_1}
= \frac{z_{AB}}{|z_{AB}|}\,\ee^{\ii\alpha}\,,
  \quad\mbox{where}\quad
  \alpha
:= \measuredangle\left(z_{AB},z_{AC}\right).    
\eeqn
Using $\cos\alpha = (\ee^{\ii\alpha}+\ee^{-\ii\alpha})/2$, from \eqref{Eq:Law_of_cosines_a} we get
\beq
  \ell_2^2
= \ell_1^2 + |z_{AB}|^2 - \ell_1\,|z_{AB}|\left(\ee^{\ii\alpha}+\ee^{-\ii\alpha}\right),
\eeq
and, after multiplication by $\ee^{\ii\alpha}$ and division by $\ell_1|z_{AB}|$, the quadratic equation
\beq
  \ee^{2\ii\alpha} - \frac{\ell_1^2-\ell_2^2+|z_{AB}|^2}{\ell_1|z_{AB}|}\,\ee^{\ii\alpha} + 1
= 0  
\eeq
with solutions
\beq
  \ee^{\ii\alpha}
= \frac{\ell_1^2-\ell_2^2+|z_{AB}|^2}{2\ell_1|z_{AB}|} \pm \sqrt{\left(\frac{\ell_1^2-\ell_2^2+|z_{AB}|^2}{2\ell_1|z_{AB}|}\right)^2 - 1}\,. 
\eeq
From \eqref{Eq:Law_of_cosines_a} it is clear that
\beq
  -1 \le \frac{\ell_1^2-\ell_2^2+|z_{AB}|^2}{2\ell_1|z_{AB}|} \le 1\,,
\eeq
and therefore the square root in
\beq
  \ee^{\ii\alpha}
= \frac{\ell_1^2-\ell_2^2+|z_{AB}|^2}{2\ell_1|z_{AB}|} \pm \ii\,\sqrt{1 - \left(\frac{\ell_1^2-\ell_2^2+|z_{AB}|^2}{2\ell_1|z_{AB}|}\right)^2} 
\eeq
is real.
From \eqref{Eq:e^{i phi_1}_with_e^{i alpha}} follows
\beqn \label{Eq:e^{i phi_1}_a_la_Stefan}
  \ee^{\ii\ph_1}
= \frac{z_{AB}}{|z_{AB}|} \left(\frac{\ell_1^2-\ell_2^2+|z_{AB}|^2}{2\ell_1|z_{AB}|} 
  \pm \ii\,\sqrt{1 - \left(\frac{\ell_1^2-\ell_2^2+|z_{AB}|^2}{2\ell_1|z_{AB}|}\right)^2}\right)  
\eeqn
(cf.\ \cite{Goessner:Planar_vector_equation}).
Substituting $h$ and $c$ into \eqref{Eq:e^{i phi_1}} gives \eqref{Eq:e^{i phi_1}_a_la_Stefan}, as can be easily verified. 

%% file: DiffGeo5_4b.tex

\subsection{Five-bar linkage} \label{Subsec:Five-bar_linkage}
\subsubsection{Basics}

\begin{SCfigure}[][ht]
\includegraphics[width=0.45\textwidth]{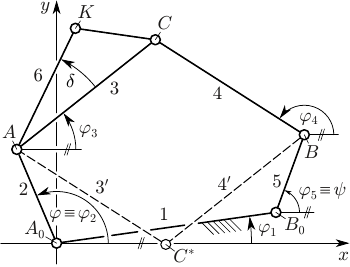}
\caption{Five-bar linkage}
\label{Abb:Five-bar_linkage01}
\end{SCfigure}

We consider the five-bar linkage\index{five-bar linkage}\index{linkage!five-bar linkage@five-bar $\sim $} in Fig.\ \ref{Abb:Five-bar_linkage01}.
It consists of five bars $1,2,3,4,5$, connected with five rotating joints at points $A_0$, $A$, $C$, $B$, $B_0$.
Bar 1 is the fixed bar (rigidly connected with the coordinate system).
The bars 2 and 5 can only rotate around their joints $A_0$ or $B_0$ in bar 1 while bars 3 and 4 (called {\em couplers}) can perform more general motions.
Bars 2 and 5 are called {\em cranks} if they perform complete revolutions. 
The five-bar linkage has degree of freedom\index{degree of freedom} two.
This means that we would have to specify two parameters (angles), e.g.\ the angles $\ph_2$ and $\ph_5$, in order to get a constrained motion.
We reduce the number of these parameters to one by assuming that the angle $\psi\equiv\ph_5$ is a $2\pi$-periodic function of the angle $\ph\equiv\ph_2$, $\psi = \psi(\ph)$.
The point $K$ is rigidly connected to bar 3.
During the motion of the mechanism, the point $K$ describes a trajectory curve $\mathcal{K}$ (not shown in Fig.\ \ref{Abb:Five-bar_linkage01}), which is called a {\em coupler curve}\index{coupler curve}. 

\subsubsection{Parametric equation (with derivatives) of the coupler curve}

The parametric equation, with parameter $\ph$, of the coupler curve $\mathcal{K}$ described by the point $K$ (see Fig.\ \ref{Abb:Five-bar_linkage01}) is simply given by
\beqn \label{Eq:z_K(phi)_5-bar_linkage}
  z_K(\ph)
= \ell_2\,\ee^{\ii\ph} + \ell_6\,\ee^{\ii(\ph_3(\ph)+\delta)}  
= \ell_2\,\ee^{\ii\ph} + \ell_6\,\ee^{\ii\delta}\,\ee^{\ii\ph_3(\ph)}\,,\quad
  0 \le \ph \le 2\pi\,.  
\eeqn
From \eqref{Eq:z_K(phi)_5-bar_linkage} one easily gets the tangent vector of the coupler curve $\mathcal{K}$ at point $K$, 
\beqn \label{Eq:z_K'(phi)}
  z_K'(\ph)
= \ell_2\,\ii\ee^{\ii\ph} + \ell_6\,\ph_3'(\ph)\,\ii\ee^{\ii\delta}\ee^{\ii\ph_3(\ph)}\,,   
\eeqn
and the derivative of the tangent vector $z_K'(\ph)$,
\beq
  z_K''(\ph)
= -\ell_2\,\ee^{\ii\ph} - \ell_6\left(\ph_3'(\ph)^2-\ii\ph_3''(\ph)\right)\ee^{\ii\delta}\ee^{\ii\ph_3(\ph)}\,.  
\eeq
$\ee^{\ii\ph_3(\ph)}$, $\ph_3'(\ph)$ and $\ph_3''(\ph)$ can be computed using the subroutine {\em DyadC} (see p.\ \pageref{Subroutine-code}) with input
\begin{gather*}
  \texttt{DyadC[l3, l4, zA, zA1, zA2, zA3, zB, zB1, zB2, zB3, vz]}\,\db,\\
  \texttt{l3} = \ell_3\,,\quad
  \texttt{l4} = \ell_4\,,\quad
  \texttt{zA} = \ell_2\,\ee^{\ii\ph}\,,\quad 
  \texttt{zA1} = \ell_2\,\ii\ee^{\ii\ph}\,,\quad
  \texttt{zA2} = -\ell_2\,\ee^{\ii\ph}\,,\quad 
  \texttt{zA3} = -\ell_2\,\ii\ee^{\ii\ph}\,,\db\\
  \texttt{zB} = \ell_1\,\ee^{\ii\ph_1} + \ell_5\,\ee^{\ii\psi(\ph)}\,,\quad
  \texttt{zB1} = \ell_5\,\psi'(\ph)\,\ii\ee^{\ii\psi(\ph)}\,,\quad
  \texttt{zB2} = -\ell_5\,\psi'(\ph)^2\,\ee^{\ii\psi(\ph)}\,,\db\\
  \texttt{zB3} = -\ell_5\,\psi'(\ph)^2\,\ii\ee^{\ii\psi(\ph)}\,,\quad
  \texttt{vz} = \pm 1\,.
\end{gather*}

\subsubsection{Area of the region bouded by the coupler curve}

To calculate the signed area $A$ of the region bounded by the coupler curve $\mathcal{K}$, we use formula
\beqn \label{Eq:Area_coupler_curve_1}
  A
= \frac{1}{2} \oint_\mathcal{K} \left[z_K(\ph),z_K'(\ph)\right] \dd\ph
\eeqn
(see Theorem \ref{Thm:Area_with_QVP} and Remark \ref{Rem:A}) with \eqref{Eq:z_K(phi)_5-bar_linkage} and \eqref{Eq:z_K'(phi)}.
Since $\psi(\ph)$ is assumed to be a $2\pi$-periodic function, \eqref{Eq:Area_coupler_curve_1} can be written as   
\beqn \label{Eq:Area_coupler_curve_2}
  A
= \frac{1}{2} \int_0^{2\pi} \left[z_K(\ph),z_K'(\ph)\right] \dd\ph\,.
\eeqn
For the numerical evaluation of this integral it is useful to simplify the quasi vector product.
We have
\begin{align} \label{Eq:[z_K,z_K']_5-bar_linkage}
  [z_K,z_K']
= {} & {\left[\ell_2\,\ee^{\ii\ph}+\ell_6\,\ee^{\ii(\ph_3+\delta)},\, \ell_2\,\ii\ee^{\ii\ph}+\ell_6\,\ph_3'\,\ii\ee^{\ii(\ph_3+\delta)}\right]}
	\nonumber\db\\[0.05cm]
= {} & {\left\langle\ell_2\,\ee^{\ii\ph}+\ell_6\,\ee^{\ii(\ph_3+\delta)},\, \ell_2\,\ee^{\ii\ph}+\ell_6\,\ph_3'\,\ee^{\ii(\ph_3+\delta)}\right\rangle}
	\quad\mbox{(see \eqref{Eq:Rechenregeln1}})\nonumber\db\\[0.05cm]	
= {} & \Big\langle\ell_2\,\ee^{\ii\ph},\, \ell_2\,\ee^{\ii\ph}\Big\rangle
	+ \left\langle\ell_2\,\ee^{\ii\ph},\, \ell_6\,\ph_3'\,\ee^{\ii(\ph_3+\delta)}\right\rangle
	+ \left\langle\ell_6\,\ee^{\ii(\ph_3+\delta)},\, \ell_2\,\ee^{\ii\ph}\right\rangle\nonumber\db\\
 &	+ \left\langle\ell_6\,\ee^{\ii(\ph_3+\delta)},\, \ell_6\,\ph_3'\,\ee^{\ii(\ph_3+\delta)}\right\rangle\nonumber\db\\[0.05cm]
= {} & \ell_2^2 + \ell_2\,\ell_6\left(\ph_3'+1\right)\left\langle\ee^{\ii\ph},\, \ee^{\ii(\ph_3+\delta)}\right\rangle + \ell_6^2\,\ph_3'
	\nonumber\db\\[0.05cm]
= {} & \ell_2^2 + \ell_2\,\ell_6\left(\ph_3'+1\right)\left\langle\ee^{\ii(\ph-\delta)},\, \ee^{\ii\ph_3}\right\rangle + \ell_6^2\,\ph_3'
	\quad\mbox{(see \eqref{Eq:Rechenregeln2})}\,.	
\end{align}
$\ee^{\ii\ph_3}$ and $\ph_3'$ can be computed with subroutine {\em DyadC}.
Since $\int_0^{2\pi} \ph_3'(\ph)\, \dd\ph = 0$, \eqref{Eq:Area_coupler_curve_2} with \eqref{Eq:[z_K,z_K']_5-bar_linkage} gives
\beqn \label{Eq:Area_coupler_curve_3}
  A
= \pi\ell_2^2 + \frac{1}{2}\,\ell_2\,\ell_6 \int_0^{2\pi} \left(\ph_3'(\ph)+1\right)
  \left\langle\ee^{\ii(\ph-\delta)},\, \ee^{\ii\ph_3(\ph)}\right\rangle	\dd\ph 
\eeqn
or alternatively
\beqn \label{Eq:Area_coupler_curve_4}
\begin{aligned}
  A
= {} & \pi\ell_2^2 + \frac{1}{2}\,\ell_2\,\ell_6 \int_0^{2\pi} \left(\ph_3'(\ph)+1\right)
  \big\{\!\cos(\ph-\delta)\cos\ph_3(\ph) + \sin(\ph-\delta)\sin\ph_3(\ph)\big\}\, \dd\ph\,, 
\end{aligned}
\eeqn
where
\beq
  \cos\ph_3(\ph)
= \Rez\left(\ee^{\ii\ph_3(\ph)}\right)\,,\qquad
  \sin\ph_3(\ph)
= \Imz\left(\ee^{\ii\ph_3(\ph)}\right)
\eeq
(see \eqref{Eq:Re_and_Im_of_exp(i*alpha)}).

\subsubsection{Example} \label{Sec:Example}

\begin{SCfigure}[][ht]
\includegraphics[width=0.55\textwidth]{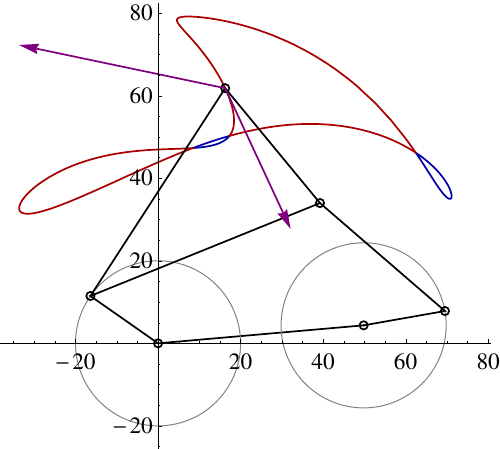}
\caption{Five-bar linkage with coupler curve $\mathcal{K}$, tangent vector $z_K'(\ph)$, and vector $z_K''(\ph)$, $\ph = 145\g$\\[0.2cm] $\ell_1 = 50$, $\ph_1 = 5\g$, $\ell_2 = 20$, $\ell_3 = 60$, $\ell_4 = 40$, $\ell_5 = 20$, sign in {\em DyadC} $= 1$, $\ell_6 = 60$, $\delta = 35\g$, $\psi(\ph) = -2\ph-\pi/3$}
\label{Abb:Five-bar_linkage02}
\end{SCfigure}

\begin{figure}[ht]
\begin{minipage}{0.48\textwidth}
  \centering
  \includegraphics[width=\textwidth]{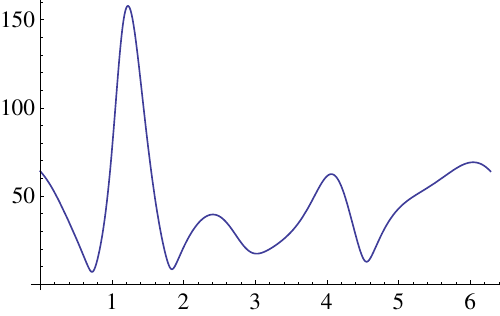}
  \caption{The function $\left|z_K'(\ph)\right|$ of the coupler curve in Fig.\ \ref{Abb:Five-bar_linkage02}}
  \label{Abb:|z_K'(phi)|}
\end{minipage}
\hfill
\begin{minipage}{0.48\textwidth}
  \includegraphics[width=\textwidth]{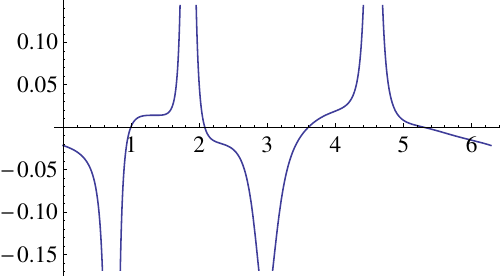}
  \caption{Curvature\index{curvature} function $\kappa(\ph)$ of the coupler curve in Fig.\ \ref{Abb:Five-bar_linkage02}}
  \label{Abb:kappa(phi)}
\end{minipage}
\end{figure}

Fig.\ \ref{Abb:Five-bar_linkage02} shows an example for a five-bar linkage with coupler curve $\mathcal{K}$, tangent vector $z_K'(\ph)$, and vector $z_K''(\ph)$.
The graph of the function $\left|z_K'(\ph)\right|$ is shown in Fig.\ \ref{Abb:|z_K'(phi)|}, and the graph of the curvature\index{curvature} function $\kappa(\ph)$ of the coupler curve in Fig.\ \ref{Abb:kappa(phi)}.
By substituting \eqref{Eq:z_K'(phi)} into \eqref{Eq:Arc_length}, numerical integration with limits $\ph = 0$ and $\ph = 2\pi$ yields the arc length $= 289.414489645\ldots$.

If the angular velocity of the crank 2 (see Fig.\ \ref{Abb:Five-bar_linkage01}) is equal to one, then $z_K'(\ph)$ is the velocity vector and $z_K''(\ph)$ the acceleration vector of point $K$.
In general, the velocity of $K$ is given by
\beq
  \dot{z}_K(\ph(t))
= \frac{\dd z_K(\ph(t))}{\dd t}
= \frac{\dd z_K(\ph(t))}{\dd \ph}\,\frac{\dd \ph}{\dd t}
= z_K'(\ph(t))\,\dot{\ph}(t)\,,
\eeq
where $t$ is the time, and $\dot{\ph}(t)$ the angular velocity of crank 2.
For the general case of the acceleration of $K$ we have
\begin{align*}
  \ddot{z}_K(\ph(t))
= {} & \frac{\dd^2 z_K(\ph(t))}{\dd t^2} 
= \frac{\dd}{\dd t}\frac{\dd z_K(\ph(t))}{\dd t}
= \frac{\dd}{\dd t}\left(z_K'(\ph(t))\,\dot{\ph}(t)\right)\db\\[0.05cm]
= {} & \frac{\dd z_K'(\ph(t))}{\dd t}\,\dot{\ph}(t) + z_K'(\ph(t))\,\frac{\dd \dot{\ph}(t)}{\dd t}
= \frac{\dd z_K'(\ph(t))}{\dd \ph}\,\frac{\dd \ph}{\dd t}\,\dot{\ph}(t) + z_K'(\ph(t))\,\ddot{\ph}(t)\db\\[0.05cm]
= {} & z_K''(\ph(t))\,\dot{\ph}(t)^2 + z_K'(\ph(t))\,\ddot{\ph}(t)\,,
\end{align*}
where $\ddot{\ph}(t)$ is the angular acceleration of crank 2.
One sees that if $\dot{\ph}(t) = \tn{const}$ and hence $\ddot{\ph}(t) = 0$, then the direction of the acceleration vector $\ddot{z}_K(\ph(t))$ is equal to the direction of the vector $z_K''(\ph(t))$.\footnote{Note that $z_K'\,\ddot{\ph}$ and $z_K''\,\dot{\ph}^2$ are not the commonly known {\em tangential acceleration}\index{acceleration!tangential acceleration@tangential $\sim $} and {\em centripetal acceleration}\index{acceleration!centripetal acceleration@centripetal $\sim $}, respectively, since $z_K''$ is in general not perpendicular to $z_K'$,}

\begin{figure}[ht]
\begin{minipage}{0.48\textwidth}
  \centering
  \includegraphics[width=\textwidth]{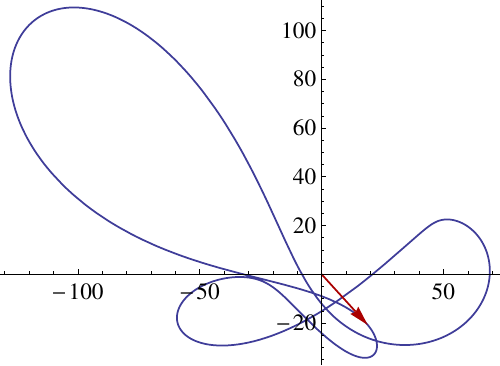}
  \caption{Hodograph (curve of the tangent vector $z_K'(\ph)$ with one end fixed) of the coupler curve in Fig.\ \ref{Abb:Five-bar_linkage02}, tangent vector for $\ph = 120\g$}
  \label{Abb:Tangent_vector_curve}
\end{minipage}
\hfill
\begin{minipage}{0.48\textwidth}
  \includegraphics[width=\textwidth]{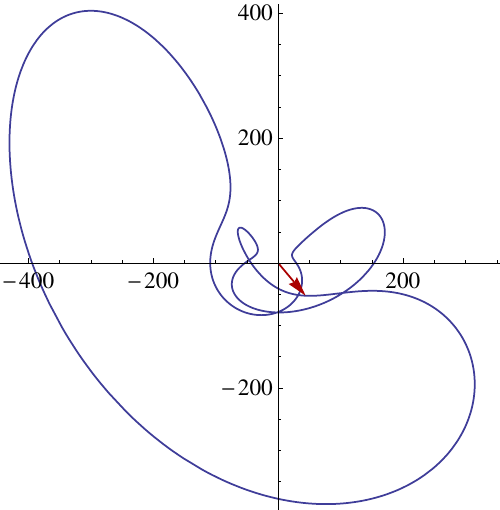}
  \caption{Curve of the vector $z_K''(\ph)$ (with one end fixed) of the coupler curve in Fig.\ \ref{Abb:Five-bar_linkage02}, vector $z_K''(\ph)$ for $\ph = 120\g$}
  \label{Abb:z_K''_vector_curve}
\end{minipage}
\end{figure}

\begin{figure}[ht]
\centering
\includegraphics[width=\textwidth]{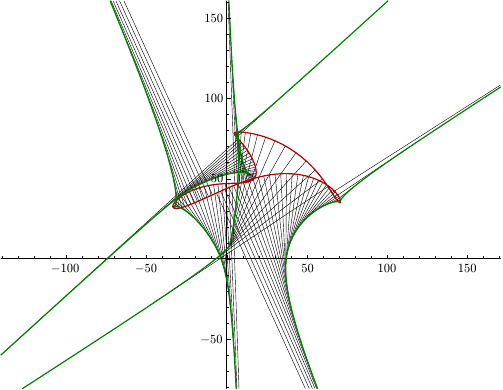}
\caption{Coupler curve (red) from Fig.\ \ref{Abb:Five-bar_linkage02} with evolute\index{evolute} (green) and line segments connecting corresponding points of both curves}
\label{Abb:Evolute_of_coupler_curve}
\end{figure}

Fig.\ \ref{Abb:Tangent_vector_curve} shows the curve of the tangent vector $z_K'(\ph)$ which is called {\em hodograph}\index{hodograph}, and Fig.\ \ref{Abb:z_K''_vector_curve} the curve of the vector $z_K''(\ph)$.
The evolute of the coupler curve is shown in Fig.\ \ref{Abb:Evolute_of_coupler_curve}.
The coupler curve has four {\em inflection points}\index{inflection point} (points at which the curvature changes sign; see Fig.\ \ref{Abb:kappa(phi)}).
Here the evolute has infinity points where it jumps from one direction to the opposite direction.

\begin{SCfigure}[][ht]
\includegraphics[width=0.48\textwidth]{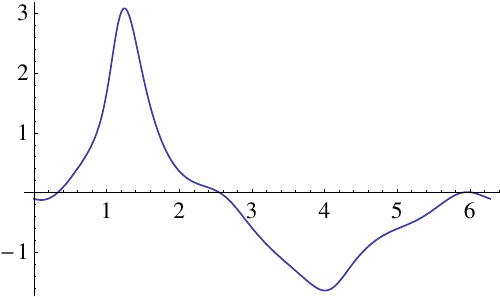}
\caption{Graph of the integrand function in \eqref{Eq:Area_coupler_curve_3} and \eqref{Eq:Area_coupler_curve_4}\\[0.2cm] Near $\ph = 6$ the function has two zeros.}
\label{Abb:Five-bar_integrand}
\end{SCfigure}

Numerical evaluation of \eqref{Eq:Area_coupler_curve_3} (and alternatively \eqref{Eq:Area_coupler_curve_4}) with {\em Mathematica}\index{Mathematica} gives
\beq
  A
\approx 984.03111500882125140  
\eeq
as algebraic sum of the signed areas of the regions enclosed by the four loops of the coupler curve $\mathcal{K}$ (see Fig.\ \ref{Abb:Five-bar_linkage02} and Remark \ref{Rem:A}).
Since the red loops are positively oriented and the blue ones are negatively oriented,\footnote{The crank 2 rotates with positive direction (cf.\ Fig.\ \ref{Abb:Five-bar_linkage01}).} the area $A$ is the sum of the areas of the red outlined regions minus the sum of the areas of the blue outlined regions.
The graph of the function to be integrated in \eqref{Eq:Area_coupler_curve_3} and \eqref{Eq:Area_coupler_curve_4} can be seen in Fig.\ \ref{Abb:Five-bar_integrand}.
The calculation of the sum of the absolute values of the partial areas is a bit more laborious.
Let us denote the areas of the right blue loop, the right red loop, the left blue loop and the left red loop respectively with $A_1$, $A_2$, $A_3$ and $A_4$.
Numerical calculation with the {\em Mathematica}\index{Mathematica} function \texttt{FindRoot}\index{FindRoot@\texttt{FindRoot}} yields that the values of the parameter $\ph$ at the right, the middle and the left self-intersection point are respectively
\begin{align*}
  \phi_1 \approx {} & \phantom{1}12.18109598202294059\g\,,\quad \phi_2 \approx \phantom{1}60.79145803836321748\g\,;\quad\mbox{(right)}\\
  \phi_3 \approx {} & 171.21734636777756458\g\,,\quad \phi_6 \approx 330.79977628841437926\g\,;\quad\mbox{(middle)}\\
  \phi_4 \approx {} & 197.44955305572769338\g\,,\quad \phi_5 \approx 322.29595150395381415\g\,.\,\quad\mbox{(left)}
\end{align*}
With
\begin{align*}
  A(\phi_a,\phi_b)
:= {} & \frac{1}{2} \left\{\ell_2^2\,(\phi_b-\phi_a) + \ell_6 \int_{\phi_a}^{\phi_b} \left(\ell_2\left(\ph_3'(\ph)+1\right)
  \left\langle\ee^{\ii(\ph-\delta)},\, \ee^{\ii\ph_3(\ph)}\right\rangle + \ell_6\,\ph_3'(\ph)\right) \dd\ph\right\} 
\end{align*}
(cf.\ \eqref{Eq:Area_coupler_curve_2} and \eqref{Eq:[z_K,z_K']_5-bar_linkage}) we get
\begin{align*}
  A_1
= {} & A(\phi_1,\phi_2)
\approx -21.73689201894431442\,,\db\\ 
  A_2
= {} & A(\phi_2,\phi_3) + A(\phi_6,2\pi+\phi_1) 
\approx 779.45336525561626079\,,\db\\ 
  A_3
= {} & A(\phi_3,\phi_4) + A(\phi_5,\phi_6) 
\approx -6.21803690630494164\,,\db\\
  A_4
= {} & A(\phi_4,\phi_5) 
\approx 232.53267867845424667\,,
\end{align*}
hence
\beq
  |A_1| + A_2 + |A_3| + A_4
\approx 1039.94097285931976352\,. 
\eeq

\clearpage

%% file: DiffGeo5_4c.tex

\subsection{Plane motion}

\begin{SCfigure}[][ht]
\includegraphics[width=0.5\textwidth]{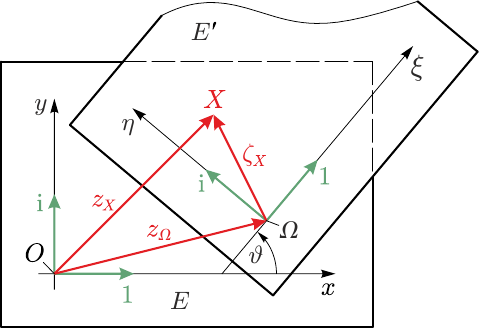}
\caption{Fixed plane $E$ and moving plane $E'$}
\label{Abb:E_and_E'}
\end{SCfigure}

In subsections \ref{Subsec:Kinematics_of_a_dyad} and \ref{Subsec:Five-bar_linkage}, we have already investigated planar motions.
This will now be deepened and extended, taking a more general viewpoint.
We consider the plane motion of a moving plane $E'$ with respect to a fixed reference plane $E$ (see Fig.\ \ref{Abb:E_and_E'}).
Let Cartesian coordinate systems $x,y$ and $\xi,\eta$ be fixed to $E$ or $E'$.   

The position of $E'$ can be defined by the position $z_\Om$ of the coordinate origin $\Om$ in the $x,y$-system and the angle $\vt$ between the positive $x$-axis and the positive $\xi$-axis (see also \cite[Fig.\ 1]{Krause&Carl}).
We consider a point $X$ fixed to $E'$.
Its position with respect to the $\xi,\eta$-system is a complex constant $\zeta_X$, while its position with respect to the $x,y$-system is given by 
\beq
  z_X
= z_\Om + \zeta_X\ee^{\ii\vt}\,.  
\eeq
If we assume $z_\Om$ and $\vt$ to be functions of a real parameter $\ph$, then     
\beqn \label{Eq:z<->zeta_a}
  z_X(\ph)
= z_\Om(\ph) + \zeta_X\ee^{\ii\vt(\ph)}  
\eeqn
is the parametric curve (trajectory) of $X$ with respect to the $x,y$-system.
\eqref{Eq:z<->zeta_a} can be written as
\beqn \label{Eq:z<->zeta_b}
  z_X(\ph)
= z_\Om(\ph) + \zeta_X\ep(\ph) 
\eeqn
with
\beqn \label{Eq:epsilon(ph)}
  \ep(\ph)
= \ee^{\ii\vt(\ph)}  
\eeqn
(see also \cite{Wunderlich:Komplexe_Zahlen}).
The $n$-th derivative of $z_X(\ph)$ with respect to $\ph$ is clearly given by
\beqn \label{Eq:z^{(n)}(phi)}
  z_X^{(n)}(\ph)
= z_\Om^{(n)}(\ph) + \zeta_X\ep^{(n)}(\ph)\,.  
\eeqn
From
\beq
  z_X^{(n)}(\ph)
= 0  
\eeq
we get
\beq
  \zeta_{P_n}(\ph)
= -\frac{z_\Om^{(n)}(\ph)}{\ep^{(n)}(\ph)}   
\eeq
as position of the so-called {\em $n$-th pole} (or {\em pole of $n$-th order}\index{pole!pole of $n$-th order@$\sim $ of $n$-th order}) $P_n$ in the $\xi,\eta$-system.
Putting this into \eqref{Eq:z<->zeta_b} yields
\beqn \label{Eq:z_{P_n}(phi)_from_Omega}
  z_{P_n}(\ph)
= z_\Om(\ph) - \frac{\ep(\ph)}{\ep^{(n)}(\ph)}\,z_\Om^{(n)}(\ph)  
\eeqn
as position of the $n$-th pole in the $x,y$-system.
Applying \eqref{Eq:z<->zeta_b} and \eqref{Eq:z^{(n)}(phi)}, \eqref{Eq:z_{P_n}(phi)_from_Omega} gives
\begin{align} \label{Eq:z_{P_n}(phi)_from_X}
  z_{P_n}(\ph)
= {} & z_X(\ph) - \zeta_X\ep(\ph) - \frac{\ep(\ph)}{\ep^{(n)}(\ph)}\left(z_X^{(n)}(\ph)-\zeta_X\ep^{(n)}(\ph)\right)\nonumber\db\\
= {} & z_X(\ph) - \frac{\ep(\ph)}{\ep^{(n)}(\ph)}\,z_X^{(n)}(\ph)\,,  
\end{align}
so the $n$-th pole can be computed from any point $z_X(\ph)$ with known $n$-th derivative.
From \eqref{Eq:z_{P_n}(phi)_from_X} it follows that 
\beqn \label{Eq:z_X^{(n)}(ph)}
  z_X^{(n)}(\ph)
= \frac{\ep^{(n)}(\ph)}{\ep(\ph)} \left(z_X(\ph)-z_{P_n}(\ph)\right). 
\eeqn
From \eqref{Eq:epsilon(ph)} one gets
\begin{align*}
  \ep'(\ph)
= {} & \ii\vt'(\ph)\ep(\ph)\,,\\[0.05cm]
  \ep''(\ph)
= {} & {\left(\ii\vt''(\ph)-\vt'(\ph)^2\right)\ep(\ph)\,,}\\[0.05cm]
  \ep'''(\ph)
= {} & {\left\{\ii\left(\vt'''(\ph)-\vt'(\ph)^3\right)-3\vt'(\ph)\vt''(\ph)\right\}\ep(\ph)\,.}      
\end{align*}
An important special case of \eqref{Eq:z_X^{(n)}(ph)} is
\beq
  z_X'(\ph)
= \frac{\ep'(\ph)}{\ep(\ph)} \left(z_X(\ph)-z_{P_1}(\ph)\right)
= \ii\vt'(\ph) \left(z_X(\ph)-z_{P_1}(\ph)\right).
\eeq
This formula shows that the vector $z_X'(\ph)$ is orthogonal to the vector $z_X(\ph)-z_{P_1}(\ph)$, and, because of
\beq
  |z_X'(\ph)|
= |\ii\vt'(\ph)\, (z_X(\ph)-z_{P_1}(\ph))|
= |\ii|\, |\vt'(\ph)|\, |z_X(\ph)-z_{P_1}(\ph)|
= |\vt'(\ph)|\, |z_X(\ph)-z_{P_1}(\ph)|\,,
\eeq
the absolute value of $z_X'(\ph)$ is proportional to the distance of point $z_X(\ph)$ from pole $z_{P_1}(\ph)$.
This means that for $\vt'(\ph) \ne 0$ we can consider the motion of $E'$ as instantaneous rotation around $P_1$.
Therefore, $P_1$ is called {\em instantaneous pole}\index{pole!instanteneous pole@instantaneous $\sim $} (also: {\em velocity pole}\index{pole!velocity pole@velocity $\sim $}).
The parametric curve 
\beq
  \zeta_{P_1}(\ph)
= -\frac{z_\Om'(\ph)}{\ep'(\ph)}
= -\frac{z_\Om'(\ph)}{\ii\vt'(\ph)\,\ee^{\ii\vt(\ph)}}
= \ii\, \frac{z_\Om'(\ph)}{\vt'(\ph)}\, \ee^{-\ii\vt(\ph)}   
\eeq
(in $E'$) is called {\em moving centrode}\index{centrode!moving centrode@moving $\sim $} while the parametric curve 
\beq
  z_{P_1}(\ph)
= z_\Om(\ph) - \frac{\ep(\ph)}{\ep'(\ph)}\,z_\Om'(\ph)
= z_\Om(\ph) - \frac{1}{\ii\vt'(\ph)}\,z_\Om'(\ph)
= z_\Om(\ph) + \ii\,\frac{z_\Om'(\ph)}{\vt'(\ph)}
\eeq
(in $E$) is called {\em fixed centrode}\index{centrode!fixed centrode@fixed $\sim $}.

From \eqref{Eq:z_{P_n}(phi)_from_Omega} one easily gets the velocity $u$ of $P_1$ (for time derivative $\frac{\dd \ph}{\dd t} = 1$),
\beqn \label{Eq:velocity_of_P_1_a}
  u(\ph)
= z_{P_1}'(\ph)
= \frac{\ep(\ph)}{\ep'(\ph)^2} \left(\ep''(\ph)z_\Om'(\ph)-\ep'(\ph)z_\Om''(\ph)\right), 
\eeqn
hence
\beq 
  u(\ph)
= \frac{\left(\vt''(\ph)+\ii\vt'(\ph)^2\right)z_\Om'(\ph)-\vt'(\ph)z_\Om''(\ph)}{\ii\vt'(\ph)^2}\,.
\eeq
\eqref{Eq:velocity_of_P_1_a} together with \eqref{Eq:z_{P_n}(phi)_from_Omega} gives 
\begin{align} \label{Eq:velocity_of_P_1_b}  
  u
= {} & \frac{\ep}{\ep'^2}\left(\ep''\,\frac{\ep'}{\ep}\left(z_\Om - z_{P_1}\right) - \ep'\,\frac{\ep''}{\ep}\left(z_\Om - z_{P_2}\right)\right)
  \nonumber\\
= {} & \frac{\ep''}{\ep'}\left(z_\Om - z_{P_1}\right) - \frac{\ep''}{\ep'}\left(z_\Om - z_{P_2}\right)\nonumber\\
= {} & \frac{\ep''}{\ep'}\left(z_{P_2} - z_{P_1}\right),   
\end{align}
thus
\beqn \label{Eq:velocity_of_P_1_c}
  u(\ph)
= \frac{\vt''(\ph)+\ii\vt'(\ph)^2}{\vt'(\ph)} \left(z_{P_2}(\ph) - z_{P_1}(\ph)\right).  
\eeqn
The instantaneous pole is independent of the parametrization which can be seen as follows: We introduce a new parameter $t$ (e.g.\ the time) with
\beqn \label{Eq:parameter_change}
  \ph = f(t),\quad f'(t) \ne 0\,,
\eeqn
and define
\beq
  \tilde{z}_X(t) := z_X(f(t))\,,\quad
  \tilde{\vt}(t) := \vt(f(t))\,.
\eeq
With this from \eqref{Eq:z_{P_n}(phi)_from_X} we get
\beqn \label{Eq:P_1_parameter_change}
\begin{aligned}
  \tilde{z}_{P_1}(t)
= {} & \tilde{z}_X(t) + \ii\,\frac{\tilde{z}_X'(t)}{\tilde{\vt}^{\mbox{$'$}}\!(t)}
= z_X(\ph) + \ii \left(\dfrac{\dd z_X(\ph)}{\dd\ph}\dfrac{\dd f(t)}{\dd t}\right)
  \!\bigg/\!\left(\frac{\dd\vt(\ph)}{\dd\ph}\frac{\dd f(t)}{\dd t}\right)\\
= {} & z_X(\ph) + \ii\,\frac{z_X'(\ph)}{\vt'(\ph)} 
= z_{P_1}(\ph)\,.  
\end{aligned}
\eeqn
Using the angle $\vt$ (see Fig.\ \ref{Abb:E_and_E'}) as parameter gives 
\beq
  z_{P_1}(\vt)
= z_X(\vt) + \ii\,\frac{\dd z_X(\vt)}{\dd\vt}
= z_X(\vt) + \ii\,z_X'(\vt)\,.  
\eeq
The (acceleration) pole $P_2$ is not independent of the parametrization, since here from
\beq 
  \tilde{z}_{P_2}(t)
= \tilde{z}_X(t) - \frac{\strut\tilde{z}_X''(t)}{\strut\ii\tilde{\vt}^{\mbox{$''$}}\!(t)-\tilde{\vt}^{\mbox{$'$}}\!(t)^2}  
\eeq
with the parameter change \eqref{Eq:parameter_change} follows
\beqn \label{Eq:tilde{z}_{P_2}(t)}
  \tilde{z}_{P_2}(t)
= z_X(\ph) - \frac{\strut z_X''(\ph)\,f'(t)^2 + z_X'(\ph)\,f''(t)}{\strut\ii\left(\vt''(\ph)\,f'(t)^2 + \vt'(\ph)\,f''(t)\right) - \vt'(\ph)^2\,f'(t)^2}\,, 
\eeqn
and this is in general different from $z_{P_2}(\ph)$. 
 
We determine the acceleration $z_{P_1}'' \equiv u'$ of $P_1$ (for time derivative $\frac{\dd \ph}{\dd t} = 1$).
From \eqref{Eq:velocity_of_P_1_b} we get
\beq
  \ep'' z_{P_1}' + \ep' z_{P_1}''
= \ep'''\left(z_{P_2}-z_{P_1}\right) + \ep''\left(z_{P_2}'-z_{P_1}'\right),  
\eeq
hence
\begin{align*}
  z_{P_1}''
= {} & \frac{\ep'''}{\ep'}\left(z_{P_2}-z_{P_1}\right) + \frac{\ep''}{\ep'}\left(z_{P_2}'-z_{P_1}'\right) - \frac{\ep''}{\ep'}\,z_{P_1}'\\ 
= {} & \frac{\ep'''}{\ep'}\left(z_{P_2}-z_{P_1}\right) + \frac{\ep''}{\ep'}\left(z_{P_2}'-2z_{P_1}'\right).
\end{align*}
For the velocity of the pole $P_2$ one finds
\begin{align*}
  z_{P_2}'
= {} & {-}\frac{\ep'}{\ep}\left(z_{P_1}-z_{P_2}\right) + \frac{\ep'''}{\ep''}\left(z_{P_3}-z_{P_2}\right)\\
= {} & {-}\frac{\ep'}{\ep}\,z_{P_1} + \left(\frac{\ep'}{\ep}-\frac{\ep'''}{\ep''}\right)z_{P_2} + \frac{\ep'''}{\ep''}\,z_{P_3}\,.    
\end{align*}
It follows
\begin{align*}
  z_{P_1}''
= {} &  {-}\frac{\ep'''}{\ep'}\,z_{P_1} + \frac{\ep'''}{\ep'}\,z_{P_2}
  + \frac{\ep''}{\ep'}\left(-\frac{\ep'}{\ep}\,z_{P_1} + \left(\frac{\ep'}{\ep}-\frac{\ep'''}{\ep''}\right)z_{P_2} + \frac{\ep'''}{\ep''}\,z_{P_3}\right)
  - 2\,\frac{\ep''^2}{\ep'^2}\left(z_{P_2}-z_{P_1}\right)\db\\[0.05cm]
= {} & {-}\left(\frac{\ep''}{\ep}-2\,\frac{\ep''^2}{\ep'^2}+\frac{\ep'''}{\ep'}\right)z_{P_1}
  + \left(\frac{\ep''}{\ep}-2\,\frac{\ep''^2}{\ep'^2}\right)z_{P_2}
  + \frac{\ep'''}{\ep'}\,z_{P_3}\db\\[0.05cm]
= {} & \left(\frac{\ep''}{\ep}-2\,\frac{\ep''^2}{\ep'^2}\right)\left(z_{P_2}-z_{P_1}\right) + \frac{\ep'''}{\ep'}\left(z_{P_3}-z_{P_1}\right).     
\end{align*}

%% file: DiffGeo5_4d.tex

\subsection{Inflection circle and curve of stationary curvature}

We are looking for the points of the moving plane $E'$ whose trajectories currently have vanishing curvature $\kappa(\ph)$.
From \eqref{Eq:Curvature} for $z'(\ph) \ne 0$ we get
\beq
  \left[z'(\ph),z''(\ph)\right] 
= 0\,.  
\eeq
Applying \eqref{Eq:z_X^{(n)}(ph)} gives
\begin{align*}
  0
= {} & {\left[\frac{\ep'(\ph)}{\ep(\ph)}\,(z-z_{P_1}(\ph)),\,\frac{\ep''(\ph)}{\ep(\ph)}\,(z-z_{P_2}(\ph))\right]}\db\\[0.05cm]
= {} & {\left[\ii\vt'(\ph)\,(z-z_{P_1}(\ph)),\,\left(\ii\vt''(\ph)-\vt'(\ph)^2\right)(z-z_{P_2}(\ph))\right]}\db\\[0.05cm]
= {} & {\left[\ii\vt'(\ph)\,(z-z_{P_1}(\ph)),\,\ii\left(\vt''(\ph)+\ii\vt'(\ph)^2\right)(z-z_{P_2}(\ph))\right]}\db\\[0.05cm]
= {} & {\left[\vt'(\ph)\,(z-z_{P_1}(\ph)),\,\left(\vt''(\ph)+\ii\vt'(\ph)^2\right)(z-z_{P_2}(\ph))\right]}\db\\[0.05cm]
= {} & \vt'(\ph)F(z)
\end{align*}
with
\beqn \label{Eq:F(z)_inflection_circle}
  F(z)
= \left[z-z_{P_1}(\ph),\,\left(\vt''(\ph)+\ii\vt'(\ph)^2\right)(z-z_{P_2}(\ph))\right].
\eeqn
We transform to an $X,Y$-coordinate system whose origin is at $z_{P_1}(\ph)$ and whose positive $X$-axis has the direction of the velocity of the instantaneous pole, $u(\ph)$.
With $Z = X + \ii Y$ and $\bar{u} = |u|\,\ee^{-\ii\arg(u)}$ we have the transformation equations
\beqn \label{Eq:z->Z_Z->z}
  Z = \left(z-z_{P_1}\right)\ee^{-\ii\arg(u)}
  \quad\mbox{or}\quad
  z = z_{P_1} + Z\,\ee^{\ii\arg(u)}\,,
\eeqn
and therefore
\begin{align*}
  \widetilde{F}(Z)
= {} & F\left(z_{P_1}+Z\ee^{\ii\arg(u)}\right)\db\\
= {} & {\left[Z\ee^{-\ii\arg(u)},\left(\vt''+\ii\vt'^2\right)(Z\ee^{\ii\arg(u)}+z_{P_1}-z_{P_2})\right]}\db\\
= {} & |Z|^2\left[1,\vt''+\ii\vt'^2\right] + \left[Z\ee^{\ii\arg(u)},\left(\vt''+\ii\vt'^2\right)\left(z_{P_1}-z_{P_2}\right)\right]\db\\ 
= {} & \vt'^2|Z|^2 + \big[Z,\left(\vt''+\ii\vt'^2\right)\left(z_{P_1}-z_{P_2}\right)\ee^{-\ii\arg(u)}\big]\db\\
= {} & \vt'^2|Z|^2 - |u|^{-1}\big[Z,\left(\vt''+\ii\vt'^2\right)\left(z_{P_2}-z_{P_1}\right)\bar{u}\big].
\end{align*}
With \eqref{Eq:velocity_of_P_1_c} follows
\beq
  \widetilde{F}(Z)
= \vt'(\ph)^2 |Z|^2 - |u(\ph)|^{-1}\big[Z,\vt'(\ph)u(\ph)\bar{u}(\ph)\big]
= \vt'(\ph)^2 |Z|^2 - \vt'(\ph)|u(\ph)|\big[Z,1\big],  
\eeq
hence
\beq
  |Z|^2 - \frac{|u(\ph)|}{\vt'(\ph)}\,[Z,1]
= 0\,,
\eeq
or in real form
\beqn \label{Eq:Inflection_circle_real}
  X^2 + Y^2 + \frac{|u(\ph)|}{\vt'(\ph)}\,Y
= 0\,.  
\eeqn
Adding $(|u|/(2\vt'))^2$ on both sides of \eqref{Eq:Inflection_circle_real} gives
\beq
  X^2 + \left(Y + \frac{|u(\ph)|}{2\vt'(\ph)}\right)^2
= \left(\frac{|u(\ph)|}{2\vt'(\ph)}\right)^2.  
\eeq
This is the equation of a circle with
\beqn \label{Eq:center_point_and_diameter}
  \mbox{center point}\; \left(0,\,-\frac{|u(\ph)|}{2\vt'(\ph)}\right)
  \;\,\mbox{and}\;\;
  \mbox{diameter}\;\, \frac{|u(\ph)|}{|\vt'(\ph)|} = \left|\frac{u(\ph)}{\vt'(\ph)}\right|
  = \left|\frac{\dd z_{P_1}(\vt)}{\dd\vt}\right|.
\eeqn
It is called {\em inflection circle}\footnote{In German: Wendekreis.}\index{inflection circle}\index{circle!inflection circle@inflection $\sim $}.

Now, we are looking for the points of the moving plane $E'$ whose trajectories have currently constant curvature $\kappa(\ph)$, hence $\kappa'(\ph) = 0$.
From \eqref{Eq:Curvature} we get
\beq
  \kappa'(\ph)
= \frac{\left[z''(\ph),z''(\ph)\right]+\left[z'(\ph),z'''(\ph)\right]}{|z'(\ph)|^3}
  - 3\,\frac{\left[z'(\ph),z''(\ph)\right]|z'(\ph)|'}{|z'(\ph)|^4}\,,  
\eeq
and with
\beq
  \left|z'(\ph)\right|'
= \left(\sqrt{\left\langle z'(\ph),z'(\ph) \right\rangle}\,\right)'
= \frac{\left\langle z''(\ph),z'(\ph) \right\rangle + \left\langle z'(\ph),z''(\ph) \right\rangle}{2\,\sqrt{\left\langle z'(\ph),z'(\ph) \right\rangle}}
= \frac{\left\langle z'(\ph),z''(\ph) \right\rangle}{\left|z'(\ph)\right|}   
\eeq
(see also the proof of Lemma \ref{Lem:T'(t)}) follows 
\beq
  \kappa'(\ph)
= \frac{\left|z'(\ph)\right|^2\left[z'(\ph),z'''(\ph)\right]-3\left[z'(\ph),z''(\ph)\right]\left\langle z'(\ph),z''(\ph) \right\rangle}
	{\left|z'(\ph)\right|^5}\,.  
\eeq
Using \eqref{Eq:z_X^{(n)}(ph)} we have found
\begin{align*}
  0
= {} & 3\left[\frac{\ep'(\ph)}{\ep(\ph)}\,(z-z_{P_1}(\ph)),\,\frac{\ep''(\ph)}{\ep(\ph)}\,(z-z_{P_2}(\ph))\right]
  \left\langle \frac{\ep'(\ph)}{\ep(\ph)}\,(z-z_{P_1}(\ph)),\,\frac{\ep''(\ph)}{\ep(\ph)}\,(z-z_{P_2}(\ph)) \right\rangle\\[0.05cm]
& - \left|\frac{\ep'(\ph)}{\ep(\ph)}\,(z-z_{P_1}(\ph))\right|^2
  \left[\frac{\ep'(\ph)}{\ep(\ph)}\,(z-z_{P_1}(\ph)),\,\frac{\ep'''(\ph)}{\ep(\ph)}\,(z-z_{P_3}(\ph))\right]
\end{align*}
as necessary condition for the points $z$ with currently constant curvature.
From this we get
\begin{align} \label{Eq:G(z)_curve_of_stationary_curvature}
  0 =
G(z)
:= {} & 3\left[\ep'\left(z-z_{P_1}\right),\ep''\left(z-z_{P_2}\right)\right]  
  \left\langle\ep'\left(z-z_{P_1}\right),\ep''\left(z-z_{P_2}\right)\right\rangle\nonumber\\
& - \vt'^2\left|z-z_{P_1}\right|^2 \left[\ep'\left(z-z_{P_1}\right),\ep'''\left(z-z_{P_3}\right)\right].
\end{align}
With the essential ingredients $u = |u|\,\ee^{\ii\arg(u)}$ and \eqref{Eq:velocity_of_P_1_b}, the transformation to the above defined $X,Y$-system gives
\begin{align*}
  0
= {} & \widetilde{G}(Z) := G\left(z_{P_1}+Z\ee^{\ii\arg(u)}\right)\db\\[0.15cm]
= {} & 3\left[\ep'Z\ee^{\ii\arg(u)},\ep''\left(Z\ee^{\ii\arg(u)}-(z_{P_2}-z_{P_1})\right)\right] 
  \left\langle\ep'Z\ee^{\ii\arg(u)},\ep''\left(Z\ee^{\ii\arg(u)}-(z_{P_2}-z_{P_1})\right)\right\rangle\\
& - \vt'^2\left|Z\ee^{\ii\arg(u)}\right|^2\left[\ep'Z\ee^{\ii\arg(u)},\ep'''\left(Z-Z_{P_3}\right)\ee^{\ii\arg(u)}\right]\db\\[0.15cm]
= {} & 3\left[\ep'\,\frac{u}{|u|}\,Z,\ep''\left(\frac{u}{|u|}\,Z-\frac{\ep'}{\ep''}\,u\right)\right] 
  \left\langle\ep'\,\frac{u}{|u|}\,Z,\ep''\left(\frac{u}{|u|}\,Z-\frac{\ep'}{\ep''}\,u\right)\right\rangle\\
& - \vt'^2|Z|^2\left[\ep'Z,\ep'''\left(Z-Z_{P_3}\right)\right]\db\\[0.15cm]
= {} & 3\left[\ep'Z,\ep''Z-\ep'|u|\right]\left\langle\ep'Z,\ep''Z-\ep'|u|\right\rangle
  - \vt'^2\,|Z|^2\left[\ep'Z,\ep'''\left(Z-Z_{P_3}\right)\right]\db\\[0.05cm]
= {} & 3\left(\left[\ep'Z,\ep''Z\right]-\left[\ep'Z,\ep'|u|\right]\right)
  \left(\left\langle\ep'Z,\ep''Z\right\rangle-\left\langle\ep'Z,\ep'|u|\right\rangle\right)\\
& - \vt'^2\,|Z|^2\left(\left[\ep'Z,\ep'''Z\right]-\left[\ep'Z,\ep'''Z_{P_3}\right]\right)\db\\[0.05cm]
= {} & 3\left(|Z|^2\left[\ep',\ep''\right]-\vt'^2|u|\left[Z,1\right]\right)
  \left(|Z|^2\left\langle\ep',\ep''\right\rangle-\vt'^2|u|\left\langle Z,1\right\rangle\right)\\
& - \vt'^2\,|Z|^2\left(|Z|^2\left[\ep',\ep'''\right]-\left[\ep'Z,\ep'''Z_{P_3}\right]\right)\db\\[0.05cm]
= {} & 3\left(|Z|^2\left[\ep',\ep''\right]+\vt'^2|u|\Imz(Z)\right)\left(|Z|^2\left\langle\ep',\ep''\right\rangle-\vt'^2|u|\Rez(Z)\right)\\
& - \vt'^2\,|Z|^2\left(|Z|^2\left[\ep',\ep'''\right]-\left[\ep'Z,\ep'''Z_{P_3}\right]\right)	
\end{align*}
with
\beq
  Z_{P_3}(\ph)
= (z_{P_3}(\ph)-z_{P_1}(\ph))\,\ee^{-\ii\arg(u(\ph))}\,.
\eeq
One easily finds
\begin{gather*}
  \left[\ep',\ep''\right]
= \vt'^3\,,\quad
  \left\langle \ep',\ep''\right\rangle
= \vt'\vt''\,,\quad
  \left[\ep',\ep'''\right]
= 3\vt'^2\vt''\,,\db\\[0.05cm]  
  \left[\ep'Z,\ep'''Z_{P_3}\right]
= \vt'\left(\vt'''-\vt'^3\right)\left[Z,Z_{P_3}\right] +3\vt'^2\vt''\left\langle Z,Z_{P_3}\right\rangle,       
\end{gather*}
and therefore
\begin{align*}
  \widetilde{G}(Z)
= {} & 3\left(\vt'^3|Z|^2+\vt'^2|u|\Imz(Z)\right)\left(\vt'\vt''|Z|^2-\vt'^2|u|\Rez(Z)\right)\\
& - \vt'^2\,|Z|^2\left(3\vt'^2\vt''|Z|^2-\vt'\left(\vt'''-\vt'^3\right)\left[Z,Z_{P_3}\right]+3\vt'^2\vt''\left\langle Z,Z_{P_3}\right\rangle\right)
  \db\\[0.05cm]
= {} & {-\vt'^2\left\{3\left(\vt'^3|u||Z|^2\Rez(Z)-\vt'\vt''|u||Z|^2\Imz(Z)+\vt'^2|u|^2\Rez(Z)\Imz(Z)\right)\right.}\\
&  \left.-\:|Z|^2 \left(\vt'\left(\vt'''-\vt'^3\right)\left[Z,Z_{P_3}\right]+3\vt'^2\vt''\left\langle Z,Z_{P_3}\right\rangle\right)\right\}\db\\[0.05cm]
= {} & {-\vt'^2\left\{|Z|^2\left(3\vt'^3|u|\Rez(Z)-3\vt'\vt''|u|\Imz(Z)-\vt'\left(\vt'''-\vt'^3\right)\left[Z,Z_{P_3}\right]
  -3\vt'^2\vt''\left\langle Z,Z_{P_3}\right\rangle\right)\right.}\\
&  \left.+\:3\vt'^2|u|^2\Rez(Z)\Imz(Z)\right\},
\end{align*}
thus
\begin{align*}
  0
= {} & \left(X^2+Y^2\right)\left\{3\vt'(\ph)^3\left|u(\ph)\right|X-3\vt'(\ph)\vt''(\ph)\left|u(\ph)\right|Y\right.\\
& \left.-\:\vt'(\ph)\left(\vt'''(\ph)-\vt'(\ph)^3\right)\left(XY_{P_3}(\ph)-YX_{P_3}(\ph)\right)\right.\\
& \left.-\:3\vt'(\ph)^2\vt''(\ph)\left(XX_{P_3}(\ph)+YY_{P_3}(\ph)\right)\right\} +\:3\vt'(\ph)^2\left|u(\ph)\right|^2XY.
\end{align*}
This shows that all points with $\kappa'(\ph) = 0$ lie on a cubic curve. 
This curve is called {\em cubic of stationary curvature}\index{cubic!cubic of stationary curvature@$\sim $ of stationary curvature}.
Since it passes through the {\em circular points at infinity}\index{circular points at infinity} {\em (cyclic points)}\index{cyclic points} with homogeneous coordinates $[1,\ii,0]$ and $[1,-\ii,0]$, it is a {\em circular cubic}\index{cubic!circular cubic@circular $\sim $}.
It has a double point at the instantaneous pole $P_1$.
If one removes all third-order terms, it can be seen that the curve touches the tangent $t$ and the normal $n$ of the {\em centrode}\index{centrode} (curve $z_{P_1}(\ph)$) at $P_1$.
(See also \cite[pp.\ 76-77]{Luck&Modler}.) 

\textsc{Ball}'s\footnote{\textsc{Robert Stawell Ball}, 1840-1913} {\em point}\index{Ball's point} $U$ is the point of the moving plane with $\kappa(\ph) = 0$ and $\kappa'(\ph) = 0$, hence the intersection point of the inflection circle and the curve of stationary curvature.
Due to $\kappa(\ph) = 0$ we have
\beqn \label{Eq:first_circle}
  \left[\ep'(z-z_{P_1}),\ep''(z-z_{P_2})\right] = 0\,,
\eeqn
and so from \eqref{Eq:G(z)_curve_of_stationary_curvature} follows
\beqn \label{Eq:zero-normal_jerk_circle}
  \left[\ep'(z-z_{P_1}),\ep'''(z-z_{P_3})\right] = 0\,,
\eeqn
therefore \textsc{Ball}'s point ist the solution $z$ of the system of equations \eqref{Eq:first_circle}, \eqref{Eq:zero-normal_jerk_circle}.
\eqref{Eq:first_circle} is the equation of the inflection circle, while \eqref{Eq:zero-normal_jerk_circle} is the equation of another circle (cf.\ \eqref{Eq:Equation_of_c_1} and \eqref{Eq:Equation_of_c_2}) called {\em zero-normal jerk circle}\index{circle!zero-normal jerk circle@zero-normal jerk $\sim $} (cf.\ \cite{Figliolini&Lanni&Tomassi}).
Apart from the {\em circular points at infinity}\index{circular points at infinity}, these two circles intersect at pole $P_1$ and \textsc{Ball}'s point $U$.

Applying the well-known calculation rules, the equation \eqref{Eq:first_circle} of the inflection circle $k_1$ can for $\vt'(\ph) \ne 0$  be written as
\begin{align*}
& z\bar{z} - \frac{1}{2}\;\overline{\left(z_{P_1}+z_{P_2} + \ii\,\frac{\vt''}{\vt'^2}\left(z_{P_1}-z_{P_2}\right)\right)}\;z
  - \frac{1}{2}\,\left(z_{P_1}+z_{P_2} + \ii\,\frac{\vt''}{\vt'^2}\left(z_{P_1}-z_{P_2}\right)\right)\,\bar{z}\\[0.1cm]
& + \frac{\left[\ep'z_{P_1},\ep''z_{P_2}\right]}{\vt'^3}
= 0\,.  
\end{align*}
Comparing with \eqref{Eq:squared_circular_equation} shows that 
\beq
  z_1(\ph)
:= \frac{1}{2}\left(z_{P_1}(\ph)+z_{P_2}(\ph) + \ii\,\frac{\vt''(\ph)}{\vt'(\ph)^2}\left(z_{P_1}(\ph)-z_{P_2}(\ph)\right)\right)  
\eeq
is the circle center point and
\beq
  r_1(\ph)^2
= |z_1(\ph)|^2 - \frac{\left[\ep'(\ph)z_{P_1}(\ph),\,\ep''(\ph)z_{P_2}(\ph)\right]}{\vt'(\ph)^3}  
\eeq
is the square of the radius.
With
\beq
  \left[\ep'z_{P_1},\,\ep''z_{P_2}\right]
= \vt'\vt''\left[z_{P_1},z_{P_2}\right] + \vt'^3\left\langle z_{P_1},z_{P_2}\right\rangle
\eeq
we get
\beq
  r_1^2
= \frac{\vt'^4+\vt''^2}{4\vt'^4} \left|z_{P_1}-z_{P_2}\right|^2.  
\eeq
From \eqref{Eq:velocity_of_P_1_c} we know that
\beq
  |u|^2
= \left|\frac{\vt''+\ii\vt'^2}{\vt'}\right|^2 \left|z_{P_1}-z_{P_2}\right|^2 
= \frac{\vt''^2+\vt'^4}{\vt'^2} \left|z_{P_1}-z_{P_2}\right|^2.
\eeq
It follows the already known result (see \eqref{Eq:center_point_and_diameter}) 
\beq
  r_1(\ph)^2 
= \frac{\left|u(\ph)\right|^2}{4\vt'(\ph)^2}\,. 
\eeq
Clearly, the poles $P_1$ and $P_2$ lie on $k_1$.
We consider the vector $z_1(\ph) - z_{P_1}(\ph)$.
It holds
\begin{align*}
  z_1 - z_{P_1}
= {} & \frac{1}{2} \left(z_{P_1} + z_{P_2} + \ii\,\frac{\vt''}{\vt'^2}\left(z_{P_1} - z_{P_2}\right)\right) - z_{P_1}  
= \frac{1}{2} \left(-(z_{P_1} - z_{P_2}) + \ii\,\frac{\vt''}{\vt'^2}\left(z_{P_1} - z_{P_2}\right)\right)\db\\  
= {} & \frac{1}{2} \left(z_{P_2} - z_{P_1} - \ii\,\frac{\vt''}{\vt'^2}\left(z_{P_2} - z_{P_1}\right)\right)
= \frac{1}{2} \left(1 - \ii\,\frac{\vt''}{\vt'^2}\right) \left(z_{P_2} - z_{P_1}\right)\db\\
= {} & \frac{1}{2}\, \frac{\vt'^2 - \ii\vt''}{\vt'^2} \left(z_{P_2} - z_{P_1}\right)
= \frac{1}{2}\,\ii\, \frac{-\ii\vt'^2 - \vt''}{\vt'^2} \left(z_{P_2} - z_{P_1}\right)\db\\
= {} & {-}\frac{1}{2}\,\ii\, \frac{\vt'' + \ii\vt'^2}{\vt'^2} \left(z_{P_2} - z_{P_1}\right)\,.
\end{align*}
Comparison with \eqref{Eq:velocity_of_P_1_c} shows that
\beq
  z_1(\ph) - z_{P_1}(\ph)
= -\frac{1}{2}\,\ii\, \frac{u(\ph)}{\vt'(\ph)}\,,
  \quad\mbox{hence}\quad
  z_1(\ph)
= z_{P_1}(\ph) - \frac{1}{2}\,\ii\, \frac{u(\ph)}{\vt'(\ph)}    
\eeq
(cf.\ \cite[pp.\ 35-36]{Blaschke&Mueller}).
The inflection circle does not depend on the parametrization which can be seen as follows:
Using the parameter change \eqref{Eq:velocity_of_P_1_c},
\beq
  \tilde{u}(t) := u(f(t))\,,\quad
  \tilde{\vt}(t) := \vt(f(t))\,.
\eeq
and \eqref{Eq:P_1_parameter_change}, we get
\begin{align*}
  \tilde{z}_1(t) 
= {} & \tilde{z}_{P_1}(t) - \frac{1}{2}\,\ii\, \frac{\tilde{u}(t)}{\tilde{\vt}^{\mbox{$'$}}\!(t)}  
= z_{P_1}(\ph) - \frac{1}{2}\,\ii \left(\frac{\dd z_{P_1}(\ph)}{\dd\ph}\,\dfrac{\dd f(t)}{\dd t}\right)
  \!\bigg/\!\left(\frac{\dd\vt(\ph)}{\dd\ph}\,\dfrac{\dd f(t)}{\dd t}\right)\\
= {} & z_{P_1}(\ph) - \frac{1}{2}\,\ii\, \frac{z_{P_1}'(\ph)}{\vt'(\ph)}
= z_{P_1}(\ph) - \frac{1}{2}\,\ii\, \frac{u(\ph)}{\vt'(\ph)}\db\\
= {} & z_{P_1}(\ph) + \left(z_1(\ph) - z_{P_1}(\ph)\right)
= z_1(\ph)\,,
\end{align*}
and so the center remains unchanged.
Because of $r_1(\ph) = |z_1(\ph)-z_{P_1}(\ph)| = |\tilde{z}_1(t)-\tilde{z}_{P_1}(t)|$, the radius also remains unchanged. 
Note that although, as shown, the pole $P_2$ is parameter-dependent, it must always lie on the parameter-independent inflection circle $k_1$.

The equation \eqref{Eq:zero-normal_jerk_circle} of the zero-normal jerk circle $k_3$ can for $\vt'(\ph) \ne 0$ be written as
\beq
  z\bar{z} - \overline{z_3(\ph)}\,z - z_3(\ph)\,\bar{z} + \left|z_3(\ph)\right|^2 - r_3(\ph)^2
= 0  
\eeq
with center point
\beq
  z_3
= \frac{1}{2}\left(z_{P_1}+z_{P_3} + \ii\,\frac{\vt'''-\vt'^3}{3\vt'\vt''}\left(z_{P_1}-z_{P_3}\right)\right)  
\eeq
and radius square
\beq
  r_3^2
= |z_3|^2 - \frac{\left[\ep'z_{P_1},\ep'''z_{P_3}\right]}{3\vt'^2\vt''}
= \frac{\vt'^6 + 9\vt'^2\vt''^2 - 2\vt'^3\vt''' + \vt'''^2}{36\vt'^2\vt''^2} \left|z_{P_1} - z_{P_3}\right|^2\,.
\eeq
Now, we will derive an explicit formula for \textsc{Ball}'s point $U$.
The elimination of $z\bar{z}$ from the non-linear equation system
\beq
\left.
\begin{array}{*{2}{c@{\:-\:}}c@{\:+\:}c@{\:-\:}c@{\:=\:}c}
  z\bar{z} & \overline{z_1(\ph)}\,z & z_1(\ph)\,\bar{z} & \left|z_1(\ph)\right|^2 & r_1(\ph)^2 & 0\,,\\[0.1cm]
  z\bar{z} & \overline{z_3(\ph)}\,z & z_3(\ph)\,\bar{z} & \left|z_3(\ph)\right|^2 & r_3(\ph)^2 & 0\,,
\end{array} 
\right\}
\eeq
gives the equation of the line through points $P_1$ and $U$ ($z_{P_1}$ and $z_U$),
\beqn \label{Eq:line_through_zP1_and_zU}
  \overline{\left(z_1(\ph)-z_3(\ph)\right)}\,z + \left(z_1(\ph)-z_3(\ph)\right)\bar{z}
= |z_1(\ph)|^2 - |z_3(\ph)|^2 - r_1(\ph)^2 + r_3(\ph)^2\,.
\eeqn 
The line through the center points $z_1(\ph)$ and $z_3(\ph)$ is given by equation
\beq
  \left[z-z_1(\ph),z-z_3(\ph)\right] = 0\,,
\eeq 
which can be written as
\beqn \label{Eq:line_through_z1_and_z2}
  \overline{\left(z_1(\ph)-z_3(\ph)\right)}\,z - \left(z_1(\ph)-z_3(\ph)\right)\bar{z}
= \overline{z_1(\ph)}\,z_3 - z_1\,\overline{z_3(\ph)}\,.   
\eeqn
The intersection point $S$ or $z_S = z_S(\ph)$ of both lines is the solution of the linear equation system \eqref{Eq:line_through_zP1_and_zU}, \eqref{Eq:line_through_z1_and_z2}, thus, using \textsc{Cramer}'s rule,
\begin{align*}
  z_S
= {} & \frac{\left|\begin{array}{@{\,}cr@{\,}} |z_1|^2-|z_3|^2-r_1^2+r_3^2 & z_1-z_3\phantom{)}\\ \bar{z}_1z_3 - z_1\bar{z}_3 & -(z_1-z_3) \end{array}
	   \right|}
	   {\left|\begin{array}{@{\,}cr@{\,}} \bar{z}_1-\bar{z}_3 & z_1-z_3\phantom{)}\\ \bar{z}_1-\bar{z}_3 & -(z_1-z_3)\, \end{array}\right|}
= \frac{\left|\begin{array}{@{\,}cr@{\,}} |z_1|^2-|z_3|^2-r_1^2+r_3^2 & 1\\ \bar{z}_1z_3 - z_1\bar{z}_3 & -1 \end{array}\right|}
	   {(\bar{z}_1-\bar{z}_3)\left|\begin{array}{@{\,}cr@{\,}} 1 & 1\\ 1 & -1 \end{array}\right|}\\[0.1cm]
= {} & \frac{|z_1|^2 - |z_3|^2 - r_1^2 + r_3^2 + \bar{z}_1z_3 - z_1\bar{z}_3}{2\left(\overline{z}_1-\overline{z}_3\right)}\,.	     
\end{align*}
\textsc{Ball}'s point can now be easily calculated by
\beq
  z_{U}(\ph)
= z_{P_1}(\ph) + 2\left(z_{S}(\ph)-z_{P_1}(\ph)\right)
= 2z_S(\ph) - z_{P_1}(\ph)\,.  
\eeq  

\begin{example}
We consider the four-bar linkage $A_0ABB_0$ in Fig.\ \ref{Abb:Curvature_cubic} with rotating joints at points $A_0$, $A \equiv \Om$, $B$, $B_0$.
The parameter $\ph$ is the (mathematically positive) angle between positive $x$-axis and vector $\vv{A_0A}$, in Fig.\ \ref{Abb:Curvature_cubic} $\ph = 345\g\cdot\pi/180\g$.
$\vt = \vt(\ph)$ is the angle between the positive $x$-axis and vector $\vv{AB}$. 

\begin{figure}[ht]
\centering
\includegraphics[width=0.65\textwidth]{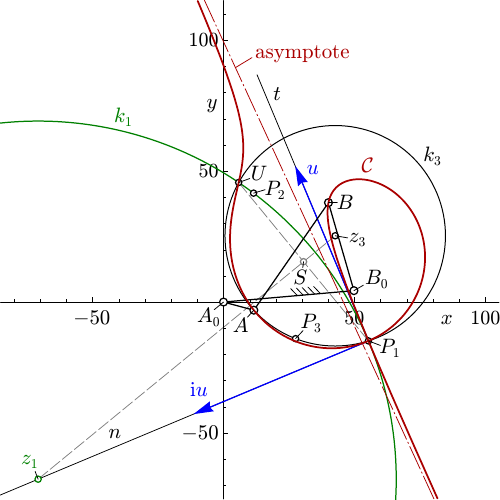}
\caption{Four-bar linkage with inflection circle $k_1$, zero-normal jerk circle\index{circle!zero-normal jerk circle@zero-normal jerk $\sim $} $k_3$, stationary curvature curve $\K$, and \textsc{Ball's} point $U$}
\label{Abb:Curvature_cubic}
\end{figure}

Calling the subroutine {\em DyadC} (see p.\ \pageref{Subroutine-code}) with
\beq
  \texttt{DyadC[l3, l4, zA, zA1, zA2, zA3, zB0, 0, 0, 0, 1]}\,,
\eeq
where
\begin{gather*}
  \texttt{l3} = \ell_3 = \overline{AB} = 50\,,\quad
  \texttt{l4} = \ell_4 = \overline{B_0B} = 35\,,\\
  \texttt{zA} = z_A(\ph) = z_\Om(\ph) = \ell_2\,\ee^{\ii\ph}
  \;\;\mbox{with}\;\;
  \ell_2 = \overline{A_0A} = 12\,,\quad 
  \texttt{zA1} = z_A'(\ph) = \ell_2\,\ii\ee^{\ii\ph}\,,\\
  \texttt{zA2} = z_A''(\ph) = -\ell_2\,\ee^{\ii\ph}\,,\quad 
  \texttt{zA3} = z_A'''(\ph) = -\ell_2\,\ii\ee^{\ii\ph}\,,\quad
  \texttt{zB0} = z_{B_0} = 50\exp(5\pi\ii/180)\,,
\end{gather*}
gives
\begin{gather*}
  \vt(\ph) \approx 0.964415 \entspr 55.2569\g\,,\quad
  \vt'(\ph) \approx -0.264825\,,\\
  \vt''(\ph) \approx -0.259041\,,\quad
  \vt'''(\ph) \approx 0.568709\,,
\end{gather*}
and consequently
\begin{gather*}
  z_{P_1} \approx 55.3600 - 14.8337\,\ii\,,\quad
  z_{P_2} \approx 11.4748 + 41.6089\,\ii\,,\\
  z_{P_3} \approx 27.5188 - 13.9759\,\ii\,,\quad
  u \approx -27.9792 + 66.8317\,\ii\,,\db\\
  z_1 \approx  -70.8208 - 67.6595\,\ii\,,\quad
  r_1 \approx 136.792\,,\quad
  z_3 \approx 42.6633 + 25.3195\,\ii\,,\quad
  r_3 \approx 42.1127\,,\\
  z_U \approx 5.7919 + 45.6660\,\ii\,.
\end{gather*}
With these values, the algebraic equation of the inflection circle $k_1$ is given by $F(z) = 0$ with $F(z)$ from \eqref{Eq:F(z)_inflection_circle}.
The algebraic equation of the curve of stationary curvature is given by \eqref{Eq:G(z)_curve_of_stationary_curvature}. $t$ and $n$ denote the tangent or normal to the curve $z_{P_1}$ at point $z_{P_1}(\ph)$.
As the trajectories of points $A$ and $B$ are respectively a circle and a circular arc, their curvature is constant and therefore $A$ and $B$ lie on the curve of constant curvature.
Since the vectors $z_A'(\ph) = \overline{A_0A}\,\ii\ee^{\ii\ph}$ and $z_A'''(\ph) = -\overline{A_0A}\,\ii\ee^{\ii\ph}$ are collinear, we have $[z_A',z_A'''] = 0$ and so from \eqref{Eq:zero-normal_jerk_circle} follows that $A$ lies on the circle $k_3$.

Fig.\ \ref{Abb:Coupler_curve_with_Ball_point} shows the four-bar linkage from Fig.\ \ref{Abb:Curvature_cubic} with added coupler point $X$ chosen so that
\beqn \label{Eq:X=U}
  \zeta_X
= \left(z_U(\ph)-z_\Om(\ph)\right) \ee^{-\ii\vt(\ph)}
  \quad\mbox{with}\quad
  \ph = 345\g\cdot\pi/180\g
\eeqn
which gives
\beq
  \zeta_X
= 36.7715 + 32.5603\,\ii\,.
\eeq
The parametric coupler curve $z_X$ is given by \eqref{Eq:z<->zeta_a}.
Because of choosing $\zeta_X$ according to \eqref{Eq:X=U}, we have $z_X(\ph) = z_U(\ph)$ in the shown position.
The coupler curve therefore has curvature $\kappa(\ph) = 0$ and derivative $\kappa'(\ph) = 0$ at point $z_X(\ph) = z_U(\ph)$, i.e. it touches its tangent $t$ very closely.\hfill\bs

\begin{SCfigure}[0.9][ht]
\includegraphics[width=0.5\textwidth]{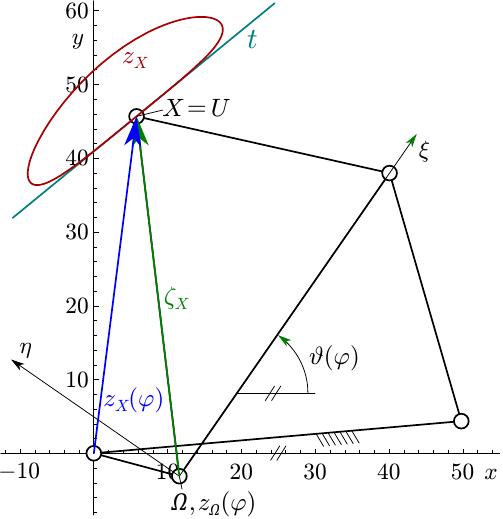}
\caption{Four-bar linkage from Fig.\ \ref{Abb:Curvature_cubic} with (part/component of) coupler curve $z_X$, \textsc{Ball}'s point $U$, and tangent $t$ at $X = U$}
\label{Abb:Coupler_curve_with_Ball_point}
\end{SCfigure}
\end{example} 
 
We now search for the geometric location of the points in the moving plane whose second derivative is perpendicular to the path curve, i.e.\ whose tangential acceleration is zero.
For these points we have
\beq
  \left\langle\ep'(z-z_{P_1}),\ep''(z-z_{P_2})\right\rangle = 0\,.
\eeq
This is the equation of a circle $k_2$ with center point
\beq
  z_2
= \frac{1}{2} \left(z_{P_1} + z_{P_2} - \ii\,\frac{\vt'^2}{\vt''}\left(z_{P_1} - z_{P_2}\right)\right)  
\eeq
and radius square
\beq
  r_2^2
= \frac{\vt'^4+\vt''^2}{4\vt''^2} \left|z_{P_1}-z_{P_2}\right|^2. 
\eeq
Clearly, the poles $P_1$ and $P_2$ lie on $k_2$.
We consider the vector $z_2(\ph) - z_{P_1}(\ph)$.
It holds
\begin{align*}
  z_2 - z_{P_1}
= {} & \frac{1}{2} \left(z_{P_1} + z_{P_2} - \ii\,\frac{\vt'^2}{\vt''}\left(z_{P_1} - z_{P_2}\right)\right) - z_{P_1}\db\\  
= {} & \frac{1}{2} \left(-(z_{P_1} - z_{P_2}) - \ii\,\frac{\vt'^2}{\vt''}\left(z_{P_1} - z_{P_2}\right)\right)\db\\  
= {} & \frac{1}{2} \left(1 + \ii\,\frac{\vt'^2}{\vt''}\right) \left(z_{P_2} - z_{P_1}\right)\db\\
= {} & \frac{1}{2} \frac{\vt'' + \ii\vt'^2}{\vt''} \left(z_{P_2} - z_{P_1}\right)\,.
\end{align*}
Comparison with \eqref{Eq:velocity_of_P_1_c} shows that
\beq
  z_2(\ph) - z_{P_1}(\ph)
= \frac{1}{2}\, \frac{\vt'(\ph)}{\vt''(\ph)}\, u(\ph)\,,
  \quad\mbox{hence}\quad
  z_2(\ph)
= z_{P_1}(\ph) + \frac{1}{2}\, \frac{\vt'(\ph)}{\vt''(\ph)}\, u(\ph)    
\eeq
(cf.\ \cite[pp.\ 36-37]{Blaschke&Mueller}).
Since $z_2$ lies on the straight line defined by the point $P_1$ ($z_{P_1}$) and the direction $u$, and the circle $k_2$ passes through $P_1$, the circle $k_2$ in $P_1$ touches the straight line defined by the point $P_1$ ($z_{P_1}$) and the direction $\ii u$.
The first line and the second line are respectively the tangent $t$ and the normal $n$ to the fixed centrode at $P_1$.
The circle $k_2$ is called {\em stationary circle}\index{circle!stationary circle@stationary $\sim $}.\footnote{In German: ``Tangentialkreis'' and also ``Gleichenkreis'' (see e.g.\ \cite[p.\ 37]{Blaschke&Mueller}).}
The circles $k_1$ and $k_2$ are called {\em \textsc{Bresse}'s\,\footnote{\textsc{Jaques Antoine Charles Bresse}, 1822-1883. According to \cite[p.\ 37]{Blaschke&Mueller}, \textsc{Bresse} considered these circles in 1853.}\,circles}.
Clearly, the radius of $k_2$ can also be written as
\beq
  r_2(\ph)
= \left|z_2(\ph) - z_{P_1}(\ph)\right|
= \frac{1}{2} \left|\frac{\vt'(\ph)}{\vt''(\ph)}\right| \left|u(\ph)\right|.   
\eeq
In contrast to the inflection circle $k_1$, the stationary circle $k_2$ depends on the parametrization.
This follows from the fact, the $k_1$ and $k_2$ intersect at $P_1$ and $P_2$, and $P_2$ is parameter-dependent (see \eqref{Eq:tilde{z}_{P_2}(t)} and also \cite[p.\ 37]{Blaschke&Mueller}).

Now, we are looking for the geometric location of the points in the moving plane whose third derivative (jerk) is perpendicular to the path curve (and the tangent vector).
For these points we have
\beq
  \left\langle\ep'(z-z_{P_1}),\ep'''(z-z_{P_3})\right\rangle = 0\,.
\eeq
This is the equation of a circle $k_4$ with center point
\beq
  z_4
= \frac{1}{2} \left(z_{P_1} + z_{P_3} - \ii\,\frac{3\vt'\vt''}{\vt'''-\vt'^3}\left(z_{P_1} - z_{P_3}\right)\right)  
\eeq
and radius square
\beq
  r_4^2
= \frac{\vt'^6+9\vt'^2\vt''^2-2\vt'^3\vt'''+\vt'''^2}{4\left(\vt'^3-\vt'''\right)^2} \left|z_{P_1}-z_{P_3}\right|^2. 
\eeq
This circle is called {\em zero-tangential jerk circle} \index{circle!zero-tangential jerk circle@zero-tangential jerk $\sim $}(cf.\ \cite{Figliolini&Lanni&Tomassi}).

\begin{example}
We consider the example four-bar linkage $A_0ABB_0$ from \cite{Goessner:Balls_point_revisited-again} (see Fig.\ \ref{Abb:Bresse_circles01}).
Let $\overline{A_0A}$ be the drive bar with angle $\ph$ and constant angular velocity $\ph'(t) = 1$, where $t$ is the time.
Calling the subroutine {\em DyadC} (see p.\ \pageref{Subroutine-code}) with
\beq
  \texttt{DyadC[l3, l4, zA, zA1, zA2, zA3, zB0, 0, 0, 0, -1]}\,,
\eeq
where
\begin{gather*}  
  \texttt{l3} = \overline{AB} = \sqrt{2}\,,\quad
  \texttt{l4} = \overline{B_0B} = 2\,,\\
  \texttt{zA} = \ee^{\ii\ph}\,,\quad 
  \texttt{zA1} = \ii\ee^{\ii\ph}\,,\quad 
  \texttt{zA2} = -\ee^{\ii\ph}\,,\quad
  \texttt{zA3} = -\ii\ee^{\ii\ph}\,,\\
  \texttt{zB0} = 3 + 2\ii\,,\quad
  \ph = \pi/2\,,
\end{gather*}
delivers 
\beq
  \exp(\ii\vt(\ph))= \frac{1+\ii}{\sqrt{2}}\,,\quad
  \vt(\ph) = \frac{\pi}{4}\,,\quad
  \vt'(\ph) = -1\,,\quad
  \vt''(\ph) = -\frac{3}{2}\,,\quad
  \vt'''(\ph) = -\frac{39}{4}\,.
\eeq
It follows
\begin{gather*}
  z_{P_1}(\ph) = 2\,\ii\,,\quad
  z_{P_2}(\ph) = -\frac{6}{13} + \frac{9}{13}\,\ii \approx -0.461538 + 0.692308\,\ii\,,\db\\
  z_{P_3}(\ph) = \frac{72}{1549} + \frac{1409}{1549}\,\ii  \approx 0.046482 + 0.909619\,\ii\,,\db\\ 
  z_1(\ph) = \frac{3}{4} + \ii\,,\;\;
  r_1 = \frac{5}{4}\,,\quad
  z_2(\ph) = -\frac{2}{3} + \frac{3}{2}\,\ii\,,\;\;
  r_2(\ph) = \frac{5}{6}\,,\db\\   
  z_3(\ph) = \frac{13}{12} + \frac{3}{2}\,\ii \approx 1.08333 + 1.50000\,\ii\,,\;\; 
  r_3(\ph) = \frac{\sqrt{205}}{12} \approx 1.19315\,,\db\\
  z_4(\ph) = -\frac{9}{35} + \frac{101}{70}\,\ii \approx -0.25714 + 1.44286\,\ii\,,\;\;
  r_4(\ph) = \frac{3}{14}\, \sqrt{\frac{41}{5}} \approx 0.613621\,,\db\\
  z_U(\ph) = \frac{51}{26} + \frac{9}{13}\,\ii \approx 1.96154 + 0.69231\,\ii\quad\mbox{(cf.\ \cite{Goessner:Balls_point_revisited-again})}\,.
\end{gather*}
This gives the situation in Fig.\ \ref{Abb:Bresse_circles01}.

\begin{figure}[ht]
\begin{minipage}{0.485\textwidth}
\includegraphics[width=\textwidth]{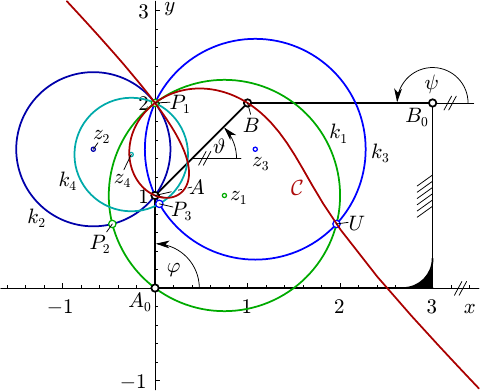}
\caption{Inflection circle $k_1$, stationary circle $k_2$, zero-normal jerk circle $k_3$, zero-tangential jerk circle $k_4$ and stationary curvature cubic $\K$ if drive point is $0$}
\label{Abb:Bresse_circles01}
\end{minipage}
\hfill
\begin{minipage}{0.485\textwidth}
\includegraphics[width=\textwidth]{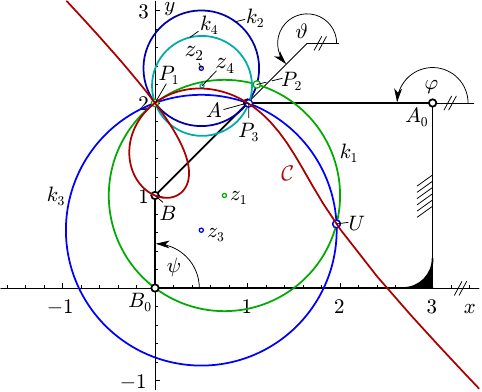}
\caption{Inflection circle $k_1$, stationary circle $k_2$, zero-normal jerk circle $k_3$, zero-tangential jerk circle $k_4$ and stationary curvature cubic $\K$ if drive point is $2 + 3\ii$}
\label{Abb:Bresse_circles02}
\end{minipage}
\end{figure}

Now, we consider the same four-bar linkage as in Fig.\ \ref{Abb:Bresse_circles01}, but we change labels (see Fig.\ \ref{Abb:Bresse_circles02}).
Let the new $\overline{A_0A}$ be the drive bar with angle $\ph$ and constant angular velocity $\ph'(t) = 1$.
Calling the subroutine {\em DyadC} with
\beq
  \texttt{DyadC[l3, l4, zA, zA1, zA2, zA3, zB0, 0, 0, 0, 1]}\,,
\eeq
where (see Fig.\ \ref{Abb:Bresse_circles02})
\begin{gather*}  
  \texttt{l3} = \overline{AB} = \sqrt{2}\,,\quad
  \texttt{l4} = \overline{B_0B} = 1\,,\\
  \texttt{zA} = 3 + 2\ii + 2\ee^{\ii\ph}\,,\quad 
  \texttt{zA1} = 2\ii\ee^{\ii\ph}\,,\quad 
  \texttt{zA2} = -2\ee^{\ii\ph}\,,\quad
  \texttt{zA3} = -2\ii\ee^{\ii\ph}\,,\\
  \texttt{zB0} = 0\,,\quad
  \ph = \pi\,,
\end{gather*}
delivers 
\beq
  \exp(\ii\vt(\ph))= -\frac{1+\ii}{\sqrt{2}}\,,\quad
  \vt(\ph) = -\frac{3\pi}{4}\,,\quad
  \vt'(\ph) = -2\,,\quad
  \vt''(\ph) = 8\,,\quad
  \vt'''(\ph) = -138\,.
\eeq
It follows
\begin{gather*}
  z_{P_2}(\ph) = \frac{11}{10} + \frac{11}{5}\,\ii = 1.1 + 2.2\,\ii\,,\quad
  z_{P_3}(\ph) = \frac{4866}{4801} + \frac{9578}{4801}\,\ii \approx 1.01354 + 1.99500\,\ii\,,\db\\
  z_2(\ph) = \frac{1}{2} + \frac{19}{8}\,\ii = 0.5 + 2.375\,\ii\,,\;\;
  r_2(\ph) = \frac{5}{8} = 0.625\,,\db\\
  z_3(\ph) = \frac{1}{2} + \frac{5}{8}\,\ii = 0.5 + 0.625\,\ii\,,\;\;
  r_3(\ph) = \frac{\sqrt{137}}{8} \approx 1.46309\,,\db\\
  z_4(\ph) = \frac{33}{65} + \frac{142}{65}\,\ii \approx 0.50769 + 2.18462\,\ii\,,\;\;
  r_4(\ph) = \frac{3\,\sqrt{137}}{65} \approx 0.540217\,.
\end{gather*}
Results for $z_{P_1}(\ph)$, $z_1(\ph)$, $r_1(\ph)$ and $z_U(\ph)$ are the same as before.
The poles $P_2$, $P_3$, the stationary circle $k_2$ and the zero-normal jerk circle $k_3$ have changed (see Fig.\ \ref{Abb:Bresse_circles02}).
While the stationary circle $k_2$ passes through point $z_A = \ii$ in Fig.\ \ref{Abb:Bresse_circles01}, it passes through point $z_A = 1 + 2\ii$ in Fig.\ \ref{Abb:Bresse_circles02}.
The reason for this is clear.

By applying \eqref{Eq:tilde{z}_{P_2}(t)}, it is also possible to get the position of $P_2$ in Fig.\ \ref{Abb:Bresse_circles02} from the situation in Fig.\ \ref{Abb:Bresse_circles01} if we consider $\ph$ there as function $f$ of (the independent variable) $\psi$, $\ph = f(\psi)$.
By choosing $X := \Om$, \eqref{Eq:tilde{z}_{P_2}(t)} gives
\beqn \label{Eq:tilde{z}_{P_2}(psi)}
  \tilde{z}_{P_2}(\psi)
= z_\Om(\ph) - \frac{\strut z_\Om''(\ph)\,f'(\psi)^2 + z_\Om'(\ph)\,f''(\psi)}
	{\strut\ii\left(\vt''(\ph)\,f'(\psi)^2 + \vt'(\ph)\,f''(\psi)\right) - \vt'(\ph)^2\,f'(\psi)^2}\,. 
\eeqn
For the inverse function $g := f^{-1}$ we have
\beq
  \psi = g(\ph)\,,\qquad
  f'(\psi) = \frac{1}{g'(\ph)}\,,\qquad
  f''(\psi) = -\frac{g''(\ph)}{g'(\ph)^3}\,,
\eeq
thus \eqref{Eq:tilde{z}_{P_2}(psi)} becomes
\beq
  \tilde{z}_{P_2}(\psi)
= z_\Om(\ph) - \frac{\dfrac{z_\Om''(\ph)}{g'(\ph)^2} - \dfrac{g''(\ph)\,z_\Om'(\ph)}{g'(\ph)^3}}
	{\ii\left(\dfrac{\vt''(\ph)}{g'(\ph)^2} - \dfrac{g''(\ph)\,\vt'(\ph)}{g'(\ph)^3}\right) - \left(\dfrac{\vt'(\ph)}{g'(\ph)}\right)^2}\,. 
\eeq
With
\begin{gather*}
  z_\Om(\pi/2) = \ii\,,\quad
  z_\Om'(\pi/2) = -1\,,\quad
  z_\Om''(\pi/2) = -\ii\,,\\[0.05cm]
  \psi = g(\pi/2) = \pi\,,\quad
  g'(\pi/2) = 1/2\,,\quad
  g''(\pi/2) = 7/4\,,\\[0.05cm]
  \vt'(\pi/2) = -1\,,\quad
  \vt''(\pi/2) = -3/2
\end{gather*}
follows
\beq
  \tilde{z}_{P_2}(\pi)
= \frac{11}{10} + \frac{11}{5}\,\ii\,,  
\eeq
which is the position of $P_2$ in Fig.\ \ref{Abb:Bresse_circles02}.\hfill\bs
\end{example}

Of course, \eqref{Eq:velocity_of_P_1_c} can be used to compute the pole $z_{P_2}$ for known velocity $u$ of $z_{P_1}$,
\beqn \label{Eq:z_{P_2}_with_u_a}
  z_{P_2}
= z_{P_1} + \frac{\vt'}{\vt''+\ii\vt'^2}\,u
= z_{P_1} + \frac{\vt'\left(\vt''-\ii\vt'^2\right)}{\vt''^2+\vt'^4}\,u\,.  
\eeqn
We also want to determine the position $Z_{P_2}$ of $P_2$ in the above used $X,Y$-coordinate system whose origin is at $z_{P_1}$ and whose $X$-axis has the direction of $u$.
With $u = |u|\,\ee^{\ii\arg(u)}$, \eqref{Eq:z_{P_2}_with_u_a} becomes
\beqn \label{Eq:z_{P_2}_with_u_b}
  z_{P_2}
= z_{P_1} + \frac{\vt'\left(\vt''-\ii\vt'^2\right)}{\vt''^2+\vt'^4}\,|u|\,\ee^{\ii\arg(u)}\,.  
\eeqn
Comparing \eqref{Eq:z_{P_2}_with_u_b} with the second transformation equation in \eqref{Eq:z->Z_Z->z} shows that
\beq
  Z_{P_2}(\ph)
= \frac{\vt'(\ph)\left(\vt''(\ph)-\ii\vt'(\ph)^2\right)}{\vt''(\ph)^2+\vt'(\ph)^4}\,|u(\ph)|
\eeq
with real and imaginary part
\beq
  X_{P_2}(\ph)
= \frac{\vt'(\ph)\,\vt''(\ph)\,|u(\ph)|}{\vt''(\ph)^2+\vt'(\ph)^4}
  \quad\mbox{and}\quad
  Y_{P_2}(\ph)
= -\frac{\vt'(\ph)^3\,|u(\ph)|}{\vt''(\ph)^2+\vt'(\ph)^4}\,,
\eeq
respectively (see also \cite[p.\ 34]{Blaschke&Mueller}).

%% file: DiffGeo5_4e.tex

\subsection{Cam mechanisms with translating flat-face follower} 
\subsubsection{Basics}

\begin{figure}[ht]
\begin{minipage}{0.56\textwidth}
\includegraphics[width=\textwidth]{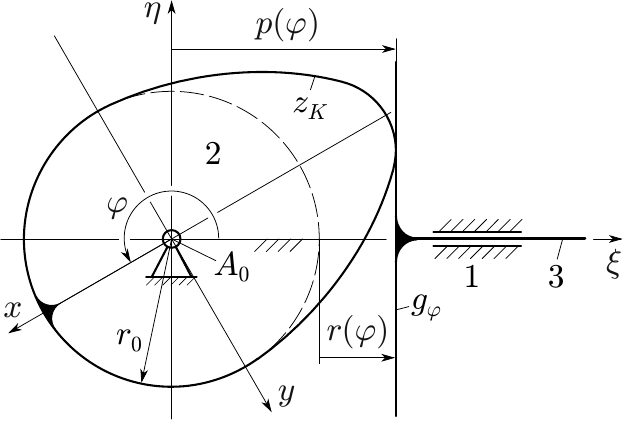}
\caption{Cam mechanism with translating flat-face follower}
\label{Abb:Flat_face_follower02}
\end{minipage}
\hfill
\begin{minipage}{0.40\textwidth}
\includegraphics[width=\textwidth]{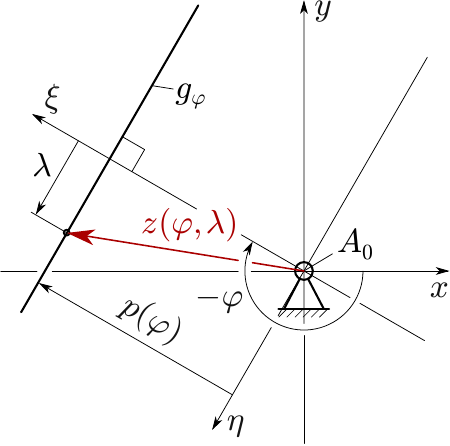}
\caption{Kinematic inversion}
\label{Abb:Flat_face_follower03}
\end{minipage}
\end{figure}

In this subsection we consider cam mechanisms consisting of a fixed frame 1 to which we assign a fixed $\xi,\eta$-coordinate system, a cam 2 with contour curve $z_K$, and a follower 3 with a flat face which is represented by the line $g_\ph$ (see Fig.\ \ref{Abb:Flat_face_follower02}).
The cam rotates around the origin $A_0$ of the $\xi,\eta$-system with angle $\ph$ (starting from the positive $\xi$-axis) in positive direction.
The follower is movable along the $\xi$-axis while the line $g_\ph$ is perpendicular to this axis and assumed to be always in contact with the contour curve $z_K$ of the cam. 
Thus the rotating motion of the cam is transformed into a linear motion of the follower.
The motion of the follower is a $2\pi$-periodic function $r(\ph)$ which is called {\em transfer function}\index{tranfer function} of the cam mechanism. 

In general, one specifies the transfer function $r(\ph)$ and uses it to determine the contour curve $z_K$ of the cam.
Without loss of generality we can choose the transfer function $r(\ph)$ in such a way that $r(\ph) \ge 0$ holds and there is at least one value of $\ph$ with $r(\ph) = 0$.
Then the chosen distance of line $g_\ph$ of the follower from $A_0$ for $r(\ph) = 0$ determines the radius $r_0$ of a circle with center at $A_0$, which is called {\em base circle}\index{base circle}, and the distance of the follower from the coordinate origin $A_0$ is given by 
\beq
  p(\ph)
= r_0 + r(\ph)\,.
\eeq
In the following subsubsection we define a transfer function $r(\ph)$ which will be used for examples afterwards.

\subsubsection{Transfer function}

Let the transfer function of the cam mechanism be given by
\beqn \label{Eq:r(phi)}
  r(\ph)
= \left\{\begin{array}{c@{\quad\;\mbox{if}\quad\;}l}
	18 f_{3,2}\!\left(\dfrac{\ph}{5\pi/6}\right) & 0 \le \ph \le \dfrac{5\pi}{6}\,,\\[0.25cm]
	18 & \dfrac{5\pi}{6} < \ph < \dfrac{7\pi}{6}\,,\\[0.25cm]
	18 \left[1-g_3\!\left(\dfrac{\ph-7\pi/6}{2\pi/3}\right)\right] & \dfrac{7\pi}{6} \le \ph \le \dfrac{11\pi}{6}\,,\\[0.3cm]
	0 & \dfrac{11\pi}{6} < \ph \le 2\pi 
  \end{array}\right.  
\eeqn
with
\beqn \label{Eq:f_{a,b}}
  \fab(x) := I(x;a+1,b+1)  
\eeqn
and
\beqn \label{Eq:g_a}
  g_a(x) := I\left(\sin^2\frac{\pi x}{2};\frac{a+1}{2},\frac{a+1}{2}\right)
\eeqn
where $I(\,\cdot\,;p,q)$ is the {\em regularized incomplete beta function}\footnote{It is available in {\em Mathematica} as \texttt{BetaRegularized[x,\,p,\,q]}.} \index{BetaRegularized@\texttt{BetaRegularized}} defined by
\beqn \label{Eq:I(x;p,q)}
  I(x;p,q)
= \frac{\int_0^x t^{p-1}(1-t)^{q-1}\, \dd t}{\int_0^1 t^{p-1}(1-t)^{q-1}\, \dd t}  
\eeqn
(see \textcite[Vol.\ 2, p.\ 346]{Gradstein&Ryshik-engl&kurz}).\footnote{For further applications of the functions \eqref{Eq:f_{a,b}} and \eqref{Eq:g_a} for cam mechanisms see \cite{Baesel:MMT2023} and \cite{Baesel:MeTrApp2023}, respectively.}

\begin{figure}[ht]
\begin{minipage}{0.48\textwidth}
\includegraphics[width=\textwidth]{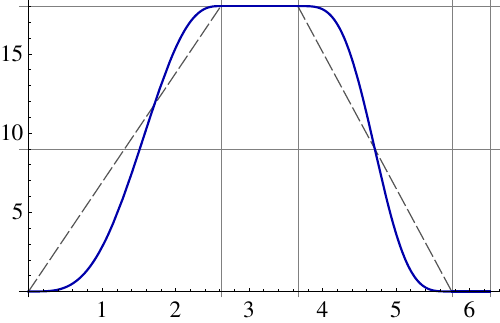}
\caption{Graph of the function $r$ defined by \eqref{Eq:r(phi)}} 
\label{Abb:r(phi)}
\end{minipage}
\hfill
\begin{minipage}{0.48\textwidth}
\includegraphics[width=\textwidth]{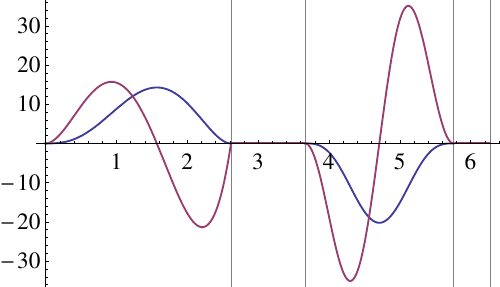}
\caption{$r'$ (blue) and $r''$ (purple) according to \eqref{Eq:r'(phi)} and \eqref{Eq:r''(phi)}}
\label{Abb:r'(phi)_and_r''(phi)}
\end{minipage}
\end{figure}

Evaluation of \eqref{Eq:f_{a,b}} with $a = 3$ and $b = 2$ yields
\beq
  f_{3,2}(x) = 15 x^4 - 24 x^5 + 10 x^6
\eeq
with derivatives
\beq
  f_{3,2}'(x) = 60 x^3 - 120 x^4 + 60 x^5\,,\qquad
  f_{3,2}''(x) = 180 x^2 - 480 x^3 + 300 x^4\,.
\eeq
Evaluation of \eqref{Eq:g_a} with $a = 3$ gives
\beq
  g_3(x) = \frac{1}{2} - \frac{9\cos(\pi x)}{16} + \frac{\cos(3\pi x)}{16}
\eeq
with derivatives
\beq
  g_3'(x) = \frac{9\pi\sin(\pi x)}{16} - \frac{3\pi\sin(3\pi x)}{16}\,,\qquad
  g_3''(x) = \frac{9\pi^2\cos(\pi x)}{16} - \frac{9\pi^2\cos(3\pi x)}{16}\,.	
\eeq
 
The graph of $r(\ph)$ is shown in Fig.\ \ref{Abb:r(phi)}.
The first and second derivative of $r(\ph)$ are
\beqn \label{Eq:r'(phi)}
  r'(\ph)
= \left\{\begin{array}{c@{\quad\;\mbox{if}\quad\;}l}
	\dfrac{18}{5\pi/6}\, f_{3,2}'\!\left(\dfrac{\ph}{5\pi/6}\right) & 0 \le \ph \le \dfrac{5\pi}{6}\,,\\[0.25cm]
	0 & \dfrac{5\pi}{6} < \ph < \dfrac{7\pi}{6}\,,\\[0.25cm]
	-\dfrac{18}{2\pi/3}\: g_3'\!\left(\dfrac{\ph-7\pi/6}{2\pi/3}\right) & \dfrac{7\pi}{6} \le \ph \le \dfrac{11\pi}{6}\,,\\[0.3cm]
	0 & \dfrac{11\pi}{6} < \ph \le 2\pi\,, 
  \end{array}\right.  
\eeqn
and
\beqn \label{Eq:r''(phi)}
  r''(\ph)
= \left\{\begin{array}{c@{\quad\;\mbox{if}\quad\;}l}
	\dfrac{18}{(5\pi/6)^2}\, f_{3,2}''\!\left(\dfrac{\ph}{5\pi/6}\right) & 0 \le \ph \le \dfrac{5\pi}{6}\,,\\[0.25cm]
	0 & \dfrac{5\pi}{6} < \ph < \dfrac{7\pi}{6}\,,\\[0.25cm]
	\dfrac{18}{(2\pi/3)^2}\: g_3''\!\left(\dfrac{\ph-7\pi/6}{2\pi/3}\right) & \dfrac{7\pi}{6} \le \ph \le \dfrac{11\pi}{6}\,,\\[0.3cm]
	0 & \dfrac{11\pi}{6} < \ph \le 2\pi\,, 
  \end{array}\right.  
\eeqn
respectively (see Fig.\ \ref{Abb:r'(phi)_and_r''(phi)}).
As explicit representation of the function $r$ (see \eqref{Eq:r(phi)}) one easily gets
\beqn \label{Eq:r(phi)_expl}
  r(\ph)
= \left\{\!\begin{array}{c@{\quad\;\mbox{if}\quad\;}l}
	\dfrac{69\,984\ph^4}{125\pi^4} - \dfrac{3\,359\,232\ph^5}{3125\pi^5} + \dfrac{1\,679\,616\ph^6}{3125\pi^6} & 0 \le \ph \le \dfrac{5\pi}{6}\,,\\[0.3cm]
	18 & \dfrac{5\pi}{6} < \ph < \dfrac{7\pi}{6}\,,\\[0.25cm]
	9 +	\dfrac{81}{8}\cos\left(\dfrac{3\ph}{2}+\dfrac{\pi}{4}\right)
	  + \dfrac{9}{8}\sin\left(\dfrac{9\ph}{2}+\dfrac{\pi}{4}\right) & \dfrac{7\pi}{6} \le \ph \le \dfrac{11\pi}{6}\,,\\[0.35cm]
	0 & \dfrac{11\pi}{6} < \ph \le 2\pi\,. 
  \end{array}\right.  
\eeqn

\subsubsection{The contour curve}

Our goal in this paragraph is to determine the contour curve $z_K$ of the cam for given transfer function $r(\ph)$ and given base circle radius $r_0$.

\begin{figure}[ht]
\begin{minipage}{0.48\textwidth}
\includegraphics[width=\textwidth]{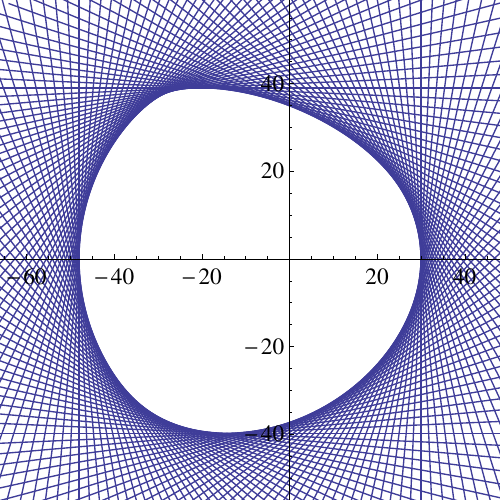}
\caption{144 lines $g_\ph$ with equation $z(\ph,\la)$} 
\label{Abb:Konturkurve_als_Huellkurve}
\end{minipage}
\hfill
\begin{minipage}{0.48\textwidth}
\includegraphics[width=\textwidth]{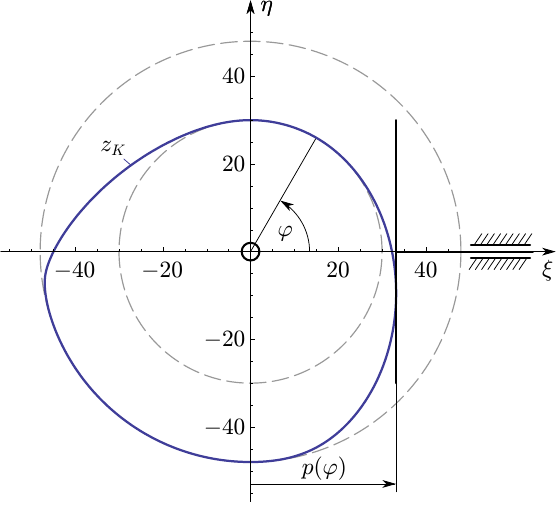}
\caption{Cam mechanism with contour curve~$z_K$ from Fig.\ \ref{Abb:Konturkurve_als_Huellkurve}}
\label{Abb:Nocken_unsymmetr}
\end{minipage}
\end{figure}

In the plane of the cam contour curve $z_K$ to be generated, we introduce an $x,y$-coordinate system with origin in $A_0$, where for $\ph = 0$ the $x$-axis coincides with the $\xi$-axis.
For the generation of the contour curve we use the method of kinematic inversion\index{kinematic inversion}\index{inversion, kinematic} (see Fig. \ref{Abb:Flat_face_follower03}). 
To do this, we fix the $x,y$-system and let the $\xi,\eta$-system with the line $g_\ph$ rotate around $A_0$ in negative (clockwise) direction.
The contour curve $z_K$ is then the envelope\index{envelope} of the family of lines $g_\ph$, $0 \le \ph < 2\pi$ (see Fig. \ref{Abb:Konturkurve_als_Huellkurve}).
A very simple parametric equation of the line $g_\ph$ in the $\xi,\eta$-system is given by
\beq
  \zeta(\ph,\la)
= p(\ph) + \la\,\ii\,,
  \quad \la\in\R\,;
\eeq
hence in the $x,y$-system by
\beqn \label{Eq:z(phi,lambda)}
  z(\ph,\la)
= [p(\ph) + \la\,\ii]\, \ee^{-\ii\ph}
= [r_0 + r(\ph) + \la\,\ii]\, \ee^{-\ii\ph}\,,  
  \quad \la\in\R\,.
\eeqn
In order to obtain the equation of the envelope, we need to determine $\la$ as a function of $\ph$.

\begin{cor} The equation of the contour curve $z_K$ is given by
\beqn \label{Eq:Huellkurvengleichung2}
  z_K(\ph) = z(\ph,\la(\ph))
= [p(\ph) - \ii\,p'(\ph)]\, \ee^{-\ii\ph}\,,\quad
  0 \le \ph < 2\pi\,,   
\eeqn
with $p(\ph) = r_0 + r(\ph)$
\end{cor}

\begin{proof}
Using Theorem \ref{Thm:Envelope} and calculation rules from \eqref{Eq:Rechenregeln2}, from \eqref{Eq:z(phi,lambda)} we get
\begin{align*}
  \left[\frac{\p z}{\p\ph}, \frac{\p z}{\p\la}\right]
= {} & \left[\left(p'(\ph) + \la - \ii p(\ph)\right)\ee^{-\ii\ph}, \ii\ee^{-\ii\ph}\right] 
= \left[p'(\ph) + \la - \ii p(\ph), \ii\right]\db\\
= {} & \left[p'(\ph)+\la,\ii\right] - \big[\ii p(\ph),\ii\big]
= \left(p'(\ph)+\la\right) \left[1,\ii\right] - p(\ph) \big[\ii,\ii\big]\\[0.05cm]
= {} & p'(\ph) + \la
= 0\,,
\end{align*}
thus
\beq
  \la = \la(\ph) = -p'(\ph)\,.
\eeq
Therefore, the equation of the contour curve (envelope) is
\beq
  z_K(\ph) = z(\ph,\la(\ph))
= [p(\ph) - \ii\,p'(\ph)]\, \ee^{-\ii\ph}\,,\quad
  0 \le \ph < 2\pi\,. \qedhere   
\eeq
\end{proof}

Note that $\zeta_K(\ph) = r_0 + r(\ph) - \ii\,r'(\ph)$ is the current contact point between the contour curve and the line $g_\ph$ in the $\xi,\eta$-coordinate system.\footnote{For the determination of the contour curve for a given transfer function using a known method to find the envelope by means of real-valued functions see \textcite[pp.\ 246-248]{Wunderlich}. Our derivation is well comprehensible in every step and due to its vectorial nature it allows to draw conclusions more easily.}

For each point of the line $g_\ph$ (see Fig. \ref{Abb:Flat_face_follower03}), as can be easily seen, the corresponding (unit) vector $z_\la$ lies on $g_\ph$. 
This means that in the point of $g_\ph$ which just generates the cam contour (thus is tangent to it) the vector $z_\ph$ must also lie on $g_\ph$. 

\begin{example} \label{Bsp:Nocken}
Fig.\ \ref{Abb:Konturkurve_als_Huellkurve} shows the lines $g_\ph$, $\ph = 2\pi k/144$, $k = 0,1,2,\ldots,143$, according to equation $z(\ph,\la)$ from \eqref{Eq:z(phi,lambda)} with $r_0 = 30$ and $r(\ph)$ from \eqref{Eq:r(phi)}.
Fig.\ \ref{Abb:Nocken_unsymmetr} shows the corresponding cam mechanism with contour curve $z_K$ according to \eqref{Eq:Huellkurvengleichung2} in position $\ph = \pi/3$, $p(\ph) = 30 + r(\pi/3)$ $= 30 + 2016/625 = 33.2256$. \phantom{}\hfill\bs     
\end{example}

\subsubsection{Perimeter and area}

Now we derive a formula for the calculation of the arc length $L(\ph_0,\ph)$ of the contour curve $z_K$ between two points with the parameter values $\ph_0$ and $\ph$. 
From \eqref{Eq:Huellkurvengleichung2} follows
\begin{align} \label{Eq:zK'(phi)}
  z_K'(\ph)
= {} & \left[p'(\ph)-\ii\,p''(\ph)\right]\ee^{-\ii\ph}
- \left[p(\ph)-\ii\,p'(\ph)\right]\ii\,\ee^{-\ii\ph}\nonumber\\[0.05cm]
= {} & \left[p'(\ph)-\ii\,p''(\ph)-\ii\,p(\ph)-p'(\ph)\right]\ee^{-\ii\ph}\nonumber\\[0.05cm]
= {} & {-}\ii\left[p(\ph)+p''(\ph)\right]\ee^{-\ii\ph}\,, 
\end{align}
hence
\beq
  L(\ph_0,\ph)
= \int_{\ph_0}^\ph \big|p(\phi)+p''(\phi)\big|\: \dd\phi\,.
\eeq
with
\beq
  p(\ph) = r_0 + r(\ph)\,,\quad 
  p'(\ph) = r'(\ph)\,,\quad
  p''(\ph) = r''(\ph)\,.
\eeq
If the envelope encloses a convex set (which must be the case for a cam with a flat-face follower) and the coordinate origin $O$ is an interior point of this set, then $p = p(\ph)$ is called the {\em support function}\index{support function} of this set with respect to $O$.
A necessary and sufficient condition for the periodic function $p$ to be the support function of a convex set is
\beqn \label{Eq:p(phi)+p''(phi)>0}
  p(\ph) + p''(\ph) > 0\,,\quad 0 \le \ph \le 2\pi\,.
  \quad\mbox{\autocite[p.\ 3]{Santalo}}
\eeqn
Thus we get
\beqn \label{Eq: L(phi_0,phi)}
  L(\ph_0,\ph)
= \int_{\ph_0}^\ph \left(p(\phi)+p''(\phi)\right) \dd\phi\,.
\eeqn
For the perimeter\index{perimeter} ($=$ total length of the convex contour curve) of the cam holds
\begin{align*}
  L
= {} & \int_0^{2\pi} \left(p(\ph) + p''(\ph)\right) \dd\ph
= \int_0^{2\pi} p(\ph)\, \dd\ph + \int_0^{2\pi} p''(\ph)\, \dd\ph\db\\[0.05cm]
= {} & \int_0^{2\pi} p(\ph)\, \dd\ph + p'(2\pi) - p'(0)
= \int_0^{2\pi} p(\ph)\, \dd\ph
\end{align*}
This is the Cauchy-Crofton\footnote{Augustin-Louis Cauchy, 1789-1857\index{Cauchy, Augustin-Louis}; Morgan Crofton, 1826-1915} formula\index{Cauchy-Crofton formula} for the perimeter of a convex figure \autocite[pp.\ 1, 12]{Blaschke:Integralgeometrie}, \autocite[slide 15]{Treibergs} -- a special case of Cauchy's surface area formula \autocite[pp.\ 295, 222]{Schneider&Weil2008}.

We calculate the length $L$ of the contour curve $z_K$ with $r(\ph)$ according to \eqref{Eq:r(phi)} and $r_0 = 30$ (see also Fig.\ \ref{Abb:Nocken_unsymmetr}):
\begin{align*}
  L
= {} & \int_0^{2\pi} p(\ph)\, \dd\ph
= \int_0^{2\pi} (r_0 + r(\ph))\, \dd\ph
=  2 \pi r_0 + \int_0^{2\pi} r(\ph)\, \dd\ph\db\\[0.05cm]
= {} & 60 \pi + \frac{129\pi}{7}
= \frac{549\pi}{7}
\approx 246.391\,.  
\end{align*} 

\begin{cor}
The oriented area of the contour curve $z_K$ with parametric equation $z_K(\ph)$ according to \eqref{Eq:Huellkurvengleichung2} is given by 
\beq
  A
= -\frac{1}{2}\int_0^{2\pi} \left(p^2(\ph)-p'^2(\ph)\right) \dd\ph\,.
\eeq
\end{cor}

\begin{proof}
With \eqref{Eq:Huellkurvengleichung2} and \eqref{Eq:zK'(phi)} we have
\begin{align*}
  [z_K(\ph),z_K'(\ph)]
= {} & \left[\left(p(\ph)-\ii\,p'(\ph)\right)\ee^{-\ii\ph},\,-\ii\left(p(\ph)+p''(\ph)\right)\ee^{-\ii\ph}\right]\\[0.05cm]  
= {} & \left[p(\ph)-\ii\,p'(\ph),\,-\ii\left(p(\ph)+p''(\ph)\right)\right]\\[0.05cm]
= {} & \left[p(\ph),\,-\ii\left(p(\ph)+p''(\ph)\right)\right] + \left[\ii\,p'(\ph),\,\ii\left(p(\ph)+p''(\ph)\right)\right]\\[0.05cm]
= {} & {-}p(\ph)\left(p(\ph)+p''(\ph)\right),
\end{align*}
hence
\begin{align} \label{Eq:Flaeche_Nocken}
  A
= -\frac{1}{2}\int_0^{2\pi} p(\ph) \left[p(\ph)+p''(\ph)\right] \dd\ph
= -\frac{1}{2}\left(\int_0^{2\pi} p^2(\ph)\:\dd\ph + \int_0^{2\pi} p(\ph)\,p''(\ph)\:\dd\ph\right).
\end{align}
Partial integration yields
\beq
  \int_0^{2\pi} p(\ph)\,p''(\ph)\:\dd\ph
= p(\ph)\,p'(\ph)\,\Big|_0^{2\pi}
- \int_0^{2\pi} p'^2(\ph)\:\dd\ph
= -\int_0^{2\pi} p'^2(\ph)\:\dd\ph\,.
\eeq
From \eqref{Eq:Flaeche_Nocken} we finally obtain
\beq
  A
= -\frac{1}{2}\int_0^{2\pi} \left[p^2(\ph)-p'^2(\ph)\right] \dd\ph
  \quad\mbox{(cf.\ \textcite[p.\ 4]{Santalo}}\,.\footnote{The negative sign of this integral (and the consequently negative area) results from the fact that the contour curve $z_K$ is negatively oriented. This means that the interior of the curve is always to the right as one moves along $z_K$ with increasing parameter $\ph$.} \qedhere
\eeq
\end{proof}

As an example, for $r(\ph)$ according to \eqref{Eq:r(phi)} and $r_0 = 30$ (see also Fig.\ \ref{Abb:Nocken_unsymmetr}) one easily gets
\beq
  |A|
= \left|\frac{1119744}{175\pi^3} - \frac{1458126909\pi}{1025024}\right|
\approx |{-}4262.65|
= 4262.65\,.
\eeq   

\subsubsection{Curvature}

\begin{cor} \label{Cor:Curvature}
The oriented curvature $\kappa(\ph)$ of the cam contour curve $z_K$ with parametric equation~\eqref{Eq:Huellkurvengleichung2} is given by
\beq
  \kappa(\ph)
= -\frac{1}{p(\ph)+p''(\ph)}\,.  
\eeq
\end{cor}

\begin{proof}
From Theorem \ref{Thm:Curvature} we know that
\beq
  \kappa(\ph)
= \frac{\left[z_K'(\ph),z_K''(\ph)\right]}{\left|z_K'(\ph)\right|^3}\,,
\eeq
and from \eqref{Eq:zK'(phi)} we have
\beq
  z_K'(\ph)
= -\ii\left(p(\ph)+p''(\ph)\right)\ee^{-\ii\ph}\,.
\eeq
It follows
\begin{align*}
  z_K''(\ph)
= {} & {-}\ii\,\Bigl\{\Bigl(p'(\ph)+p'''(\ph)\Bigr)\ee^{-\ii\ph} - \ii\,\Bigl(p(\ph)+p''(\ph)\Bigr)\ee^{-\ii\ph}\Bigr\}\\ 
= {} & {-}\Bigl\{\Bigl(p(\ph)+p''(\ph)\Bigr) + \ii\,\Bigl(p'(\ph)+p'''(\ph)\Bigr)\Bigr\}\,\ee^{-\ii\ph}\,, 
\end{align*}
thus, using the rotation rule from \eqref{Eq:Rechenregeln2}, and \eqref{Eq:[z_1,z_2]_mit_Im},
\begin{align*}
  [z_K'(\ph),z_K''(\ph)]
= {} & \left[-\ii\left(p(\ph)+p''(\ph)\right)\ee^{-\ii\ph},\, -\!\left\{\left(p(\ph)+p''(\ph)\right)+\ii\left(p'(\ph)+p'''(\ph)\right)\right\}\ee^{-\ii\ph}
  \right]\\[0.05cm]  
= {} & \left[\ii\left(p(\ph)+p''(\ph)\right),\, \left(p(\ph)+p''(\ph)\right)+\ii\left(p'(\ph)+p'''(\ph)\right)\right]\db\\[0.05cm]
= {} & \Imz\Bigl\{-\ii\,\Bigl(p(\ph)+p''(\ph)\Bigr) \Bigl(\left(p(\ph)+p''(\ph)\right)+\ii\left(p'(\ph)+p'''(\ph)\right)\Bigr)\Bigr\}\db\\[0.05cm]
= {} & {-}\left(p(\ph)+p''(\ph)\right)^2.
\end{align*}
With $|z_K'(\ph)| = |p(\ph)+p''(\ph)| = p(\ph)+p''(\ph) > 0$ the assertion follows.
\end{proof}

According to Corollary \ref{Cor:Curvature}, $\kappa(\ph)$ is always negative in the $x,y$-coordinate system with the assumed positive rotation direction of the cam.
This means that the tangent vector $z_K'(\ph)$ always turns to the right.
One obtains the result of Corollary \ref{Cor:Curvature} also by means of Lemma \ref{Lem:T'(t)} (see Remark \ref{Bem:kappa(phi)_with_Frenet}).

\begin{remark} \label{Bem:kappa(phi)_with_Frenet}
According to Lemma \ref{Lem:T'(t)} we have
\beq
  \kappa(\ph)
= \frac{T'(\ph)}{\ii\,T(\ph)\,|z_K'(\ph)|}\,.    
\eeq
$T(\ph)$ is independent of the transfer function, as can be easily seen from the Figures \ref{Abb:Flat_face_follower02} and \ref{Abb:Flat_face_follower03}, and we have
\beq
  T(\ph) = -\ii\,\ee^{-\ii\ph}\,.  
\eeq
With $T'(\ph) = -\ee^{-\ii\ph}$ and, see \eqref{Eq: L(phi_0,phi)}, $|z_K'(\ph)| = p(\ph) + p''(\ph)$ we obtain
\beq
  \kappa(\ph)
= \frac{-\ee^{-\ii\ph}}{-\ii^2\,\ee^{-\ii\ph}\,\left(p(\ph)+p''(\ph)\right)}
= -\frac{1}{p(\ph)+p''(\ph)}\,,
\eeq
which is the result of Corollary \ref{Cor:Curvature}.
\end{remark}

The curvature radii $\rh$ of the contour curve are of interest, among other things, for determining the pressures (stresses) between the flat-face follower and the cam.
The curvature radius\index{curvature radius} $\rh$ at point $\ph$ of the contour curve $z_K$ is given by
\beqn \label{Eq:Curvature_radius}
  \rh(\ph)
= \frac{1}{|\kappa(\ph)|}
= \left|p(\ph) + p''(\ph)\right|
= \left|r_0 + r(\ph) + r''(\ph)\right|.
\eeqn
One obtains \eqref{Eq:Curvature_radius} also directly from the definition
\beq
  \rh(\ph)
= \frac{\dd s(\ph)}{\dd \ph}  
\eeq
of the curvature radius, where $\dd s$ is the differential of the arc length $s$ and $\dd\ph$ the corresponding differential of the tangent angle $\ph$.
From \eqref{Eq:Arc_length} we know that
\beq
  \frac{\dd s(\ph)}{\dd \ph}
= \big|z_K'(\ph)\big|\,,  
\eeq
and from \eqref{Eq:zK'(phi)} it follows that
\beq
  \big|z_K'(\ph)\big|
= \left|p(\ph)+p''(\ph)\right|,   
\eeq
and so we have got \eqref{Eq:Curvature_radius} again.

\begin{figure}[H]
\begin{minipage}{0.48\textwidth}
\includegraphics[width=\textwidth]{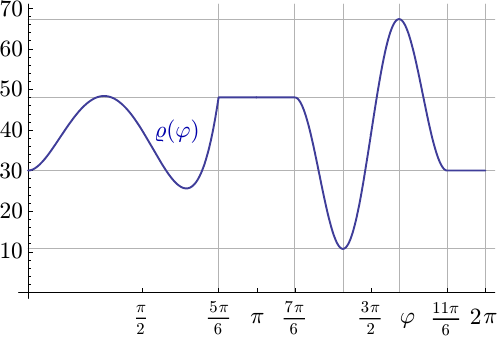}
\caption{Curvature radius function $\rh(\ph)$ of the curve $z_K$ in Figs.\ \ref{Abb:Konturkurve_als_Huellkurve} and \ref{Abb:Nocken_unsymmetr}} 
\label{Abb:Curvature_radius02}
\end{minipage}
\hfill
\begin{minipage}{0.48\textwidth}
\includegraphics[width=\textwidth]{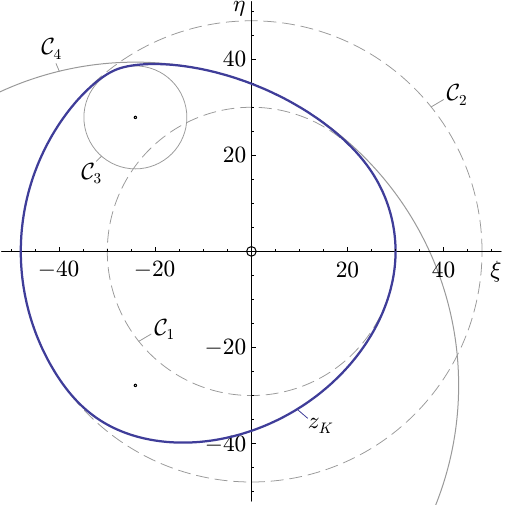}
\caption{Curvature circles}
\label{Abb:Curvature_circles}
\end{minipage}
\end{figure}

\begin{example} 
Fig.\ \ref{Abb:Curvature_radius02} shows the curvature radius function for the contour curve $z_K$ with transfer function $r(\ph)$ according to \eqref{Eq:r(phi)} and base circle radius $r_0 = 30$.
For $5\pi/6 \le \ph \le 7\pi/6$ we have $\rh(\ph) = 30 + 18 = 48$.
$\rh$ takes its absolute minimum $\approx 10.7168$ at point $\ph \approx 4.32691$ ($247.914\g$), and its absolute maximum $\approx 67.2832$ at point $\ph \approx 5.09787$ ($292.086\g$). 
$\rh$ is equal to the base circle radius $30$ for $11\pi/6 \le \ph \le 2\pi$. 
$\mathcal{C}_1$ in Fig.\ \ref{Abb:Curvature_circles} is the base circle, $\mathcal{C}_2$ the curvature circle for $5\pi/6 \le \ph \le 7\pi/6$, $\mathcal{C}_3$ is the smallest curvature circle, and $\mathcal{C}_4$ the largest. \hfill\bs
\end{example}
 
For $\rh(\ph) = 0$ a cusp\index{cusp} of $z_K$ occurs.
If the expression $r_0 + r(\ph) + r''(\ph)$ would change its sign, which does not happen with convex curve as explained, then {\em undercut}\index{undercut} occurs, which practically means that $z_K$ has a loop\index{loop} with cusps, consequently the cam is not generated completely and so the transfer function $r(\ph)$ is not fully realized.

\begin{figure}[ht]
\begin{minipage}{0.48\textwidth}
\includegraphics[width=\textwidth]{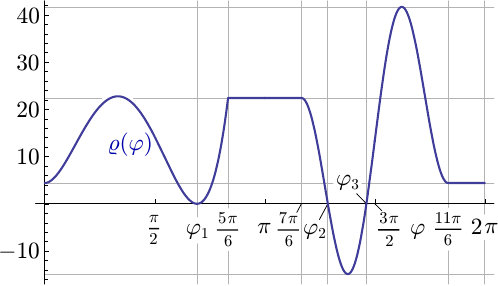}
\caption{Curvature radius function $\rh(\ph)$ of the curve $z_K$ in Figs.\ \ref{Abb:Kurve_als_Huellkurve} and \ref{Abb:Nocken_Unterschnitt}}
\label{Abb:Curvature_radius03}
\end{minipage}
\hfill
\begin{minipage}{0.48\textwidth}
\includegraphics[width=\textwidth]{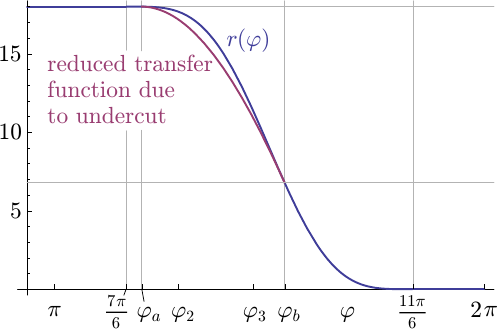}
\caption{Intended transfer function $r$, and reduced transfer function due to undercut}
\label{Abb:Nockengetriebe_Hubverlust}
\end{minipage}
\end{figure}

\begin{figure}[ht]
\begin{minipage}{0.48\textwidth}
\includegraphics[width=0.96\textwidth]{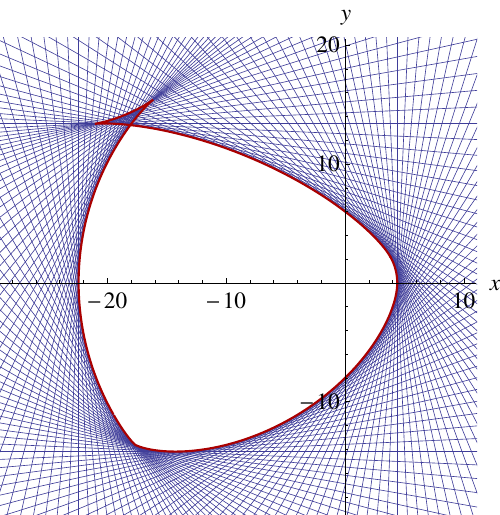}
\caption{Lines generating the curve $z_K$ in Abb.\ \ref{Abb:Nocken_Unterschnitt}} 
\label{Abb:Kurve_als_Huellkurve}
\end{minipage}
\hfill
\begin{minipage}{0.48\textwidth}
\includegraphics[width=0.96\textwidth]{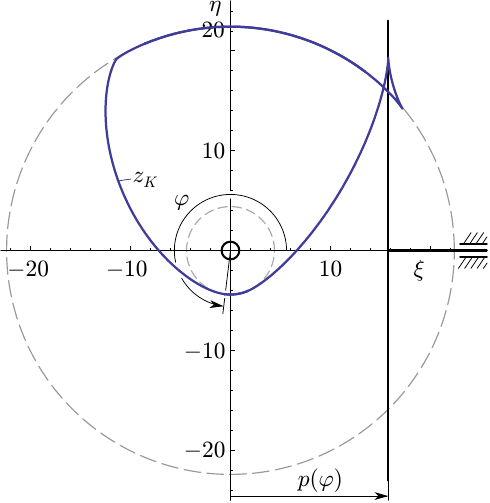}
\caption{Curve $z_K$ with loop and cusps}
\label{Abb:Nocken_Unterschnitt}
\end{minipage}
\end{figure}

\begin{example}
We consider the curvature radius function $\rh(\ph)$ in Fig.\ \ref{Abb:Curvature_radius02} (with $r(\ph)$ according to \eqref{Eq:r(phi)} and $r_0 = 30$).
From \eqref{Eq:Curvature_radius} we see that the necessary condition for a local extremum is $\rh'(\ph) = r'(\ph) + r'''(\ph) = 0$.
Numerical calculation with the {\em Mathematica} function \texttt{FindRoot}\index{FindRoot@\texttt{FindRoot}} shows that the minimum between $\ph = \pi/2$ and $\ph = 5\pi/6$ occurs at point $\ph_1 \approx 2.17476$ $(124.605\g)$.
By choosing $r_0 = 4.4022287\ldots \approx 4.40223$ as base circle radius we get the radius of curvature at this point to be zero (see Fig.\ \ref{Abb:Curvature_radius03}), so that $z_K(\ph_1) = -17.6385 - 13.6288\,\ii$ is a cusp (see Fig.\ \ref{Abb:Kurve_als_Huellkurve}).
The contour curve has a loop with a self-intersection point and two cusps.
From the singularity condition $z_K'(\ph) = 0$ (or simply from $\rh(\ph) = 0$) with the help of \texttt{FindRoot}\index{FindRoot@\texttt{FindRoot}} we get the parameter values $\ph_2 \approx 4.04221$ ($231.601\g$) and $\ph_3 \approx 4.59237$ ($263.123\g$) for these cusps.\footnote{Note that we also allow negative radii of curvature here by using formula \eqref{Eq:Curvature_radius} without the absolute value bars, hence $\rh(\ph) = r_0 + r(\ph) + r''(\ph)$.}
One finds that $\ph_a \approx 3.78000$ ($216.578\g$) and $\ph_b \approx 4.82395$ ($276.392\g$) are the parameter values of the self-intersection point $z_K(\ph_a) = z_K(\ph_b) \approx -18.0484 + 13.2662\,\ii$ (see Fig.\ \ref{Abb:Kurve_als_Huellkurve}).  
Since the loop practically does not exist (see the contour in Fig.\ \ref{Abb:Kurve_als_Huellkurve}) due to manufacturing and functionality, the transfer function $r(\ph)$ cannot be completely realized.
Thus in the interval $\ph_a \le \ph \le \ph_b$ the transfer function $r(\ph)$ is replaced by the reduced transfer function $\Rez[z_K(\ph_a)\,\ee^{\ii\ph}]$  (see Fig.\ \ref{Abb:Nockengetriebe_Hubverlust}).\\  
\phantom{} \hfill\bs
\end{example}

%% file: DiffGeo5_4f.tex

\subsection{Cam mechanisms with swinging flat-face follower}
\subsubsection{Basics}

\begin{SCfigure}[][ht]
\includegraphics[width=0.56\textwidth]{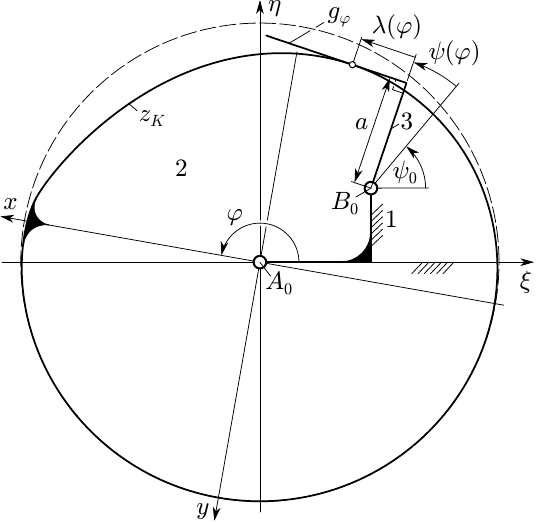}
\caption{Cam mechanism with swinging flat-face follower}
\label{Abb:Swinging_flat-face_follower01}
\end{SCfigure}

In this subsection we consider cam mechanisms consisting of a fixed frame 1 to which we assign a $\xi,\eta$-coordinate system with origin $A_0$, a cam 2 with contour curve $z_K$, and a swinging follower 3 with a flat face which is represented by the line $g_\ph$ (see Fig.\ \ref{Abb:Swinging_flat-face_follower01}).
The cam rotates around $A_0$ with angle $\ph$ in positive direction.
The rotating motion of the cam is transformed into a swinging motion of the follower around point $B_0$.

The $2\pi$-periodic function $\psi(\ph)$ is the {\em transfer function}\index{tranfer function} of this cam mechanism type.
$\psi(\ph)$ is uniquely determined if we choose the constant angle $\psi_0$ so that $\psi(\ph) \ge 0$ (or alternatively $\psi(\ph) \le 0$) holds and there is at least one value of $\ph$ with $\psi(\ph) = 0$. 

\subsubsection{The contour curve}

In this paragraph we determine the contour curve $z_K$ for given transfer function $\psi(\ph)$, and given parameters $\zeta_{B_0}$, $\psi_0$ and $a$ (see Fig.\ \ref{Abb:Swinging_flat-face_follower01}).

In the plane of the contour curve $z_K$ to be generated, we introduce an $x,y$-system whose origin is also located at $A_0$, and the rotation angle $\ph$ of the cam is the angle between the positive $\xi$-axis and the positive $x$-axis. 

A parametric equation of the line $g_\ph$ in the $\xi,\eta$-system is given by
\begin{align*}
  \zeta(\ph,\la)
= {} & \zeta_{B_0} + a\ee^{\ii(\psi_0+\psi(\ph))} + \la\ii\ee^{\ii(\psi_0+\psi(\ph))}\\
= {} & \zeta_{B_0} + \left(a + \la\ii\right)\ee^{\ii(\psi_0+\psi(\ph))}\,,\quad \la\in\R\,.  
\end{align*}
The corresponding parametric equation of $g_\ph$ in the $x,y$-system is
\beqn \label{Eq:z(phi,lambda)_SFFF} 
  z(\ph,\la)
= \zeta(\ph,\la)\,\ee^{-\ii\ph}
= \zeta_{B_0}\ee^{-\ii\ph} + \left(a + \la\ii\right)\ee^{\ii(\psi_0+\psi(\ph)-\ph)}\,,\quad \la\in\R\,.  
\eeqn
The contour curve $z_K$ of the cam is the envelope of the two-parameter family $z(\ph,\la)$ of lines $g_\ph$.

\begin{cor}
The parametric equation of the contour curve $z_K$ is given by
\beqn \label{Eq:zK(phi)_SFFF}
  z_K(\ph)
= z(\ph,\la(\ph))
= \zeta_{B_0}\ee^{-\ii\ph} + \left(a + \ii\la(\ph)\right)\ee^{\ii(\psi_0+\psi(\ph)-\ph)}
\eeqn
with
\beqn \label{Eq:lambda(phi)_SFFF}
  \la(\ph)
= \frac{\left[\zeta_{B_0},\ee^{\ii(\psi_0+\psi(\ph))}\right]}{1-\psi'(\ph)}\,.  
\eeqn
\end{cor}

\begin{proof}
We apply Theorem \ref{Thm:Envelope}.
For the partial derivatives of $z(\ph,\la)$ in \eqref{Eq:z(phi,lambda)_SFFF} we have
\begin{align*}
  \frac{\p z}{\p \ph}
= {} & {-}\zeta_{B_0}\ii\ee^{-\ii\ph} + \left(a+\la\ii\right)\ii\left(\psi'(\ph)-1\right)\ee^{\ii(\psi_0+\psi(\ph)-\ph)}\,,\db\\
  \frac{\p z}{\p \la}
= {} & \ii\ee^{\ii(\psi_0+\psi(\ph)-\ph)}\,.  
\end{align*}
Using rules from \eqref{Eq:Rechenregeln1} and \eqref{Eq:Rechenregeln2}, we get
\begin{align*}
  0
= {} & {\left[\frac{\p z}{\p \ph},\frac{\p z}{\p \la}\right]}\db\\
= {} & {-}\left[\zeta_{B_0},\ee^{\ii(\psi_0+\psi(\ph))}\right] + \left[\left(a+\la\ii\right)\left(\psi'(\ph)-1\right),1\right]\db\\  
= {} & {-}\left[\zeta_{B_0},\ee^{\ii(\psi_0+\psi(\ph))}\right] + \la\left(\psi'(\ph)-1\right)\left[\ii,1\right]\db\\
= {} & {-}\left[\zeta_{B_0},\ee^{\ii(\psi_0+\psi(\ph))}\right] + \la\left(1-\psi'(\ph)\right),
\end{align*}
hence
\beq
  \la
= \la(\ph)
= \frac{\left[\zeta_{B_0},\ee^{\ii(\psi_0+\psi(\ph))}\right]}{1-\psi'(\ph)}\,. \qedhere
\eeq 
\end{proof}

From \eqref{Eq:zK(phi)_SFFF} we get the tangent vector
\beqn \label{Eq:zK'(phi)_SFFF}
  z_K'(\ph)
= -\zeta_{B_0}\ii\ee^{-\ii\ph} + \left\{\left(1-\psi'(\ph)\right)\la(\ph) - \ii\left(a\left(1-\psi'(\ph)\right)-\la'(\ph)\right)\right\}
  \ee^{\ii(\psi_0+\psi(\ph)-\ph)} 
\eeqn
with $\la(\ph)$ according to \eqref{Eq:lambda(phi)_SFFF}, and
\begin{align} \label{Eq:lambda'(phi)_SFFF}
  \la'(\ph)
= {} & \frac{\left[\zeta_{B_0},\ii\psi'(\ph)\ee^{\ii(\psi_0+\psi(\ph))}\right]}{1-\psi'(\ph)}
  + \frac{\left[\zeta_{B_0},\ee^{\ii(\psi_0+\psi(\ph))}\right]}{(1-\psi'(\ph))^2}\,\psi''(\ph)\nonumber\db\\
= {} & \frac{\psi'(\ph)}{1-\psi'(\ph)} \left\langle\zeta_{B_0},\ee^{\ii(\psi_0+\psi(\ph))}\right\rangle
  + \frac{\psi''(\ph)}{(1-\psi'(\ph))^2} \left[\zeta_{B_0},\ee^{\ii(\psi_0+\psi(\ph))}\right].
\end{align}
From \eqref{Eq:zK'(phi)_SFFF}, a straightforward calculation yields
\begin{align*}
  z_K''
= {} & {-}\zeta_{B_0}\ee^{-\ii\ph} + \Big\{2\left(1-\psi'\right)\la' - \psi''\la - a(1-\psi')^2\\
& \hspace{2.7cm}\left. +\;\ii\left(\la'' - \left(1-\psi'\right)^2\la + a\psi''\right)\right\} \ee^{\ii(\psi_0+\psi-\ph)}\,.  
\end{align*}
It remains to calculate $\la''(\ph)$.
From \eqref{Eq:lambda'(phi)_SFFF} we get
\begin{align} \label{Eq:lambda''(phi)_SFFF}
  \la''
= {} & \left(\frac{\psi''}{1-\psi'} + \frac{\psi'\psi''}{(1-\psi')^2}\right) \left\langle\zeta_{B_0},\ee^{\ii(\psi_0+\psi)}\right\rangle
  + \frac{\psi'}{1-\psi'} \left\langle\zeta_{B_0},\ii\psi'\ee^{\ii(\psi_0+\psi)}\right\rangle\nonumber\\
& + \left(\frac{\psi'''}{(1-\psi')^2} + \frac{2\psi''^2}{(1-\psi')^3}\right) \left[\zeta_{B_0},\ee^{\ii(\psi_0+\psi)}\right] 
  + \frac{\psi''}{(1-\psi')^2} \left[\zeta_{B_0},\ii\psi'\ee^{\ii(\psi_0+\psi)}\right]\nonumber\db\\[0.1cm] 
= {} & \left(\frac{\psi''}{1-\psi'} + \frac{\psi'\psi''}{(1-\psi')^2}\right) \left\langle\zeta_{B_0},\ee^{\ii(\psi_0+\psi)}\right\rangle
  - \frac{\psi'^2}{1-\psi'} \left[\zeta_{B_0},\ee^{\ii(\psi_0+\psi)}\right]\nonumber\\
& + \left(\frac{\psi'''}{(1-\psi')^2} + \frac{2\psi''^2}{(1-\psi')^3}\right) \left[\zeta_{B_0},\ee^{\ii(\psi_0+\psi)}\right] 
  + \frac{\psi'\psi''}{(1-\psi')^2} \left\langle\zeta_{B_0},\ee^{\ii(\psi_0+\psi)}\right\rangle\nonumber\db\\[0.1cm]
= {} & \frac{\left(1+\psi'\right)\psi''}{(1-\psi')^2} \left\langle\zeta_{B_0},\ee^{\ii(\psi_0+\psi)}\right\rangle
  + \left(\frac{\psi'''}{(1-\psi')^2} + \frac{2\psi''^2}{(1-\psi')^3} - \frac{\psi'^2}{1-\psi'}\right) \left[\zeta_{B_0},\ee^{\ii(\psi_0+\psi)}\right].  
\end{align}

\subsubsection{Example} 

\begin{figure}[ht]
\begin{minipage}{0.48\textwidth}
\includegraphics[width=\textwidth]{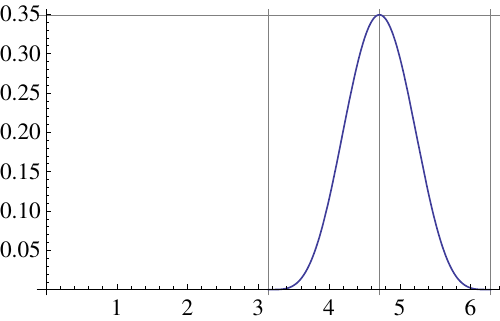}
\caption{Graph of the function $\psi$ defined by \eqref{Eq:psi(phi)_SFFF}} 
\label{Abb:psi(phi)_SFFF}
\end{minipage}
\hfill
\begin{minipage}{0.48\textwidth}
\includegraphics[width=\textwidth]{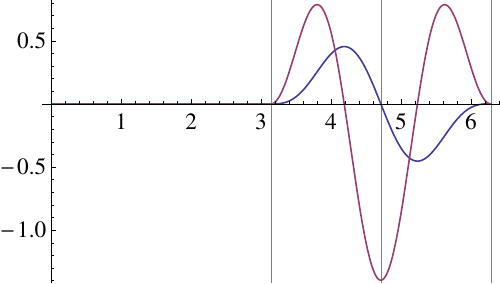}
\caption{$\psi'$ (blue) and $\psi''$ (purple) according to \eqref{Eq:psi'(phi)_SFFF} and \eqref{Eq:psi''(phi)_SFFF}}
\label{Abb:psi'(phi)_and_psi''(phi)_SFFF}
\end{minipage}
\end{figure}

Let the transfer function of the cam mechanism be given by
\beqn \label{Eq:psi(phi)_SFFF}
  \psi(\ph)
= \left\{\begin{array}{c@{\quad\;\mbox{if}\quad\;}l}
	0 & 0 \le \ph < \pi\,,\\[0.15cm]
	\dfrac{\pi}{9}\,f\!\left(\dfrac{\ph-\pi}{\pi/2}\right) & \pi \le \ph \le 2\pi
  \end{array}\right.  
\eeqn
with
\beqn \label{Eq:f(x)_SFFF}
  f(x) = \sin^4\frac{\pi x}{2}\,.  
\eeqn
The graph of $\psi(\ph)$ is shown in Fig.\ \ref{Abb:psi(phi)_SFFF}.
Using
\begin{align*}
  \sin^4 x
= {} & \left(1-\cos^2\right)^2
= 1 - 2\cos^2 x + \left(\cos^2 x\right)^2\db\\
= {} & 1 - 2\cdot\frac{1}{2}\left(1+\cos(2x)\right) + \frac{1}{4}\left(1+\cos(2x)\right)^2\db\\
= {} & {-}\cos(2x) + \frac{1}{4}\left(1+2\cos(2x)+\cos^2(2x)\right)\db\\
= {} & {-}\cos(2x) + \frac{1}{4} + \frac{1}{2}\cos(2x) + \frac{1}{4}\cdot\frac{1}{2}\left(1+\cos(4x)\right)\\
= {} & \frac{3}{8} - \frac{1}{2}\cos(2x) + \frac{1}{8}\cos(4x)\,,
\end{align*}
we can write \eqref{Eq:f(x)_SFFF} as
\beq
  f(x)
= \frac{3}{8} - \frac{1}{2}\cos(\pi x) + \frac{1}{8}\cos(2\pi x)\,.  
\eeq
For the first and second derivative we get 
\begin{align*}
  f'(x)
= {} & \frac{\pi}{2}\sin(\pi x) - \frac{\pi}{4}\sin(2\pi x)\,\\ 
  f''(x)
= {} & \frac{\pi^2}{2}\cos(\pi x) - \frac{\pi^2}{2}\cos(2\pi x)\,,
\end{align*}
hence
\beqn \label{Eq:psi'(phi)_SFFF}
  \psi'(\ph)
= \left\{\begin{array}{c@{\quad\;\mbox{if}\quad\;}l}
	0 & 0 \le \ph < \pi\,,\\[0.15cm]
	\dfrac{\pi/9}{\pi/2}\,f'\!\left(\dfrac{\ph-\pi}{\pi/2}\right) & \pi \le \ph \le 2\pi
  \end{array}\right.  
\eeqn
and
\beqn \label{Eq:psi''(phi)_SFFF}
  \psi''(\ph)
= \left\{\begin{array}{c@{\quad\;\mbox{if}\quad\;}l}
	0 & 0 \le \ph < \pi\,,\\[0.15cm]
	\dfrac{\pi/9}{(\pi/2)^2}\,f''\!\left(\dfrac{\ph-\pi}{\pi/2}\right) & \pi \le \ph \le 2\pi
  \end{array}\right.  
\eeqn
(see Fig.\ \ref{Abb:psi'(phi)_and_psi''(phi)_SFFF}).\footnote{The function \eqref{Eq:psi(phi)_SFFF} is also used in \cite[pp.\ 309-311]{Pennestri&Cera} as transfer function for an example of a cam mechanism with a swinging flat-face follower.}

\begin{SCfigure}[1.1][ht]
\includegraphics[width=0.45\textwidth]{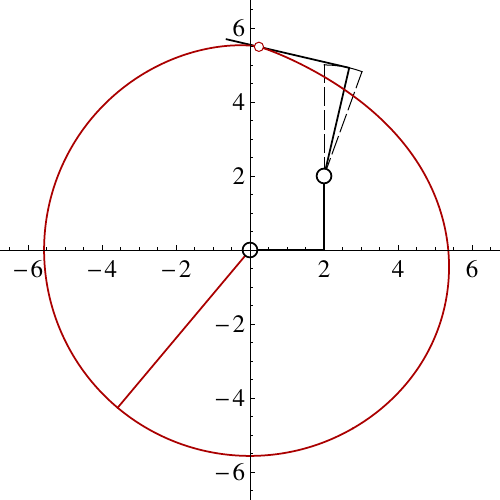}
\caption{Cam mechanism with swinging flat-face follower\\[0.2cm] Position: $\ph = 230\g$, $\psi(\ph) = 6.88725\g$\\ (in radians: $\ph \approx 4.01426$, $\psi(\ph) = 0.120205)$}
\label{Abb:Swinging_flat-face_follower02}
\end{SCfigure}

Choosing
\beq
  \zeta_{B_0} = 2 + 2\ii\,,\qquad \psi_0 = 70\g\cdot\pi/180\g\,,\qquad a = 3
\eeq
gives the mechanism in Fig.\ \ref{Abb:Swinging_flat-face_follower02}.

\begin{figure}[ht!]
\begin{minipage}{0.48\textwidth}
\includegraphics[width=\textwidth]{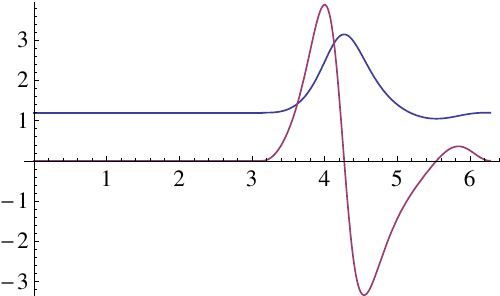} 
\caption{$\la(\ph)$ (blue) and $\la'(\ph)$ (purple)} 
\label{Abb:lambda_and_lambda1_SFFF}
\end{minipage}
\hfill
\begin{minipage}{0.48\textwidth}
\includegraphics[width=\textwidth]{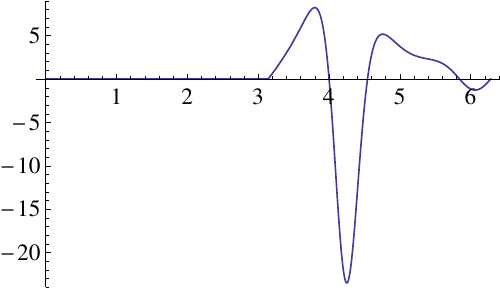} 
\caption{$\la''(\ph)$}
\label{Abb:lambda2_SFFF}
\end{minipage}
\end{figure}

The Figs.\ \ref{Abb:lambda_and_lambda1_SFFF} and \ref{Abb:lambda2_SFFF} show the graphs of the functions $\la(\ph)$, $\la'(\ph)$ and $\la''(\ph)$ according to \eqref{Eq:lambda(phi)_SFFF}, \eqref{Eq:lambda'(phi)_SFFF} and \eqref{Eq:lambda''(phi)_SFFF}, respectively.  

\begin{figure}[ht!]
\begin{minipage}{0.48\textwidth}
\includegraphics[width=\textwidth]{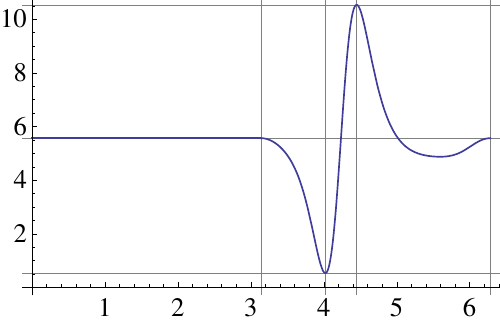}
\caption{Curvature radius $|\kappa(\ph)|^{-1}$} 
\label{Abb:Curvature_radius_SFFF01}
\end{minipage}
\hfill
\begin{minipage}{0.48\textwidth}
\includegraphics[width=\textwidth]{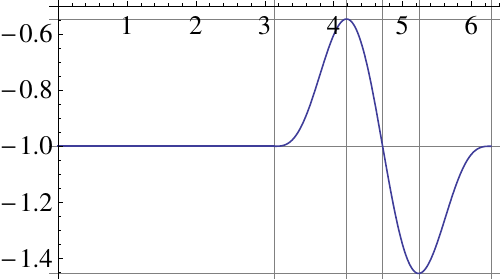}
\caption{$\kappa(\ph) \cdot |z_K'(\ph)|$}
\label{Abb:kappa_times_abs(zK1)_SFFF01}
\end{minipage}
\end{figure}

Fig.\ \ref{Abb:Curvature_radius_SFFF01} shows the graph of the curvature radius function while Fig.\ \ref{Abb:kappa_times_abs(zK1)_SFFF01} shows the graph of the function $\kappa(\ph) \cdot |z_K'(\ph)|$.
One can see that the graph in Fig.\ \ref{Abb:kappa_times_abs(zK1)_SFFF01} looks similar to that of $\psi'(\ph)$ in Fig.\ \ref{Abb:psi'(phi)_and_psi''(phi)_SFFF}; both graphs are point-symmetrical in the interval $[\pi,2\pi]$.
Clearly, we have
\beqn \label{Eq:Total_curvature_special_case}
  \int_0^{2\pi} \kappa(\ph) \left|z_K'(\ph)\right| \dd\ph
= -2\pi\,.  
\eeqn
This is special case of the total curvature\index{curvature!total curvature@total $\sim $} formula
\beq
  \oint \kappa\, \dd s
= 2\pi   
\eeq
for a simple closed curve\index{curve!simple closed curve@simple closed $\sim $} (closed curve without self-intersections, Jordan curve\index{curve!Jordan curve@Jordan $\sim $}) with arc length element $\dd s$ (see \cite[p.\ 16]{Voss}).\footnote{Of course, the negative sign in \eqref{Eq:Total_curvature_special_case} results from the negative orientation of our contour curve.}

%% file: DiffGeo5_4g.tex

\subsection{Cam mechanisms with pivoted roller follower} \label{Sec:Pivoted_follower}
\subsubsection{Basics}

\begin{SCfigure}[][ht]
  \includegraphics[width=0.48\textwidth]{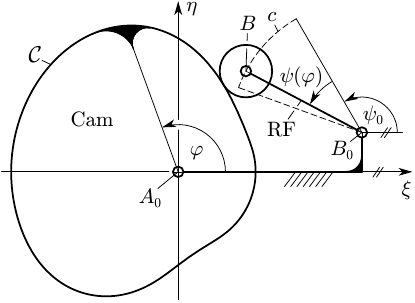}
  \caption{Cam mechanism with roller follower}
  \label{Abb:Roller_follower} 
\end{SCfigure}

In this subsection we consider cam mechanisms with roller follower (see Fig.\ \ref{Abb:Roller_follower}).
The cam rotates with angle $\ph$ around a fixed point $A_0$.
The roller follower RF makes an oscillating motion around a fixed (pivot) point $B_0$, where the roller center point $B$ describes a circular arc $c$ around $B_0$.
The angle of the oscillating motion is the angle $\psi_0 + \psi(\ph)$, $\psi_0 = \tn{const}$, between the positive $\xi$-axis and the oriented line segment $B_0B$.
The {\em transfer function} of this mechanism is the function $\psi(\ph)$, $0 \le \ph \le 2\pi$.   

In the following, we will use a several times the transfer function
\beqn \label{Eq:psi(phi)}
  \psi(\ph)
= \left\{\begin{array}{c@{\quad\;\mbox{if}\quad\;}l}
	\dfrac{40\g\cdot\pi}{180\g}\, f_{3,2}\!\left(\dfrac{\ph}{5\pi/6}\right) & 0 \le \ph \le \dfrac{5\pi}{6}\,,\\[0.3cm]
	\dfrac{40\g\cdot\pi}{180\g} & \dfrac{5\pi}{6} < \ph < \dfrac{7\pi}{6}\,,\\[0.3cm]
	\dfrac{40\g\cdot\pi}{180\g} \left[1-g_3\!\left(\dfrac{\ph-7\pi/6}{2\pi/3}\right)\right] & \dfrac{7\pi}{6} \le \ph \le \dfrac{11\pi}{6}\,,\\[0.3cm]
	0 & \dfrac{11\pi}{6} < \ph \le 2\pi\,, 
  \end{array}\right.  
\eeqn
which is similar to \eqref{Eq:r(phi)} (see there for $f_{3,2}$ and $g_3$). 

\subsubsection{The contour curve}

\begin{figure}[ht]
\begin{minipage}{0.43\textwidth}
  \includegraphics[width=\textwidth]{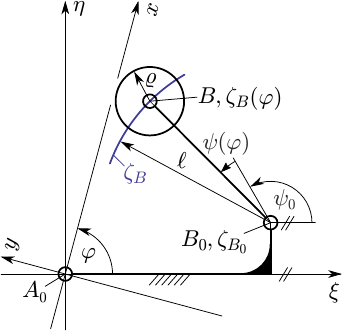}
  \caption{\small Parameters of the cam mechanism, and circular arc $\zeta_B$}
  \label{Abb:Kurvengetriebe_zetaB}
\end{minipage}
\hfill
\begin{minipage}{0.53\textwidth}
  \centering
  \includegraphics[width=\textwidth]{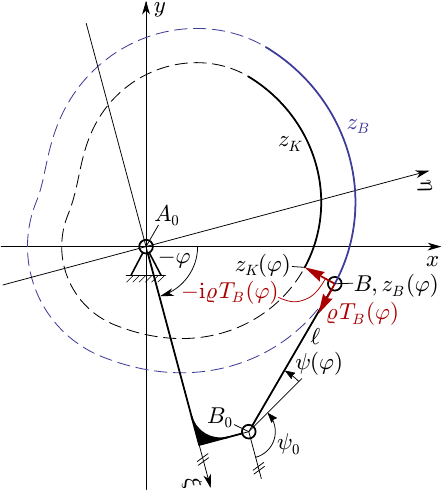}
  \caption{\small Generation of center point curve $z_B$ and contour curve $z_K$}
  \label{Abb:Kurvengetriebe_zB_und_zK}
\end{minipage}
\end{figure}

We introduce a frame-fixed $\xi,\eta$-coordinate system such that the origin lies at the pivot point $A_0$ (see Fig. \ref{Abb:Kurvengetriebe_zetaB}).
The roller center point $B$ describes a circular arc $\zeta_B$ with parametric equation
\beqn \label{Eq:zetaB(phi)}
  \zeta_B(\ph)
= \zeta_{B_0} + \ell\,\ee^{\ii[\psi_0+\psi(\ph)]}\,,\quad
  0 \le \ph < 2\pi\,,\footnote{Note that the complex-valued function\index{function!complex-valued function@complex-valued $\sim $} $\zeta_B(\ph)$, for example, can also be regarded as a transfer function of the cam mechanism; $\zeta_B(\ph)$ describes a position (the position of the point $B$) as a function of the input angle $\ph$. In this sense the $k$-th derivative $\zeta_B^{(k)}(\ph)$ is a $k$-th order transfer function.}
\eeqn
where $\psi_0$ is a constant which can be freely chosen within certain limits, and $\psi(\ph)$ is the transfer function.

Again, we let an $x,y$-coordinate system rotate around the point $A_0$ in positive direction ($=$ assumed rotation direction of the cam) with angle $\ph$ between the positive $\xi$-axis and the positive $x$-axis. 
With this, the path $z_B$ of the roller center $B$ (roller center curve\index{roller center curve}) can now be determined in the $x,y$-coordinate system.
For this we use the principle of {\em kinematic inversion}\index{inversion, kinematic}, and assume the $x,y$-system to be the fixed reference system.
The $\xi,\eta$-system (frame) then rotates with rotation angle $-\ph$ around $A_0$ (see Fig. \ref{Abb:Kurvengetriebe_zB_und_zK}).
With \eqref{Eq:zetaB(phi)} we get the parametric equation
\beq
  z_B(\ph)
= \zeta_B(\ph)\,\ee^{-\ii\ph} 
= \left(\zeta_{B_0} + \ell\,\ee^{\ii[\psi_0+\psi(\ph)]}\right)\ee^{-\ii\ph}\\
= \zeta_{B_0}\,\ee^{-\ii\ph} + \ell\,\ee^{\ii[\psi_0+\psi(\ph)-\ph]} 
\eeq
with which $z_B$ is a negatively oriented curve\index{curve!negatively oriented curve@negatively oriented $\sim $}.
The tangent unit vector is given by
\beq
  T_B(\ph)
= \frac{z_B'(\ph)}{|z_B'(\ph)|}
\eeq
with
\beqn \label{Eq:z_B'(phi)}
  z_B'(\ph)
= -\ii\,\zeta_{B_0}\,\ee^{-\ii\ph} + \ii\,\ell\,(\psi'(\ph)-1)\,\ee^{\ii[\psi_0+\psi(\ph)-\ph]}
\eeqn
and
\beqn \label{Eq:|z_B'(phi)|}
\begin{aligned}
& \hspace{-0.55cm} \left|z_B'(\ph)\right|\\[0.05cm]
= {} & |\ii| \left|\ell\left(\psi'(\ph)-1\right)\ee^{\ii[\psi_0+\psi(\ph)]} - \zeta_{B_0}\right| \left|\ee^{-\ii\ph}\right|
= \left|\ell\left(\psi'(\ph)-1\right)\ee^{\ii[\psi_0+\psi(\ph)]} - \zeta_{B_0}\right|\\[0.05cm]
= {} & \sqrt{\left[\ell\,(\psi'(\ph)-1)\cos(\psi_0+\psi(\ph))-\xi_{B_0}\right]^2 + \left[\ell\,(\psi'(\ph)-1)\sin(\psi_0+\psi(\ph))-\eta_{B_0}\right]^2}\;.  
\end{aligned}
\eeqn
With the chosen roller radius $\rh$, the contour curve $z_K$ is obtained as {\em inner parallel curve}\index{curve!inner parallel curve@inner parallel $\sim $} of $z_B$ by means of
\beqn \label{Eq:z_K(phi)}
  z_K(\ph) = z_B(\ph) - \ii\rh T_B(\ph)\,.
\eeqn
If one replaces the minus sign in \eqref{Eq:z_K(phi)} by the plus sign, then one receives as $z_K$ obviously the {\em outer parallel curve}\index{curve!outer parallel curve@outer parallel $\sim $} of $z_B$.

\subsubsection{Curvature}

According to Eq.\ \eqref{Eq:Curvature} in Theorem \ref{Thm:Curvature} we can calculate the oriented curvature\index{curvature!oriented curvature@oriented $\sim $} $\kappa_B(\ph)$ at a point~$\ph$ of the center point curve $z_B$ by means of
\beq
  \kappa_B(\ph)
= \frac{\left[z_B'(\ph),z_B''(\ph)\right]}{|z_B'(\ph)|^3}
\eeq
with $z_B'(\ph)$ from \eqref{Eq:z_B'(phi)} and
\beq
  z_B''(\ph)
= -\zeta_{B_0}\,\ee^{-\ii\ph} - \ell \left[\left(\psi'(\ph)-1\right)^2 - \ii\psi''(\ph)\right] \ee^{\ii[\psi_0+\psi(\ph)-\ph]}\,.\footnote{Changing the orientation of $z_B(\ph)$ changes the sign of the curvature $\kappa_B(\ph)$ due to the reversal of the direction of the tangent vector $z_B'(\ph)$. The ``acceleration'' vector $z_B''(\ph)$ remains the same.}    
\eeq
Therefore, the osculating circle $\mathcal{C}_B(\ph)$ to the center point curve $z_B$ at point $\ph$ has radius 
\beq
  \rh_B(\ph) = \frac{1}{|\kappa_B(\ph)|}
\eeq
and, according to Eq.\ \eqref{Eq:Center_of_curvature} in Theorem \ref{Thm:Curvature}, center point
\beq
  \tilde{z}_B(\ph) = z_B(\ph) + \ii\, \frac{1}{\kappa_B(\ph)}\, \frac{z_B'(\ph)}{|z_B'(\ph)|}\,.
\eeq
Since the contour curve $z_K$ is the parallel curve of $z_B$ with distance $\rh$ ($=$ roller radius), the osculating circle $\mathcal{C}_K(\ph)$ to $z_K$ at point~$\ph$ has radius
\beqn \label{Eq:rho_K(phi)}
  \rh_K(\ph) = \left|\frac{1}{\kappa_B(\ph)} + \rh\right|
\eeqn
and center point $\tilde{z}_K(\ph) = \tilde{z}_B(\ph)$.

In contrast to cam mechanisms with flat-face follower, the contour curve $z_K$ does not have to be completely convex\index{convex}, but can also have concave\index{concave} sections.
With our negative orientation of $z_B$, this curve is convex\index{convex} in a point $\ph$ if $\kappa_B(\ph) < 0$, and concave\index{concave} if $\kappa_B(\ph) > 0$.
The analogous applies to the curve $z_K$.

From \eqref{Eq:rho_K(phi)} one sees that $\rh_K(\ph) = \rh_B(\ph) - \rh$ if $\kappa_B(\ph) < 0$, and $\rh_K(\ph) = \rh_B(\ph) + \rh$ if $\kappa_B(\ph) > 0$.

\begin{SCfigure}[][ht]
\includegraphics[width=0.53\textwidth]{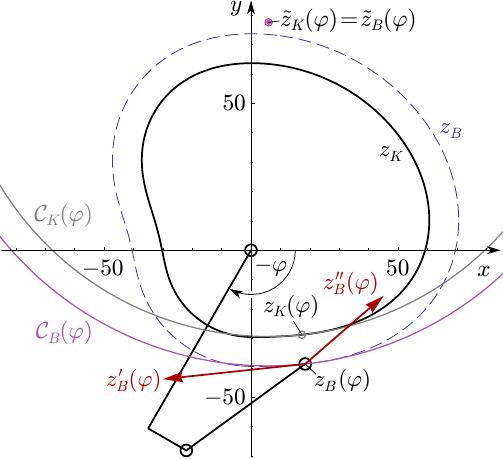}
\caption{Center point curve $z_B$ with a curvature circle $\mathcal{C}_B(\ph)$, and contour curve~$z_K$ with a curvature circle $\mathcal{C}_K(\ph)$\\[0.2cm]
Transfer function $\psi$ according to \eqref{Eq:psi(phi)}, $\xi_{B_0} = 70$, $\eta_{B_0} = 15$, $\ell = 50$, $\rh = 10$, $\psi_0 = 120\g$\\[0.2cm]
Position $\ph = 120\g$}
\label{Abb:Kurveng_mit_Kruemmungskreisen}
\end{SCfigure} 

\begin{SCfigure}[][ht]
\includegraphics[width=0.53\textwidth]{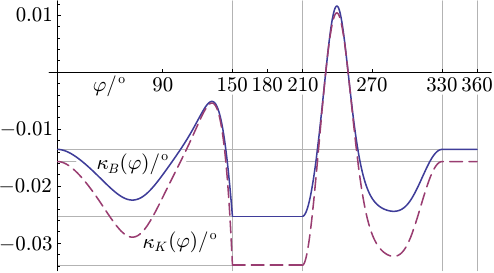}
\caption{Curvature functions $\kappa_B(\ph)$ and $\kappa_K(\ph)$ for the curves $z_B$ and $z_K$, respectively, in Fig.\ \ref{Abb:Kurveng_mit_Kruemmungskreisen}}
\label{Abb:Curvature_function}
\end{SCfigure} 

\begin{example} \label{Bsp:Cam_with_curvature_circles}
An example for a cam with curvature circles $\mathcal{C}_B(\ph)$ and $\mathcal{C}_K(\ph)$ is shown in Fig.~\ref{Abb:Kurveng_mit_Kruemmungskreisen}.
At point $z_B(\ph)$, the directional angle between the vectors $z_B'(\ph)$ and $z_B''(\ph)$ is greater than $180\g$, from \eqref{Eq:QVP_with_sine} it follows that the outer product $[z_B'(\ph),z_B''(\ph)]$ is negative, therefore the curvature $\kappa_B(\ph)$ is negative and the tangent vector $z_B'(\ph)$ rotates in negative direction.
Fig.\ \ref{Abb:Curvature_function} shows the curvature functions of the curves of $z_B$ and $z_K$. \hfill\bs
\end{example}

\begin{figure}[ht]
\begin{minipage}{0.22\textwidth}
  \includegraphics[width=0.8\textwidth]{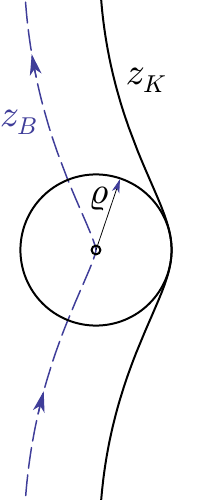}
  \caption{\small Curve $z_B$ with cusp}
  \label{Abb:zB_with_cusp}
\end{minipage}
\hfill
\begin{minipage}{0.22\textwidth}
  \centering
  \includegraphics[width=0.96\textwidth]{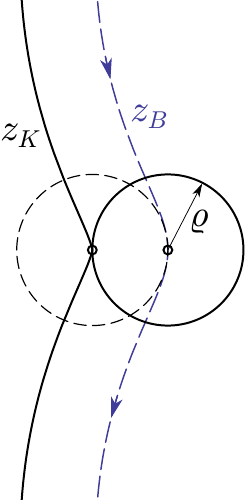}
  \caption{\small Curve $z_K$ with cusp}
  \label{Abb:zK_with_cusp}
\end{minipage}
\hfill
\begin{minipage}{0.22\textwidth}
  \centering
  \includegraphics[width=0.9\textwidth]{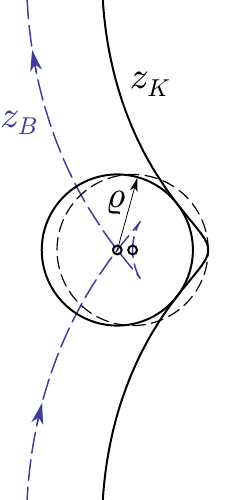}
  \caption{\small Curve $z_B$ with fishtail}
  \label{Abb:zB_with_fishtail}
\end{minipage}
\hfill
\begin{minipage}{0.22\textwidth}
  \centering
  \includegraphics[width=1.01\textwidth]{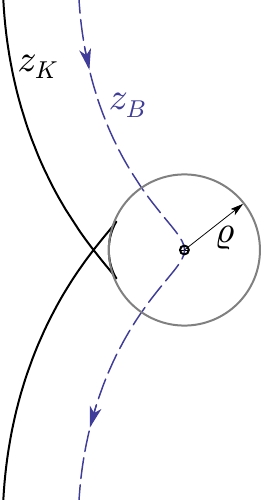}
  \caption{\small Curve $z_K$ with fishtail}
  \label{Abb:zK_with_fishtail}
\end{minipage}
\end{figure}

We consider some cases where cusps\index{cusp} or loops\index{loop} of center point curve $z_B$ or contour curve $z_K$ occur:\mynobreakpar
\begin{itemize}[leftmargin=0.5cm]
\item Fig.\ \ref{Abb:zB_with_cusp}: The center point curve $z_B$ has a cusp\index{cusp}, but no loop\index{loop}.
At the point of the cusp, the curvature radius $\rh_B$ is equal to zero.
Consequently, according to \eqref{Eq:rho_K(phi)}, the curvature radius~$\rh_K$ is equal to the roller radius $\rh$.
The curvature $\kappa_K$ at this point is positive.
\item Fig.\ \ref{Abb:zK_with_cusp}: In the position shown is $\kappa_B = -\rho^{-1}$, otherwise $\kappa_B > -\rho^{-1}$.
Therefore, in the shown position holds $\rh_K = |\kappa_B^{-1} + \rh| = |{-}\rh + \rh| = 0$ and $z_K$ has a cusp, but no loop. 
\item Fig.\ \ref{Abb:zB_with_fishtail}: The center point curve $z_B$ has a loop (with two cusps) in the form of a fishtail.
One sees that the smallest curvature radius of $z_K$ is smaller than the roller radius $\rh$.
Consequently, the roller cannot roll on the entire curve $z_K$.
The theoretical roller position is drawn dashed and the practical position with a solid line.
\item Fig.\ \ref{Abb:zK_with_fishtail}: The smallest curvature radius of $z_B$ is smaller than the roller radius $\rh$.
This results in $z_K$ having a loop in the form of a fishtail, which of course is of no practical use. 
\end{itemize}

The fact that $z_B$ and $z_K$ are mutually parallel curves results in the duality visible in Figures \ref{Abb:zB_with_cusp} and \ref{Abb:zK_with_cusp} as well as in Figures \ref{Abb:zB_with_fishtail} and \ref{Abb:zK_with_fishtail}.

For practical requirements it is necessary in any case that for convex segments of $z_K$ the condition $\rh_B > \rh$ holds, and for concave segments of $z_K$ the condition $\rh_K > \rh$.

The curvature functions of the cams in Fig.\ \ref{Abb:P_and_F-cam} are shown in Fig.\ \ref{Abb:kappa_P_und_F}.
 
\subsubsection{Length and area}

\begin{thm} \label{Thm:L(z_K)}
The length $L(z_K)$ of the inner parallel curve $z_K$ at distance $\rh$ from $z_B$ is given by
\beq
  L(z_K)
= \int_0^{2\pi} \left|1 + \rh\,\kappa_B(\ph)\right| \left|z_B'(\ph)\right| \dd\ph\,.  
\eeq
\end{thm}

\begin{proof}
According to \eqref{Eq:Arc_length}, the length of the contour curve $z_K$ is given by
\beq
  L(z_K)
= \int_0^{2\pi} \left|z_K'(\ph)\right| \dd\ph\,.
\eeq
With \eqref{Eq:z_K(phi)} we get
\beq
  L(z_K)
= \int_0^{2\pi} \left|z_B'(\ph) - \ii\,\rh\,T_B'(\ph)\right| \dd\ph\,,   
\eeq
and with the result of Lemma \ref{Lem:T'(t)}
\begin{align*}
  L(z_K)
= {} & \int_0^{2\pi} \left|z_B'(\ph) + \rh\,\kappa_B(\ph) \left|z_B'(\ph)\right| T_B(\ph)\right| \dd\ph\db\\[0.05cm]
= {} & \int_0^{2\pi} \left|z_B'(\ph) + \rh\,\kappa_B(\ph) \left|z_B'(\ph)\right| \frac{z_B'(\ph)}{|z_B'(\ph)|}\right| \dd\ph\db\\[0.05cm]   
= {} & \int_0^{2\pi} \left|z_B'(\ph) + \rh\,\kappa_B(\ph)\,z_B'(\ph)\right| \dd\ph
= \int_0^{2\pi} \left|\left(1 + \rh\,\kappa_B(\ph)\right) z_B'(\ph)\right| \dd\ph\\[0.05cm]
= {} & \int_0^{2\pi} \left|1 + \rh\,\kappa_B(\ph)\right| \left|z_B'(\ph)\right| \dd\ph\,. \qedhere
\end{align*}
\end{proof}

We consider the integral from Theorem \ref{Thm:L(z_K)}
\beq
  L(z_K)
= \int_0^{2\pi} \left|1 + \rh\,\kappa_B(\ph)\right| \left|z_B'(\ph)\right| \dd\ph\,.  
\eeq
If $1 + \rh\,\kappa_B(\ph) \ge 0$ for $0 \le \ph < 2\pi$,
which can obviously be written equivalently as
\beqn \label{Eq:kappa_B(phi)>=-1/rho}
  \kappa_B(\ph) \ge -\frac{1}{\rh} \quad\mbox{for}\quad 0 \le \ph < 2\pi\,,
\eeqn
then we have
\begin{align*}
  L(z_K)
= {} & \int_0^{2\pi} \left(1 + \rh\,\kappa_B(\ph)\right) \left|z_B'(\ph)\right| \dd\ph
= \int_0^{2\pi} \left|z_B'(\ph)\right| \dd\ph + \rh \int_0^{2\pi} \kappa_B(\ph) \left|z_B'(\ph)\right| \dd\ph\db\\[0.05cm]   
= {} & L(z_B) + \rh \int_0^{2\pi} \kappa_B(\ph) \left|z_B'(\ph)\right| \dd\ph\,,
\end{align*}
where $L(z_B)$ is the length of $z_B$.
It remains to evaluate the last integral.
With the definition of the curvature we have
\beq
  \kappa_B
= \frac{\dd \tau_B}{\dd s_B}
= \frac{\dd \tau_B}{\dd \ph}\, \frac{\dd \ph}{\dd s_B}\,,
  \quad\mbox{hence}\quad
  \kappa_B(\ph)
= \frac{\tau_B'(\ph)}{s_B'(\ph)}\,.
\eeq
Using
\beq
  s_B'(\ph)
= \sqrt{\left(x_B'(\ph)\right)^2 + \left(y_B'(\ph)\right)^2}
= \sqrt{z_B'(\ph)\,\overline{z_B'(\ph)}}
= \left|z_B'(\ph)\right|
\eeq      
we get
\beq
  \tau_B'(\ph)
= \kappa_B(\ph) \left|z_B'(\ph)\right|\,. 
\eeq
and therefore
\beqn \label{Eq:Total_curvature}
  \int_0^{2\pi} \kappa_B(\ph) \left|z_B'(\ph)\right| \dd\ph
= \int_0^{2\pi} \tau_B'(\ph)\, \dd\ph
= \tau_B(2\pi) - \tau_B(0)
= -2\pi\,.\footnote{The negative sign follows from the negative orientation of $z_B$. In total, the tangent vector rotates by the angle $-2\pi$.}  
\eeqn
So we finally have
\beqn \label{Eq:L(z_K)_with_L(z_B)}
  L(z_K) = L(z_B) - 2\pi\rh
\eeqn
if the condition \eqref{Eq:kappa_B(phi)>=-1/rho} is fulfilled.

We denote with $A(z_B)$ and $A(z_K)$ the area of the region bounded by the curve $z_B$ and the curve $z_K$, respectively. 

\begin{thm} \label{Thm:A(z_K)}
The signed area $A(z_K)$ (with negative sign) bounded by the inner parallel curve $z_K$ at distance $\rh$ from $z_B$ is given by
\beq
  A(z_K) = A(z_B) + \rh L(z_B) - \pi \rh^2\,,
\eeq
where $A(z_B)$ is the signed area (with negative sign) of the region bounded by $z_B$, 
\beq 
  A(z_B) = \frac{1}{2} \int_0^{2\pi} \left[z_B(\ph),z_B'(\ph)\right] \dd\ph\,,
\eeq
and $L(z_B)$ is the length of $z_B$,
\beq
  L(z_B)
= \int_0^{2\pi} \left|z_B'(\ph)\right| \dd\ph\,.
\eeq
\end{thm}

\begin{proof}
According to Theorem \ref{Thm:Area_with_QVP}, the area $A(z_K)$ is given by
\beq
  A(z_K)
= \frac{1}{2} \int_0^{2\pi} \left[z_K(\ph),z_K'(\ph)\right] \dd\ph\,.
\eeq
Taking into account \eqref{Eq:z_K(phi)}, for the outer product in the integrand we have
\begin{align} \label{Eq:[z_K,z_K']}
  \left[z_K,z_K'\right]
= {} & \left[z_B-\ii\rh T_B,z_B'-\ii\rh T_B'\right]\nonumber\db\\[0.05cm]
= {} & \left[z_B,z_B'\right] - \rh\left[z_B,\ii T_B'\right] - \rh\left[\ii T_B,z_B'\right] + \rh^2\left[\ii T_B,\ii T_B'\right]\nonumber\db\\[0.05cm]
= {} & \left[z_B,z_B'\right] - \rh\left[z_B,\ii T_B'\right] + \rh\left[z_B',\ii T_B\right] + \rh^2\left[T_B,T_B'\right].
\end{align}
For the third outer product in \eqref{Eq:[z_K,z_K']} we get
\beqn \label{Eq:[z_B',iT_B]}
  \left[z_B',\ii T_B\right]
= \left|z_B'\right| \cdot \left|\ii T_B\right| \cdot \sin\measuredangle\left(z_B',\ii T_B\right)
= \left|z_B'\right| \sin(\pi/2)
= \left|z_B'\right|.  
\eeqn
Applying the product rule $[z_1(\ph),z_2(\ph)]' = [z_1'(\ph),z_2(\ph)] + [z_1(\ph),z_2'(\ph)]$ and \eqref{Eq:[z_B',iT_B]}, for the second outer product in \eqref{Eq:[z_K,z_K']} we get
\beqn \label{Eq:[z_B,iT_B']}
  \left[z_B,\ii T_B'\right]
= \left[z_B,\ii T_B\right]' - \left[z_B',\ii T_B\right]  
= \left[z_B,\ii T_B\right]' - \left|z_B'\right|.
\eeqn
Using Lemma \ref{Lem:T'(t)}, for the last outer product in \eqref{Eq:[z_K,z_K']} one finds
\beqn \label{Eq:T_B,T_B'}
  \left[T_B,T_B'\right]
= \left[T_B,\kappa_B\left|z_B'\right| \ii\, T_B\right]
= \kappa_B \left|z_B'\right| \left[T_B,\ii T_B'\right]
= \kappa_B \left|z_B'\right|. 
\eeqn
Substituting \eqref{Eq:[z_B',iT_B]}, \eqref{Eq:[z_B,iT_B']} and \eqref{Eq:T_B,T_B'} into \eqref{Eq:[z_K,z_K']} yields
\begin{align*}
  \left[z_K,z_K'\right]
= {} & \left[z_B,z_B'\right] - \rh\left(\left[z_B,\ii T_B\right]' - \left|z_B'\right|\right) + \rh\left|z_B'\right| 
  + \rh^2\kappa_B \left|z_B'\right|\\[0.05cm]
= {} & \left[z_B,z_B'\right] - \rh\left[z_B,\ii T_B\right]' + 2\rh\left|z_B'\right| + \rh^2\kappa_B \left|z_B'\right|.
\end{align*}
Applying Theorem \ref{Thm:Area_with_QVP}, Eq.\ \eqref{Eq:Arc_length} and Eq.\ \eqref{Eq:Total_curvature}, we finally get
\begin{align*}
  A(z_K)
= {} & \frac{1}{2} \int_0^{2\pi} \left[z_K,z_K'\right] \dd\ph\db\\[0.05cm]
= {} & \frac{1}{2} \int_0^{2\pi} \left(\left[z_B,z_B'\right] - \rh\left[z_B,\ii T_B\right]' + 2\rh\left|z_B'\right| + \rh^2\kappa_B \left|z_B'\right|
  \right) \dd\ph\db\\[0.05cm]
= {} & \frac{1}{2} \int_0^{2\pi} \left[z_B,z_B'\right] \dd\ph - \frac{1}{2}\,\rh \int_0^{2\pi} \left[z_B,\ii T_B\right]'\, \dd\ph
  + \rh \int_0^{2\pi} \left|z_B'\right| \dd\ph + \frac{1}{2}\,\rh^2 \int_0^{2\pi} \kappa_B \left|z_B'\right| \dd\ph\db\\[0.05cm]
= {} & A(z_B) - \frac{1}{2}\,\rh\left(\left[z_B(2\pi),\ii T_B(2\pi)\right]-\left[z_B(0),\ii T_B(0)\right]\right)
  + \rh L(z_B) - \pi\rh^2\db\\[0.05cm]
= {} & A(z_B) + \rh L(z_B) - \pi\rh^2\,. \qedhere  
\end{align*}
\end{proof}

\begin{example}
We calculate the length $L(z_K)$ and the area $A(z_K)$ for the contour curve $z_K$ in Fig.\ \ref{Abb:Kurveng_mit_Kruemmungskreisen}.
From Fig.\ \ref{Abb:Curvature_function} it can be seen that with $\rh=10$ the inequality \eqref{Eq:kappa_B(phi)>=-1/rho} is satisfied.
Therefore we can use \eqref{Eq:L(z_K)_with_L(z_B)}.
For $L(z_B)$ we have
\beq
  L(z_B)
= \int_0^{2\pi} \left|z_B'(\ph)\right| \dd\ph
= \left(\int_0^{5\pi/6}+\int_{5\pi/6}^{7\pi/6}+\int_{7\pi/6}^{11\pi/6}+\int_{11\pi/6}^{2\pi}\right) \left|z_B'(\ph)\right| \dd\ph\,.
\eeq
Numerical integration gives $L(z_B) = 367.5036483978\ldots \approx 367.5036$ and thus
\beq
  L(z_K) \approx 367.5036 - 2\pi \cdot 10 = 304.6718\,.
\eeq
According to Theorem \ref{Thm:A(z_K)} the signed area $A(z_K)$ can be calculated with
\beq
  A(z_K) = A(z_B) + \rh L(z_B) - \pi \rh^2\,.
\eeq
For the signed area $A(z_B)$ we have
\begin{align*}
  A(z_B)
= {} & \frac{1}{2} \int_0^{2\pi} \left[z_B(\ph),z_B'(\ph)\right] \dd\ph\db\\[0.05cm]
= {} & \frac{1}{2} \left(\int_0^{5\pi/6}+\int_{5\pi/6}^{7\pi/6}+\int_{7\pi/6}^{11\pi/6}+\int_{11\pi/6}^{2\pi}\right) \left[z_B(\ph),z_B'(\ph)\right]
	   \dd\ph\,.
\end{align*}
Numerical integration yields $A(z_B) \approx -10589.1488$.
With this we get
\beq
  A(z_K) \approx -10589.1488 + 3675.03648 - \pi \cdot 10^2 = -7228.2716\,.
\eeq
The actual area is of course the absolute value of $A(z_K)$. \hfill\bs  
\end{example}
 
\subsubsection{Transmission angle (function)}

{\bf Transmission angle:}\index{transmission angle} 
The {\em transmission angle} $\mu$ according to Alt (Hermann Alt) is the acute angle between the tangent of the path curve of the roller midpoint $B$ with respect to the frame and the tangent of the path curve of $B$ with respect to the cam \autocite[p. 251]{Luck&Modler}.  

\begin{SCfigure}[][ht]
\includegraphics[width=0.5\textwidth]{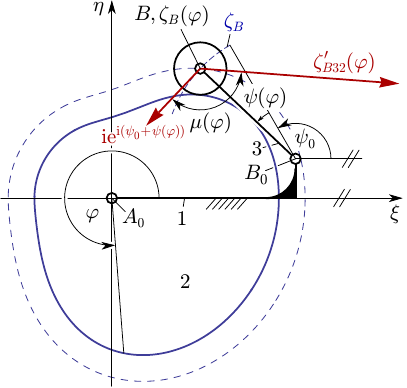}
\caption{Cam mechanism from Fig.\ \ref{Abb:Kurveng_mit_Kruemmungskreisen} with transmission angle $\mu = \mu(\ph)$.\\[0.2cm]
Position $\ph = 274.4098\g$ with\\ $\max\mu \approx \mu(274.4098\g) = 129.1068\g$\\[0.2cm]
For better visibility, the vector $\ii\ee^{\ii(\psi_0+\psi(\ph))}$ was stretched by the factor 30.}
\label{Abb:Kurvengetriebe_mu}
\end{SCfigure}

The force transmission between the cam and the roller is normal to the roller centre curve (normal force $\vv*{F}{n}$).
This normal force $\vv*{F}{n}$ can be decomposed into a tangential force component $\vv*{F}{t}$ in the direction of the path curve of $B$ with respect to the frame and into a radial force component~$\vv*{F}{r}$ in the direction of the line $B_0B$.
The smaller the angle $\alpha$ {\em (pressure angle)}\index{pressure angle} between $\vv*{F}{t}$ and $\vv*{F}{n}$, the better is the force transmission from the cam to the roller and the follower.
Angle $\alpha$ and transmission angle $\mu$ add up to $90\g$.
A particularly good force transmission can therefore be expected with a transmission angle close to $90\g$.
Note: Especially for high-speed cam mechanisms with high forces (mass or inertia forces, etc.), a design using methods of machine dynamics is essential (see e.\,g.\ \textcite{Dresig&Vulfson}).
Here we define the transmission angle (function) without the restriction to the interval $[0,\pi/2]$ by means of: 

\begin{defin}
The {\em transmission angle function}\index{transmission angle function} $\mu = \mu(\ph)$ (see Fig.\ \ref{Abb:Kurvengetriebe_mu}) is the angle $\mu$, depending on $\ph$, between the tangent unit vector $\ii\ee^{\ii(\psi_0+\psi(\ph))}$ of the curve $\zeta_B$ at the point $\ph$, and the tangent vector
\beq
  \zeta_{B32}'(\ph)
= z_B'(\ph)\,\ee^{\ii\ph}  
\eeq
with $z_B'(\ph)$ according to \eqref{Eq:z_B'(phi)}.\footnote{Clearly, $\zeta_{B32}'(\ph)$ is the tangent vector of the curve $z_B$ (see Fig.\ \ref{Abb:Kurvengetriebe_zB_und_zK}) rotated by the angle $\ph$ around point $A_0$ at the point $\ph$. It is also, and therefore the notation, the first order transfer function of the roller midpoint $B$, as point of the follower 3, with respect to the cam 2.}
\end{defin}

So using the inner product (scalar product), from \eqref{Eq:Angle_with_inner_product} we immediately get
\begin{align} \label{Eq:Transmission_angle_with_complex_terms}
  \mu(\ph)
= {} & \measuredangle\left(\ii\ee^{\ii(\psi_0+\psi(\ph))},\, \zeta_{B32}'(\ph)\right)  
= \arccos\frac{\left\langle \ii\ee^{\ii(\psi_0+\psi(\ph))},\, \zeta_{B32}'(\ph) \right\rangle}{\big|\ii\ee^{\ii(\psi_0+\psi(\ph))}\big| \left|\zeta_{B32}'(\ph)\right|}\db\nonumber\\[0.05cm]
= {} & \arccos\frac{\left\langle \ii\ee^{\ii(\psi_0+\psi(\ph))},\, \zeta_{B32}'(\ph) \right\rangle}{\left|\zeta_{B32}'(\ph)\right|}\,.  
\end{align}
Eq.\ \eqref{Eq:Transmission_angle_with_complex_terms} assigns to each angle $\ph\in[0,2\pi)$ an angle $\mu(\ph)$ uniquely determined in the interval $[0,\pi]$.

\begin{thm}
The transmission angle function \eqref{Eq:Transmission_angle_with_complex_terms} is  in real form given by
\beq
  \mu(\ph)
= \arccos\left(\frac
	{\ell\left(\psi'(\ph)-1\right) - \xi_{B_0}\cos(\psi_0+\psi(\ph)) - \eta_{B_0}\sin(\psi_0+\psi(\ph))}
	{|z_B'(\ph)|}\right)  
\eeq
with $|z_B'(\ph)|$ according to \eqref{Eq:|z_B'(phi)|}.
\end{thm}

\begin{proof}
Taking into account that
\begin{align*}
  \zeta_{B32}'(\ph)
= {} & \left(-\ii\,\zeta_{B_0}\,\ee^{-\ii\ph} + \ii\,\ell\,(\psi'(\ph)-1)\,\ee^{\ii(\psi_0+\psi(\ph)-\ph)}\right) \ee^{\ii\ph}\\[0.05cm]
= {} & \ii\,\ell\,(\psi'(\ph)-1)\,\ee^{\ii(\psi_0+\psi(\ph))} - \ii\,\zeta_{B_0}
\end{align*}
we get
\begin{align*}
  \left\langle \ii\ee^{\ii(\psi_0+\psi(\ph))},\, \zeta_{B32}'(\ph) \right\rangle
= {} & \left\langle \ii\ee^{\ii(\psi_0+\psi(\ph))},\, \ii\,\ell\,(\psi'(\ph)-1)\,\ee^{\ii(\psi_0+\psi(\ph))} - \ii\,\zeta_{B_0} \right\rangle\db\\[0.05cm]
= {} & \left\langle \ee^{\ii(\psi_0+\psi(\ph))},\, \ell\,(\psi'(\ph)-1)\,\ee^{\ii(\psi_0+\psi(\ph))} - \zeta_{B_0} \right\rangle\db\\[0.05cm]   
= {} & \left\langle \ee^{\ii(\psi_0+\psi(\ph))},\, \ell\,(\psi'(\ph)-1)\,\ee^{\ii(\psi_0+\psi(\ph))} \right\rangle
	- \left\langle \ee^{\ii(\psi_0+\psi(\ph))},\, \zeta_{B_0} \right\rangle\db\\[0.05cm]
= {} & \left\langle 1,\, \ell\,(\psi'(\ph)-1) \right\rangle
	- \left\langle \ee^{\ii(\psi_0+\psi(\ph))},\, \zeta_{B_0} \right\rangle\db\\[0.05cm]
= {} & \ell\left(\psi'(\ph)-1\right) - \xi_{B_0}\cos(\psi_0+\psi(\ph)) - \eta_{B_0}\sin(\psi_0+\psi(\ph))\,.
\end{align*}
From \eqref{Eq:Transmission_angle_with_complex_terms} with $|\zeta_{B32}'(\ph)| = \left|z_B'(\ph)\,\ee^{\ii\ph}\right| = \left|z_B'(\ph)\right|$ and $|z_B'(\ph)|$ according to \eqref{Eq:|z_B'(phi)|} the result of the theorem follows.
\end{proof}

For the design of cam mechanisms with pivoted roller follower, experience has shown that the transmission angle function should be within the following intervals:
\begin{itemize}[leftmargin=0.5cm]
\item $45\g \le \mu(\ph) \le 135\g$, $0\g \le \ph \le 360\g$, if for the input angular velocity $\dot{\ph}$ the inequality $\dot{\ph} \leq 2\pi \cdot 30\tn{min}^{-1}$ holds, 
\item $60\g \le \mu(\ph) \le 120\g$, $0\g \le \ph \le 360\g$, if $\dot{\ph} > 2\pi \cdot 30\tn{min}^{-1}$ holds \autocite[p.~252]{Luck&Modler}, \autocite[p. 17]{VDI2142-1-engl}.
\end{itemize}

\begin{SCfigure}[][ht]
\includegraphics[width=0.5\textwidth]{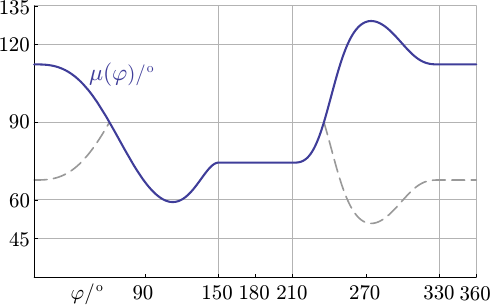}
\caption{Transmission function $\mu$ for the cam mechanism in Fig.\ \ref{Abb:Kurveng_mit_Kruemmungskreisen}}
\label{Abb:Kurvengetriebe_mu(phi)}
\end{SCfigure}

Fig.\ \ref{Abb:Kurvengetriebe_mu(phi)} shows the transmission angle function $\mu$ for the mechanism in Fig.\ \ref{Abb:Kurvengetriebe_mu}.
It can be seen that the mechanism is only recommended for cam's angular velocity $\dot{\ph} \leq 2\pi \cdot 30\tn{min}^{-1}$.
The graph modified with the dashed line is the graph of the conventional transfer angle function with values $\mu(\ph)$, $0 \le \mu(\ph) \le 90\g$.
(For this function the above conditions are $\mu(\ph) \ge 45\g$ and $\mu(\ph) \ge 60\g$, respectively.) 

\begin{figure}[H]
\begin{minipage}{0.45\textwidth}
\centering
\includegraphics[width=\textwidth]{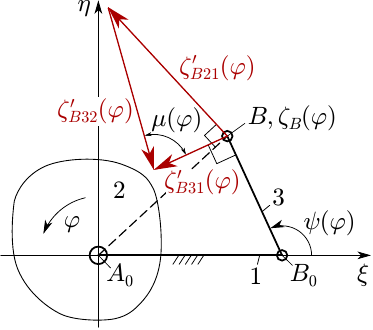}
\caption{Vectors (transfer functions of first order) $\zeta_{B21}'(\ph)$, $\zeta_{B31}'(\ph)$, $\zeta_{B32}'(\ph)$, and transmission angle $\mu(\ph)$}
\label{Abb:Kurvengetriebe_mit_mu01}
\end{minipage}
\hfill
\begin{minipage}{0.45\textwidth}
\centering
\includegraphics[width=\textwidth]{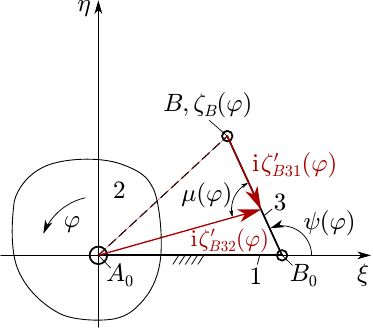}
\caption{Vectors $\zeta_{B31}'(\ph)$, $\zeta_{B32}'(\ph)$, each with angle $90\g$ rotated in the rotation sense of the cam, and transmission angle $\mu(\ph)$}
\label{Abb:Kurvengetriebe_mit_mu02}
\end{minipage}

\vspace{0.5cm}

\begin{minipage}{0.45\textwidth}
\centering
\includegraphics[width=\textwidth]{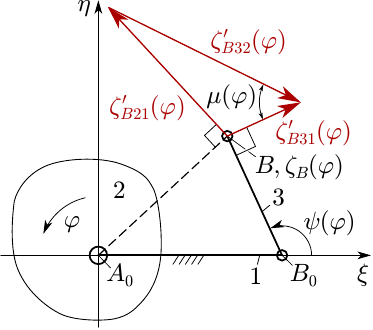}
\caption{Vectors (transfer functions of the first order) $\zeta_{B21}'(\ph)$, $\zeta_{B31}'(\ph)$, $\zeta_{B32}'(\ph)$, and transmission angle $\mu(\ph)$}
\label{Abb:Kurvengetriebe_mit_mu03a}
\end{minipage}
\hfill
\begin{minipage}{0.45\textwidth}
\centering
\includegraphics[width=\textwidth]{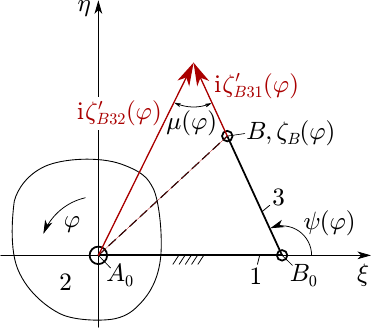}
\caption{Vectors $\zeta_{B31}'(\ph)$, $\zeta_{B32}'(\ph)$, each with angle $90\g$ rotated in the rotation sense of the cam, and transmission angle $\mu(\ph)$}
\label{Abb:Kurvengetriebe_mit_mu03b}
\end{minipage}
\end{figure}

In order to derive the {\em hodograph method}\index{hodograph method} described in Section \ref{Sec:A0-regions}, we make the following usual preliminary consideration (see \textcite[pp. 374-375]{Volmer:Lehrbuch}, \textcite[pp. 253-254]{Luck&Modler}), but we consistently use complex-valued transfer functions instead of velocities, drawing scales and velocity scales.
The zero-order transfer function of the point $B$ regarded as point of the cam 2 with respect to the frame 1 is given by $\zeta_{B21}(\ph) = \overline{A_0B}\,\ee^{\ii(\ph_0+\ph)}$ with $\ph_0 = \tn{const}$.
For the corresponding first-order transfer function holds
\beq
  \zeta_{B21}'(\ph) 
= \ii\,\overline{A_0B}\,\ee^{\ii(\ph_0+\ph)}\,, 
\eeq
hence
\beq
  \left|\zeta_{B21}'(\ph)\right| 
= \big|\ii\,\overline{A_0B}\,\ee^{\ii(\ph_0+\ph)}\big|
= |\ii| \cdot \overline{A_0B} \cdot |\ee^{\ii(\ph_0+\ph)}\big|
= \overline{A_0B}\,.
\eeq
From
\beq
  \zeta_{B12}'(\ph) + \zeta_{B23}'(\ph) + \zeta_{B31}'(\ph) = 0
  \qquad\mbox{and}\qquad
  \zeta_{Bki}'(\ph) = -\zeta_{Bik}'(\ph)
\eeq
(see \textcite[p.\ 111]{Luck&Modler}) follows
\beq
  -\zeta_{B21}'(\ph) - \zeta_{B32}'(\ph) + \zeta_{B31}'(\ph) = 0
  \quad\Longrightarrow\quad
  \zeta_{B32}'(\ph) = \zeta_{B31}'(\ph) - \zeta_{B21}'(\ph)\,;
\eeq
see Fig.\ \ref{Abb:Kurvengetriebe_mit_mu01}.

Rotation of the vectors $\zeta_{B31}'(\ph)$ and $\zeta_{B32}'(\ph)$ by $90\g$ in the rotational direction of the cam results in the situation in Fig.\ \ref{Abb:Kurvengetriebe_mit_mu02}, which allows a particularly simple determination of the transmission angle $\mu(\ph)$, for which only the vector $\zeta_{B31}'(\ph)$ needs to be known:
Drawing a line from the cam pivot point $A_0$ to the arrowhead of the vector $\zeta_{B31}'(\ph)$ attached to point $B$ and rotated by $90\g$ in the rotational direction of the cam $\zeta_{B31}'(\ph)$, we get $\mu(\ph)$ immediately as the acute angle between this line and the line $B_0B$.
$\zeta_{B31}'(\ph)$ follows from
\beq
  \zeta_{B31}(\ph)
= \zeta_B(\ph)
= \ell_1 + \overline{B_0 B}\, \ee^{\ii\psi(\ph)}   
\eeq
(see Fig.\ \ref{Abb:Kurvengetriebe_zetaB}) to
\beq
  \zeta_{B31}'(\ph)
= \zeta_B'(\ph)  
= \overline{B_0 B}\, \ii\, \psi'(\ph)\, \ee^{\ii\psi(\ph)}\,.  
\eeq
With positive rotational direction of the cam and negative rotational direction of the follower, the situation in Fig.\ \ref{Abb:Kurvengetriebe_mit_mu03a} results.
By rotating $\zeta_{B31}'(\ph)$ and $\zeta_{B32}'(\ph)$ by $90\g$ in the rotational direction of the cam, we obtain the situation in Fig.\ \ref{Abb:Kurvengetriebe_mit_mu03b}.
So, again, only a line from $A_0$ to the arrowhead of the vector $\ii \zeta_{B31}'(\ph)$ attached to $B$ has to be drawn to be able to measure the transmission angle $\mu(\ph)$. 

\subsubsection[$A_0$-regions (hodograph method)\index{hodograph method}]{$\boldsymbol{A_0}$-regions (hodograph method)}
\index{hodograph method}\index{A0-region@$A_0$-region}
\label{Sec:A0-regions}

In the following we need this definition (see \cite[p.\ 73]{Volmer:Kurvengetriebe}, \cite[p.\ 372]{Volmer:Lehrbuch}, \cite[pp.\ 252-253]{Luck&Modler}):

\begin{defin}
For given direction of cam rotation, a {\em P-cam mechanism}\index{P-cam mechanism} (P $=$ centri\underline{p}etal) is one where the distance between the cam pivot point $A_0$ and and the roller center point $B$ decreases when the roller follower rotates around $B_0$ in the same direction as the cam; a {\em F-cam mechanism}\index{F-cam mechanism} (F $=$ centri\underline{f}ugal) is one where this distance increases.\footnote{Note that this definition actually depends on the cam rotation direction.}
\end{defin}

\begin{SCfigure}[][ht]
\includegraphics[width=0.7\textwidth]{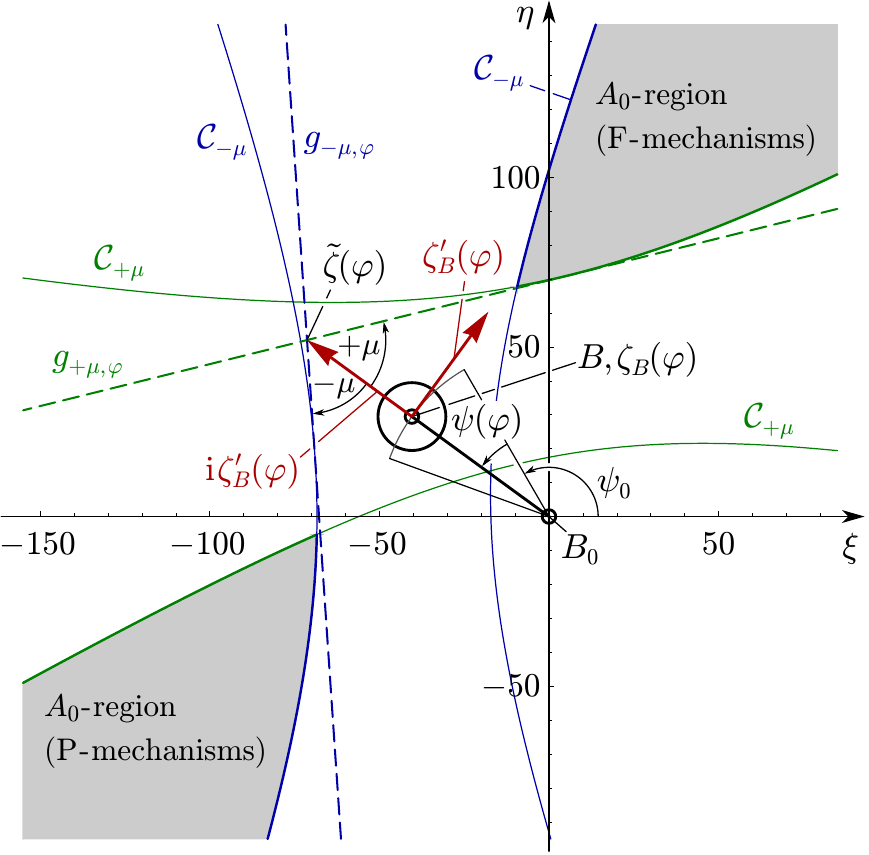}
\caption{Determining the $A_0$-regions by means of the hodograph method for $\mu = 50\g$.\\[0.2cm]
Position:\\
\hspace*{0.5cm} $\ph = 265\g$,\\
\hspace*{0.5cm} $\psi(265\g) \approx 23.8935\g$.\\[0.2cm]
During the rise, the lower branch of $\mathcal{C}_{+\mu}$ and the right branch of $\mathcal{C}_{-\mu}$ are generated; during the return the upper branch of $\mathcal{C}_{+\mu}$ and the left branch of $\mathcal{C}_{-\mu}$.}
\label{Abb:A0-Bereiche}
\end{SCfigure}

\begin{SCfigure}[][ht]
\includegraphics[width=0.7\textwidth]{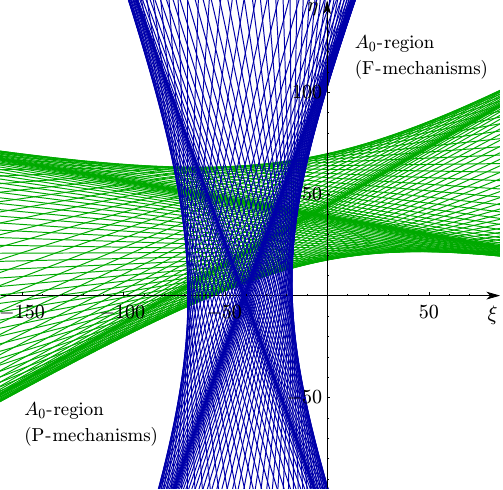}
\caption{Families of lines $g_{+\mu,\ph}$ (green) and $g_{-\mu,\ph}$ (blue)}
\label{Abb:Families_of_lines}
\end{SCfigure}

From the relationship between the vector $\zeta_B'(\ph) \equiv \zeta_{B31}'(\ph)$ rotated in the rotational direction of the cam by $90\g$ and attached to $B$, the transmission angle $\mu(\ph)$, and the cam pivot point $A_0$, shown in Fig.\ \ref{Abb:Kurvengetriebe_mit_mu02} and Fig.\ \ref{Abb:Kurvengetriebe_mit_mu03b}, the so-called {\em hodograph method}\index{hodograph method} (see \textcite[pp. 375-379]{Volmer:Lehrbuch}, \textcite[pp.\ 255-259]{Luck&Modler}) results.
In this method the desired/required (not to fall below) minimum value $\mu > 0$ of the transmission angle is given; with this, the regions in which the cam pivot point $A_0$ has to be located are determined (see Fig.\ \ref{Abb:A0-Bereiche}).
The angles $+\mu$ and $-\mu$ are attached to the vector $\ii\zeta_{B}'(\ph)$ in its arrowhead $\tilde{\zeta}(\ph)$.
(We assume in Fig.\ \ref{Abb:A0-Bereiche} mathematical positive rotational direction of the cam. Since $A_0$ is not yet known, the follower pivot point $B_0$ was placed in the coordinate origin.)
This gives the two straight lines $g_{+\mu,\ph}$ and $g_{-\mu,\ph}$.
Since $\ph$ is varying, $\psi(\ph)$, $\zeta_{B}'(\ph)$ and $\ii\zeta_{B}'(\ph)$ are also varying; $g_{+\mu,\ph}$ and $g_{-\mu,\ph}$ are now families of straight lines\index{family of lines} (see Fig.\ \ref{Abb:Families_of_lines}).
The envelopes\index{envelope} $\mathcal{C}_{+\mu}$ and $\mathcal{C}_{-\mu}$ of these families bound the two $A_0$-regions.
Depending on the region in which the point $A_0$ is chosen, one obtains a P-cam mechanism\index{P-cam mechanism} or a F-cam mechanism\index{F-cam mechanism}.
If the rotation direction of the cam is negative, then the $A_0$-regions change their roles. 

\begin{thm}
The equations of the envelopes $\mathcal{C}_{\pm\mu}$ are given by
\beq
  \zeta_{\pm\mu}(\ph,\la_{\pm\mu}(\ph))
= \zeta_{B_0} + \left[\ell\left(1-\psi'(\ph)\right) + \la_{\pm\mu}(\ph)\,\ee^{\pm\ii\mu}\right] \ee^{\ii(\psi_0+\psi(\ph))}
\eeq
with
\beqn \label{Eq:lambda_with_pm_only}
  \la_{\pm\mu}(\ph)
= -\ell\left(\left(1-\psi'(\ph)\right)\cos(\pm\mu) + \frac{\psi''(\ph)}{\psi'(\ph)}\,\sin(\pm\mu)\right).  
\eeqn
\end{thm}

\begin{proof}
The equation for point $B$ is
\beq
  \zeta_B(\ph) = \zeta_{B0} + \ell\, \ee^{\ii(\psi_0+\psi(\ph))}\,.\footnote{Here, unlike Fig.\ \ref{Abb:A0-Bereiche}, the point $B_0$ can be in general position.}
\eeq
From this follows
\beq
  \zeta_{B}'(\ph) = \ii\,\ell\,\psi'(\ph)\,\ee^{\ii(\psi_0+\psi(\ph))}\,.
\eeq
We rotate this vector in the mathematical positive rotational direction of the cam by $90\g$ and get
\beq
  \ii\,\zeta_{B}'(\ph) = -\ell\,\psi'(\ph)\,\ee^{\ii(\psi_0+\psi(\ph))}\,.
\eeq
(Note that the vectors $\ell\,\ee^{\ii(\psi_0+\psi(\ph))}$ and $\ii\,\zeta_{B}'(\ph)$ have the same direction for $\psi'(\ph) < 0$ as shown in Fig.~\ref{Abb:A0-Bereiche}.)
Thus the equation for the point of the follower (or its extension), through which pass the lines $g_{+\mu,\ph}$ und $g_{-\mu,\ph}$ is as follows:
\beq
  \tilde{\zeta}(\ph)
= \zeta_B(\ph) + \ii\,\zeta_B'(\ph)  
= \zeta_{B_0} + \ell \left(1 - \psi'(\ph)\right) \ee^{\ii(\psi_0+\psi(\ph))}\,.
\eeq
For given minimum value $\mu > 0$ of the transmission angle, the equations of the lines $g_{\pm\mu,\ph}$ are
\beq
  \zeta_{\pm\mu}(\ph,\la)
= \tilde{\zeta}(\ph) + \la\, \ee^{\ii(\psi_0+\psi(\ph)\pm\mu)}\,,\quad \la \in \R\,,  
\eeq
hence
\beq
  \zeta_{\pm\mu}(\ph,\la)
= \zeta_{B_0} + \left[\ell\left(1-\psi'(\ph)\right) + \la\,\ee^{\pm\ii\mu}\right] \ee^{\ii(\psi_0+\psi(\ph))}
\eeq
(cf.\ \eqref{Eq:z(phi,lambda)}).

Now we determine the envelopes\index{envelope} $\mathcal{C}_{\pm\mu}$ of the families of lines\index{family of lines} $g_{\pm\mu,\ph}$, $0 \le \ph < 2\pi$.
Using Theorem \ref{Thm:Envelope} and Eq.\ \eqref{Eq:[z_1,z_2]_mit_Im}, we know that
\beq
  0
= \left[\frac{\p\zeta}{\p\ph}, \frac{\p\zeta}{\p\la}\right]
= -\Imz\left(\tfrac{\p\zeta}{\p\ph}\: \overline{\left(\tfrac{\p\zeta}{\p\la}\right)}\right)
= -\Imz\left(\frac{\p\zeta}{\p\ph}\, \frac{\p\overline{\zeta}}{\p\la}\right).   
\eeq
(So again the value $\la = \la(\ph)$ has to be determined for the point where the straight line is tangent to the (still unknown) envelope.)
We have
\begin{align*}
  \frac{\p\zeta_{\pm\mu}(\ph,\la)}{\p\ph}
= {} & {-}\ell\, \psi''(\ph)\, \ee^{\ii(\psi_0+\psi(\ph))}
	+ \left[\ell\left(1-\psi'(\ph)\right) + \la_{\pm\mu}\,\ee^{\pm\ii\mu}\right] \ii\,\psi'(\ph)\,\ee^{\ii(\psi_0+\psi(\ph))}\\
= {} & \left[-\ell\,\psi''(\ph) + \ii\,\ell\,\psi'(\ph)\left(1-\psi'(\ph)\right)
	+ \ii\,\la_{\pm\mu}\,\psi'(\ph)\,\ee^{\pm\ii\mu}\right] \ee^{\ii(\psi_0+\psi(\ph))}\,,\db\\[0.3cm]
  \overline{\zeta_{\pm\mu}(\ph,\la)}
= {} & \overline{\zeta_{B_0}} + \ell\left(1-\psi'(\ph)\right)\ee^{-\ii(\psi_0+\psi(\ph))} + \la_{\pm\mu}\,\ee^{\mp\ii\mu}\,\ee^{-\ii(\psi_0+\psi(\ph))}\,,
  \\[0.3cm]
  \frac{\p\overline{\zeta_{\pm\mu}(\ph,\la)}}{\p\la}
= {} & \ee^{\mp\ii\mu}\,\ee^{-\ii(\psi_0+\psi(\ph))}\,,     	 
\end{align*}
hence
\begin{align*}
  \Imz\left(\frac{\p\zeta_{\pm\mu}}{\p\ph}\, \frac{\p\overline{\zeta_{\pm\mu}}}{\p\la}\right)
= {} & \Imz\left[-\ell\,\psi''(\ph)\,\ee^{\mp\ii\mu} + \ii\,\ell\,\psi'(\ph)\left(1-\psi'(\ph)\right)\,\ee^{\mp\ii\mu}
	+ \ii\,\la\,\psi'(\ph)\right]\\[-1mm]
= {} & {\pm}\ell\,\psi''(\ph)\,\sin\mu + \ell\,\psi'(\ph)\left(1-\psi'(\ph)\right)\cos\mu + \la\,\psi'(\ph)\\[2mm]
= {} & 0\,.
\end{align*}
From this follows
\beqn \label{Eq:lambda_mit_pm_und_mp}
  \la_{\pm\mu}
= \la_{\pm\mu}(\ph)  
= \ell\left(\left(\psi'(\ph)-1\right)\cos\mu \mp \frac{\psi''(\ph)}{\psi'(\ph)}\,\sin\mu\right).  
\eeqn
The minus sign of ``$\mp$'' in \eqref{Eq:lambda_mit_pm_und_mp} is for $+\mu$, hence for the envelope $\mathcal{C}_{+\mu}$ of the family of lines $g_{+\mu,\ph}$, while the plus sign of ``$\mp$'' is for $-\mu$, hence for the envelope $\mathcal{C}_{-\mu}$ of the family of lines $g_{-\mu,\ph}$.
Thus, \eqref{Eq:lambda_mit_pm_und_mp} can be written as \eqref{Eq:lambda_with_pm_only}.
\end{proof}

In \cite{Ji&Manna}, the pressure angle is used as criterion, whereby the case of a translating roller follower is also taken into account.

\begin{SCfigure}[][ht]
\includegraphics[width=0.58\textwidth]{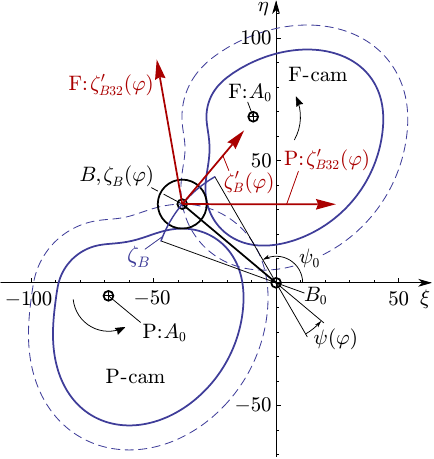}
\caption{Minimum size cams for $\mu = 50\g$ (Position $\ph = 270\g$, $\psi(270\g) = 20\g$)}
\label{Abb:P_and_F-cam}
\end{SCfigure}

\begin{figure}[ht]
\begin{minipage}[b]{0.48\textwidth}
\includegraphics[width=\textwidth]{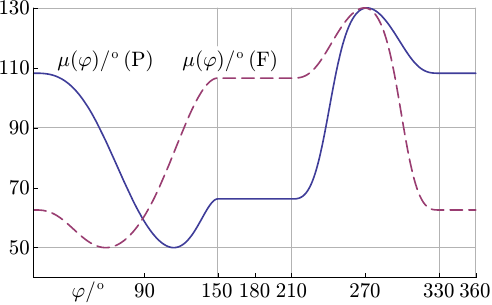}
\caption{Transmission angle functions for the P-mechanism and the F-mechanism in Fig.~\ref{Abb:P_and_F-cam}}
\label{Abb:mu_P_und_F}
\end{minipage}
\hfill
\begin{minipage}[b]{0.48\textwidth}
\includegraphics[width=\textwidth]{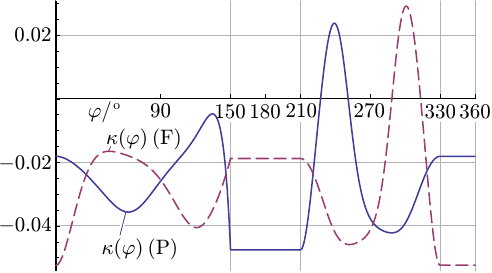}
\caption{Curvature functions for the contour curves of the P-cam and the F-cam in Fig.~\ref{Abb:P_and_F-cam}}
\label{Abb:kappa_P_und_F}
\end{minipage}
\end{figure}

\begin{example}
An example for the determination of the $A_0$-regions is shown in Fig.\ \ref{Abb:A0-Bereiche}.
The transfer function $\psi$ is given by \eqref{Eq:psi(phi)}.
Parameters used are $\ell = 50$ and $\psi_0 = 120\g$ as in Figures \ref{Abb:Kurveng_mit_Kruemmungskreisen} and \ref{Abb:Kurvengetriebe_mu}.
Differently from Fig.\ \ref{Abb:Kurvengetriebe_mu}, we have placed the pivot point $B_0$ ($\zeta_{B_0}$) of the follower in the coordinate origin. 
Each of the envelopes $\mathcal{C}_{+\mu}$ and $\mathcal{C}_{-\mu}$ consists of two branches.
The values $\ph_1$ and $\ph_2$ of $\ph$ at an intersection point of the curves $\mathcal{C}_{+\mu}$ and $\mathcal{C}_{-\mu}$ are solutions of the nonlinear systems of equations
\beq
\left.\begin{array}{c@{\;=\;}l}
  \Rez\left[\zeta_{+\mu}(\ph_1,\la_{+\mu}(\ph_1))\right] & \Rez\left[\zeta_{-\mu}(\ph_2,\la_{-\mu}(\ph_2))\right],\\[0.15cm]
  \Imz\left[\zeta_{+\mu}(\ph_1,\la_{+\mu}(\ph_1))\right] & \Imz\left[\zeta_{-\mu}(\ph_2,\la_{-\mu}(\ph_2))\right].  
\end{array}\right\}
\eeq
It can be solved for suitable initial values by means of Newton's approximation method\index{Newton's approximation method} or e.\,g. in {\em Mathematica}\index{Mathematica} with the function \texttt{FindRoot}\index{FindRoot@\texttt{FindRoot}}.  
With the solution $\ph_1 \approx 114.142\g$, $\ph_2 \approx 271.485\g$ we get the intersection point $\xi \approx -68.4388$, $\eta \approx -5.3116$, and with the solution $\ph_1 \approx 269.281\g$, $\ph_2 \approx 58.7963\g$ the intersection point $\xi \approx -9.34509$, $\eta \approx 67.8047$.
If we take the first intersection point as pivot point $A_0$, we get the P-cam mechanism with smallest cam for $50\g \le \mu(\ph) \le 130\g$ (see Fig.\ \ref{Abb:P_and_F-cam}).
Putting $A_0$ into the second intersection point, gives us the F-cam mechanism with smallest cam for this interval (see also Fig.\ \ref{Abb:P_and_F-cam}).
The maximum radii of the P-cam and F-cam are 55.1931 and 53.1475, respectively.
The transmission angle functions $\mu(\ph)$ for both mechanisms are shown in Fig.\ \ref{Abb:mu_P_und_F} (cf.\ Fig.\ \ref{Abb:Kurvengetriebe_mu(phi)}), the curvature functions of the contour curves in Fig.\ \ref{Abb:kappa_P_und_F} (cf.\ Fig.\ \ref{Abb:Curvature_function}). \hfill\bs
\end{example}

\clearpage

%% file: DiffGeo5_5.tex

\section{Application in the field of machine elements: polygon profiles}
\subsection{Basics}

\begin{SCfigure}[][ht]
\includegraphics[width=0.4\textwidth]{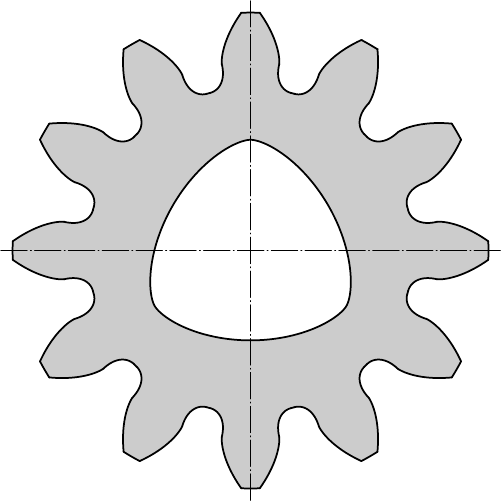}
\caption{Gear with polygon profile (P3G profile)}
\label{Abb:Gear_with_polygon_profile}
\end{SCfigure}

In mechanical engineering, the connection of a rotating machine part\index{machine part} (hub\index{hub}), e.g.\ a gear\index{gear}, to a shaft is an important and frequent task.
For this purpose, so-called polygon profiles\index{profile!polygon profile@polygon $\sim $}\footnote{For the derivation of the polygon profiles from K-profiles\index{K-profile} (epitrochoids\index{epitrochoid}) see \cite{Musyl55}; for K-profiles see \cite{Musyl46}.} can be used as form-fit connections (for an example see Fig.\ \ref{Abb:Gear_with_polygon_profile}).
These profiles are not polygons in the mathematical sense of the term.
The German industrial standards \textcite{DIN32711-1-engl} and \textcite{DIN32712-1-engl} define the so-called P3G profiles or P4C profiles, the calculation and dimensioning is subject of the standards DIN 32711-2 and DIN 32712-2 (see also \cite{Filemon}, \cite[pp.\ 147, 150-151]{Haberhauer18}, \cite[pp.\ 529-536]{Steinhilper&Sauer_Band_1}, \cite{Ziaei07a-engl}, \cite{Ziaei07b-engl}, and, for further references, \cite[p.\ 441]{Martini&Montejano&Oliveros}).

From \textcite{DIN32711-1-engl} and \textcite{DIN32712-1-engl} we generalize the following definition (see also \cite{HEXAGON_PnG}, \cite{HEXAGON_PnC}):

\begin{defin}
A parametric curve given by the real parametric representation
\beqn \label{Eq:Profile_in_kartesischen_Koordinaten}
\left.
\begin{array}{*{3}{c@{\;}}l@{\;}c@{\;}l}
  x(\ph) & = & [R-e\cos(n\ph)]\cos\ph & - & ne\sin(n\ph)\sin\ph\,,\\[0.1cm]
  y(\ph) & = & [R-e\cos(n\ph)]\sin\ph & + & ne\sin(n\ph)\cos\ph\,,
\end{array}
\right\} \quad 0 \le \ph \le 2\pi\,,
\eeqn
with $n\in\N$ and real parameters $R > 0$ and $e > 0$ is referred to as P$n$ {\em curve}.

A compact convex subset of $\R^2$ is called P$n$G {\em profile} if its boundary is given by the complete curve \eqref{Eq:Profile_in_kartesischen_Koordinaten}; it is called P$n$C {\em profile} if its boundary is alternately piecewise defined by \eqref{Eq:Profile_in_kartesischen_Koordinaten} and arcs of a circle with center at coordinate origin.\footnote{For the shaft we consider the material, for the hub, as in Fig.\ \ref{Abb:Gear_with_polygon_profile}, the hole. Not all profiles according to this definition are technically feasible.}
\end{defin}

\begin{figure}[h!t]
\begin{minipage}[t]{0.45\textwidth}
  \includegraphics[width=\textwidth]{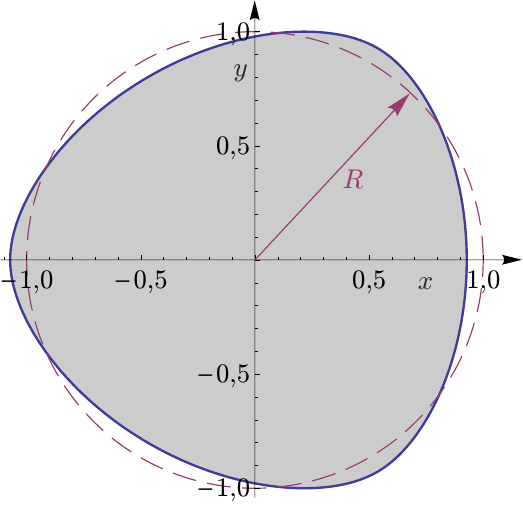}
  \caption{\small P3 curve (blue) and P3G profile (gray shaded) with $R = 1$ and $e = 0.072$\\[0.2cm]}
  \label{Abb:P3G-Profil}
\end{minipage}
\hfill
\begin{minipage}[t]{0.45\textwidth}
  \includegraphics[width=\textwidth]{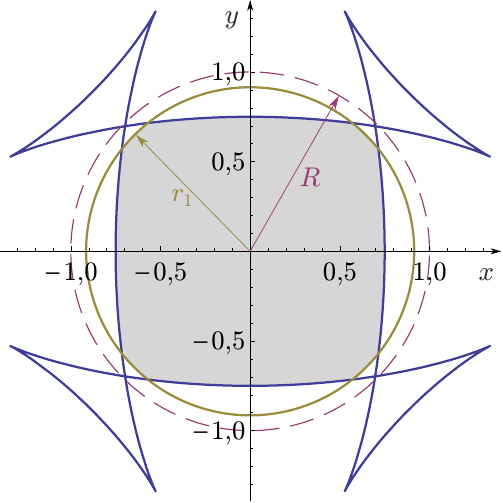}
  \caption{\small P4 curve (blue) and P4C profile (gray shaded) with $R = 1$, $e = 0.25$ and $r_1 = 0.914$}
  \label{Abb:P4C-Profil}
\end{minipage}
%
\begin{minipage}[t]{0.45\textwidth}
  \includegraphics[width=\textwidth]{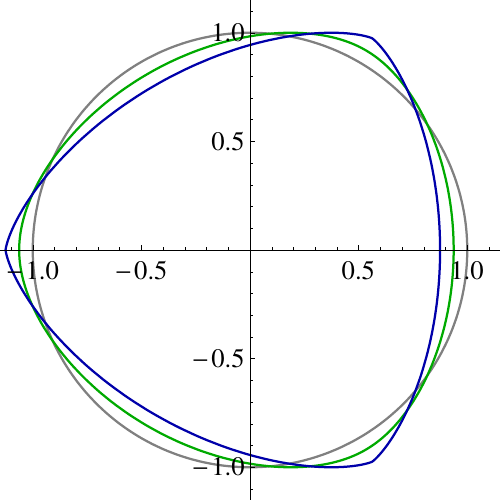}
  \caption{\small P3 curves with $R = 1$ and $e = 1/16$, $1/8$, and unit circle}
  \label{Abb:P3G-Profile}
\end{minipage}
\hfill
\begin{minipage}[t]{0.45\textwidth}
  \includegraphics[width=\textwidth]{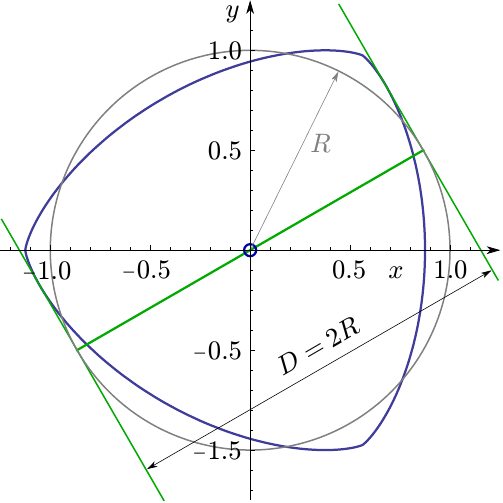}
  \caption{\small P3G profile with constant width (diameter) $D = 2R$}
  \label{Abb:P3G_Gleichdick}
\end{minipage}
\end{figure}

Fig.\ \ref{Abb:P3G-Profil} shows an example of a P3 curve with corresponding P3G profile\index{profile!P3G profile@P3G $\sim $}.
Fig.\ \ref{Abb:P4C-Profil} shows an example of a P4 curve and a P4C profile\index{profile!P4C profile@P4C $\sim $} bounded piecewise by parts of the P4 curve and arcs of the circle of radius $r_1$.
Fig.\ \ref{Abb:P3G-Profile} shows P3 curves with different values of the eccentricity $e$.

The parametric representation \eqref{Eq:Profile_in_kartesischen_Koordinaten} can in complex form be written as
\begin{align} \label{Eq:z(phi)_with_cos_and_sin}
  z(\ph)
= {} & \left(R - e\cos(n\ph)\right)\ee^{\ii\ph} + ne\sin(n\ph) \left(\ii\cos\ph - \sin\ph\right)\nonumber\db\\[0.02cm] 
= {} & \left(R - e\cos(n\ph)\right)\ee^{\ii\ph} + \ii ne\sin(n\ph) \left(\cos\ph + \ii\sin\ph\right)\nonumber\db\\[0.02cm]
= {} & \left(R - e\cos(n\ph)\right)\ee^{\ii\ph} + \ii ne\sin(n\ph)\,\ee^{\ii\ph}\nonumber\\[0.02cm]
= {} & \ee^{\ii\ph} \left(R - e\cos(n\ph) + \ii ne\sin(n\ph)\right)
\end{align}
with derivative
\beqn \label{Eq:z'(phi)_with_cos_and_sin}
  z'(\ph)
= \ii\ee^{\ii\ph} \left(R + \left(n^2-1\right)e\cos(n\ph)\right).  
\eeqn
Clearly, a P$n$ curve can only be the boundary of a P$n$G profile if it has no self-intersections; the curve must therefore be free of singularities apart from possible cusps.
$z'(\ph)$ can only be equal to zero if the expression in the outer brackets vanishes.
A P$n$G profile is therefore only obtained if this expression is for $0 \le \ph < 2\pi$ greater than or equal to zero.
Because of $-1 \le \cos(n\ph) \le 1$, we must have
\beqn \label{Eq:Condition_for_no_self-intersections_1}
  R - \left(n^2-1\right)e \ge 0\,,\quad n = 1,2,3,\ldots,
\eeqn
hence
\beqn \label{Eq:Condition_for_no_self-intersections_2}
  e \le \frac{R}{n^2-1}\,,\quad n = 2,3,4,\ldots.
\eeqn
If $R = \left(n^2-1\right)e$, then the boundary curve of the P$n$G profile has cusps.

\subsection{Kinematic generation of the profile curves}

\begin{figure}[ht]
\begin{minipage}[t]{0.48\textwidth}
  \centering	
  \includegraphics[width=\textwidth]{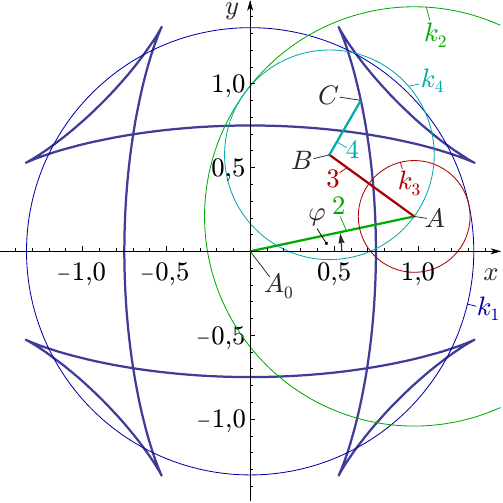}
  \caption{\small Kinematic generation of the P4 curve in Fig.\ \ref{Abb:P4C-Profil}}
  \label{Abb:Polygon_kinematisch01}
\end{minipage}
\hfill
\begin{minipage}[t]{0.48\textwidth}
  \includegraphics[width=\textwidth]{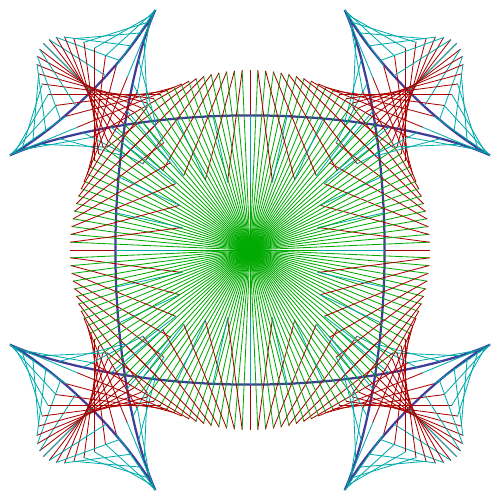}
  \caption{\small Motion sequence for the mechanism from Fig.\ \ref{Abb:Polygon_kinematisch01} with 144 positions, $\Delta\ph = 2.5\g$}
  \label{Abb:Polygon_kinematisch02}
\end{minipage}
\end{figure}

From \eqref{Eq:z(phi)_with_cos_and_sin} we get
\begin{align*}
  z(\ph)
= {} & \ee^{\ii\ph} \left(R - e\,\frac{\ee^{\ii n\ph} + \ee^{-\ii n\ph}}{2} + ne\,\frac{\ee^{\ii n\ph} - \ee^{-\ii n\ph}}{2}\right)\db\\[0.05cm]
= {} & \ee^{\ii\ph} \left(R - \frac{\left(n+1\right)e}{2}\,\ee^{-\ii n\ph} + \frac{\left(n-1\right)e}{2}\,\ee^{\ii n\ph}\right)\db\\[0.05cm]
= {} & R\,\ee^{\ii\ph} - \frac{\left(n+1\right)e}{2}\,\ee^{-\ii(n-1)\ph} + \frac{\left(n-1\right)e}{2}\,\ee^{\ii(n+1)\ph}\db\\[0.05cm]
= {} & R\,\ee^{\ii\ph} + \frac{\left(n+1\right)e}{2}\,\ee^{\ii[\pi-(n-1)\ph]} + \frac{\left(n-1\right)e}{2}\,\ee^{\ii(n+1)\ph}\,,
\end{align*}
hence
\beqn \label{Eq:z(phi)_for_profile}
  z(\ph)
= \ell_2\,\ee^{\ii\ph} + \ell_3\,\ee^{\ii[\pi-(n-1)\ph]} + \ell_4\,\ee^{\ii(n+1)\ph} 
\eeqn
with
\beq
  \ell_2 = R = \overline{A_0A}\,,\quad
  \ell_3 = \frac{1}{2}\,(n+1)\,e = \overline{AB}\,,\quad
  \ell_4 = \frac{1}{2}\,(n-1)\,e = \overline{BC}\,.\footnote{See also \cite[Eq.\ (5)]{Schoenwandt}.}
\eeq
The line segment $\overline{A_0A}$ can be regarded as a crank (part 2) rotating around the coordinate origin $A_0$ with angle $\ph$ (see Fig.\ \ref{Abb:Polygon_kinematisch01}).
At point $A$, part 2 is connected to the line segment $AB$ (part 3) with a rotary joint.
Point $B$ is the rotary joint for the connection of $AB$ and $BC$ (part 4).
Point $C$ describes the profile curve.
Part 1 is the frame.
The motion of parts 2, 3 and 4 is shown in Fig.~\ref{Abb:Polygon_kinematisch02}.

We now determine the radii of two pairs of circles connected to the parts, whose rolling on each other gives the rotation angles of the parts.
The circles can then be used as pitch or rolling circles of gears.
The angular velocity $\om_{21}$ of part 2 with respect to the frame 1 is $\dot\ph$, where the point denotes the time derivative.
For the angular velocities $\om_{31}$ and $\om_{41}$ of parts elements 3 and 4 with respect to 1, the following applies
\beq
  \om_{31} = -(n-1)\dot\ph = -(n-1)\om_{21}\,,\qquad
  \om_{41} = (n+1)\dot\ph = (n+1)\om_{21}\,.
\eeq
The circle $k_1$ with center $A_0$ and radius $\rh_1$ of the first pair belongs to frame 1.
The circle $k_3$ with center $A$ and radius $\rh_3$ of this pair belongs to coupler 3.
For the angular velocities $\om_{32}$ and $\om_{12}$ of parts 3 or 1 with respect to the crank 2, the following relationship applies
\beq
  \frac{\rh_1}{\rh_3} = \frac{\om_{32}}{\om_{12}}\,.
\eeq
We have
\begin{gather*}
  \om_{12} + \om_{23} + \om_{31} = 0
  \quad\Longrightarrow\quad
  -\om_{21} - \om_{32} + \om_{31} = 0 \\[0.05cm]
  \Longrightarrow\quad
  \om_{32} = \om_{31} - \om_{21} = -(n-1)\om_{21} - \om_{21} = - n\om_{21}\,,
\end{gather*}     
from which
\beq
  \frac{\rh_1}{\rh_3} = \frac{-n\om_{21}}{-\om_{21}} = n
\eeq
follows.
Furthermore, we have
\beq
  \ell_2
= \rh_1 - \rh_3 = n\rh_3 - \rh_3 = \rh_3(n-1)
\eeq
and consequently
\beq
  \rh_3
= \frac{1}{n-1}\,\ell_2
= \frac{1}{n-1}\,R
  \qquad\mbox{and}\qquad
  \rh_1
= \frac{n}{n-1}\,\ell_2
= \frac{n}{n-1}\,R\,.    
\eeq
The circle $k_2$ with center $A$ and radius $\rh_2$ of the second pair belongs to crank 2.
The circle $k_4$ with center $B$ and radius $\rh_4$ of this pair belongs to coupler 4.
It applies
\beq
  \frac{\rh_2}{\rh_4} = \frac{\om_{43}}{\om_{23}}\,.
\eeq
With
\begin{gather*}
  \om_{13} + \om_{34} + \om_{41} = 0
  \quad\Longrightarrow\quad
  -\om_{31} - \om_{43} + \om_{41} = 0\\[0.05cm]
  \Longrightarrow\quad
  \om_{43} = \om_{41} - \om_{31} = (n+1)\om_{21} + (n-1)\om_{21} = 2n\om_{21}
\end{gather*}
and $\om_{23} = n\om_{21}$ follows
\beq
  \frac{\rh_2}{\rh_4} = \frac{2n\om_{21}}{n\om_{21}} = 2\,.
\eeq
Now we have
\beq
  \ell_3
= \rh_2 - \rh_4
= 2\rh_4 - \rh_4
= \rh_4\,,
\eeq
thus
\beqn \label{Eq:rho_4_and_rho_2}
  \rh_4 = \ell_3 = \frac{1}{2}\,(n+1)\,e
  \qquad\mbox{and}\qquad
  \rh_2 = 2\ell_3 = (n+1)\,e\,.
\eeqn

\subsection{Area, width und length}

Now, we compute the area $A$ enclosed by a profile curve.
From Eq.\ \eqref{Eq:A} in Theorem \ref{Thm:Area_with_QVP} we know that
\beq
  A
= \frac{1}{2} \oint \left[z(\ph),z'(\ph)\right] \dd\ph
= \frac{1}{2} \int_0^{2\pi} \left[z(\ph),z'(\ph)\right] \dd\ph\,,  
\eeq
and with \eqref{Eq:z(phi)_with_cos_and_sin} and \eqref{Eq:z'(phi)_with_cos_and_sin} we get
\begin{align*}
  \left[z(\ph),z'(\ph)\right]
= {} & \left[\ee^{\ii\ph}\left(R - e\cos(n\ph) + \ii en\sin(n\ph)\right),\;
	\ii\ee^{\ii\ph}\left(R + \left(n^2-1\right)e\cos(n\ph)\right)\right]\db\\[0.05cm]  
= {} & \left[R - e\cos(n\ph) + \ii en\sin(n\ph),\;
	\ii\left(R + \left(n^2-1\right)e\cos(n\ph)\right)\right]\db\\[0.05cm]
= {} & \left\langle R - e\cos(n\ph) + \ii en\sin(n\ph),\;
	R + \left(n^2-1\right)e\cos(n\ph)\right\rangle\db\\[0.05cm]
= {} & \left(R  - e\cos(n\ph)\right)\left(R + \left(n^2-1\right)e\cos(n\ph)\right)\db\\[0.05cm]
= {} & R^2 + \left(n^2-1\right)Re\cos(n\ph)	- Re\cos(n\ph) - \left(n^2-1\right)e^2\cos^2(n\ph)\,.
\end{align*}
It follows that
\begin{align*}
  A
= {} & \frac{1}{2}\,R^2 \int_0^{2\pi} \dd\ph - \frac{1}{2}\left(n^2-1\right) e^2 \int_0^{2\pi} \cos^2(n\ph)\, \dd\ph\db\\[0.05cm]
= {} & \pi R^2 - \frac{1}{4}\left(n^2-1\right) e^2 \int_0^{2\pi} \left(1+\cos(2n\ph)\right) \dd\ph\,,
\end{align*}
thus
\beqn \label{Eq:A_profile_curve}
  A
= \pi R^2 - \frac{1}{2}\,\pi\left(n^2-1\right) e^2\,.
\eeqn
Note that if the profile curve has self-intersections (see Figs.\ \ref{Abb:Polygon_kinematisch01} and \ref{Abb:Polygon_kinematisch02}), then \eqref{Eq:A_profile_curve} is the algebraic sum of the signed areas of all loops (see Remark \ref{Rem:A} and Fig.\ \ref{Abb:Area_as_sum}).

\eqref{Eq:z'(phi)_with_cos_and_sin} immediately shows that the tangent vector $z'(\ph)$ is always orthogonal to crank 2 (with unit vector $\ee^{\ii\ph}$).
Therefore, we obtain the support function as projection of $z(\ph)$ (see \eqref{Eq:z(phi)_with_cos_and_sin}) onto the direction of the vector $\ee^{\ii\ph}$ (see \eqref{Eq:Projection_with_inner_product} and Fig.\ \ref{Abb:Projection}):
\begin{align} \label{Eq:a(phi)_profile_curve}
  a(\ph)
= {} & \left\langle z(\ph),\, \ee^{\ii\ph} \right\rangle\nonumber\db\\[0.05cm]
= {} & \left\langle \ee^{\ii\ph} \left(R - e\cos(n\ph) + \ii en\sin(n\ph)\right),\, \ee^{\ii\ph} \right\rangle\nonumber\db\\[0.05cm]  
= {} & \left\langle R - e\cos(n\ph) + \ii en\sin(n\ph),\, 1 \right\rangle\nonumber\db\\[0.05cm]
= {} & R - e\cos(n\ph)\,.
\end{align}

\begin{SCfigure}[][ht]
\includegraphics[width=0.5\textwidth]{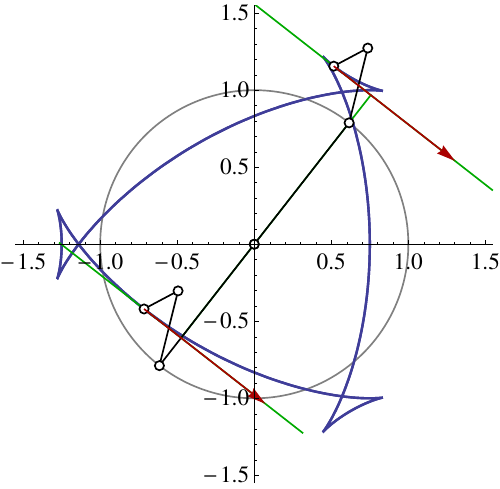}
\caption{P3 curve ($R=1$, $e=1/4$) with gen\-e\-rat\-ing mechanisms and unit tangent vectors\index{vector!tangent vector@tangent $\sim $}\\[0.2cm] positions: $\ph = 52\g$ and $\ph = 232\g$}
\label{Abb:Non-convex_P3_curve}
\end{SCfigure}

We determine the width function $w$ of a P$n$G profile from \eqref{Eq:a(phi)_profile_curve}.
For even $n = 2,4,6,\ldots$ we have
\begin{align*}
  w(\ph)
= {} & a(\ph) + a(\ph+\pi)\db\\[0.05cm]
= {} & R - e\cos(n\ph) + R - e\cos(n(\ph+\pi))\db\\[0.05cm]
= {} & R - e\cos(n\ph) + R - e\cos(n\ph+n\pi)\db\\[0.05cm]
= {} & 2R - e\cos(n\ph) - e\cos(n\ph)\db\\[0.05cm]
= {} & 2R - 2e\cos(n\ph)\db\\[0.05cm]
= {} & 2(R - e\cos(n\ph))\db\\[0.05cm]
= {} & 2 a(\ph)\,.
\end{align*}
For odd $n = 1,3,5,7,\ldots$ we have
\begin{align*}
  w(\ph)
= {} & a(\ph) + a(\ph+\pi)\db\\[0.05cm]
= {} & R - e\cos(n\ph) + R - e\cos(n(\ph+\pi))\db\\[0.05cm]  
= {} & R - e\cos(n\ph) + R - e\cos(n\ph+n\pi)\db\\[0.05cm]
= {} & 2R - e\cos(n\ph) + e\cos(n\ph)\\[0.05cm]
= {} & 2R\,.
\end{align*}
All P$n$G profiles with odd $n$ are therefore figures of constant width (for an example see Fig.\ \ref{Abb:P3G_Gleichdick}). 
In a certain sense, even a non-convex P$n$ curve can be regarded as a curve of constant width: the distance between the parallel tangents to the curve at corresponding points $z(\ph)$ and $z(\ph+\pi)$ is always equal to $2R$.
For an example see Fig.\ \ref{Abb:Non-convex_P3_curve}.
The positions of the bars in the two positions show that the distance between the tangents is actually equal to $2R = 2$. 

Now, we derive a formula for the length $L$ of the boundary of a P$n$G profile.
From \eqref{Eq:Arc_length} with \eqref{Eq:z'(phi)_with_cos_and_sin} we get
\begin{align*}
  L
= {} & \int_0^{2\pi} \left|z'(\ph)\right| \dd\ph\db\\[0.05cm]
= {} & \int_0^{2\pi} \left|\ii\ee^{\ii\ph} \left(R + \left(n^2-1\right)e\cos(n\ph)\right)\right| \dd\ph\db\\[0.05cm]  
= {} & \int_0^{2\pi} \left|R + \left(n^2-1\right)e\cos(n\ph)\right| \dd\ph\,.
\end{align*}
The expression between the absolute value bars becomes minimal for $\cos(n\ph) = -1$.
In order for the said expression to always be greater than or equal to zero, it must be $R - \left(n^2-1\right)e \ge 0$.
This is the case for P$n$G profiles (see \eqref{Eq:Condition_for_no_self-intersections_1}), and so we have
\beq
  L
= R \int_0^{2\pi} \dd\ph + \left(n^2-1\right)e \int_0^{2\pi} \cos(n\ph)\, \dd\ph
= 2\pi R\,.  
\eeq
This result is not surprising for odd $n$, because $L = 2\pi R$ holds for every figure of constant width $2R$ (see \cite[p.\ 8]{Santalo}); interestingly, however, it is also valid for P$n$G profiles with even $n$.

\begin{figure}[ht]
\begin{minipage}[t]{0.45\textwidth}
  \centering	
  \includegraphics[width=1\textwidth]{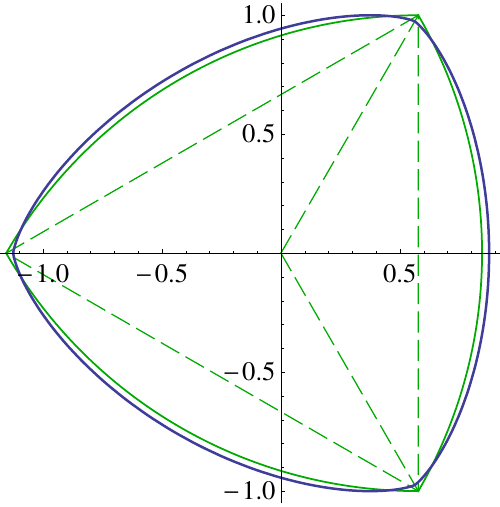}
  \caption{\small \textsc{Reuleaux} triangle (green solid line) and P3 curve ($R = 1$, $e = 1/8$, blue) with same constant width $D = 2R = 2$)}
  \label{Abb:Reuleaux_and_P3_curve01}
\end{minipage}
\hfill
\begin{minipage}[t]{0.45\textwidth}
  \centering
  \includegraphics[width=1\textwidth]{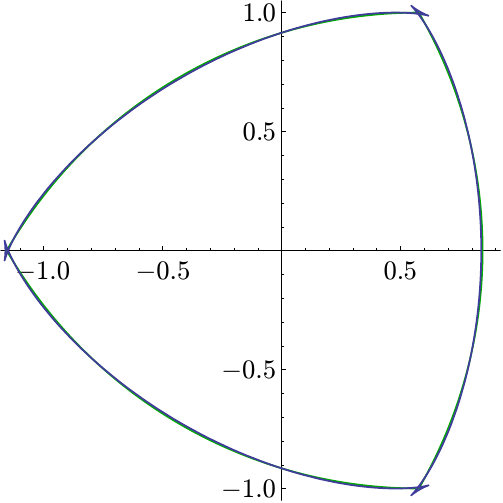}
  \caption{\small \textsc{Reuleaux} triangle (green) with $D = 2$, and P3 curve (blue) with $R = 1$ and $e \approx 0.160201$}
  \label{Abb:Reuleaux_and_P3_curve02}
\end{minipage}
\end{figure}

Now, we compare the area enclosed by a P$n$ curve, $n = 3,5,7,\ldots$, with that of a {\em \textsc{Reuleaux} triangle}\index{Reuleaux triangle} (see Fig.\ \ref{Abb:Reuleaux_and_P3_curve01}) having same constant width $D = 2R$.
For the area $A^*$ of the \textsc{Reuleaux} triangle we have
\beq
  A^*
= 3\times\frac{1}{6}\,\pi D^2 - 2\,\frac{\sqrt{3}}{4}\,D^2
= \frac{1}{2} \left(\pi - \sqrt{3}\right) D^2  
\eeq
(see also \cite{wiki:Reuleaux_Dreieck_russian}).
From this formula we get
\beq
  A^*
= \frac{1}{2} \left(\pi - \sqrt{3}\right) (2R)^2
= 2\pi R^2 -2\,\sqrt{3}\,R^2\,,
\eeq
hence
\beqn \label{Eq:A_Reuleaux_triangle}
  A^*
= \pi R^2 - \left(2\,\sqrt{3} - \pi\right) R^2\,. 
\eeqn

We derive a condition for $A^* < A$ with $A^*$ and $A$ according to \eqref{Eq:A_Reuleaux_triangle} and \eqref{Eq:A_profile_curve}, respectively.
From
\beq
  \pi R^2 - \left(2\,\sqrt{3} - \pi\right) R^2 < \pi R^2 - \frac{1}{2}\,\pi\left(n^2-1\right) e^2 
\eeq
we get
\beq
  \pi\left(n^2-1\right) e^2 < 2 \left(2\,\sqrt{3} - \pi\right) R^2\,, 
\eeq
hence
\beqn \label{Eq:(e/R)^2-1}
  \frac{\pi\left(n^2-1\right)}{2\left(2\,\sqrt{3}-\pi\right)} < \left(\frac{R}{e}\right)^2.
\eeqn
The condition \eqref{Eq:Condition_for_no_self-intersections_1} for a P$n$G profile can be written as
\beqn \label{Eq:(e/R)^2-2}
  \left(n^2-1\right)^2 \le \left(\frac{R}{e}\right)^2.
\eeqn
For the terms on the left-hand sides of the inequalities \eqref{Eq:(e/R)^2-1} and \eqref{Eq:(e/R)^2-2}, a simple calculation shows that
\beqn \label{Eq:Left-hand_sides}
  \frac{\pi\left(n^2-1\right)}{2\left(2\,\sqrt{3}-\pi\right)} < \left(n^2-1\right)^2
  \quad\mbox{if}\quad
  n \ge 3\,.
\eeqn
If $n = 1$, then \eqref{Eq:(e/R)^2-1} gives $0 < (R/e)^2$, hence the area $A = \pi R^2$ of the set enclosed by the P1 curve (with $R$ and $e$) is greater than the area of the \textsc{Reuleaux} triangle with same width $2R$, and \eqref{Eq:(e/R)^2-2} gives $0 \le (R/e)^2$, hence the enclosed set is a P1G profile.\footnote{A P1 curve is a circle of radius $R$ with center at point $(-e,0)$.}
From \eqref{Eq:Left-hand_sides} together with the special case $n = 1$ we see that the area of the \textsc{Reuleaux} triangle is smaller than that of any P$n$G profile with odd $n \ge 1$ and the same constant width.
This is special case of the \textsc{Blaschke}–\textsc{Lebesgue} theorem\index{theorem!Blaschke-Lebesgue theorem@Blaschke-Lebesgue $\sim $} which states that the \textsc{Reuleaux} triangle has the least area of all plane sets of given constant width \cite{Blaschke:Konvexe_Bereiche}, \cite[p.\ 323]{Martini&Montejano&Oliveros}.
Fig.\ \ref{Abb:Reuleaux_and_P3_curve01} shows the P3G profile with smallest possible value of $e$ and the Reuleaux triangle with the same width.
It can be seen that the area of the Reuleaux triangle is smaller than that of the P3G profile.
Fig.\ \ref{Abb:Reuleaux_and_P3_curve02} shows the P3 curve (\underline{no} P3G profile) with $R = 1$ which encloses the same area as the \textsc{Reuleaux}-Polygon with width $D = 2R =2$.
(The three small loops of this P3 curve have ``negative'' area.)
As can be seen, this P3 curve gives a very good approximation to the \textsc{Reuleaux} triangle especially when removing the small loops.
However, these small loops are necessary if the P3 curve is to be a curve of constant width in the extended sense given in the above explanations for Fig.\ \ref{Abb:Non-convex_P3_curve}.   

\textsc{Rabinowitz} gives in \cite{Rabinowitz} (see also \cite[pp.\ 111-113]{Martini&Montejano&Oliveros}) an example for a polynomial equation defining a curve {\em (\textsc{Rabinowitz} curve)}\index{curve!Rabinowitz curve@Rabinowitz $\sim $} of constant width.
\textsc{Rabinowitz} deduces this polynomial equation from the parametric equations
\beq
\begin{array}{c@{\;=\;}*{5}{c@{\;}}}
  x(\ph) & 9\cos\ph & + & 2\cos(2\ph) & - & \cos(4\ph)\,,\\[0.05cm]
  y(\ph) & 9\sin\ph & - & 2\sin(2\ph) & - & \sin(4\ph)\,. 
\end{array}
\eeq 
We can write this equations in complex form as
\beq
  z(\ph)
= 9\ee^{\ii\ph} - 2\ee^{\ii(\pi-2\ph)} - \ee^{4\ii\ph}\,.  
\eeq
The reflection with respect to the imaginary axis gives the parametric equation
\beq
  \tilde{z}(\ph)
= -9\ee^{-\ii\ph} + 2\ee^{-\ii(\pi-2\ph)} + \ee^{-4\ii\ph}\,.  
\eeq
We substitute $\ph \rightarrow \pi-\ph$ and get
\begin{align} \label{Eq:Rabinowitz_as_special_P3_curve}
  \hat{z}(\ph)
:= {} & \tilde{z}(\pi-\ph)\nonumber\\[0.05cm]
= {} & {-}9\ee^{-\ii(\pi-\ph)} + 2\ee^{-\ii(\pi-2(\pi-\ph))} + \ee^{-4\ii(\pi-\ph)}\nonumber\db\\[0.05cm]  
= {} & {-}9\ee^{-\ii\pi}\ee^{\ii\ph} + 2\ee^{-\ii(-\pi+2\ph)} + \ee^{\ii(4\ph-4\pi)}\nonumber\\[0.05cm]
= {} & 9\ee^{\ii\ph} + 2\ee^{\ii(\pi-2\ph)} + \ee^{4\ii\ph}\,. 
\end{align}
This is \eqref{Eq:z(phi)_for_profile} with $R = 9$, $n = 3$ and $e = 1$, and consequently the \textsc{Rabinowitz} curve is a special P3 curve (see Fig.\ \ref{Abb:Rabinowitz_curve}).
Substituting $-x$ for $x$ in the polynomial equation of \textsc{Rabinowitz} gives the polynomial equation
\begin{align*}
720^3 = {} & \left(x^2+y^2\right)^4 - 45\left(x^2+y^2\right)^3 - 41283\left(x^2+y^2\right)^2 + 7950960\left(x^2+y^2\right)\\[0.05cm]
& + 16\left(x^2-3y^2\right)^3 + 48\left(x^2+y^2\right)\left(x^2-3y^2\right)^2\\[0.05cm]
& - \left(x^2-3y^2\right)x\left(16\left(x^2+y^2\right)^2-5544\left(x^2+y^2\right)+266382\right)   
\end{align*}
of the P3 curve defined by equation \eqref{Eq:Rabinowitz_as_special_P3_curve}.

\subsection{Manufacturing}

\begin{figure}[ht!]
\begin{minipage}[t]{0.48\textwidth}
  \centering
  \includegraphics[width=0.9\textwidth]{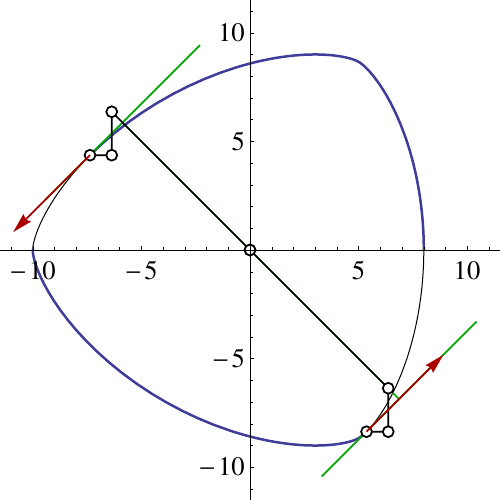}
  \caption{\small \textsc{Rabinowitz} curve as special P3 curve with generating mechanisms according to \eqref{Eq:Rabinowitz_as_special_P3_curve} and tangent vectors normalized to length 5; positions $\ph = 135\g$ and $\ph = 315\g$}
\label{Abb:Rabinowitz_curve}
\end{minipage}
\hfill
\begin{minipage}[t]{0.48\textwidth}
  \centering
  \includegraphics[width=0.9\textwidth]{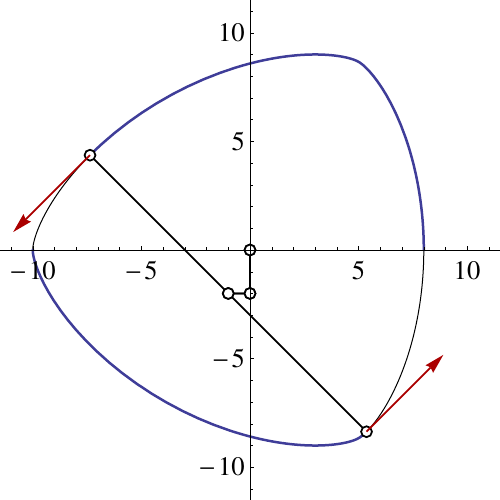}
  \caption{\small Same curve as in Fig.\ \ref{Abb:Rabinowitz_curve}, but with changed order of the bars of the generating mechanism; positions $\ph = 135\g$ and $\ph = 315\g$}
\label{Abb:Changed_order_of_bars}
\end{minipage}
\end{figure}

The three terms on the right-hand side of \eqref{Eq:z(phi)_for_profile} can be arranged in six different orders, which therefore also applies to the bars of the generating mechanism.
Using the order
\beqn \label{Eq:z(phi)_for_profile_order:bca}
  z(\ph)
= \ell_3\,\ee^{\ii[\pi-(n-1)\ph]} + \ell_4\,\ee^{\ii(n+1)\ph} + \ell_2\,\ee^{\ii\ph}  
\eeqn
gives in the special case of the \textsc{Rabinowitz} curve the mechanisms in Fig.\ \ref{Abb:Changed_order_of_bars}.
(In this subsection, we will always use the \textsc{Rabinowitz} curve in the examples.)

\begin{figure}[ht!]
\begin{minipage}[t]{0.48\textwidth}
  \centering
  \includegraphics[width=1\textwidth]{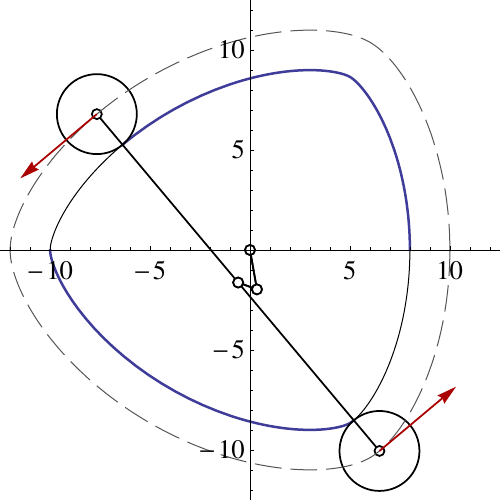}
  \caption{\small Manufacturing an outer contour (shaft) using circular tools; positions $\ph = 130\g$ and $\ph = 310\g$}
  \label{Abb:Manufacturing_outer_contour}
\end{minipage}
\hfill
\begin{minipage}[t]{0.48\textwidth}
  \centering
  \includegraphics[width=1\textwidth]{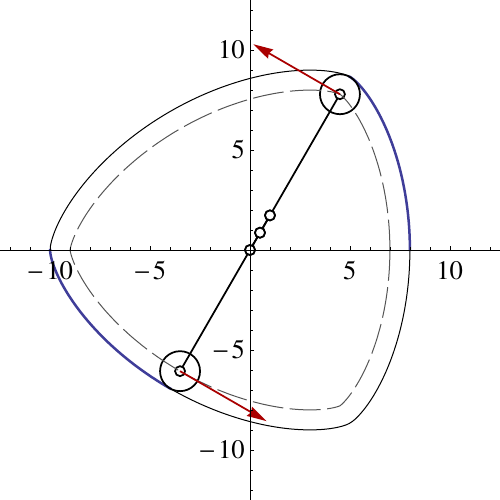}
  \caption{\small Manufacturing an inner contour (hub) using circular tools; positions $\ph = 60\g$ and $\ph = 240\g$}
  \label{Abb:Manufacturing_inner_contour}
\end{minipage}
\end{figure}

As already shown, the tangent vectors (here: normalized to length 5) are orthogonal to the bars of length $\ell_2 = \overline{A_0A} = R$ (in Fig.\ \ref{Abb:Changed_order_of_bars}: $\ell_2 = 9$).
This makes it possible to manufacture a profile with a circular tool (milling cutter\index{milling cutter} or grinding disk\index{grinding disk}) whose radius we denote by $r$.
The workpiece is fixed, and the path curve of the tool center point is the parallel curve with distance $r$ to the contour to be produced.
Fig.\ \ref{Abb:Manufacturing_outer_contour} shows a principle of manufacuring the outer contour (shaft\index{shaft}) using a tool of radius $r = 2$, and Fig.\ \ref{Abb:Manufacturing_inner_contour} shows this principle used for manufacturing the inner contour of the hub\index{hub} using a tool of radius $r = 1$.
For this inner contour, $r = 1$ is the maximum possible tool radius because, as can be easily calculated, 1 is the smallest radius of curvature of the curve.
On the other hand, the tool radius $r$ for the outer contour can be chosen arbitrarily large. 

\begin{figure}[ht]
\begin{minipage}[t]{0.48\textwidth}
  \includegraphics[width=\textwidth]{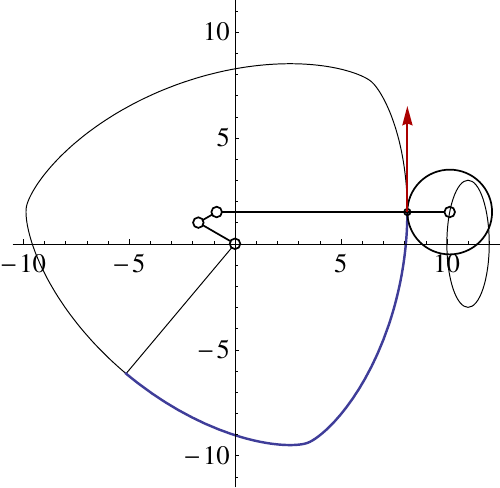}
  \caption{\small Manufacturing the outer contour of a rotating workpiece (shaft); position $\ph = 130\g$ (kinematic inversion\index{kinematic inversion} of the mechanism in Fig.\ \ref{Abb:Manufacturing_outer_contour})}
  \label{Abb:Rabinowitz_rotating_bca}
\end{minipage}
\hfill
\begin{minipage}[t]{0.48\textwidth}
  \includegraphics[width=\textwidth]{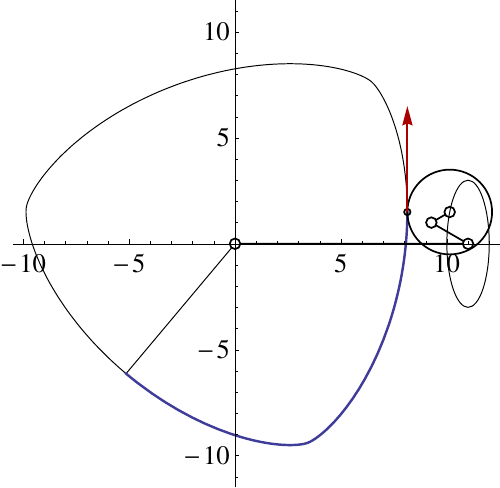}
  \caption{\small Manufacturing the outer contour of a rotating workpiece (shaft); position $\ph = 130\g$}
  \label{Abb:Rabinowitz_rotating_abc}
\end{minipage}
\end{figure}

Multiplying \eqref{Eq:z(phi)_for_profile_order:bca} by $\ee^{-\ii\ph}$ (kinematic inversion\index{kinematic inversion}) gives
\beqn \label{Eq:Rotating_curve_1}
  z(\ph)\,\ee^{-\ii\ph}
= \ell_3\,\ee^{\ii(\pi-n\ph)} + \ell_4\,\ee^{\ii n\ph} + \ell_2\,.
\eeqn
This equation means that the curve (the workpiece) now rotates, while the bar with length $\ell_2$ only carries out a translational motion (see Fig.\ \ref{Abb:Rabinowitz_rotating_bca}).
The vector shown is the tangent vector normalized to length 5 to the curve at its currently generated point.
We consider the path curve of the tool center point.
The equation of this curve is given by the right-hand side of \eqref{Eq:Rotating_curve_1}, hence, changing the order of the summands, using $\ell_2$, $\ell_3$ and $\ell_4$ from \eqref{Eq:z(phi)_for_profile}, and denoting, as above, by $r$ the tool radius, we have
\begin{align*}
  \ell_2 + r + \ell_3\,\ee^{\ii(\pi-n\ph)} + \ell_4\,\ee^{\ii n\ph}
= {} & \ell_2 + r - \ell_3\,\ee^{-\ii n\ph} + \ell_4\,\ee^{\ii n\ph}\db\\[0.05cm]
= {} & R + r - \frac{1}{2}\,(n+1)\,e\,\ee^{-\ii n\ph} + \frac{1}{2}\,(n-1)\,e\,\ee^{\ii n\ph}\db\\[0.05cm]
= {} & R + r - e\,\frac{\ee^{\ii n\ph}+\ee^{-\ii n\ph}}{2} + ne\,\frac{\ee^{\ii n\ph}-\ee^{-\ii n\ph}}{2}     
\end{align*}
and so
\beqn \label{Eq:Ellipse}
  \ell_2 + r + \ell_3\,\ee^{\ii(\pi-n\ph)} + \ell_4\,\ee^{\ii n\ph}
= R + r - e\cos(n\ph) + \ii ne\sin(n\ph)\,.  
\eeqn

This is the parametric equation of an ellipse (cf.\ \eqref{Eq:z(phi)_with_cos_and_sin}) passing through $n$ times in the mathematically negative sense with center of gravity at point $R + r$, and starting from point $R + r - e$.
The length of the minor semi-axis is equal to $e$, that of the major one is equal to $ne$.

Writing \eqref{Eq:Rotating_curve_1} as
\beq
  z(\ph)\,\ee^{-\ii\ph}
= \ell_2 + \ell_3\,\ee^{\ii(\pi-n\ph)} + \ell_4\,\ee^{\ii n\ph}
\eeq
gives the mechanism in Fig.\ \ref{Abb:Rabinowitz_rotating_abc}.
The parametric equation of the tool center point is again given by \eqref{Eq:Ellipse}.
Compared to the variant in Fig.\ \ref{Abb:Rabinowitz_rotating_bca}, the variant in Fig.\ \ref{Abb:Rabinowitz_rotating_abc} has the advantage that one part less has to be moved.
The variant in Fig.\ \ref{Abb:Rabinowitz_rotating_abc} is essentially the kinematic inversion of the variant in Fig.\ \ref{Abb:Rabinowitz_curve}. 
This can be seen from \eqref{Eq:z'(phi)_with_cos_and_sin}: If we replace $\ell_2$ by $\ell_2 + r$, and hence $R$ by $R + r$, then the tangent vector only changes its length, but not its direction.
The not shown tangent vector of the tool center point in Fig.\ \ref{Abb:Rabinowitz_rotating_abc} with respect to the rotating curve has the same direction as the shown tangent vector to the curve at its currently generated point.
Of course, this must be the case, as the tool center point curve with respect to the rotating workpiece is the parallel curve of the contour curve to be generated.      
The variant in Fig.\ \ref{Abb:Rabinowitz_rotating_abc} is essentially the solution in \cite[Fig.\ 9]{Durrenbach}, \cite[Fig.\ 2]{Filemon}.
There, too, the tool center point describes the aforementioned ellipse (with respect to the fixed part/coordinate system).

\begin{figure}[ht]
\begin{minipage}[t]{0.48\textwidth}
  \centering
  \includegraphics[width=0.86\textwidth]{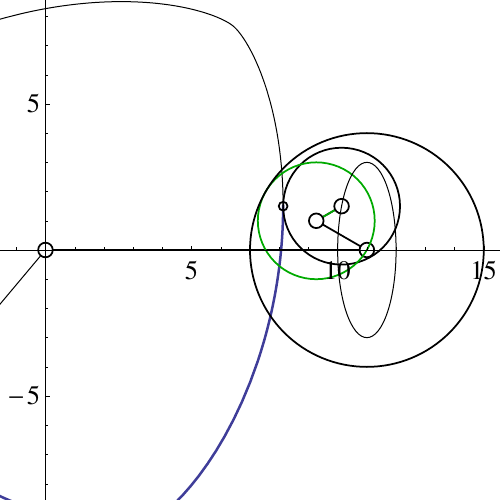}
  \caption{\small Epicyclic gear train for generating the motion of the tool center point in Fig.\ \ref{Abb:Rabinowitz_rotating_abc}; tool radius $r = 2$, position $\ph = 130\g$\\}
  \label{Abb:Epicyclic_gear_train_1}
\end{minipage}
\hfill
\begin{minipage}[t]{0.48\textwidth}
  \centering
  \includegraphics[width=0.86\textwidth]{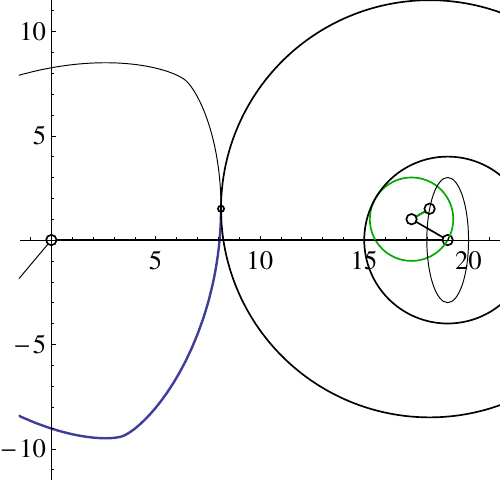}
  \caption{\small Epicyclic gear train as in Fig.\ \ref{Abb:Epicyclic_gear_train_1}; tool radius $r = 10$, position $\ph = 130\g$}
  \label{Abb:Epicyclic_gear_train_2}
\end{minipage}
\end{figure}

From \eqref{Eq:rho_4_and_rho_2} it is clear that the elliptical path of the tool center point can be generated by rolling a circle (green) of radius $\rh_4 = (n+1)\,e/2$ within a fixed circle of double radius $\rh_2$ (see Fig.\ \ref{Abb:Epicyclic_gear_train_1}), and consequently the ellipse\index{ellipse} is a special hypotrochoid\index{hypotrochoid}. 
These circles can be the pitch circles of an external and an internal gear wheel, hence the mechanical solution is an epicyclic gear train\index{epicyclic gear train}.
It is a coincidence that the tool and the rolling gear have the same radius $r = 2$ in Fig.\ \ref{Abb:Epicyclic_gear_train_1}.
Fig.\ \ref{Abb:Epicyclic_gear_train_2} shows a variant with tool radius $r = 10$.  

We can write \eqref{Eq:Rotating_curve_1} as
\beq
  z(\ph)\,\ee^{-\ii\ph} - \ell_3\,\ee^{\ii(\pi-n\ph)} - \ell_4\,\ee^{\ii n\ph}
= \ell_2 
\eeq
and consequently as
\beq
  \ell_3\,\ee^{-\ii n\ph} + \ell_4\,\ee^{\ii(\pi+n\ph)} + z(\ph)\,\ee^{-\ii\ph}
= \ell_2\,.
\eeq

\begin{SCfigure}[][ht]
\includegraphics[width=0.45\textwidth]{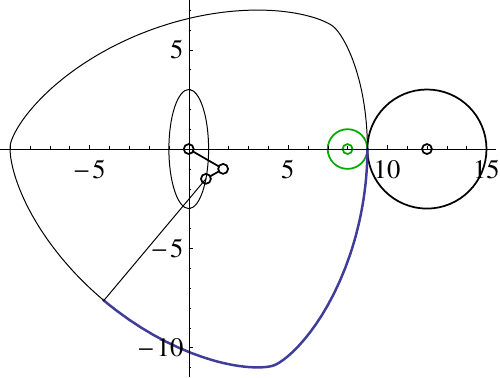}
\caption{Variant with fixed tool center points; position $\ph = 130\g$}
\label{Abb:Fixed_tool_center_points}
\end{SCfigure}

This gives the solution in Fig.\ \ref{Abb:Fixed_tool_center_points}, in which possible tools with fixed center points for manufacturing the outer or inner contour are shown (see also \cite[Fig.\ 26]{Filemon}).
This variant is essentially a kinematic inversion of the solution in Fig.\ \ref{Abb:Rabinowitz_rotating_bca} and has the advantage that an existing stationary grinding wheel can be used, and the motion of the workpiece can be realized by an additional device; however, the motion of the workpiece is more difficult to realize than, for example, with the variant in Figs.\ \ref{Abb:Epicyclic_gear_train_1} and \ref{Abb:Epicyclic_gear_train_2}. 

\begin{figure}[ht]
\begin{minipage}[t]{0.45\textwidth}
  \includegraphics[width=\textwidth]{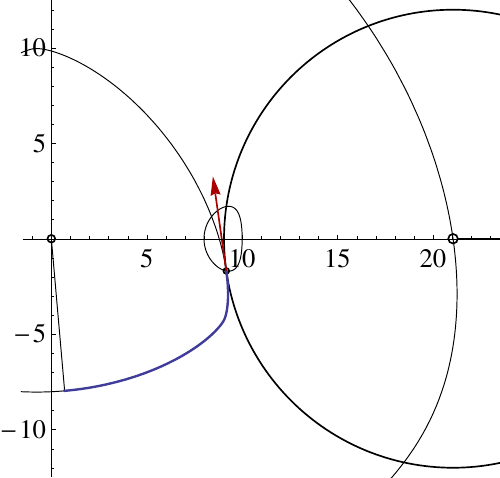}
  \caption{\small Tool with linear motion along the $x$-axis; tool radius $r = 12$, position $\ph = 85\g$}
  \label{Abb:Rabinowitz_linear_1}
\end{minipage}
\hfill
\begin{minipage}[t]{0.45\textwidth}
  \includegraphics[width=\textwidth]{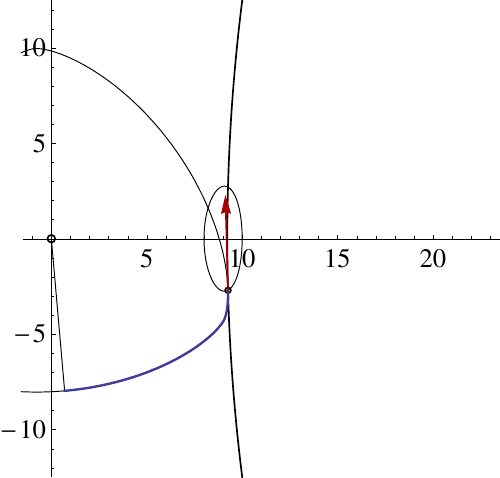}
  \caption{\small Tool with linear motion along the $x$-axis; tool radius $r = 100$, position $\ph = 85\g$}
  \label{Abb:Rabinowitz_linear_2}
\end{minipage}
\end{figure}

The production variants discussed up to now require special machine tools.
Nowadays, flexible CNC-machines\footnote{computer numerically controlled machines} are preferred for manufacturing P$n$ profiles.
For our purposes, a CNC-machine in which the tool center can perform a controlled rectilinear motion is sufficient.
We assume that the tool center point moves along the (horizontal) $x$-axis.
Then the current $x$-coordinate of the tool center point is the right intersection between the parallel curve at distance $r$ to the P$n$ curve rotated by angle $-\ph$ and the $x$-axis (see Fig.\ \ref{Abb:Rabinowitz_linear_1}).

\begin{figure}[ht]
\begin{minipage}{0.48\textwidth}
  \includegraphics[width=\textwidth]{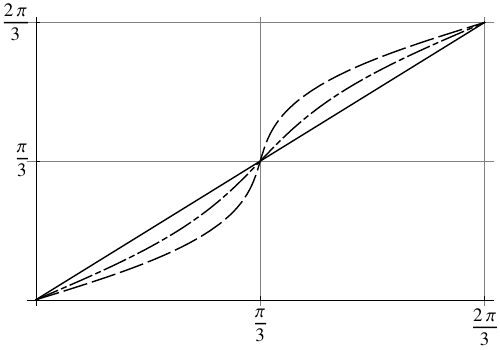}
  \caption{\small Graphs of the functions $f_1$ (dashed), $f_{12}$ (dot-dashed), and $f_\infty$ (solid)}
  \label{Abb:Rabinowitz_linear_angle}
\end{minipage}
\hfill
\begin{minipage}{0.48\textwidth}
  \includegraphics[width=\textwidth]{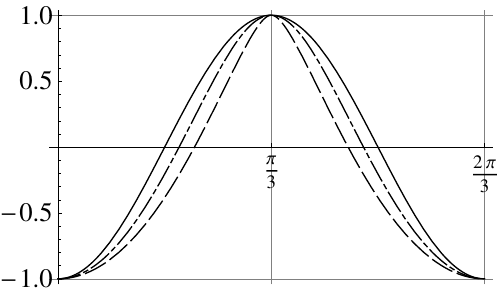}
  \caption{\small Graphs of the function $\widetilde{x}_1$ (dashed), $\widetilde{x}_{12}$ (dot-dashed), and $\widetilde{x}_\infty$ (solid)}
  \label{Abb:Rabinowitz_linear_motion}
\end{minipage}
\end{figure}

In order to get this intersection, we must first determine the value $\psi$ of the parameter of the parametric curve $z$ at the intersection point, for which we define in {\em Mathematica}\index{Mathematica}, using \texttt{FindRoot}\index{FindRoot@\texttt{FindRoot}}, the function $f_r(\ph)$ by means of
\begin{verbatim}
    f[\[CurlyPhi]_, r_] := \[Psi] /. 
        FindRoot[Im[(z[\[Psi]] - I r z1[\[Psi]]/Abs[z1[\[Psi]]])*
        Exp[-I \[CurlyPhi]]] == 0, {\[Psi], \[CurlyPhi]}]
\end{verbatim}
where \texttt{z1} is the first derivative of $z$ (see \eqref{Eq:z'(phi)_with_cos_and_sin}).
The $x$-coordinate function $x_r(\ph)$ of the tool center point is then calculated using
\begin{verbatim}
  Re[(z[f[\[CurlyPhi], r]] - 
     I r z1[f[\[CurlyPhi], r]]/Abs[z1[f[\[CurlyPhi], r]]])*Exp[-I \[CurlyPhi]]]
\end{verbatim}
For example, the graphs of the functions $f_1$ and $f_{12}$ are shown in Fig.\ \ref{Abb:Rabinowitz_linear_angle}.
The corresponding graphs of the shifted functions
\beq
  \widetilde{x}_r(\ph) := x_r(\ph) - R - r
\eeq
for $r = 1$ and $r = 12$ are shown in Fig.\ \ref{Abb:Rabinowitz_linear_motion}.

For $r \rightarrow \infty$, the contact point between tool and workpiece tends towards the contact point in traditional manufacturing with rotating workpiece (see Figs.\ \ref{Abb:Rabinowitz_rotating_abc}, \ref{Abb:Epicyclic_gear_train_1}, \ref{Abb:Epicyclic_gear_train_2}), hence
\beq
  f_\infty(\ph) := \lim_{r\rightarrow\infty} f_r(\ph) = \ph
\eeq
(see Fig.\ \ref{Abb:Rabinowitz_linear_angle}).
For $r \rightarrow \infty$, the path curve of the contact point tends to the ellipse with parametric equation
\beq
  R - e\cos(n\ph) + \ii ne\sin(n\ph)  
\eeq
(see Figs.\ \ref{Abb:Rabinowitz_linear_1}, \ref{Abb:Rabinowitz_linear_2} and \ref{Abb:Rabinowitz_rotating_bca}, \ref{Abb:Rabinowitz_rotating_abc}, and cf.\ \eqref{Eq:z(phi)_with_cos_and_sin}, \eqref{Eq:Ellipse}).
Note that $R - e\cos(n\ph)$ is the support function of the P$n$ curve (see \eqref{Eq:a(phi)_profile_curve}).
From \eqref{Eq:Ellipse} it immediately follows that
\beq
 \widetilde{x}_\infty(\ph) := \lim_{r\rightarrow\infty} \widetilde{x}_r(\ph) = - e\cos(n\ph)
\eeq
(see Fig.\ \ref{Abb:Rabinowitz_linear_motion})
         
In \textcite{DIN3689-1-engl}, the use of hypotrochoids\index{hypotrochoid} as profile curves is standardized, which takes us back to the starting point of the development in \textcite{Musyl46} (see also \cite{Ziaei07b-engl}), where, however, epitrochoids were used.
Of course, CNC-machines allow the production of almost arbitrary profile curves; however, the resulting flexibility also makes the practical application much more complex and difficult to handle.

%% file: DiffGeo5_L.tex

\addcontentsline{toc}{section}{References}
\printbibliography

\bigskip
{\small\textsc{Uwe Bäsel, Leipzig University of Applied Sciences (HTWK Leipzig), Faculty of Engineering, Germany}, \texttt{uwe.baesel@htwk-leipzig.de}}

\addcontentsline{toc}{section}{Index}
\printindex